\documentclass{article}
\usepackage{graphicx} 
\usepackage{amsthm}
\usepackage[english]{babel}
\usepackage[letterpaper,top=2cm,bottom=2cm,left=3cm,right=3cm,marginparwidth=1.75cm]{geometry}

\newtheorem{remark}{Remark}[section]

\newtheorem{theorem}{Theorem}[section]
\newtheorem{lemma}[theorem]{Lemma}
\newtheorem{corollary}[theorem]{Corollary}
\newtheorem{proposition}[theorem]{Proposition}

\usepackage{amsmath}
\usepackage{amssymb}
\usepackage{graphicx}
\usepackage{subfigure}
\usepackage{float}
\usepackage{caption}
\usepackage{wrapfig}
\usepackage{longtable}
\usepackage{listings}
\usepackage{todonotes}
\usepackage{makecell}
\usepackage{amsmath,amssymb}
\usepackage{algorithm}
\usepackage{algpseudocode}
\usepackage[most]{tcolorbox}
\usepackage{tabularx}
\usepackage{enumitem}
\usepackage{mathrsfs}
\usepackage{lmodern}

\newcommand{\dInt}[1]{{#1\in[1,p]_\mathbb{N}}}

\usepackage[colorlinks=true, allcolors=blue]{hyperref}

\usepackage{enumitem}
\numberwithin{equation}{section}
\usepackage{authblk}
\title{Propagation failure for bistable differential equations on periodic trees}
\author[1]{Hermen Jan Hupkes \thanks{\tt hhupkes@math.leidenuniv.nl}}
\author[2]{Petr Stehl\'{\i}k\thanks{\tt pstehlik@fav.zcu.cz}}
\author[2]{Vladim\'{\i}r \v{S}v\'{\i}gler \thanks{\tt sviglerv@fav.zcu.cz}}
\author[1]{Jialing Yu \thanks{\tt j.yu@math.leidenuniv.nl}}
\affil[1]{\small Mathematisch Instituut, Universiteit Leiden, P.O. Box 9512, 2300 RA Leiden, The Netherlands}
\affil[2]{\small Department of Mathematics and NTIS, Faculty of Applied Sciences, University of West Bohemia,\authorcr Univerzitn\'\i~8, 306 14 Plze\v{n}\\ Czech Republic}

\begin{document}
\maketitle

%%%%%%%%%%%%%%%%%%%%%%%%%%%%%%%%%%%%%%%%%%%%%%%%%%
%%%%%%%%%%%%%%%%%%%% ABSTRACT %%%%%%%%%%%%%%%%%%%%
%%%%%%%%%%%%%%%%%%%%%%%%%%%%%%%%%%%%%%%%%%%%%%%%%%
\begin{abstract}
    We study traveling fronts for bistable reaction-diffusion equations on periodic heterogeneous trees. We use comparison principles to estimate the pinning and propagation regions for both cubic and general bistabilities. Our analytical results and numerical simulations
    show that the interesting  propagation reversal effects found in homogeneous trees become more pronounced in the heterogeneous setting, paving the way towards understanding more general graphs. % and lead to non-smooth pinning boundaries. 
    As the diffusion strength increases, the waves end up moving in the same direction irrespective of the bistability. This direction is determined solely by the underlying structure of the heterogeneous tree. Our techniques also lead to sharper results for the homogeneous setting.
\end{abstract}



%%%%%%%%%%%%%%%%%%%%%%%%%%%%%%%%%%%%%%%%%%%%%%%%%%%%%%%%%%%%%%%%%%%
%%%%%%%%%%%%%%%%%%%%%%%%%%%%%%%%%%%%%%%%%%%%%%%%%%%%%%%%%%%%%%%%%%%
% INTRODUCTION
%%%%%%%%%%%%%%%%%%%%%%%%%%%%%%%%%%%%%%%%%%%%%%%%%%%%%%%%%%%%%%%%%%%
%%%%%%%%%%%%%%%%%%%%%%%%%%%%%%%%%%%%%%%%%%%%%%%%%%%%%%%%%%%%%%%%%%%

\section{Introduction}\label{sec:intro}

\begin{figure}[t]
    \centering
    \includegraphics[width=1\linewidth]{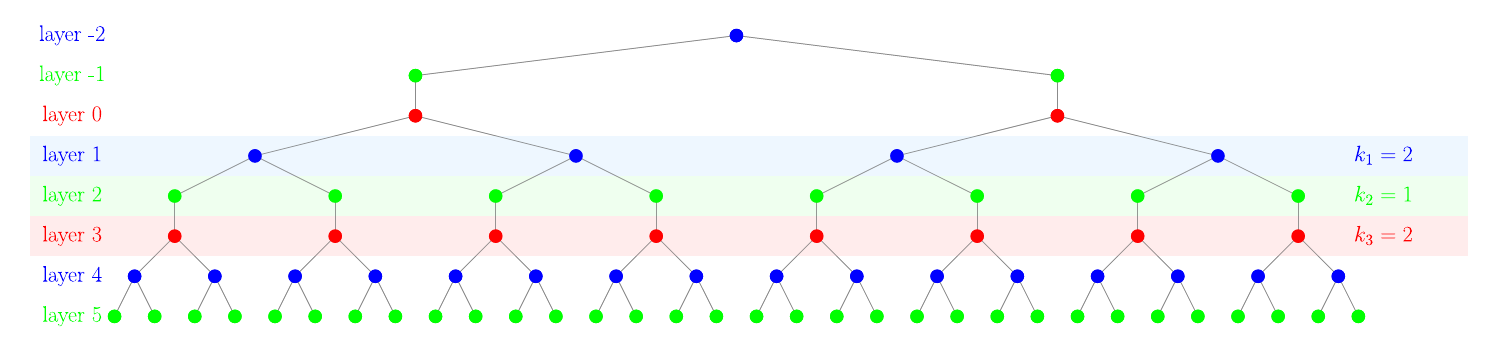}
    \caption{A part of a biinfinite periodic tree with $\mathbf{k}=\langle 2,1,2\rangle$. We focus on layer solutions, where the values $u_i(t)$  are identical for all vertices in the layer $i\in\mathbb{Z}$.}
    \label{fig:periodic_tree}
\end{figure}

In this paper we investigate traveling fronts for the periodically heterogeneous bistable reaction-diffusion-advection lattice differential equation 
\begin{equation}\label{eq:LDE}
    \dot{u}_i=d(k_iu_{i+1}-(k_i+1)u_i+u_{i-1})+g(u_i;a),\qquad i\in\mathbb{Z},\ t>0.
\end{equation}
Here $d>0$ is a diffusion parameter and the reaction term \(g(u;a)\) is assumed to be bistable, for example the standard cubic nonlinearity
\begin{equation}\label{eq:cubic}
    g_\mathrm{cub}(u;a)=u(1-u)(u-a),
    \qquad
    a\in(0,1).
\end{equation}
The spatial heterogeneity is encoded by a $p$-periodic positive coefficient sequence
\[
    \mathbf k=(k_i)_{i\in\mathbb Z} = \langle k_1,k_2,\ldots,k_p\rangle,
\]
where the angle brackets indicate the biinfinite repetition of the fundamental block $(k_1,k_2\ldots, k_p)$ of length $p\in\mathbb{N}$, i.e., $k_{i+p}=k_i$. When the coefficients \(k_i\) are natural numbers, the sequence \(\mathbf k\) can be interpreted as a periodically varying branching factor of a biinfinite tree. The results developed here extend those of \cite{hupkes2023propagation} from the homogeneous setting (i.e., $p=1$ and $\mathbf{k}=\langle k \rangle$, $k\in\mathbb{N}$) to periodically heterogeneous media with $p>1$, $\mathbf{k}= \langle k_1,k_2,\ldots,k_p\rangle$, see Figures~\ref{fig:periodic_tree} and~\ref{fig:trees_examples}. The solution of the LDE~\eqref{eq:LDE} is then the so-called layer solution on the given tree, a solution that has the same value on each tree layer, see Figure~\ref{fig:periodic_tree} for illustration. Alternatively, when the coefficients $k_i>0$ are reals, we can interpret them as advection coefficients, since the equation~\eqref{eq:LDE} can be rearranged into
\[
\dot{u}_i=d(u_{i+1}-2 u_i+u_{i-1}) + d (k_i-1) (u_{i+1}-u_i) + g(u_i;a).
\]

Our analysis concerns traveling waves of the form
\[
    u_i(t)=\Phi_i(i-ct),
    \qquad
    \Phi_{i+p}=\Phi_i,
\]
whose $p$ phase-dependent profiles $\boldsymbol{\Phi}=\langle\Phi_1,\ldots,\Phi_p\rangle$ connect the two stable states,
\[
    \lim_{\xi\to-\infty}\Phi_i(\xi)=0,
    \qquad
    \lim_{\xi\to+\infty}\Phi_i(\xi)=1,\qquad \dInt{i}.
\]
The central question is how the detuning parameter \(a\in(0,1)\), the diffusion strength \(d>0\), and the periodic branching (advection) coefficients \(\mathbf k\) determine whether the front advances, retreats, or remains pinned.

\begin{figure}[t]
        \centering
        \subfigure[$\mathbf{k}=\langle 1 \rangle$]{
            \includegraphics[width=0.45\linewidth]{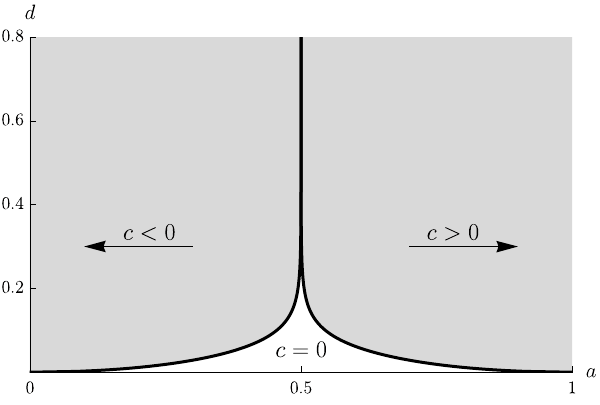}
        }
        \subfigure[$\mathbf{k}=\langle 2 \rangle$]{
            \includegraphics[width=0.45\linewidth]{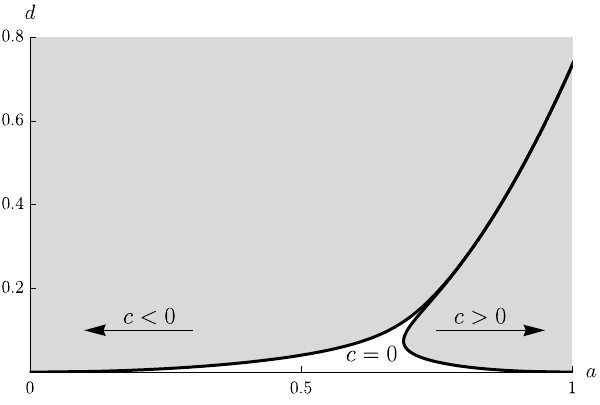}
        }
        \subfigure[$\mathbf{k}=\langle 2,1 \rangle$]{
            \includegraphics[width=0.45\linewidth]{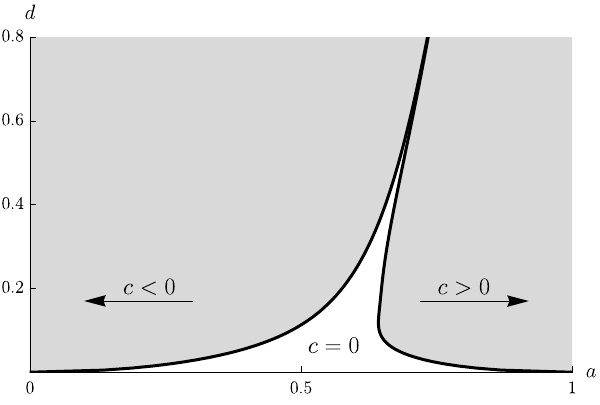}
        }
        \subfigure[$\mathbf{k}=\langle 2,1,1 \rangle$]{
            \includegraphics[width=0.45\linewidth]{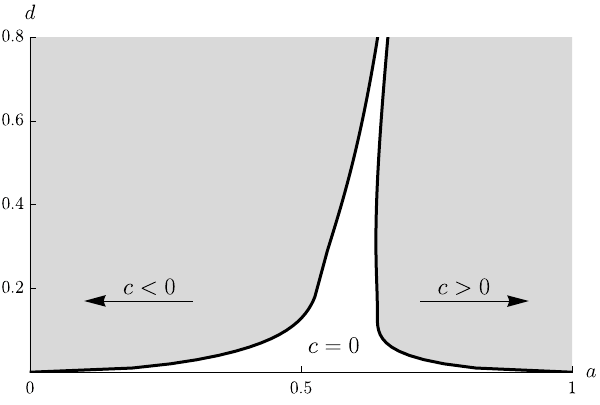}
        }
        \caption{Propagation failure regions for (a) the integer lattice, (b) the homogeneous binary tree, (c)--(d) periodically branching trees (cf. Figure~\ref{fig:trees_examples}). The propagation failure region (white) depicts values $(a,d)$ for which the waves do not move. In cases (b)--(d) we observe propagation reversal, i.e., by increasing $d>0$ for a fixed $a\in(0,1)$ the wave is pinned $c=0$ for small $d$, then moves right $c>0$ then is pinned $c=0$ again and eventually moves left $c<0$.}  \label{fig:4_ad_regions}
\end{figure}

\begin{figure}[t]
    \centering
        \subfigure[]{
            \includegraphics[height=0.15\textheight]{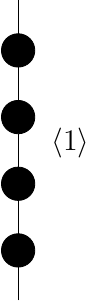}
            \label{fig:tree-1}
        }
        \subfigure[]{
            \includegraphics[height=0.15\textheight]{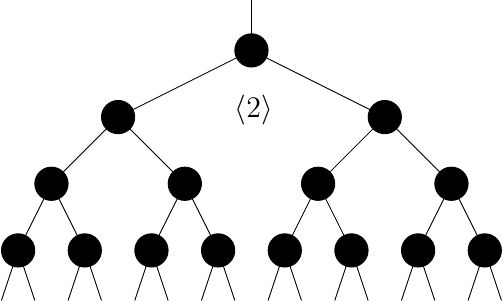}
            \label{fig:tree-2}
        }
        \subfigure[]{
            \includegraphics[height=0.15\textheight]{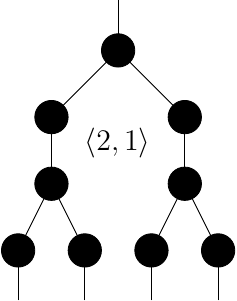}
            \label{fig:tree-2x1}
        }
        \subfigure[]{
            \includegraphics[height=0.15\textheight]{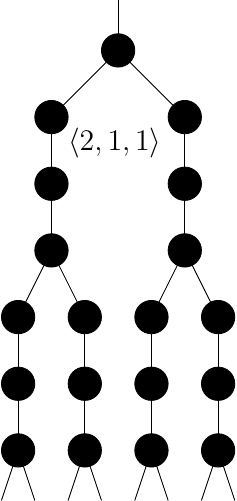}
            \label{fig:tree-2x1x1}
        }
        \subfigure[]{
            \includegraphics[height=0.15\textheight]{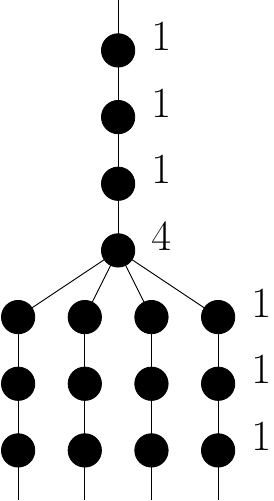}
            \label{fig:tree-4-aper}
        }
    \caption{From left to right: Parts of a homogeneous unary tree with $\mathbf{k}=\langle 1 \rangle$, homogeneous binary tree with $\mathbf{k}=\langle 2 \rangle$, a heterogeneous periodic tree with $\mathbf{k}=\langle 2,1 \rangle$, and a heterogeneous periodic tree with $\mathbf{k}=\langle 2,1,1 \rangle$. The rightmost illustration (e) depicts part of an aperiodic tree with a single four-fold branching. Consult Figure~\ref{fig:4_ad_regions} for respective pinning and traveling regions of traveling waves on structures (a)--(d).} 
    \label{fig:trees_examples}
\end{figure}

Our principal conclusions are summarized as follows (and illustrated in Figure~\ref{fig:4_ad_regions}):
\begin{enumerate}
    \item For a general bistable nonlinearity \(g\), wave pinning occurs in a nonempty analytical region in which the speed \(c=0\); see Theorem~\ref{thm:periodic-pinning}.

    \item For the cubic nonlinearity $g_\mathrm{cub}$, an improved analytical wave pinning region is derived by a separate approach that exploits the explicit form of $g_\mathrm{cub}$; see Theorem~\ref{thm:pinning}.
    
    \item For the cubic nonlinearity $g_\mathrm{cub}$, an explicit criterion for pinning in the small-diffusion regime is obtained; see Corollary~\ref{cor:envelope}.

    \item For a general bistable nonlinearity and all $a\in(0,1)$ we have $c<0$ whenever $d$ is sufficiently large and the product
    \begin{equation}\label{eq:khat}
        \widehat{k} := \prod_{i=1}^p k_i
    \end{equation}
    satisfies $\widehat{k}>1$, see Theorem~\ref{t:left-travel}. 
\end{enumerate}

Taken together, these results demonstrate that the three principal phenomena that were identified in \cite{hupkes2023propagation} for the spatially homogeneous problem, namely:
\begin{enumerate}[label=(\alph*)]
\item wave pinning - the inability of waves to move (i.e., $c=0$) for small diffusion $d$, 
\item diffusion-driven propagation reversal - the change of wave direction (i.e., the sign of $c$) solely by the decrease/increase of diffusion $d$,
\item strictly negative wave speed in the strong-diffusion regime $d\gg0$,
\end{enumerate}
 also occur in periodically heterogeneous media. See Figure~\ref{fig:4_ad_regions} where the wave pinning occurs in all four cases, while the propagation reversal and negative wave speed for large $d$ arise in panels (b)--(d). The underlying comparison and sub-solution mechanisms remain closely related to those in the homogeneous setting, but the alignment of the front with the periodic coefficients introduces phase-dependent constraints that substantially affect the propagation and pinning regions. 

\paragraph{Background and motivation.}
Front propagation in spatially discrete bistable systems is strongly influenced by both the coupling strength and the underlying geometry of discrete media. Already on one-dimensional lattices, spatial discreteness may prevent a front from propagating and produce pinning over nontrivial parameter regimes; see, for example,
\cite{keener1987propagation}. When the integer lattice, Figure~\ref{fig:trees_examples}(a), is replaced by a homogeneous tree, Figure~\ref{fig:trees_examples}(b), the geometry introduces an additional directional effect and its asymmetry leads to propagation, pinning and retreating mechanisms. More precisely, for homogeneous $k$-ary trees, Figure~\ref{fig:trees_examples}(b), these effects are investigated in \cite{kouvaris2012traveling} and, more recently, in \cite{hupkes2023propagation}. In particular, the latter work establishes rigorous sufficient conditions for pinning and propagation by means of the comparison principle. Together with \cite{hupkes2024propagation}, which analyses the large-diffusion regime through spectral convergence and provides an asymptotic description of propagation reversal, these results demonstrate that the branching factor enters directly into the balance between diffusion and bistable reaction that determines the motion of the interface.

Geometrically speaking, periodically branching trees (see Figures~\ref{fig:periodic_tree} and~\ref{fig:trees_examples}) provide a natural next step. They allow the branching geometry to vary from layer to layer while retaining translation invariance over one complete period. A general theory of discrete periodic bistable media including the existence, uniqueness and stability of pulsating traveling waves is developed in \cite{chen2008traveling}. Hence, periodic models incorporate nontrivial spatial heterogeneity and maintain the structural regularity needed for a rigorous analytical treatment of wave propagation. The influence of this heterogeneity cannot in general be reduced to one single
effective branching parameter, as would be the case for a homogeneous tree characterized by a constant branching factor $k$. A similar phenomenon occurs in the monostable setting, where periodic heterogeneity changes critical diffusion rates even when the mean degree is preserved; see \cite{hoffman2016invasion}.

Heuristically, the periodic branching setting provides more than a generalization of the homogeneous tree, it offers a richer local description of heterogeneous branching and therefore goes further towards capturing the propagation mechanisms that arise in randomly branching structures. More precisely, considering sufficiently small $d$, where front propagation is strongly influenced by the branching geometry near the interface, this locality allows the homogeneous $k$-ary tree theory to be applied in irregular trees and sparse random networks, see \cite{kouvaris2012traveling}. The periodic setting suggests a refinement of this idea. Instead of representing the local environment by a single branching factor $k$, a finite fundamental block $(k_1,\ldots,k_p)$ can be retained to periodically form a $p$-periodic branching sequence $\mathbf{k}=\langle k_1,\ldots,k_p\rangle$ and a tractable reference problem. This analysis may yield a more detailed local description of propagation and pinning in irregularly branching structures.

A second motivation comes from propagation through isolated branching points, Figure~\ref{fig:trees_examples}(e). Here, an isolated branching point refers to a graph that is locally chain-like away from a single junction: all vertices have degree two, except for one vertex whose degree is at least three. Such a junction can therefore be viewed as a localized geometric heterogeneity embedded in a homogeneous medium. Localized changes in the underlying spatial geometry are known to have a significant influence on front propagation and may, in particular, lead to pinning or propagation failure. For example, in continuous nerve models, changes in axon geometry, and branching in particular, may prevent an excitation wave from propagating further \cite{pauwelussen1981nerve}. For bistable equations on star and tree-like metric graphs, the number and arrangement of incoming and outgoing branches can determine whether a front is propagating or blocked by a standing stationary state, see \cite{jimbo2019entire,jimbo2024front}. In discrete star networks, the degree of the central node, together with the coupling strength, governs transitions between spreading, pinning and retraction, see \cite{kouvaris2017stationary}. More examples that are related to front propagation within spatially heterogeneous discrete bistable media refer to one-dimensional Nagumo chains with spatially heterogeneous coupling or diffusion \cite{morfu2002wave,humphries2011front} and two-dimensional lattices containing geometric obstacles \cite{hoffman2017entire}.

More closely related to an isolated change in connectivity, Yanagita \cite{yanagita2007input} considers a network of discrete FitzHugh--Nagumo excitable elements arranged in three one-dimensional cables joined at a single Y-junction. By direct numerical simulations, the author shows that the pulse may split at the junction and propagate into both outgoing cables, become blocked precisely at the junction, or decay before reaching it. The periodic framework considered here suggests another route towards the isolated-junction problem. Indeed, consider the sparse periodic branching pattern
\[
    \mathbf{k}=\langle k_1,1,\ldots,1 \rangle,
\]
for which branching layers are separated by \(p-1\) homogeneous one-dimensional layers. As \(p\) increases, the whole setting approximates to a structure with one isolated junction connecting one branch to \(k_1>1\) branches. The phase of a pulsating front then acquires a geometric interpretation as its position relative to the branching layer: phases far from the junction correspond to propagation along a homogeneous chain, whereas the phases close to the branching point describe the interaction of the front with the local heterogeneity in connectivity.

For \(k_1>3\), this leads naturally to a class of higher-degree star-like junctions that is structurally different from both metric star graphs considered in \cite{jimbo2019entire,jimbo2024front} and the finite bistable star networks studied in \cite{kouvaris2017stationary}. The large-period regime therefore suggests a possible connection between periodically branching trees and bistable dynamics through isolated high-degree junctions, which is a special heterogeneous setting. This interaction would be an interesting direction for future work.

\paragraph{Comparison principle.}

The comparison principle is a standard tool in the analysis of bistable lattice equations and has also been used to determine propagation and pinning regimes on $k$-ary trees; see, for example, \cite{chen2008traveling,hupkes2023propagation}. For the periodic tree problem \eqref{eq:LDE}, our propagation and pinning results rely on two distinct types of arguments that involve the comparison principle. The first involves the characterization of subsets in sequence space that are invariant under the flow of \eqref{eq:LDE}. The second involves the construction of explicit sub- and super-solutions that either block or force propagation.

The invariant subsets we construct are inspired by the approach developed by Keener \cite{keener1987propagation} and generalize
the invariant intervals analyzed in \cite{hupkes2023propagation} for homogeneous trees. In particular, our sequence subsets  incorporate conditions on at least $L \ge p$ neighboring sites in order to account for the $p$-periodic heterogeneity. Keeping track of the various dependencies, we 
obtain a region $\mathcal P^-(\mathbf k)$ in $(a,d)$ space where $c\leq0$ and a region $\mathcal P^+(\mathbf k)$ where $c\geq0$; their intersection $\mathcal P^0(\mathbf k)$ therefore guarantees $c=0$ and is part of the pinning region. As in  \cite{hupkes2023propagation,keener1987propagation}
the upper boundary of this region is a graph $\big(a, d(a) \big)$. However, the novelty of our construction is the free parameter $L \ge p$ which can be used to significantly enlarge $\mathcal P^0(\mathbf k)$, also for the previously considered homogeneous settings.

Our first explicit sub-solution construction is stationary and tailored specifically to the cubic bistable nonlinearity. As above, it provides
regions $\mathcal D^\pm(\mathbf k)$ in the $(a,d)$ plane where $c\leq0$ respectively $c \ge 0$ holds. However, the intersection $\mathcal D^0(\mathbf k)$ where pinning is guaranteed is significantly larger than $\mathcal P^0(\mathbf k)$ above. A crucial feature is that the boundary is no longer a graph, allowing us to closely track the bends in the actual pinning region.

This construction is based on the exponential behavior of the wave tails established by the linearization at the stable states \cite{mallet1999global}. In the periodic setting, these tails are modulated by periodic factors that are characterized by the theory developed in \cite{chen2008traveling}. Our main contribution here is that we add second-order correction terms to neutralize residuals with the incorrect sign, heavily exploiting the polynomial nature of the cubic;
see \S\ref{sec:construction} and \S\ref{sec:subsolution}.
The resulting sub-solutions are particularly effective for small $d$, where traveling profiles remain steep, in accordance with the classical weak-coupling framework for discrete bistable equations \cite{keener1987propagation}. 
The reflection symmetry of the cubic allows us to limit our computations to the $c \le 0$ case.

Our second sub-solution construction 
involves explicit \emph{moving} profiles, which
provide the \emph{strict} inequality $c < 0$
for sufficiently large $d \gg 1$ under the condition $\widehat{k}=\prod_{i=1}^p k_i >1$, see \S\ref{sec:large-d}. Considering general bistable families again, the comparison principle is applied directly to the mixed functional differential equations
that describe the shapes and speed of traveling waves; see Corollary~\ref{cor:comparison-MFDE}.
The complication is that one needs to work with continuous profiles defined on $\mathbb{R}$ instead of doubly-infinite sequences. Therefore, we work exclusively in the $d \gg 1$ regime and the threshold obtained this way cannot be considered optimal. Our approach is to exploit the directional asymmetry in the graph Laplacian of heterogeneous trees to construct a family of sub-solution profiles for which the Laplacian leaves a nonnegative residual. This can be used to outweigh the potential negative contributions stemming from the nonlinearity. The upside is that we only require the nonlinear term $g$ to be of bistable nature and sufficiently smooth in \S\ref{sec:large-d}.

\paragraph{Weak coupling.}
The small-diffusion regime is closely connected to the classical weak-coupling mechanism for discrete bistable reaction diffusion systems. When the coupling is weak, the reaction term dominates the diffusive interaction at individual lattice sites, and transition fronts become strongly localized near an interface. As a consequence, most lattice sites remain close to the stable states, while only a small number of sites near the interface determine the leading-order balance. This viewpoint has been used extensively to analyze propagation failure in discrete reaction diffusion systems. In particular, Keener \cite{keener1987propagation} uses upper and lower solutions together with perturbation arguments to distinguish the weak- and strong-coupling
regimes. In \cite{erneux1993propagating}, the authors develop an asymptotic construction of traveling fronts in the limit $a\to0^+$ and $d(a)\to0^+$, obtaining an
approximation of the critical coupling for propagation failure. More recent related work
on the discrete Nagumo equation
includes \cite{bustamante2025propagation},
where the propagation-failure region is approximated using asymptotic and modulation methods, 
together with \cite{hupkes2019bichromatic} where asymptotic expansions are used to understand the appearance and disappearance of so-called bichromatic waves that connect 2-periodic equilibria.

We employ this weak-coupling framework as the starting point for our small-diffusion analysis. 
In particular, we study the parameter space near the corner $(a,d) = (1, 0)$ and extract
the relevant scalings of each of the variables
that underpin the construction of the stationary
sub-solutions in \S\ref{sec:subsolution},
which include the periodic exponential tails.
This leads to an explicit description
of the boundary of our pinning region $\mathcal{D}^0(\mathbf k)$ in terms of tangency conditions that mark the disappearance of admissible variables.

Based on these scalings, we expand the boundary curve and its associated parameters in powers of $\sqrt d$; see \S\ref{sec:estimate}. These expansions depend critically on the positioning of the interface where the two exponential tails are tied together relative to the periodicity,
which provides a significant layer of complexity that is absent in earlier work. These expansions are highly accurate for small $d$, but robust tracking at intermediate $d$ requires a numerical continuation of the explicit tangency conditions discussed above.

\paragraph{Organization}
The remainder of the paper is organized as follows. In~\S\ref{sec:main_results}, we introduce the mathematical setting, fix the notation, and formulate the principal results. In~\S\ref{sec:preliminaries}, we explore symmetries implied by the cubic nonlinearity and develop some consequences of comparison principle needed throughout the paper. In~\S\ref{sec:pinned-waves}, we treat general bistable nonlinearities and derive an analytical pinning region $\mathcal P^0$ through the construction of phase-dependent invariant rectangular regions, which adapts the idea from~\cite{hupkes2023propagation} to a general periodic setting. The next three sections focus on the cubic nonlinearity $g_\mathrm{cub}$. Inspired by \cite{chen2008traveling}, we construct the exponential tails used in the stationary sub-solution construction in~\S\ref{sec:construction}. In~\S\ref{sec:subsolution} we show how these tails can be incorporated into complete sub-solutions and verify the existence. The conclusions of sub-solution construction are brought together in~\S\ref{sec:speed-sign}, where the comparison principle converts the existence of stationary sub-solutions into propagation direction $c\leq 0$ on region $\mathcal D^-$. We then apply the symmetry transformation of the system to derive $c\geq 0$. The intersection of these two regions consequently yields the pinning region for the cubic nonlinearity, throughout which $c=0$. In~\S\ref{sec:estimate}, we derive explicit asymptotic boundary curves by means of the approximation strategy for the analytical pinning region in the small-diffusion regime. In~\S\ref{sec:large-d}, we turn to the strong diffusion regime and establish a criterion that guarantees strictly negative wave speed, \(c<0\).  We introduce a finite-core refinement in~\S\ref{sec:numerics-sharp} that allows us to sharpen our pinning regions even further and outline the numerical procedures underpinning our results.

\paragraph{Acknowledgments}
HJH acknowledges support from the Netherlands Organisation for Scientific Research (NWO) through grant 10.61686/NMHBE05119. JY is partially supported by the China Scholarship Council (No. 202206250014).

The authors acknowledge the use of large language models (LLMs) for language editing. Moreover, the proofs of certain technical statements in~\S\ref{sec:estimate} and~\S\ref{sec:large-d} have been obtained with the help of Aristotle and Claude. The authors take full responsibility for the correctness of all statements, proofs, and conclusions.

%%%%%%%%%%%%%%%%%%%%%%%%%%%%%%%%%%%%%%%%%%%%%%%%%%%%%%%%%%%%%%%%%%%
%%%%%%%%%%%%%%%%%%%%%%%%%%%%%%%%%%%%%%%%%%%%%%%%%%%%%%%%%%%%%%%%%%%
% MAIN RESULTS
%%%%%%%%%%%%%%%%%%%%%%%%%%%%%%%%%%%%%%%%%%%%%%%%%%%%%%%%%%%%%%%%%%%
%%%%%%%%%%%%%%%%%%%%%%%%%%%%%%%%%%%%%%%%%%%%%%%%%%%%%%%%%%%%%%%%%%%

\section{Main Results}
\label{sec:main_results}

Let $p \in \mathbb{N}$, $p \ge 2$ be a period, $\mathbf{k} = \langle k_1, k_2, \ldots, k_p\rangle$ be a positive $p$-periodic sequence of branching factors. We consider the bistable reaction diffusion equation on periodic trees~\eqref{eq:LDE} in a slightly more general form by allowing the diffusion to be also $p$-periodic. Besides the scalar diffusion (scaling) parameter $d>0$ we also introduce a $p$-periodic sequence of positive diffusion weights $\boldsymbol{\omega}=\langle \omega_{1/2}, \omega_{3/2}, \ldots, \omega_{p-1/2} \rangle$. The fractional notation emphasizes the fact that all the edges between two layers share the same diffusion. For example, the edges between the first and second layer (and also $np+1$-st and $np+2$-nd layer) have diffusion weight $\omega_{3/2}$. Consequently, we rewrite~\eqref{eq:LDE} as
\begin{equation} \label{e:main}
    \dot{u}_{np+j}(t) = d \Delta_{j} [\mathbf{u}(t)]_{n p + j}  + g(u_{np+j}(t);a), \qquad t>0, \quad n \in \mathbb{Z}, \quad \dInt{j},
\end{equation}
where the advection-diffusion operator $\Delta_{j}$ is given by
\[ 
    \Delta_{j}[\mathbf{u}(t)]_{i } := \omega_{j+1/2} k_j (u_{i+1}(t) - u_i(t))
    + \omega_{j-1/2} \big( u_{i-1}(t) - u_i(t) \big)
     ,\qquad i=np+j.     
\]

Besides the general form of results below, the periodicity of diffusion weights $\boldsymbol{\omega}$ allows us to consider only waves with $c\le0$ and obtain the results for $c\ge0$ using a duality argument, see \S\ref{sec:duality}. Throughout the entire paper we fix the period $p\in\mathbb{N}$ and the positive and $p$-periodic sequences $\mathbf{k}, \boldsymbol{\omega}$.

We require the nonlinearity $g$ to satisfy the following standard bistability assumption. 
\begin{enumerate}[label=(Hg), ref=Hg]
    \item \label{Hg} The map $(u, a) \mapsto g(u; a)$ is $C^2$-smooth on $\mathbb R \times [0, 1]$. For every $a\in(0,1)$, we have
$$g(0;a) = g(a;a) = g(1;a) = 0,$$
$$g'(0;a) <0,\quad g'(1;a) < 0, \quad g'(a;a) > 0.$$ 
In addition, the function $g$ satisfies the inequalities
\[
    g(v;a)>0\text{ for }v\in(-\infty,0)\cup(a,1),\qquad g(v;a)<0\text{ for }v\in(0,a)\cup(1,\infty). 
\]
\end{enumerate}

We are interested in traveling wave solutions $\mathbf{u} = \mathbf{u}^{\rm tw}$ of \eqref{e:main}. Such solutions consist of $p$ periodically repeating profiles \begin{equation}\label{eq:wave-ansatz}
    \boldsymbol{\Phi}(\xi) = \langle\Phi_j(\xi)\rangle, \qquad
    u^{\rm tw}_{n p + j}(t)= \Phi_j(n p + j - c t),
\end{equation}
where $c\in\mathbb R$ is the (common) wave speed. The corresponding profiles satisfy the mixed functional differential equation (MFDE)
\begin{align} \label{e:main-mfde}
    - c \Phi'_j(\xi) = d \, \mathfrak{D}_{\boldsymbol{\omega},\mathbf{k}}[\boldsymbol{\Phi}(\, \cdot \,)]_j(\xi) + g(\Phi_j(\xi);a), \quad \xi \in \mathbb{R} \quad \dInt{j}, 
\end{align}
with $2p$ boundary conditions
\begin{align} \label{e:main-mfde-bc}
    \Phi_j(-\infty) = 0, \qquad \Phi_j(+\infty) = 1, \qquad \dInt{j},
\end{align}
where the diffusion operator $\mathfrak{D}_{\boldsymbol{\omega},\mathbf{k}}$ is defined by
\begin{align} \label{e:lapl-cont}
    \mathfrak{D}_{\boldsymbol{\omega},\mathbf{k}} [\boldsymbol{\Phi}(\, \cdot \,)]_j(\xi) := \omega_{j+1/2} k_j \Phi_{j+1}(\xi+1) - (\omega_{j+1/2} k_j+ \omega_{j-1/2}) \Phi_j(\xi) + \omega_{j-1/2} \Phi_{j-1}(\xi-1). 
\end{align}
We build our theory on the following existence result, which we prove in {\S}\ref{sec:existence}.

\begin{proposition}[Existence of pulsating traveling waves]
\label{prop:periodic-wave-existence}
 Let $(a,d)\in\mathcal{H}:=(0,1)\times(0,\infty)$. Then there exist a speed $c=c(a,d;\boldsymbol{\omega},\mathbf k)\in\mathbb R$ and $p$ non-decreasing profiles
\[
    \boldsymbol{\Phi}(\xi)=\langle\Phi_j(\xi)\rangle,
\]
that satisfy the MFDE \eqref{e:main-mfde} with the boundary conditions \eqref{e:main-mfde-bc}. Moreover, the wave speed $c$ is uniquely determined. If $c\neq0$, then each $\Phi_j$ is $C^1$, strictly increasing, and unique up to translation.
\end{proposition}

For bistable families that admit the \emph{reflection identity} (e.g., the cubic \eqref{eq:cubic})
\begin{equation}\label{eq:mr:reflection}
    g(1-u;a)=-g(u;1-a),
\end{equation}
there is a useful duality principle that can be exploited.
Given two strictly positive $p$-periodic sequences $\boldsymbol{\omega}=\langle \omega_{1/2}, \omega_{3/2}, \ldots, \omega_{p-1/2} \rangle$ and $\mathbf k=\langle k_1, k_2, \ldots, k_p\rangle$, 
this principle revolves around the transformed positive \(p\)-periodic sequences
\[
    \widetilde{\boldsymbol{\omega}}=(\widetilde{\omega}_{j-1/2})_{j\in\mathbb Z} =\langle \widetilde{\omega}_{1/2},\ldots,\widetilde{\omega}_{p-1/2} \rangle,
    \qquad
    \widetilde{\mathbf k}=(\widetilde{k}_j)_{j\in\mathbb Z}=\langle \widetilde{k}_1,\ldots,\widetilde{k}_p \rangle,
\]
defined through the identities
\begin{equation}
\label{eq:mr:refl:dual:transf}
    \widetilde \omega_{j-1/2}:=\omega_{-j+1/2}k_{-j}=\omega_{p-j+1/2}k_{p-j},
    \qquad
    \widetilde k_j:=\frac1{k_{-j-1}} = \frac1{k_{p-j-1}}.
\end{equation}
Under this transformation the wave speed reverses precisely, providing the periodic analogue of \cite[Lemma~2.2]{hupkes2023propagation}, which we prove in {\S}\ref{sec:duality}. 
\begin{lemma}[Duality]
\label{lem:mr:reflection}
Assume that the reflection identity \eqref{eq:mr:reflection} holds.  Then for any $a \in (0,1)$ and $d>0$ we have
\begin{equation}
\label{eq:mr:refl:c}
    c(a,d; \boldsymbol{\omega}, \mathbf{k} ) = - c \big( 1 - a,d;\widetilde{\boldsymbol{\omega}}, \widetilde{\mathbf k} \big) .
\end{equation}
\end{lemma}

\subsection{Pinned waves via trapped sequences}
\label{sec:pinned-waves-main}

\begin{figure}[t]
    \centering
    \includegraphics[width=0.9\linewidth]{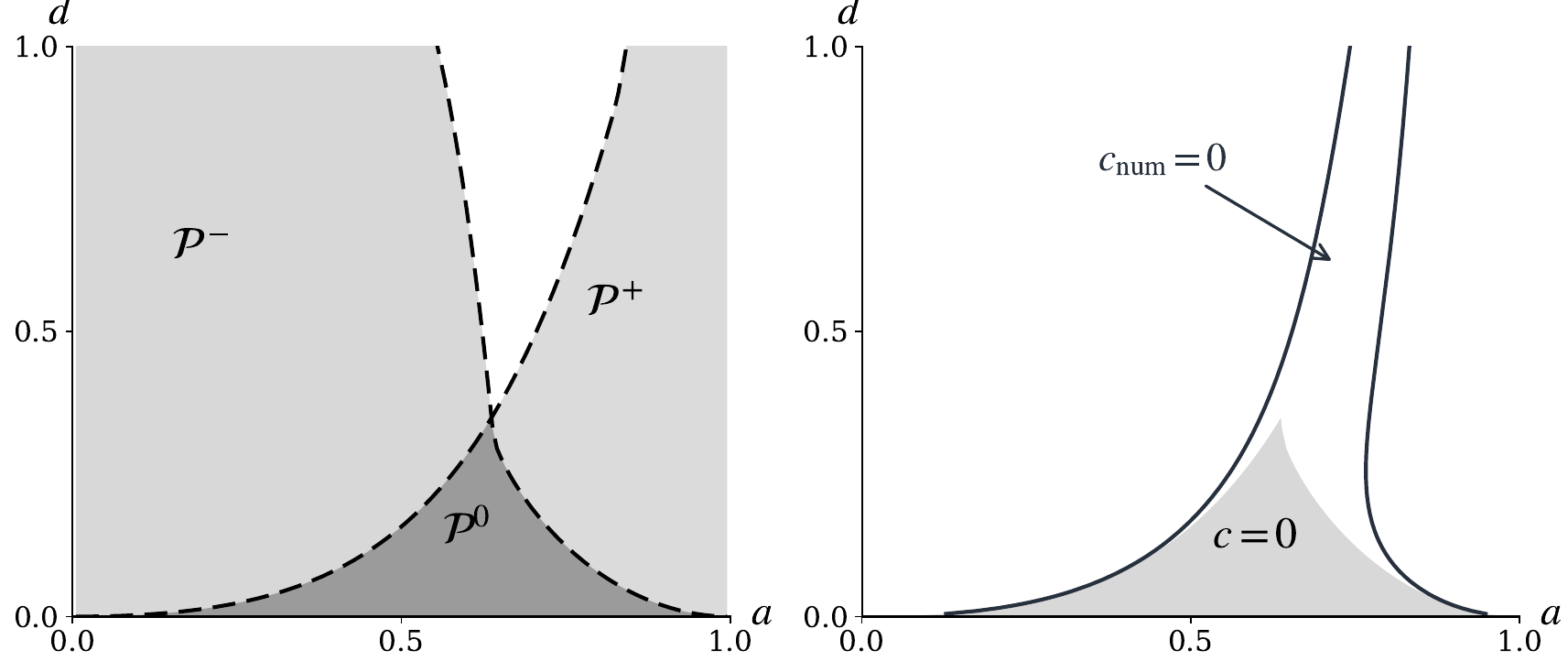}
    \caption{Illustration of Theorem~\ref{thm:periodic-pinning} 
    with $L=p=4$
    for the bistable family $g_\mathrm{exp}$ defined in \eqref{eq:mr:def:g:exp}, with \(\mathbf k=\langle1,2,2,2\rangle\) and $\boldsymbol{\omega}=\langle 1\rangle$. Left: the  regions \(\mathcal P^\pm(\boldsymbol{\omega},\mathbf k)\) and the theoretical pinning region \(\mathcal P^0(\boldsymbol{\omega},\mathbf k)\).
    Right: comparison of the latter with the numerically computed propagation boundaries.}
    \label{fig:general}
\end{figure}

Our first result identifies a nonempty region in the \((a,d)\)-plane in which the periodic traveling wave is pinned, that is, $c(a,d;\boldsymbol{\omega},\mathbf{k})=0$; see Figure~\ref{fig:general}. We achieve this by showing that there exist two regions with nonempty intersection, one with $c\le 0$ and the other with $c\ge 0$. 
Their intersection corresponds to $c=0$. 

 We fix an integer \(L\ge p\). For a position $s\in\mathbb Z$ and a vector \(\mathbf x =(x_1,\ldots,x_L)\in(a,1)^L\) we introduce a \emph{non-periodic} extension by writing $x_0:=0$ and $x_{L+1}:=x_{L+1-p}$, which allows us to define the quantities
\begin{equation}
\label{eq:B-}
    B_{j,s}^-(\mathbf x;\boldsymbol{\omega},\mathbf{k}):=k_{s+j}\omega_{s+j+1/2}
    (x_j - x_{j+1})
    +\omega_{s+j-1/2}(x_j -x_{j-1}),
    \qquad 1\le j\le L.
\end{equation}
In addition, we define the sequence space 
\begin{equation}
\label{eq:mr:seq:trap}
    \mathcal{S}_s^-(\mathbf x)
    = \{ \mathbf{u}: \mathbb Z \to [0,1] : \mathbf{u} \hbox{ is p-monotone and } x_j\le u_{j+s} \le 1 \hbox{ for all } 1\le j\le L\}.
\end{equation}
The key idea follows from the fact that if $x_{j_*} = u_{j_*+s}$ in~\eqref{eq:mr:seq:trap} for  $j_* \in \{1, \ldots, L \}$ , then by construction we have
\begin{equation}
    \dot{u}_{j_*+s} \ge g(u_{j_*+s} ; a) - d B^-_{j_*,s}(\mathbf x; \boldsymbol{\omega},\mathbf{k}) = g(x_{j_*} ; a)  - d B^-_{j_*,s}(\mathbf x; \boldsymbol{\omega},\mathbf{k}).
\end{equation}
Since \eqref{Hg} implies $g(x_{j_*} ; a) > 0$, we can hence guarantee the \emph{trapping condition} $\dot{u}_{j_* + s} > 0$ by restricting the size of $d$ if $B^-_{j_*,s}(\mathbf x; \boldsymbol{\omega},\mathbf{k}) > 0$. In particular, by requiring that
\begin{equation}
    0 < d <  \frac{g(x_j ; a)}{B^-_{j,s}(\mathbf x ;\boldsymbol{\omega},\mathbf{k})} \qquad \qquad \hbox{ for all } j \hbox{ with } B^-_{j,s}(\mathbf x; \boldsymbol{\omega},\mathbf{k}) > 0,
\end{equation}
we can show that $\mathcal{S}_s^-(\mathbf x)$ is invariant under the flow of \eqref{e:main} and forces $c \le 0$.

Utilizing the freedom to choose the phase \(\dInt{s}\) and the vector \(\mathbf x\) in an optimal fashion, we now define the  diffusion threshold 
\begin{equation}
\label{eq:d-minus-periodic}
    d^-(a; \boldsymbol{\omega},\mathbf k)
    :=
    \max_{\dInt{s}}
    \sup_{\mathbf x\in(a,1)^L}
    \min_{\{j:B_{j,s}^-(\mathbf x;\boldsymbol{\omega},\mathbf k)>0\}}
     \dfrac{g(x_j;a)}
    {B_{j,s}^-(\mathbf x; \boldsymbol{\omega},\mathbf{k})},
\end{equation}
noting that the minimum over an empty set is understood to be $\infty$. In particular, \eqref{Hg} directly implies $d^-(a; \boldsymbol{\omega},\mathbf k) > 0$.

Similarly, we extend any \(\mathbf y=(y_1,\ldots,y_{L})\in(0,a)^L\) non-periodically by writing
$y_0:=y_p$, $y_{L+1}:=1$, 
which allows us to define the quantities
\begin{equation}
\label{eq:B+}
    B_{j,s}^+(\mathbf y; \boldsymbol{\omega},\mathbf{k}):=
k_{s+j}\omega_{s+j+1/2}
    (y_j - y_{j+1})
    +\omega_{s+j-1/2}(y_j -y_{j-1}),
    \qquad 1\le j\le L,
\end{equation}
together with the associated diffusion threshold 
\begin{equation}
\label{eq:d-plus-periodic}
    d^+(a; \boldsymbol{\omega},\mathbf{k})
    :=
    \max_{\dInt{s}}
    \sup_{\mathbf y\in(0,a)^L}
    \min_{\{j:B_{j,s}^+(\mathbf y;\boldsymbol{\omega},\mathbf k)<0\}}
    \dfrac{g(y_j;a)}
    {B_{j,s}^+(\mathbf y; \boldsymbol{\omega},\mathbf{k})} > 0,
\end{equation}
again writing $\infty$ for the minimum over an empty set. The corresponding sequence space is now given by
\begin{equation}
\label{eq:mr:seq:trap-y}
    \mathcal{S}_s^+(\mathbf y)
    = \{ \mathbf{u}: \mathbb Z \to [0,1] : \mathbf{u} \hbox{ is p-monotone and } 0\le u_{j+s} \le y_j \hbox{ for all } 1\le j\le L\}.
\end{equation}

We then introduce two nonempty regions
\begin{equation}
    \mathcal P^-(\boldsymbol{\omega},\mathbf{k}):=\left\{(a,d)\in\mathcal H:0<d<d^-(a;\boldsymbol{\omega},\mathbf{k})
    \right\},
\end{equation}
\begin{equation}
    \mathcal P^+(\boldsymbol{\omega},\mathbf{k}):=\left\{(a,d)\in\mathcal H:0<d<d^+(a;\boldsymbol{\omega},\mathbf{k})
    \right\},
\end{equation}
and define their intersection by
\begin{equation}
\label{eq:pinning-region-periodic}
    \mathcal P^0(\boldsymbol{\omega},\mathbf{k})
    :=
    \left\{
        (a,d)\in\mathcal H:
        0<d<
        \min\left\{
            d^-(a;\boldsymbol{\omega},\mathbf{k}),
            d^+(a;\boldsymbol{\omega},\mathbf{k})
        \right\}
    \right\} = \mathcal{P}^-(\boldsymbol{\omega},\mathbf{k}) \cap \mathcal{P}^+(\boldsymbol{\omega},\mathbf{k}).
\end{equation}

\begin{theorem}[{Propagation regions, general $g$, {\S}\ref{sec:pinned-waves}}]
\label{thm:periodic-pinning}

Let $(a,d)\in\mathcal H$, then the following statements hold.
\begin{enumerate}[label=(\roman*), ref=(\roman*)]
    \item \label{thm:periodic-pinning-1} We have $\mathcal P^-(\boldsymbol{\omega},\mathbf{k})\neq \varnothing$,  $\mathcal P^+(\boldsymbol{\omega},\mathbf{k})\neq \varnothing$ and $\mathcal P^0(\boldsymbol{\omega},\mathbf k)\neq \varnothing$. 
    \item \label{thm:periodic-pinning-2} If $(a,d)\in\mathcal P^-(\boldsymbol{\omega},\mathbf{k})$, we have $c(a,d;\boldsymbol{\omega},\mathbf{k})\le0$.
    \item \label{thm:periodic-pinning-3} If $(a,d)\in\mathcal P^+(\boldsymbol{\omega},\mathbf{k})$, we have $c(a,d;\boldsymbol{\omega},\mathbf{k})\ge0$.
    \item \label{thm:periodic-pinning-4} If $(a,d)\in\mathcal P^0(\boldsymbol{\omega},\mathbf{k})$, we have $c(a,d;\boldsymbol{\omega},\mathbf{k})=0$.
\end{enumerate}
\end{theorem}

\begin{figure}[t]
    \centering
    \includegraphics[width=0.9\linewidth]{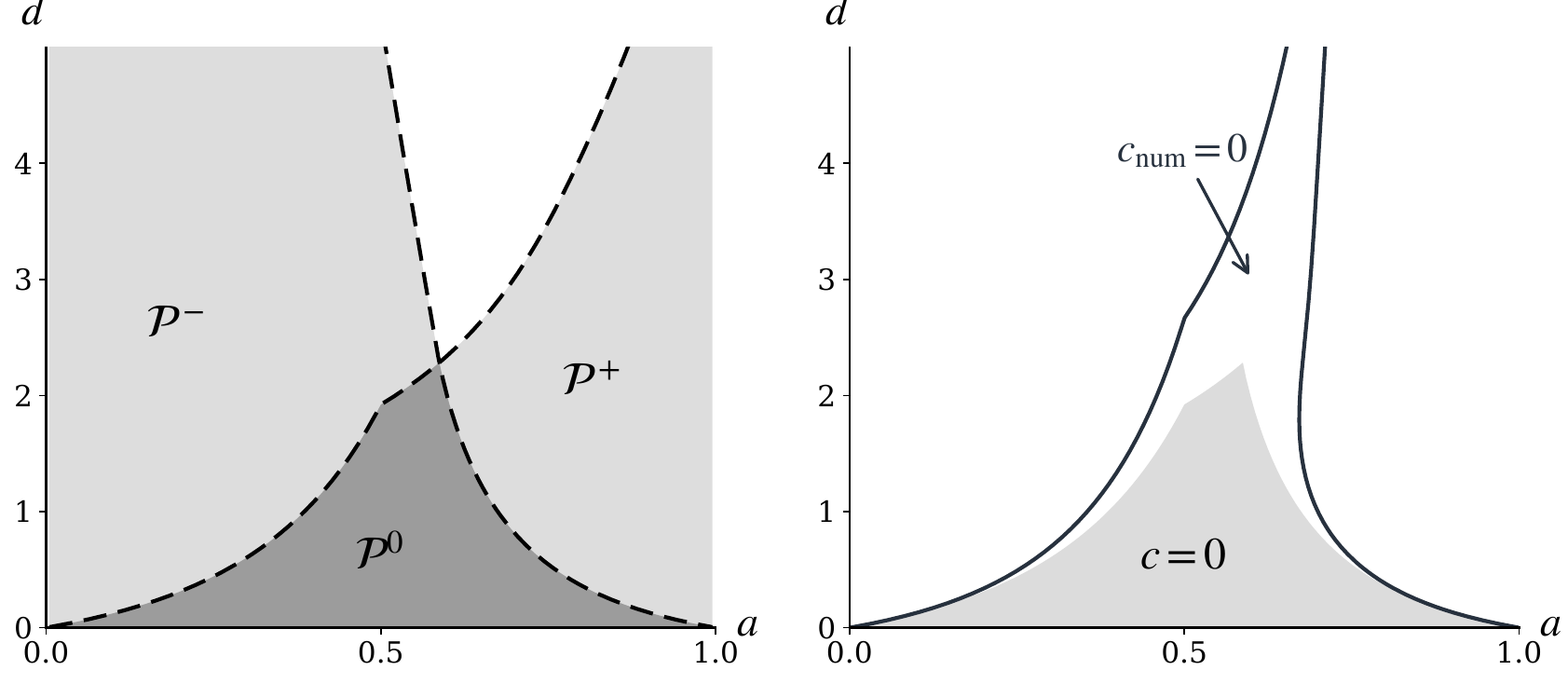}
    \caption{Illustration of Theorem~\ref{thm:periodic-pinning}
    with $L =4$ and $p = 2$
    for the bistable family $g_\mathrm{sin}(u;a)$ defined in~\eqref{eq:sine:bistability}, 
    with  \(\mathbf k=\langle1,2\rangle\) and $\boldsymbol{\omega}=\langle 1 \rangle$.
    Left: the  regions \(\mathcal P^\pm(\boldsymbol{\omega},\mathbf k)\) and the theoretical pinning region \(\mathcal P^0(\boldsymbol{\omega},\mathbf k)\).
    Right: comparison of the latter with the numerically computed propagation boundaries.
    }
    \label{fig:general:2}
\end{figure}

Figure~\ref{fig:general} illustrates Theorem~\ref{thm:periodic-pinning} for the bistable family 
\begin{equation}
\label{eq:mr:def:g:exp}
    g_\mathrm{exp}(u;a)=e^{2u}u(1-u)(u-a),
\end{equation} while Figure~\ref{fig:general:2} provides theoretical and numerical propagation regions for the bistable family 

\begin{equation}\label{eq:sine:bistability}
    {g}_{\sin}(u;a)=\begin{cases}\displaystyle
    -\min\left\lbrace\dfrac{a}{1-a},1\right\rbrace \cdot \sin\left(\dfrac{\pi u}{a}\right) &\text{ for } u\in[0,a],\\
    -\min\left\lbrace 1,\dfrac{1-a}{a}\right\rbrace \cdot \sin\left(\dfrac{\pi (u-1)}{1-a}\right) &\text{ for } u\in(a,1].
    \end{cases}
\end{equation}

\begin{remark}
Restricting the threshold vectors in \eqref{eq:d-minus-periodic}
and \eqref{eq:d-plus-periodic}
to the constant choices
\[
    \mathbf x=(x,\ldots,x),
    \qquad
    \mathbf y=(y,\ldots,y),
    \qquad x\in(a,1),\quad y\in(0,a),
\]
and writing $\omega_{\min}:=\min_{s\dInt{}}\omega_{s+1/2}$, $(k\omega)_{\min}:=\min_{s\dInt{}}k_s\omega_{s+1/2}$ leads directly to the lower bounds
\begin{equation}
\label{eq:mr:seq:d:pm:c}
\begin{aligned}
    d^-(a;\boldsymbol{\omega},\mathbf{k})& \ge & d_{\mathrm c}^-(a;\boldsymbol{\omega},\mathbf{k})
    &: =
    \frac1{\omega_{\min}}\sup_{x\in(a,1)}\frac{g(x;a)}{x},
    \\
    d^+(a;\boldsymbol{\omega},\mathbf{k})& \ge &
    d_{\mathrm c}^+(a;\boldsymbol{\omega},\mathbf{k})
    &: =
    \frac1{(k\omega)_{\min}}
    \sup_{y\in(0,a)}\frac{-g(y;a)}{1-y}.
\end{aligned}
\end{equation}
In the homogeneous case $\boldsymbol{\omega} = \langle1\rangle$ and \(\mathbf k = \langle k \rangle\) the lower bounds reduce
to the expressions in \cite[Eq.~(2.10)]{hupkes2023propagation}, while for $k=1$ the bounds \cite[Theorem~2.8]{keener1987propagation} are recovered. 
\end{remark}

\begin{figure}[t]
    \centering
    \subfigure[$\mathbf{k}=\langle 1 \rangle$]{
    \includegraphics[width=0.31\linewidth]{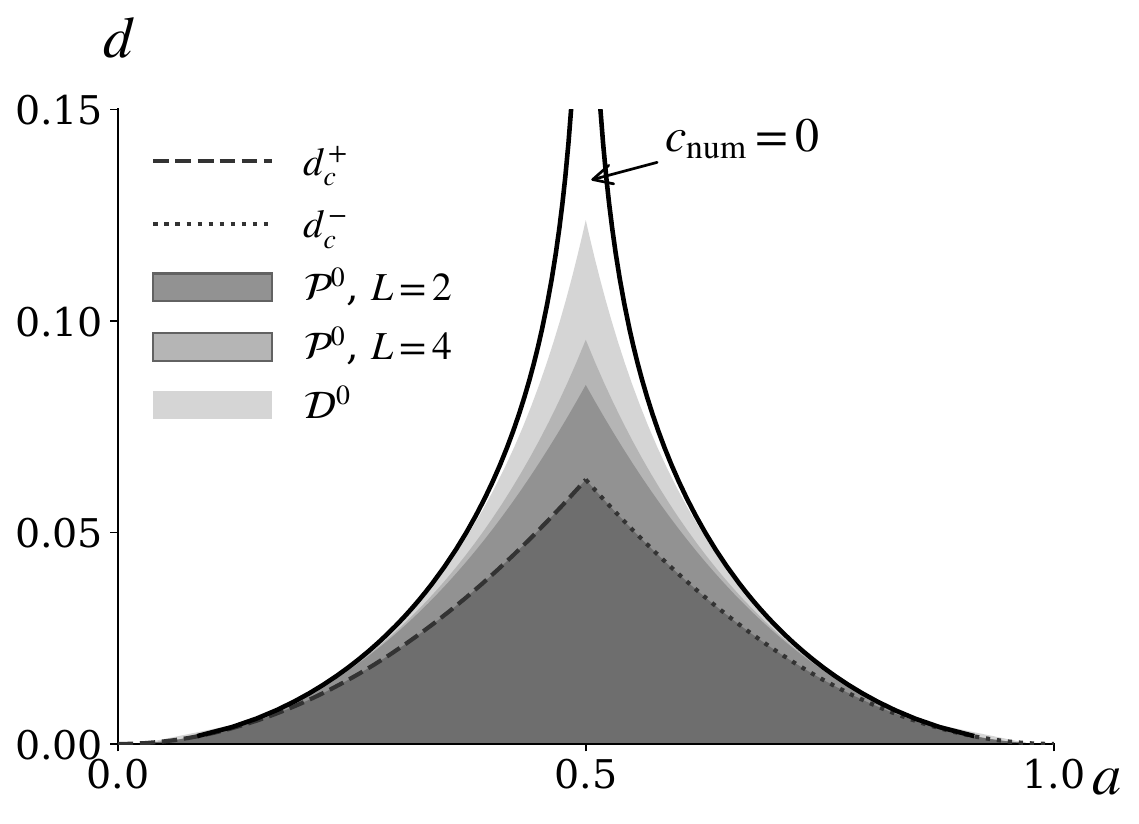}
    }
    \subfigure[$\mathbf{k}=\langle 1,2 \rangle$]{
    \includegraphics[width=0.31\linewidth]{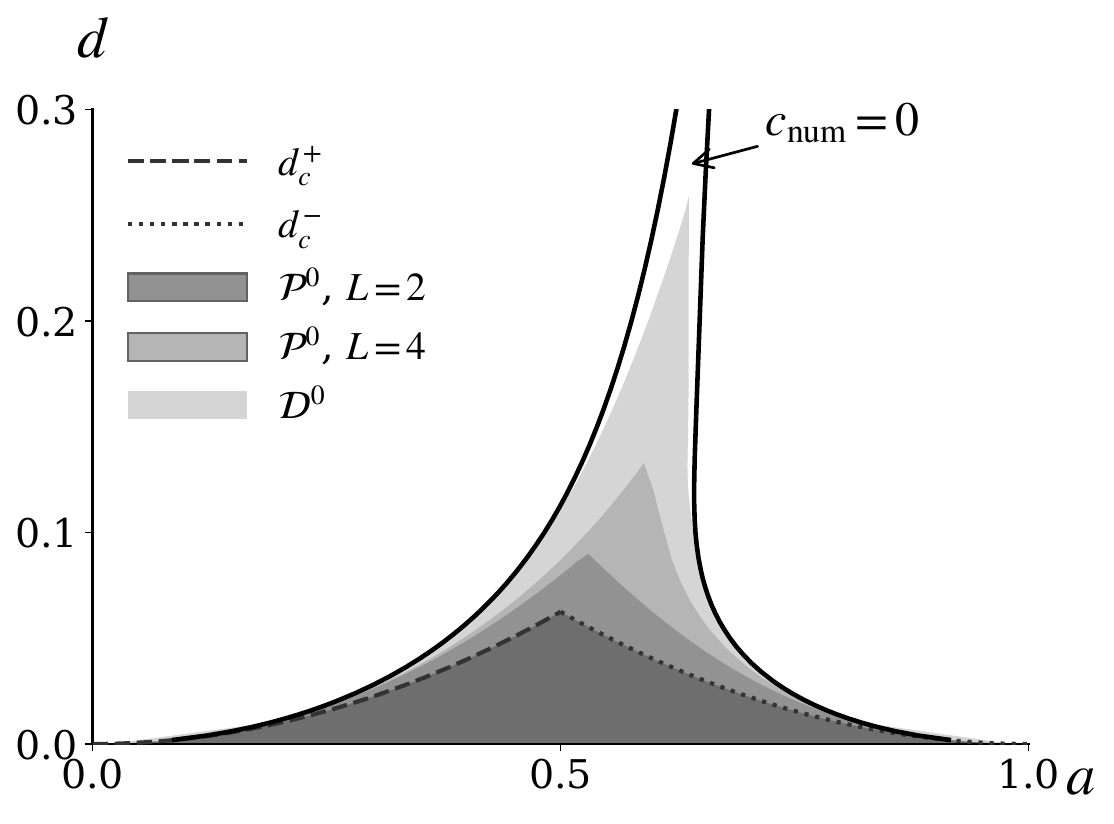}
    }
    \subfigure[$\mathbf{k}=\langle 3 \rangle$]{
    \includegraphics[width=0.31\linewidth]{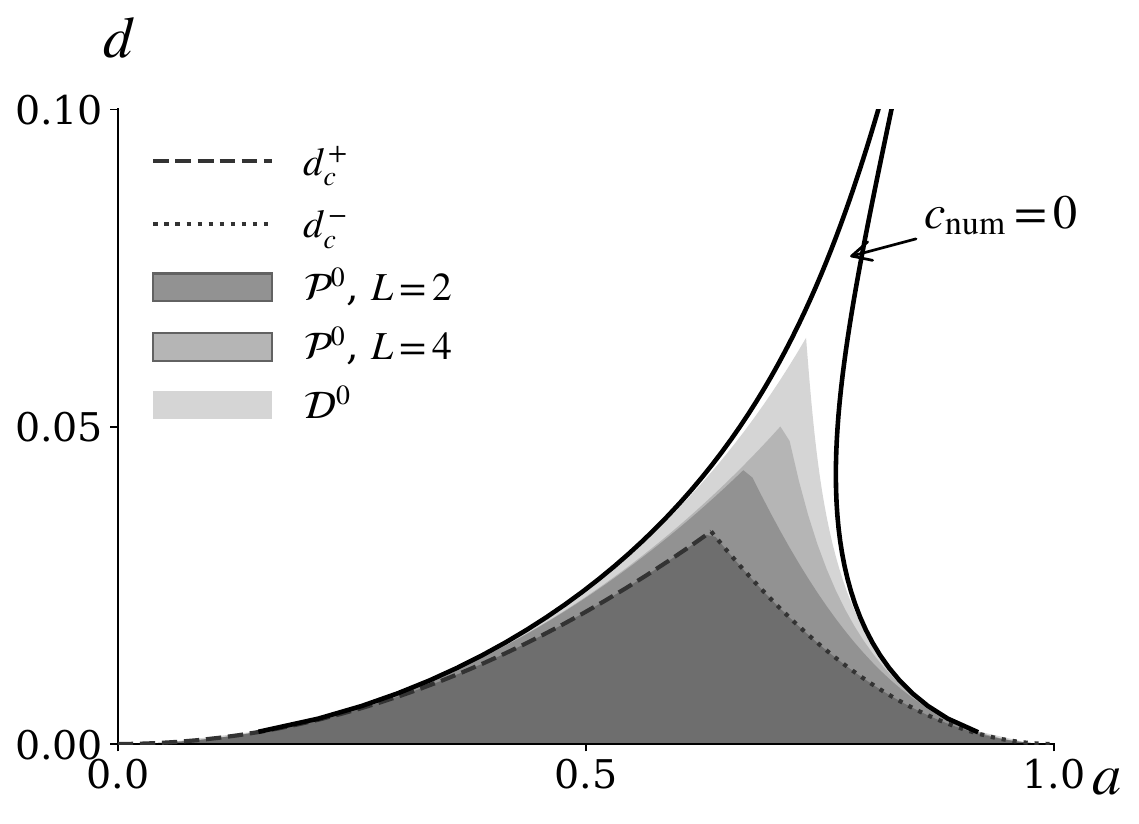}   
    }
    \caption{Comparison of various theoretical pinning regions %$P^0(\boldsymbol{\omega},\mathbf{k})$ 
    for the cubic nonlinearity~\eqref{eq:cubic},
    with %the periodic medium with \(\mathbf k=\langle1,2\rangle\) and 
    $\boldsymbol{\omega}=\langle 1 \rangle$. The darkest gray region represents the constant-threshold bounds $d^\pm_c$ defined in \eqref{eq:mr:seq:d:pm:c}. Two medium gray regions 
    depict $\mathcal P^0(\boldsymbol{\omega},\mathbf{k})$
    with the choices \(L=p=2\) and \(L=2p=4\), respectively. The light gray region is $\mathcal D^0(\boldsymbol{\omega},\mathbf{k})$ from Theorem \ref{thm:pinning}\ref{thm:pinning-3}.\label{fig:pinning-cp}}    
\end{figure}

\begin{remark}
The window width \(L \ge p\) provides a further parameter that can be used to trade off accuracy against complexity. Figure~\ref{fig:pinning-cp} illustrates the improvements that can be achieved by passing from homogeneous to heterogeneous $\mathbf x$ and $\mathbf y$ and subsequently increasing $L$ beyond the minimal value $p$. 
\end{remark}

\subsection{Pinning region via stationary barriers for cubic nonlinearity}
\label{sec:pinning-region-for-cubic}
\begin{figure}[t]
    \centering
    \includegraphics[width=0.9\linewidth]{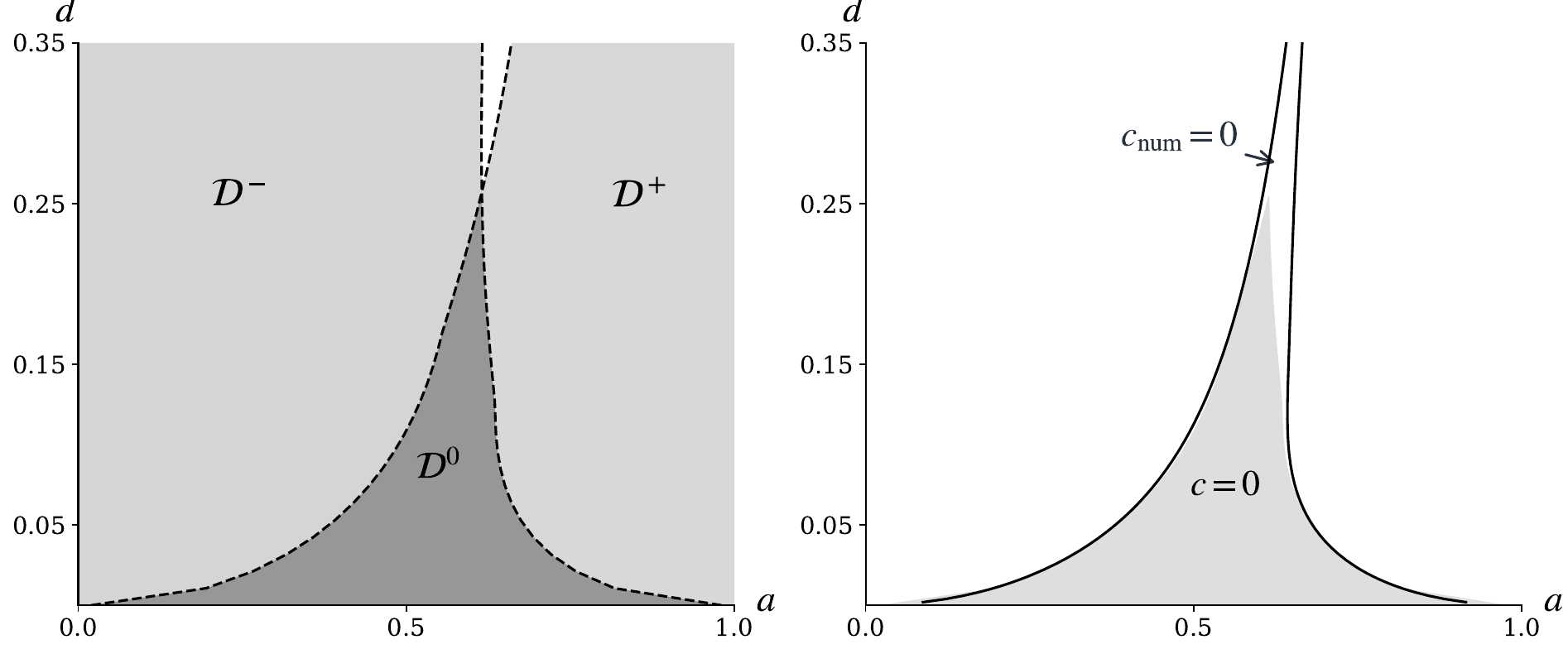}
    \caption{Illustration of Theorem~\ref{thm:pinning}
    for the cubic nonlinearity,
    with \(\mathbf{k}=\langle 1,2\rangle \) and $\boldsymbol{\omega}=\langle 1 \rangle$.
    Left: the regions \(\mathcal D^\pm(\boldsymbol{\omega},\mathbf{k})\) and the theoretical pinning region \(\mathcal D^0(\boldsymbol{\omega},\mathbf{k})\). Right: comparison of the latter 
    with the numerically computed propagation boundaries.}
    \label{fig:pin-region-k12}
\end{figure}

\begin{figure}[t]
    \centering
    \subfigure[$\mathbf k=\langle1,2\rangle$]{
    \includegraphics[width=0.31\linewidth]{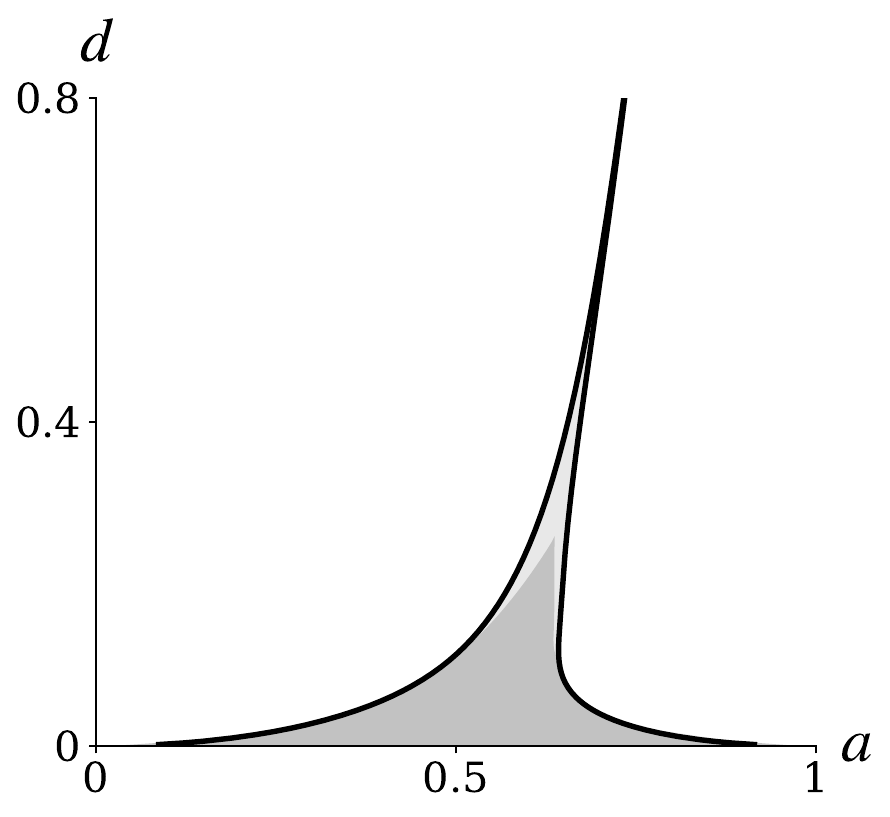}
    }\subfigure[$\mathbf k=\langle1,2,2,2\rangle$]{
    \includegraphics[width=0.31\linewidth]{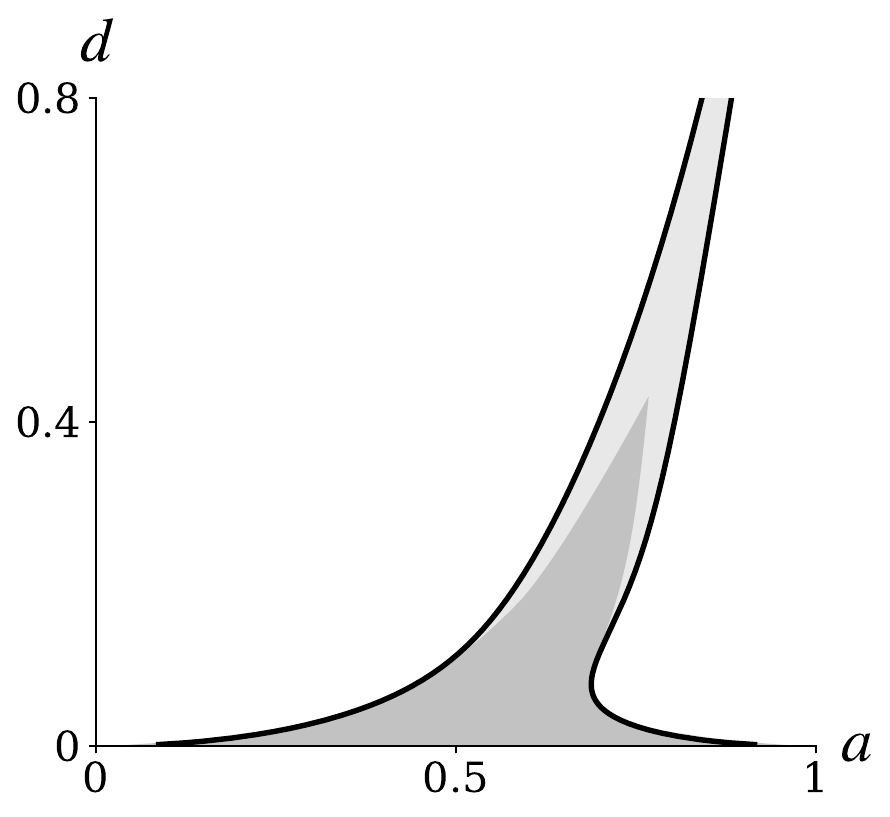}
    }\subfigure[$\mathbf k=\langle1,2,2,2,2,2\rangle$]{
    \includegraphics[width=0.31\linewidth]{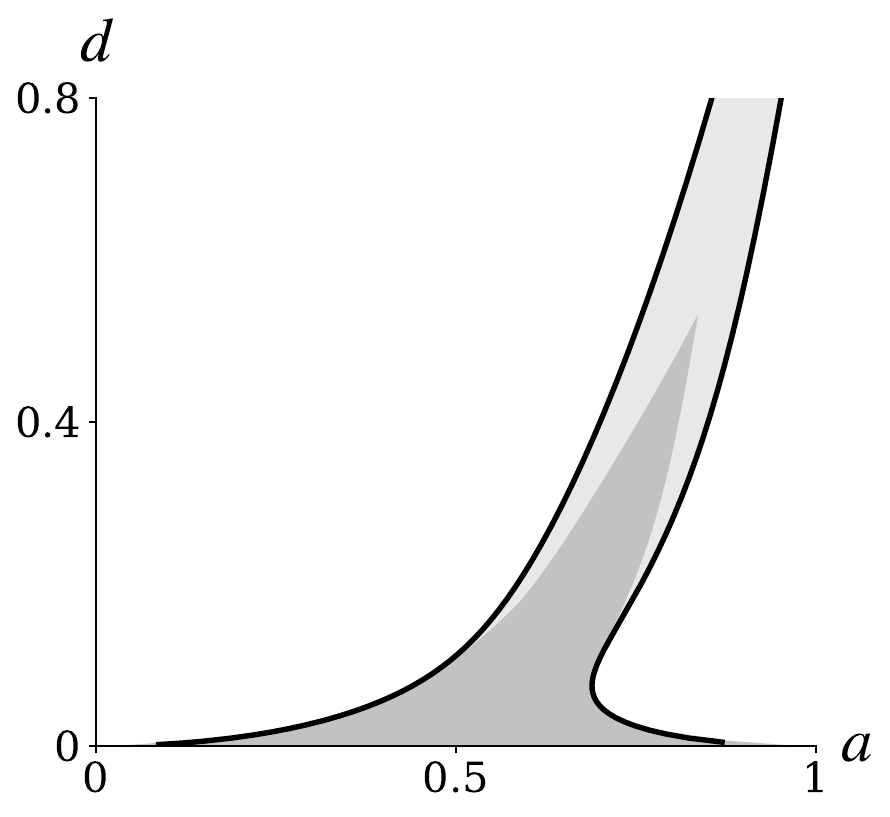}
    }
    \caption{Comparison of several pinning regions for the cubic nonlinearity \eqref{eq:cubic}, with diffusion strength $\boldsymbol{\omega}=\langle 1 \rangle$, and three periodic branching patterns: (a) \(\mathbf k=\langle1,2\rangle\), (b) \(\mathbf k=\langle1,2,2,2\rangle\), and (c) \(\mathbf k=\langle1,2,2,2,2,2\rangle\). The black curves represent the numerically computed boundaries of the actual pinning region. The dark gray regions show the theoretical pinning regions $\mathcal D^0(\boldsymbol{\omega},\mathbf k)$ obtained from the exponential narrow-interface sub-solution construction \eqref{eq:mr:def:u:i:s} and Theorem~\ref{thm:pinning}. The light gray regions show the theoretical pinning regions $\mathcal D_{\rm fc}^0(\boldsymbol{\omega},\mathbf k; \mathcal{N})$ detected by the finite-core refinement \eqref{eq:fc:def:fc:0}, optimized over the core lengths \(\mathcal N=\{0,1,2,10,20\}\).}
    \label{fig:comparison_region}
\end{figure}

In this section we establish a sharper sufficient pinning region for the cubic bistable nonlinearity~\eqref{eq:cubic}; see Figure~\ref{fig:pin-region-k12}. As above, this is achieved by identifying parametric regions in $(a,d)$ where \(c\leq0\) and \(c\geq0\), respectively. 
The key difference to \S\ref{sec:pinned-waves-main} is that - unlike \eqref{eq:pinning-region-periodic} - the boundaries of our regions are not described by functions of $a$. This allows us to follow the fold of the pinning boundary (and thus capture the propagation reversal phenomenon as well). Indeed, our numerical results  indicate that the boundaries of our regions follow the corresponding edges of the numerically observed pinning region rather closely over a broad range of parameters; see Figures~\ref{fig:pin-region-k12} and~\ref{fig:comparison_region}.

Our strategy is to fix a phase $\dInt{s}$ and construct a stationary sub-solution for 
\eqref{e:main} by connecting two stationary one-sided sub-solution tails, defined on $i 
\le s$ and $i 
\ge s+1$ respectively. 
This construction then yields matching equalities at the interface points $\{s, s+1\}$ that lead naturally to conditions on  the system parameters $(a, d,\boldsymbol{\omega},\mathbf k)$ that guarantee $c \le 0$. We then use the duality from Lemma~\ref{lem:mr:reflection} to obtain results for $c \ge 0$. 
To be more concrete, we  pick two strictly positive parameters 
$(v,w)$ and introduce the sequence $\big(\underline u_i^{(s)}(v,w)\big)_{i \in \mathbb Z}$ by writing
\begin{equation}
\label{eq:mr:def:u:i:s}
    \underline u_i^{(s)}(v,w)=
    \begin{cases}
        v\mathcal L_{i;s}, & i\le s,\\[1mm]
        1-w\mathcal R_{i;s+1}
        +C_iw^2\mathcal R_{i;s+1}^{\,2},
        & i\ge s+1 .
    \end{cases}    
\end{equation}
Here $(\mathcal L_{i;s})_{i \le s+1}$ is 
the unique strictly positive solution to the linearization of \eqref{e:main} around the spatially homogeneous equilibrium $u =0$
that decays exponentially as $i \to - \infty$ and admits the normalization $\mathcal L_{s;s} =1 $.
Defined analogously for the linearization of \eqref{e:main} around $u =1$, the sequence $\mathcal R_{i,s+1}$ decays exponentially as $i \to + \infty$ and is normalized to have $\mathcal R_{s+1;s+1} = 1$. The positive $p$-periodic coefficients $\mathbf C=(C_i)_{i\in\mathbb Z}=\langle C_1,\dots,C_p\rangle$ provide a second order correction to the right tail to ensure that $\underline{\mathbf{u}}^{(s)}$ is a sub-solution for the full nonlinear system in the regime $i \ge s+2$. The coefficients $C_i$ depend on the system parameters $(a,d,\boldsymbol{\omega},\mathbf k)$ and are constructed in \S\ref{sec:construction}.

The normalization conditions $\mathcal L_{s;s} =1 $ and $\mathcal R_{s+1;s+1} = 1$ yield
\[
    \underline u_s^{(s)}=v,
    \qquad
    \underline u_{s+1}^{(s)}
    =1-w+C_{s+1}w^2.
\]
Thus $v \in (0,1)$ directly provides the value at the last left-tail site, whereas $w \in (0, w_s^*)$ controls the amplitude of the corrected right tail. Here $w_s^* > 1$ is an explicit constant that guarantees $\underline{u}^{(s)} > 0$.

To ensure that the sequence $\underline{\mathbf{u}}^{(s)}(v,w)$ from \eqref{eq:mr:def:u:i:s} is a global stationary sub-solution that can be used to show that $c \le 0$, there are several conditions on $(v,w)$ that need to be imposed. We perform this analysis in \S\ref{sec:subsolution},
where we construct explicit polynomials $\theta_s$, $h_s$ and $V_s$ of order $2$, $3$ and $6$ so that every pair $(v,w)$ in 
the set
\[
\mathcal A_s^-(a,d;\boldsymbol{\omega},\mathbf k)
= \{ (v, w) \in (0,1) \times (0, w_s^*): h_s(v) \le \theta_s(w) \hbox{ and } V_s(w) \le v \}
\]
generates a suitable sub-solution.
Since this set can be empty in general, we introduce the corresponding feasible parameter region 
\begin{equation}
\label{eq:Ds-minus}
\mathcal D_s^-(\boldsymbol{\omega},\mathbf k):=\{(a,d)\in\mathcal H:\mathcal A_s^-(a, d; \boldsymbol{\omega}, \mathbf k)\neq\varnothing\}.
\end{equation}
Thus a parameter pair \((a,d)\) belongs to \(\mathcal D_s^-(\boldsymbol{\omega},\mathbf k)\) precisely when there exists at least one choice of the interface amplitudes \((v,w)\) for which the phase-\(s\) construction yields a stationary sub-solution,
which then implies $c(a,d) \le 0$.

At this point, we exploit the remaining freedom to choose the  phase $\dInt{s}$ where the matching is performed in an optimal fashion. In particular, we write
\begin{equation}
\label{eq:D-minus}
\mathcal D^-(\boldsymbol{\omega},\mathbf k):=\bigcup_{\dInt{s}}\mathcal D_s^-(\boldsymbol{\omega},\mathbf k),
\end{equation}
for the union of the regions where our construction establishes $c(a, d) \le 0$.

Turning to the opposite case $c \ge 0$, we apply the transformation \eqref{eq:mr:refl:dual:transf}  and the speed reversal \eqref{eq:mr:refl:c} to define the sets
\begin{equation}
\label{eq:Ds-plus}
\begin{array}{lcl}
\mathcal D_s^+(\boldsymbol{\omega},\mathbf k) &:= &\{(a,d)\in\mathcal H:\mathcal A_s^-(1-a,d;\widetilde{\boldsymbol{\omega}},\widetilde{\mathbf k})\neq\varnothing\},
\\[0.2cm]
\mathcal{D}^+(\boldsymbol{\omega},\mathbf k) & :=  & \bigcup_{\dInt{s}}\mathcal D_s^+(\boldsymbol{\omega},\mathbf k),
\end{array}
\end{equation}
for which $c \ge 0$.  Accordingly, the intersection $\mathcal{D}^+(\boldsymbol{\omega},\mathbf k) \cap \mathcal{D}^-(\boldsymbol{\omega},\mathbf k)$ - which we show to be non-empty -
determines a parameter region in which the waves are pinned. 

\begin{theorem}[{Propagation regions, cubic $g$, {\S}\ref{sec:construction},\ref{sec:subsolution},\ref{sec:speed-sign}}]
\label{thm:pinning}
For the cubic nonlinearity~\eqref{eq:cubic}, the following statements hold.

\begin{enumerate}[label=(\roman*), ref=(\roman*)]
\item \label{thm:pinning-1} $\mathcal D^-(\boldsymbol{\omega},\mathbf k)\neq\varnothing$ and $\mathcal D^+(\boldsymbol{\omega},\mathbf k)\neq\varnothing$.
\item \label{thm:pinning-2} If $(a,d)\in\mathcal D^-(\boldsymbol{\omega},\mathbf k)$, we have $c(a,d;\boldsymbol{\omega},\mathbf k)\leq0$. And if $(a,d)\in\mathcal D^+(\boldsymbol{\omega},\mathbf k)$, we have $c(a,d;\boldsymbol{\omega},\mathbf k)\geq0$.
\item \label{thm:pinning-3} $\mathcal D^0(\boldsymbol{\omega},\mathbf k):=\mathcal D^-(\boldsymbol{\omega},\mathbf k)\cap\mathcal D^+(\boldsymbol{\omega},\mathbf k)\neq\varnothing$, and $c(a,d;\boldsymbol{\omega},\mathbf k)=0$ for $(a,d)\in\mathcal D^0(\boldsymbol{\omega},\mathbf k)$.
\item \label{thm:pinning-4} For any $(a,d) \in \overline{\mathcal D^-(\boldsymbol{\omega},\mathbf k)}\cap \mathcal H$, we have $c(a',d) \le 0$ for all $0<a'<a$. Similarly, for any $(a,d) \in \overline{\mathcal D^+(\boldsymbol{\omega},\mathbf k)}\cap \mathcal H$, we have $c(a',d)\ge0$ for all $a<a'<1$.

\end{enumerate}
\end{theorem}

By Theorem~\ref{thm:pinning}\ref{thm:pinning-4}, at fixed $d$ we have $c\le0$ to the left of the right boundary of $\mathcal D^-(\boldsymbol{\omega},\mathbf k)$ and $c\ge0$ to the right of the left boundary of $\mathcal D^+(\boldsymbol{\omega},\mathbf k)$. The overlap of these two regions therefore guarantees pinning.
We pursue this idea in {\S}\ref{sec:estimate}, where we  construct explicit curves $(d,a^\pm(d))$ that track the boundary of the regions $\mathcal{D}^\pm$ for $0 < d \ll 1$. These curves are obtained by expanding the free variables $(v,w)$ in~\eqref{eq:mr:def:u:i:s} and the detuning parameter $a$ simultaneously in powers of $\sqrt{d}$. These expansions can be truncated at any desired order and yield high precision results up to intermediate values of $d$. 

We emphasize that in the homogeneous regime $
\boldsymbol{\omega} = \langle 1 \rangle$ and $\mathbf k = \langle k \rangle$ our results actually improve the 
prior thresholds $d^\pm_c$ defined in \eqref{eq:mr:seq:d:pm:c} which were obtained in \cite{hupkes2023propagation,keener1987propagation}. This can already be seen at the level of the expressions themselves: \(d^-_c(a)\) does not respond to changes in the branching factor \(k\),
whereas the region \(\mathcal D^0(\boldsymbol\omega,\mathbf k)\) obtained from our stationary exponential-tail construction incorporates the branching coefficients through the
periodic tails \(\mathcal L\) and \(\mathcal R\), the correction coefficients \(\mathbf C\), and the two interface residuals. The resulting improvements for the regions where pinning can be guaranteed are illustrated in Figure~\ref{fig:pinning-cp}.

The estimate of the pinning region \(\mathcal D^0(\boldsymbol{\omega},\mathbf k)\) obtained above is based on the two-parameter family \eqref{eq:mr:def:u:i:s}. It is therefore natural to ask whether additional free parameters could further improve the size of the region where we can guarantee pinning. We explore this in \S\ref{sec:numerics-sharp}, where we enlarge the matching region from the two adjacent sites $\{s, s+1\}$ to $N+2$ sites $\{s, \dots,s+N+1\}$ by introducing free profile values at $N$ sites $\{s+1,\dots,s+N\}$; see \eqref{eq:fc-profile}. This additional flexibility  indeed leads to an improved approximation of the pinning boundary; see Figure~\ref{fig:comparison_region}.

\subsection{Large diffusion propagation criteria}
Returning to general bistable families $g$, we  impose the condition $\widehat{k}>1$ (see \eqref{eq:khat}) on the LDE~\eqref{e:main} and show that there exists a critical diffusion threshold $d^*(a;\boldsymbol{\omega},\mathbf{k})$ 
so that $c(a, d; \boldsymbol{\omega},\mathbf k) < 0$ for all $d>d^*(a;\boldsymbol{\omega},\mathbf{k})$. This threshold can be bounded  uniformly for $a \in (0,1)$. In particular, for sufficiently large $d$ the waves of Proposition \ref{prop:periodic-wave-existence} all travel to the left.

The core idea behind our construction is that we can find explicit superharmonic functions  for the graph Laplacian 
$\mathfrak{D}_{\boldsymbol{\omega},\mathbf{k}}$.
In particular, we introduce the set
\begin{equation} \label{e:omega}
\Omega_\mathbf{k} := \left\{ \mathbf{m} \in \mathbb{R}^p \, | \, 0 < \mathbf{m} < \mathbf{k}, \; \widehat{m}>1 \right\}
\end{equation}
which is non-empty on account of the asymmetry of the diffusion operator encoded by the condition $\widehat{k}>1$. For any $\mathbf{m} \in \Omega_{\mathbf k}$ we then
consider the function $\boldsymbol\chi_\mathbf{m}(\xi;\boldsymbol\omega) = \langle \chi_{j,\mathbf{m}}(\xi;\boldsymbol\omega) \rangle_\dInt{j}$ given by
\begin{equation}
    \chi_{j,\mathbf{m}}(\xi;\boldsymbol{\omega}) = 
    \alpha_{j,\mathbf{m}}(\boldsymbol\omega) \beta_{j,\mathbf{m}}(\xi),
\end{equation}
where $\alpha_{j,\mathbf{m}}(\boldsymbol{\omega}) > 0$ and $\beta_{j,\mathbf{m}}(\xi) > 0 $ are given by
\begin{align} \label{e:alpha-intro}
        \alpha_{j,\mathbf{m}}(\boldsymbol\omega) &:=\frac{1}{p} \left(\sum_{n=1}^{p} \frac{\widehat{\omega}}{\omega_{j+n-1/2}} \prod_{i=1}^{p-n} m_{j-i}  \right), \\
        \beta_{j,\mathbf{m}}(\xi) &:= \left( \prod_{i=1}^{p} \sqrt[p]{m_{j-i}^{i-p}} \right) \sqrt[p]{\widehat{{m}}}^{-\xi}. \label{e:beta_intro}
\end{align}
The exponential nature of $\boldsymbol\chi(\, \cdot \,;\boldsymbol\omega, \mathbf m)$ then leads to the computation
\[
-\mu_{j,\mathbf{m}}(\boldsymbol{\omega},\mathbf{k}) \mathfrak{D}_{\boldsymbol{\omega},\mathbf{k}}[\boldsymbol\chi_\mathbf{m}(\, \cdot \,; \boldsymbol\omega)]_j(\xi) = \frac{\widehat{\omega}}{p}
\chi_{j,\mathbf{m}}(\xi; \boldsymbol{\omega}),
\]
where the expression $\mu_{j,\mathbf{m}}(\boldsymbol{\omega},\mathbf{k})$ is defined by 
\begin{equation}
\mu_{j,\mathbf{m}}(\boldsymbol{\omega},\mathbf{k}):=  \frac{m_j \alpha_{j,\mathbf{m}}(\boldsymbol\omega)}{(\widehat{{m}}-1)(k_j-m_j)} > 0 \label{e:mu_j-intro}
\end{equation}
using \eqref{e:omega} for the strict positivity.

For $d \gg 1$, this sign condition can be used to ensure
that $1-\boldsymbol\chi_\mathbf{m}$
is a strict sub-solution for the MFDE~\eqref{e:main-mfde} in the sense 
\[
-d \mathfrak{D}_{\boldsymbol{\omega},\mathbf{k}}[1-\boldsymbol{\chi}_\mathbf{m}(\, \cdot \,;\boldsymbol{\omega})]_j(\xi) - g\big(1-\chi_{j,\mathbf{m}}(\xi;\boldsymbol{\omega}); a\big) < 0, \qquad \qquad \xi \in \mathbb{R}, \quad \dInt{j}. 
\]
Rearranging, using the linearity of $\mathcal{D}_{\boldsymbol\omega,\mathbf{k}}$ together with $\mathcal{D}_{\boldsymbol\omega,\mathbf{k}}[1] = 0$, we see that it suffices for $d$ to satisfy
\begin{equation} \label{e:d-star-prelim}
d > \frac{p}{\widehat{\omega}} \mu_{j,\mathbf{m}}(\boldsymbol{\omega},\mathbf{k}) \left( - \frac{g\big(1-\chi_{j,\mathbf{m}}(\xi;\boldsymbol{\omega}) ; a \big) } {\chi_{j,\mathbf{m}}(\xi; \boldsymbol{\omega}) } \right), \qquad \qquad \xi \in \mathbb{R}, \quad \dInt{j}. 
\end{equation}
We now take the liberty of choosing $\mathbf{m} \in \Omega_\mathbf{k}$ in order for the threshold $d^*(a;\boldsymbol\omega,\mathbf{k})$ to be as low as possible. To this end, we define
\begin{align}\label{e:mu}  
        \mu(\boldsymbol{\omega},\mathbf{k}) &:=  \min_{\mathbf{m} \in \Omega_\mathbf{k}} \max_{\dInt{j}} \mu_{j,\mathbf{m}}(\boldsymbol{\omega},\mathbf{k}) , \\
        \label{e:nu}
        \nu(a) &:= \max\limits_{s \in [1-a,1]} -\frac{g(1-s;a)}{s},
 \end{align}
 where $\mu(\boldsymbol{\omega},\mathbf{k})$ minimizes the first factor in~\eqref{e:d-star-prelim} and $\nu(a)$ handles the worst-case scenario for the second factor in~\eqref{e:d-star-prelim}. Note that each $\mu_{j,\mathbf{m}}(\boldsymbol\omega,\mathbf{k})$ is continuous for $\mathbf{m} \in \mathrm{int} \, \Omega_\mathbf{k}$ and $\max\limits_\dInt{j} \mu_{j,\mathbf{m}}(\boldsymbol\omega,\mathbf{k})$ diverges to $+\infty$ for $\mathbf{m} \to \partial\Omega_\mathbf{k}$. The function $\mu(\boldsymbol\omega,\mathbf{k})$ is thus well-defined. The sign conditions in $\eqref{Hg}$ ensure that $\nu(a)$ is strictly positive, while the smoothness condition provides a uniform bound for this quantity on $a \in (0,1)$.
 
We emphasize that the construction above cannot be used directly, since $\boldsymbol{\chi}$ is unbounded and our comparison principles require bounded sub-solutions. Nevertheless, the result below shows that the resulting lower bound for $d$ suffices to obtain left-traveling waves for the full problem.

\begin{theorem}[Large diffusion and $c<0$, \S\ref{sec:large-d}] \label{t:left-travel}
    Let $\mathbf{k}$ be such that $\widehat{k} > 1$. Then for any 
    \begin{equation}\label{e:dStar}
        d > d^*(a;\boldsymbol{\omega},\mathbf{k}) := 
        \frac{p}{\widehat\omega} \mu(\boldsymbol\omega,\mathbf{k}) \nu(a),
    \end{equation}
    we have $c(a,d;\boldsymbol{\omega},\mathbf{k}) < 0$.    
\end{theorem}

\begin{figure}[t]
    \centering
    \subfigure[cubic bistability $g_\mathrm{cub}$, \eqref{eq:cubic}]{
    \includegraphics[width=0.45\linewidth]{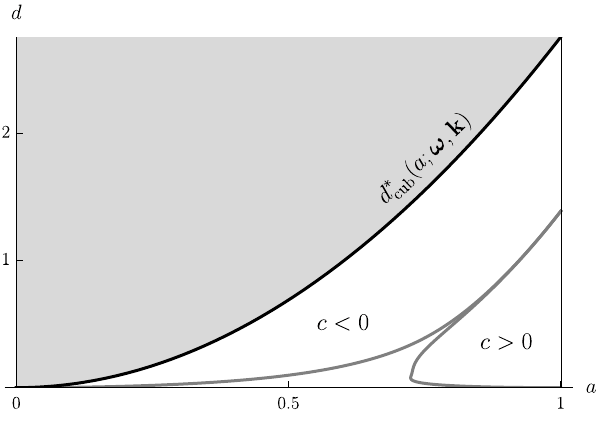}
    \label{fig:d-star-cub}
    }\subfigure[exponential bistability $g_\mathrm{exp}$, \eqref{eq:mr:def:g:exp}]{
    \includegraphics[width=0.45\linewidth]{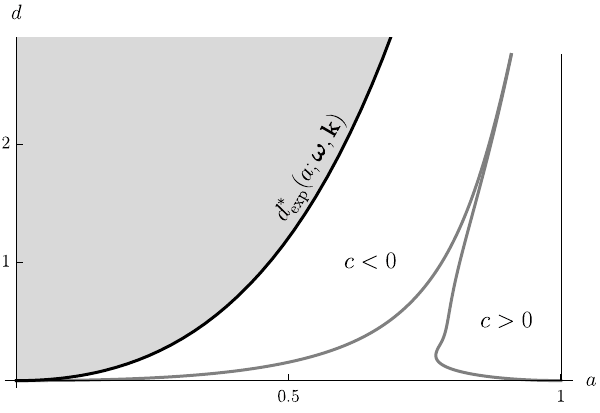}
    \label{fig:d-star-exp}
    }
    \caption{Illustration of Theorem~\ref{t:left-travel} for the cubic (left) and exponential (right) bistable families. We compare the region $d > d^*(a;\boldsymbol\omega,\mathbf{k})$ where $c < 0$ is guaranteed with its numerical counterpart, considering the diffusion strength $\boldsymbol\omega =\langle 1 \rangle$ and the branching pattern $\mathbf{k}=\langle 1,3 \rangle$.
    } 
    \label{fig:d-star}
\end{figure}

Figure~\ref{fig:d-star} illustrates Theorem~\ref{t:left-travel} for the cubic~\eqref{eq:cubic} and exponential~\eqref{eq:mr:def:g:exp} bistable families. We provide the details and the rigorous proof of Theorem~\ref{t:left-travel} in \S\ref{sec:large-d}. Moreover, in \S\ref{sec:large-d-analysis} we analyze the properties of the $p$-dimensional constrained optimization problem for $\mu(\boldsymbol\omega,\mathbf{k})$ in~\eqref{e:mu}. We design a strategy to obtain estimates for the optima and outline a numerical treatment for the problem. We discuss the properties of $\nu(a)$ there too which mostly follow from \eqref{Hg}. 

%%%%%%%%%%%%%%%%%%%%%%%%%%%%%%%%%%%%%%%%%%%%%%%%%%%%%%%%%%%%%%%%%%%
%%%%%%%%%%%%%%%%%%%%%%%%%%%%%%%%%%%%%%%%%%%%%%%%%%%%%%%%%%%%%%%%%%%
% PRELIMINARIES
%%%%%%%%%%%%%%%%%%%%%%%%%%%%%%%%%%%%%%%%%%%%%%%%%%%%%%%%%%%%%%%%%%%
%%%%%%%%%%%%%%%%%%%%%%%%%%%%%%%%%%%%%%%%%%%%%%%%%%%%%%%%%%%%%%%%%%%

\section{Preliminaries}\label{sec:preliminaries}

In this section we state and prove several useful preliminary results. We start in {\S}\ref{sec:existence} by proving the existence of wave profiles and their speed. In {\S}\ref{sec:duality}, we then show how to exploit the structure generated by nonlinearities with reflection symmetry to construct a parameter transformation that reverses the wave speed. In {\S}\ref{sec:comparison_principle} we formulate a comparison principle for our main LDE \eqref{e:main}, which we then exploit to derive sign conditions on the wave speed.

\subsection{Existence of pulsating waves}
\label{sec:existence}
Our main existence result follows from minor modifications of those from~\cite{chen2008traveling}.
\begin{proof}[Proof of Proposition~\ref{prop:periodic-wave-existence}]
Existence follows by adapting the argument in \cite[proof of Theorem~7]{chen2008traveling}. Although our coupling is not necessarily of the special form considered there, the spatial-ordering property used in their construction remains valid for arbitrary positive nearest-neighbor coefficients and a spatially independent reaction term. Indeed, write
\[
    A_i:=d k_i\omega_{i+1/2}>0,
    \qquad B_i:=d\omega_{i-1/2}>0,
    \qquad z_i(t):=u_{i+1}(t)-u_i(t).
\]
Subtracting the equations at adjacent sites gives
\[
    \dot z_i=A_{i+1}z_{i+1}+B_i z_{i-1}
    +\bigl(q_i(t)-A_i-B_{i+1}\bigr)z_i,
    \qquad
    q_i(t):=\int_0^1\partial_u g\bigl(u_i(t)+\theta z_i(t);a\bigr)\,d\theta.
\]
This system is cooperative, so spatially increasing initial data remain spatially increasing. The same difference calculation shows that the finite-interval integral operators used in \cite[\S8.2]{chen2008traveling} preserve the ordering $\Phi_i(\xi)\le\Phi_{i+1}(\xi+1)$ when their auxiliary constant is chosen sufficiently large. Their construction therefore yields non-decreasing periodic profiles with spatially constant limits. Assumption~\eqref{Hg} excludes the intermediate limiting state by the argument in \cite[proof of Theorem~7]{chen2008traveling}, giving the connection from $0$ to $1$.

Since $0$ and $1$ are stable equilibria, \cite[Theorem~3]{chen2008traveling} further implies that the wave speed is uniquely determined and, when $c\neq0$, the traveling wave is unique up to a time translation.
\end{proof}

\subsection{Reflection symmetry}\label{sec:duality}
Our goal here is to establish Lemma \ref{lem:mr:reflection}. To this end, we recall the transformed parameters $\widetilde{\boldsymbol{\omega}},\widetilde{\mathbf{k}}$  from  \eqref{eq:mr:refl:dual:transf} and consider the associated dual problem
\begin{equation}\label{eq:prlm:dual}
    \begin{cases}
    \dot{\widetilde{u}}_{np+j}(t) = d \widetilde{\Delta}_{j} [\widetilde{u}(t)]_{n p + j}  + g(\widetilde{u}_{np+j}(t);1-a), \qquad t>0, \quad n \in \mathbb{Z}, \quad \dInt{j},\\

    \widetilde{\Delta}_{j}[\mathbf{u}(t)]_{i } := \widetilde{\omega}_{j+1/2} \widetilde{k}_j (\widetilde{u}_{i+1}(t) - \widetilde{u}_i(t))
    + \widetilde{\omega}_{j-1/2} \big( \widetilde{u}_{i-1}(t) - \widetilde{u}_i(t) \big)
     ,\qquad i=np+j.
     \end{cases}
\end{equation}
We define a transformed wave by writing
\begin{equation}\label{eq:prlm:transf:wave}
    \widetilde{\Phi}_j(\xi)=1-{\Phi}_{-j}(-\xi), \qquad \qquad \widetilde{c} = -c.
\end{equation}
By construction, we note that $\widetilde{\boldsymbol{\Phi}}$ satisfies the boundary conditions \eqref{e:main-mfde-bc} if and only if they are satisfied by $\boldsymbol{\Phi}$. Recalling the MFDE~\eqref{e:main-mfde}, the profile $\boldsymbol{\Phi}$ satisfies
\begin{equation}
    \label{eq:mr:wave:mfde:new}
\begin{aligned}
-c\Phi_j'(\xi)
&=d\omega_{j+1/2}k_j\bigl(\Phi_{j+1}(\xi+1)-\Phi_j(\xi)\bigr)\\
&\quad+d\omega_{j-1/2}\bigl(\Phi_{j-1}(\xi-1)-\Phi_j(\xi)\bigr)
+g(\Phi_j(\xi);a).
\end{aligned}
\end{equation}

\begin{lemma}
\label{lem:prlm:symmetry:dual}
Assume that the bistable family $g$ satisfies the reflection symmetry~\eqref{eq:mr:reflection} and pick $(a,d)\in\mathcal{H}$. If $u_{n p + j}(t)= \Phi_j(n p + j - c t)$ is a solution to the LDE~\eqref{e:main}, then $\widetilde{u}_{n p + j}(t) =\widetilde{\Phi}_{j}(np + j - \widetilde{c} t) $ solves the dual problem~\eqref{eq:prlm:dual}.
\end{lemma}

\begin{proof}
In order to exploit the reflection symmetry,
we apply the transformations $j\mapsto-j$ and $\xi\mapsto-\xi$ 
to \eqref{eq:mr:wave:mfde:new} and multiply the equation by $-1$
to find
\[
\begin{aligned}
c\Phi_{-j}'(-\xi)
&=-d\omega_{-j+1/2}k_{-j}
\bigl(\Phi_{1-j}(1-\xi)-\Phi_{-j}(-\xi)\bigr)\\
&\quad-d\omega_{-j-1/2}
\bigl(\Phi_{-j-1}(-\xi-1)-\Phi_{-j}(-\xi)\bigr)
-g(\Phi_{-j}(-\xi);a).
\end{aligned}
\]
    In particular, the transformed wave \eqref{eq:prlm:transf:wave} 
    satisfies
\[
\begin{aligned}
-\widetilde c\,\widetilde\Phi_j'(\xi)
&=d\omega_{-j+1/2}k_{-j}
\bigl(\widetilde\Phi_{j-1}(\xi-1)-\widetilde\Phi_j(\xi)\bigr)\\
&\quad+d\omega_{-j-1/2}
\bigl(\widetilde\Phi_{j+1}(\xi+1)-\widetilde\Phi_j(\xi)\bigr)
-g(1-\widetilde\Phi_j(\xi);a).
\end{aligned}
\]
    Exploiting the reflection symmetry~\eqref{eq:mr:reflection} and recalling the transformation \eqref{eq:mr:refl:dual:transf}, we obtain
\[
\begin{aligned}
-\widetilde c\,\widetilde\Phi_j'(\xi)
&=d\widetilde\omega_{j+1/2}\widetilde k_j
\bigl(\widetilde\Phi_{j+1}(\xi+1)-\widetilde\Phi_j(\xi)\bigr)\\
&\quad+d\widetilde\omega_{j-1/2}
\bigl(\widetilde\Phi_{j-1}(\xi-1)-\widetilde\Phi_j(\xi)\bigr)
+g(\widetilde\Phi_j(\xi);1-a).
\end{aligned}
\]
    In particular, the pair $(\widetilde{c}, \widetilde{\boldsymbol{\Phi}})$
    satisfies the traveling wave MFDE~\eqref{eq:mr:wave:mfde:new} associated
    to the dual problem \eqref{eq:prlm:dual},
    which is hence satisfied by $\tilde{u}$.
\end{proof}

\begin{proof}[Proof of Lemma \ref{lem:mr:reflection}]
Consider the traveling wave $(c, \boldsymbol{\Phi})$ 
given by Proposition \ref{prop:periodic-wave-existence}
together with its associated pair $(\widetilde{c}, \widetilde{\boldsymbol{\Phi}})$ defined by \eqref{eq:prlm:transf:wave}.  
Lemma \ref{lem:prlm:symmetry:dual} shows that this pair is a traveling wave solution to the dual problem 
\eqref{eq:prlm:dual}. By construction, it also satisfies the boundary 
conditions \eqref{e:main-mfde-bc}.
The uniqueness of the wave speed
guaranteed by Proposition \ref{prop:periodic-wave-existence}
for the dual system now
implies that 
\[
    c(1-a, d; \widetilde{\boldsymbol{\omega}}, \widetilde{\mathbf k} ) = \widetilde{c} = - c(a, d; \boldsymbol{\omega},\mathbf k ),
\]
i.e., the equality~\eqref{eq:mr:refl:c}, as desired.
\end{proof}

\subsection{Comparison principle}\label{sec:comparison_principle}
The main tool used in this paper is the comparison principle for the periodic lattice dynamical system~\eqref{e:main}. We state it first in a full-lattice form adapted to the periodic tree equation. We then derive a corollary of the comparison principle showing how the existence of a stationary sub-solution can be used to obtain the non-positivity of the wave speed $c(a,d;\boldsymbol{\omega},\mathbf k).$

For a fixed sequence $\mathbf{u}=(u_i)_{i \in \mathbb{Z}}$ we define the stationary residual
\begin{equation}
\label{eq:residual}
F_i(\mathbf{u}):=d\omega_{i+1/2}k_i\left(u_{i+1}-u_i\right) + d\omega_{i-1/2} \left(u_{i-1}-u_i\right)+g(u_i;a).
\end{equation}
In addition, for a time-dependent sequence \(\mathbf{w}(t) =(w_i(t))_{i \in \mathbb{Z}}\)
we define the dynamic residual
\[
\mathcal N[\mathbf{w}]_i(t)
:=
\dot w_i(t)-F_i(\mathbf{w}(t)).
\]

\begin{proposition}[Comparison principle]\label{prop:comparison}
Let $\underline{\mathbf{u}},\overline{\mathbf{u}} \in C^1\big([0, \infty), \ell^\infty(\mathbb Z) \big)$
be such that
%$\underline{\mathbf{u}}(\cdot)=(\underline u_i(\cdot))_{i\in\mathbb Z}$ and $\bar{\mathbf{u}}(\cdot)=(\bar{u}_i(\cdot))_{i\in\mathbb Z}$ are bounded sequences. If
\[
\mathcal N[\underline{\mathbf{u}}]_i(t)\le0,
\qquad
\mathcal N[\overline{\mathbf{u}}]_i(t)\ge0,
\]
for all $i\in\mathbb Z$ and $t \ge 0$. Suppose furthermore that 
%and for some $t_0$, it holds
\[
\underline u_i(0)\le\overline u_i(0),
\qquad i\in\mathbb Z.
\]
Then
\[
\underline u_i(t)\le\overline u_i(t),
\qquad
i\in\mathbb Z,\quad t\ge 0. %t_0.
\]
\end{proposition}

\begin{proof}
The statement is the reformulation of \cite[Lemma 1]{chen2008traveling} with $j=\infty$.
\end{proof}

The first corollary implies that the LDE~\eqref{e:main} preserves $p$-monotonicity. 

\begin{lemma}[$p$-monotonicity invariance]\label{lemma:comparison-principle-monotone}
Let $p\in\mathbb{N}$ and \(\mathbf{u}^0=(u_i^0)_{i\in\mathbb Z}\in\ell^\infty(\mathbb Z)\)
be $p$-increasing 
\begin{equation}
\label{eq:comparison-p-monotone-0}
    u_i^0\leq u_{i+p}^0,
    \qquad i\in\mathbb Z.
\end{equation}
Then the solution $\mathbf{u}(t)$ to \eqref{e:main} with $\mathbf{u}(0)=\mathbf{u}^0$ is also $p$-increasing
\begin{equation}
\label{eq:comparison-p-monotone-t}
    u_i(t)\leq u_{i+p}(t),
    \qquad i\in\mathbb Z, \quad t>0.
\end{equation}
\end{lemma}

\begin{proof}
Following the proof of \cite[Lemma 3.2]{hupkes2023propagation}, we define $v_i(t):=u_{i+p}(t)$, which is also a solution to \eqref{e:main}. Moreover,
\[
    u_i(0)\leq v_i(0)=u_{i+p}(0),
\]
and the application of Proposition~\ref{prop:comparison} yields the desired $p$-monotonicity~\eqref{eq:comparison-p-monotone-t}.
\end{proof}

The next consequence implies that the existence of a sub-solution prevents the front solution from traveling to the right.

\begin{corollary}[Stationary sub-solution criterion]
\label{cor:speed}
Consider a traveling wave $(c, \boldsymbol{\Phi})$  
as constructed in
Proposition~\ref{prop:periodic-wave-existence} and
write $\mathbf{u}^{\mathrm{tw}}$ for the associated solution to \eqref{e:main}
defined via \eqref{eq:wave-ansatz}.
Suppose that there exists a sequence
$\underline{\mathbf{u}}=(\underline u_i)_{i\in\mathbb Z}$ such that
\[
0<\underline u_i<1,
\qquad
F_i(\underline u)\ge0,
\qquad i\in\mathbb Z,
\]
and
\[
\lim_{i\to-\infty}\underline u_i=0,
\qquad
\lim_{i\to+\infty}\underline u_i=1.
\]
Assume, in addition, that there exists an integer $N\ge1$ such that
\[
\underline u_i\le u^{\mathrm{tw}}_i(0),
\qquad |i|\ge N.
\]
Then $c\le0$.
\end{corollary}

\begin{proof}
Arguing by contradiction, suppose that $c>0$. For each fixed $|i|<N$, one has $u^{\mathrm{tw}}_i(t)\to1$ as $t\to-\infty$. Since $\underline u_i<1$, there exists $t_*<0$ such that
\[
\underline u_i<u^{\mathrm{tw}}_i(t_*),
\qquad |i|<N.
\]
By the monotonicity of the profiles in
Proposition~\ref{prop:periodic-wave-existence}, $u^{\mathrm{tw}}_i(t_*)>u^{\mathrm{tw}}_i(0)$ for $i\in\mathbb Z$. Consequently, $\underline u_i<u^{\mathrm{tw}}_i(t_*)$ for $i\in\mathbb Z$. Applying the comparison principle from Proposition~\ref{prop:comparison}, we have
\[
\underline u_i\le u^{\mathrm{tw}}_i(t),
\qquad i\in\mathbb Z,\quad t\ge t_*.
\]
Fix any $i_0\in\mathbb Z$. Since $c>0$, we have $i_0-ct\to-\infty$ as $t\to+\infty$ and hence
$$0<\underline u_{i_0}\le u^{\mathrm{tw}}_{i_0}(t) \to 0$$
which is impossible. Therefore $c\le0$.
\end{proof}

Finally, we formulate a comparison principle directly for the traveling MFDE \eqref{e:main-mfde}. To this end, we define the MFDE residual
\begin{equation}\label{e:cal-I}
    \boldsymbol{\mathcal{I}}_{a,d,\boldsymbol{\omega},\mathbf{k}}[c,\boldsymbol{\Phi}]_j(\xi) :=  -c \Phi'_j(\xi) - d \mathfrak{D}_{\boldsymbol{\omega}, \mathbf k}[\boldsymbol{\Phi}]_j(\xi)- g(\Phi_j(\xi);a) . 
\end{equation}

\begin{corollary}[Comparison principle for MFDE]
\label{cor:comparison-MFDE}
    Let $(a,d)\in\mathcal H$ and the pair $(c,\boldsymbol\Phi)$
    be a solution to~\eqref{e:main-mfde},~\eqref{e:main-mfde-bc}.
    Assume that there exists a constant $\bar c<0$
    (resp. $\bar c>0$) and bounded profiles
    $\boldsymbol\Upsilon=\langle\Upsilon_j\rangle$,
    $\Upsilon_j\in C^1(\mathbb R)$, that satisfy the
    elementwise inequalities
    \begin{enumerate}
        \item
        $\sup\limits_{\xi\in\mathbb R}\boldsymbol\Upsilon(\xi)>0$
        (resp.
        $\inf\limits_{\xi\in\mathbb R}\boldsymbol\Upsilon(\xi)<1$),
        \item
        $\boldsymbol\Upsilon(\xi)\le\boldsymbol\Phi(\xi)$
        (resp. $\boldsymbol\Upsilon(\xi)\ge\boldsymbol\Phi(\xi)$)
        for all $\xi\in\mathbb R$,
        \item
        $\boldsymbol{\mathcal I}_{a,d,\boldsymbol\omega,\mathbf k}
        [\bar c,\boldsymbol\Upsilon](\xi)\le0$
        (resp.
        $\boldsymbol{\mathcal I}_{a,d,\boldsymbol\omega,\mathbf k}
        [\bar c,\boldsymbol\Upsilon](\xi)\ge0$)
        for all $\xi\in\mathbb R$.
    \end{enumerate}
    Then we have $c\le\bar c<0$
    (resp. $c\ge\bar c>0$).
\end{corollary}

\begin{proof}
The proof closely follows the proof of~\cite[Corollary~3.3]{hupkes2023propagation} which covers the case $p = 1$. We prove the case $c \le \overline{c}$. The reverse inequality $c \ge \overline{c}$ can be shown similarly. 

Consider $\boldsymbol{\Upsilon}(\xi) \le \boldsymbol{\Phi}(x)$ for all $\xi \in \mathbb{R}$ elementwise. We define 
\begin{align*}
v_{np +j}(t) &:= \Upsilon_j(np+j-\overline{c}t), \\
u_{np +j}(t) &:= \Phi_j(np+j-ct).
\end{align*}
By construction, these time-dependent sequences satisfy the assumptions of Proposition~\ref{prop:comparison} and thus we have 
\begin{equation} \label{e:cor-compar}
   \Upsilon_j(np +j -\overline{c}t) \le \Phi_j(np + j-ct), \qquad n \in \mathbb{Z}, \quad \dInt{j}.
\end{equation}
We are now set up to proceed by contradiction. Assume $c> \overline{c}$. Without loss of generality, pick $\dInt{j}$. There exists $M > 0$ and $\xi_0 \in \mathbb{R}$ such that $\Upsilon_j(\xi_0) = M$. For the sake of simplicity, shift the profiles so that $\xi_0:=0$. Since $\lim_{\xi \to -\infty}\Phi_j(\xi) = 0$ (due to the boundary condition~\eqref{e:main-mfde-bc}), we can find an integer $n_1<0$ such that $\Phi_j(\xi)<M/2$ for all $\xi\le n_1 p +j <0 $. Now define $t_1>0$ by 
\[
    n_1 p + j = -(c- \overline{c})t_1. 
\]
Next, pick $t_2 \ge t_1 >0$ and $n_2 \in \mathbb{Z}$ such that $ \overline{c}t_2 = n_2 p +j \in \mathbb{Z}$ holds. Putting all the above together we have 
\[
v_{n_2 p +j}(t_2) = \Upsilon_j(n_2 p +j - \overline{c}t_2) = \Upsilon_j(0) = M.
\]
On the other hand, 
\[
n_2 p + j -c t_2 = (\overline{c} - c) t_2 \le (\overline{c}-c)t_1 = n_1 p + j
\]
which implies $u_{n_2 p +j}(t_2) = \Phi_j(n_2 p + j - c t_2) < M/2$ which contradicts~\eqref{e:cor-compar}.
\end{proof}

%%%%%%%%%%%%%%%%%%%%%%%%%%%%%%%%%%%%%%%%%%%%%%%%%%%%%%%%%%%%%%%%%%%
%%%%%%%%%%%%%%%%%%%%%%%%%%%%%%%%%%%%%%%%%%%%%%%%%%%%%%%%%%%%%%%%%%%
% PINNED PERIODIC WAVES
%%%%%%%%%%%%%%%%%%%%%%%%%%%%%%%%%%%%%%%%%%%%%%%%%%%%%%%%%%%%%%%%%%%
%%%%%%%%%%%%%%%%%%%%%%%%%%%%%%%%%%%%%%%%%%%%%%%%%%%%%%%%%%%%%%%%%%%

\section{Pinned periodic monotonic waves}
\label{sec:pinned-waves}
\begin{figure}[t]
    \centering
    \includegraphics[width=0.8\linewidth]{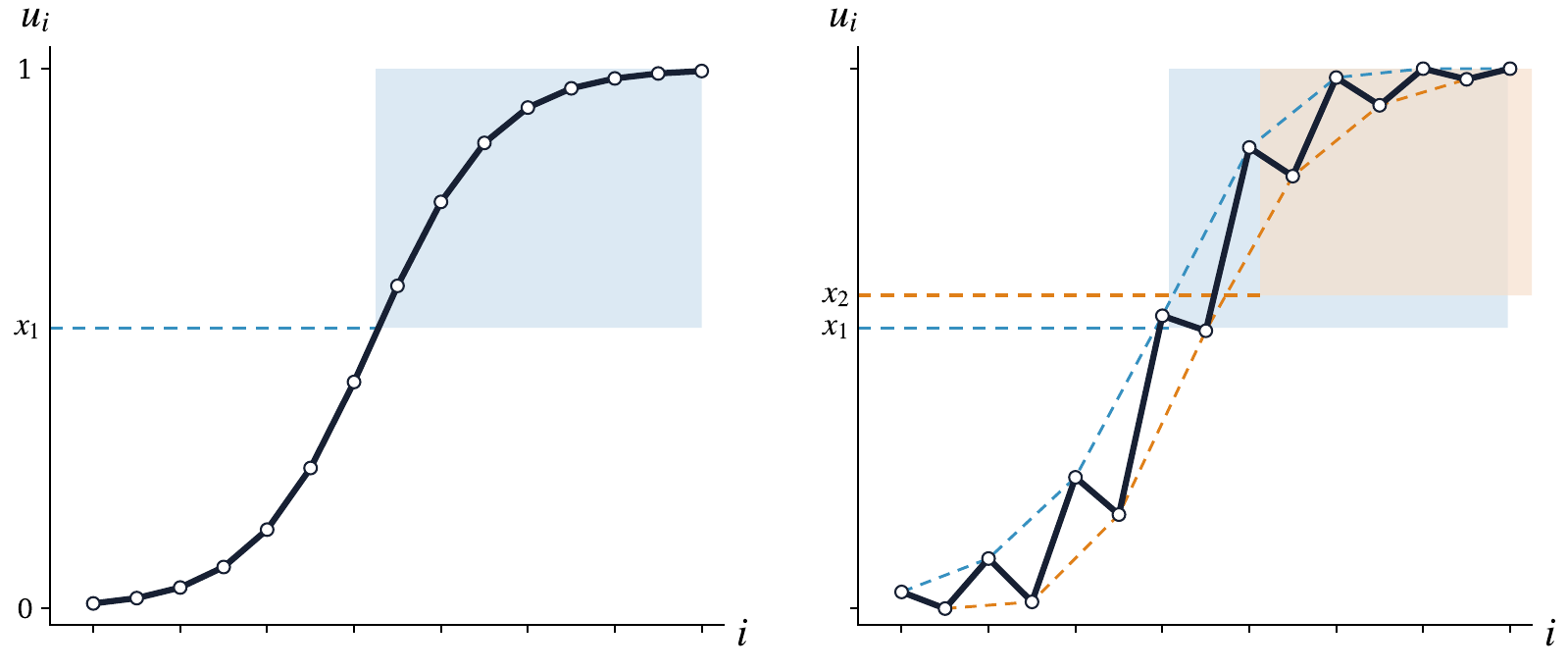}
    \caption{Illustration of the constraints defining \(\mathcal S^-_{s}(\mathbf x)\) in \eqref{eq:mr:seq:trap}. For \(p=L=1\) (left), a single threshold \(x_1\) describes the trapped interval \([x_1,1]\), indicated by the blue shading. For \(p=L=2\) (right), the two consecutive constrained sites have separate lower thresholds \(x_1\) and \(x_2\), one for each phase. The black curve connects the lattice values in spatial order, while the blue and orange dashed curves connect the two phase subsequences in the right panel. Both panels illustrate the case \(L=p\).}
    \label{fig:invariant}
\end{figure}
In this section, we establish Theorem \ref{thm:periodic-pinning} inspired by the invariant interval approach developed by Keener for the spatially homogeneous setting \cite{keener1987propagation}. 
Indeed, in the special case $L = p = 1$ the interval $[x_1,1]$ that appears in the definition \eqref{eq:mr:seq:trap} of our sequence space $\mathcal{S}_s^-$ corresponds precisely to the invariant interval constructed by Keener. In \cite[\S4]{hupkes2023propagation} this construction was adapted to homogeneous trees and our generalization here covers the periodic setting. We emphasize, however, that the freedom to pick $L \ge p$ in the thresholds $d^\pm$, \eqref{eq:d-minus-periodic} and~\eqref{eq:d-plus-periodic}, is new here and sharpens the prior bounds from \cite{keener1987propagation,hupkes2023propagation}.

We adapt this argument to the periodic setting in order to establish conditions that guarantee the invariance of the sets
\(\mathcal S_s^-(\mathbf x)\) and \(\mathcal S_s^+(\mathbf y)\) defined in \eqref{eq:mr:seq:trap} under the flow of \eqref{e:main}; see Figure~\ref{fig:invariant}. In particular, Lemma \ref{lemma:comparison-principle-monotone} shows that $p$-monotonicity is invariant, which allows us to restrict our attention to solutions of \eqref{e:main} that satisfy $0\le u_i(t)\le1$ and $u_i(t)\le u_{i+p}(t)$ for all $i\in\mathbb Z$ and $t\ge0$. 

For $(a,d)\in\mathcal P^-(\boldsymbol{\omega},\mathbf k)$, we show that there exist a phase $\dInt{s}$ and a vector $\mathbf x\in(a,1)^L$ such that, for each fixed $n\in\mathbb Z$, the
simultaneous bounds $u_{np+s+j}(0)>x_j$ for all $1\le j\le L$ imply $u_{np+s+j}(t)>x_j$ for all $1\le j\le L$ and $t\ge0$. 
Since the traveling-wave is also $p$-monotonic,
this prevents the low state from invading the high state. Propagation in the positive direction therefore cannot occur, and hence $c(a,d;\boldsymbol{\omega},\mathbf k)\leq0$. Similarly, \((a,d)\in\mathcal P^+(\boldsymbol{\omega},\mathbf k)\) blocks propagation to the left and consequently yields $c(a,d;\boldsymbol{\omega},\mathbf k)\geq0$.

Recall
the definitions 
\eqref{eq:d-minus-periodic}
and \eqref{eq:d-plus-periodic}
for $d^\pm$, together with 
\eqref{eq:B-} and \eqref{eq:B+}
for $B^\pm$. As the first step, 
we spell out the consequences
of the inequalities $ 0 < d < d^\pm$.

\begin{lemma}
\label{lemma:B-}
    Pick $(a,d)\in\mathcal H$. If $0<d<d^-(a;\boldsymbol{\omega},\mathbf k)$, then there exists a pair $(s,\mathbf x)\in[1,p]_{\mathbb N}\times(a,1)^L$ such that
    \begin{equation}
    \label{eq:ineq-B-}
        g(x_j;a)-d B_{j,s}^-(\mathbf x;\boldsymbol{\omega},\mathbf k)>0, \qquad 1\le j\le L.
    \end{equation}
\end{lemma}

\begin{proof}
Since \(d<d^-(a;\boldsymbol{\omega},\mathbf k)\), it follows from \eqref{eq:d-minus-periodic} that there exist \(s^0\dInt{}\) and \(\mathbf x^0\in(a,1)^L\) such that
\[
    d<
    \min_{\{j:B_{j,s^0}^-(\mathbf x^0;\boldsymbol{\omega},\mathbf k)>0\}}
    \frac{g(x_j^0;a)}
    {B_{j,s^0}^-(\mathbf x^0;\boldsymbol{\omega},\mathbf k)}.
\]
For any $1\le j\le L$, if \(B_{j,s^0}^-(\mathbf x^0;\boldsymbol{\omega},\mathbf k)>0\), multiplication by this positive quantity gives
\[
    dB_{j,s^0}^-(\mathbf x^0;\boldsymbol{\omega},\mathbf k)<g(x_j^0;a).
\]
On the other hand, if \(B_{j,s^0}^-(\mathbf x^0;\boldsymbol{\omega},\mathbf k)\le 0\), then \(x_j^0\in(a,1)\) gives \(g(x_j^0;a)>0\). Together this yields \eqref{eq:ineq-B-}.
\end{proof}

The same result can be stated and proved analogously for $d^+$.

\begin{lemma}
\label{lemma:B+}
    Pick $(a,d)\in\mathcal H$.
    If $0<d<d^+(a;\boldsymbol{\omega},\mathbf k)$, then there exists a pair $(s,\mathbf y)\in[1,p]_{\mathbb N}\times(0,a)^L$ such that
    \begin{equation}
       g(y_j;a)- d B_{j,s}^+(\mathbf y;\boldsymbol{\omega},\mathbf k)<0, \qquad 1\le j\le L.
    \end{equation}
\end{lemma}

Next, we show the sequence spaces $S_{s}^\pm$ are invariant under the flow of \eqref{e:main}.

\begin{lemma}
\label{lemma:invariant-rectangles}
Pick $(a,d)\in\mathcal H$ and fix an integer $L\ge p$. 
Then the following claims hold for the solution $\mathbf{u}(t)$ of \eqref{e:main}.
\begin{enumerate}[label=(\roman*), ref=(\roman*)]
\item \label{lemma:invariant-rectangles-1} If $0<d<d^-(a;\boldsymbol{\omega},\mathbf k)$, then, there exist a phase $\dInt{s}$ and a vector $\mathbf x\in(a,1)^L$ such that, for every $n\in\mathbb Z$, $u(0)\in \mathcal S_{np+s}^-(\mathbf x)$ implies $\mathbf{u}(t)\in \mathcal S_{np+s}^-(\mathbf x)$ for all $t>0$.
\item \label{lemma:invariant-rectangles-2} If $0<d<d^+(a;\boldsymbol{\omega},\mathbf k)$, then, there exist a phase $\dInt{s}$ and a vector $\mathbf y\in(0,a)^L$ such that, for every $n\in\mathbb Z$, $u(0)\in \mathcal S_{np+s}^+(\mathbf y)$ implies $\mathbf{u}(t)\in \mathcal S_{np+s}^+(\mathbf y)$ for all $t>0$.
\end{enumerate}
\end{lemma}

\begin{proof}
We first prove \ref{lemma:invariant-rectangles-1}. By Lemma~\ref{lemma:B-}, choose a phase $\dInt{s}$ and a vector $\mathbf x\in(a,1)^L$ such that $g(x_j;a)-dB_{j,s}^-(\mathbf x;\boldsymbol{\omega},\mathbf k)>0$ for $1\le j\le L$. Fix \(n\in\mathbb Z\) and suppose that $u(0)\in\mathcal S_{np+s}^-(\mathbf x)$. By Proposition~\ref{prop:comparison}, we have \(0\le u_i(t)\le1\) for all \(i\in\mathbb Z\) and \(t\ge0\). Moreover, Lemma~\ref{lemma:comparison-principle-monotone} ensures that $u_i(t)\le u_{i+p}(t)$ for $i\in\mathbb Z$ and $t\ge0$.

By contradiction, assume that \(\mathbf{u}(t)\) leaves \(\mathcal S_{np+s}^-(\mathbf x)\) for the first time at \(\tau>0\).  Then there exists \(1\le j_*\le L\) such that $u_{np+s+j_*}(\tau)=x_{j_*}$ and $\dot u_{np+s+j_*}(\tau)\leq0$, together with $u_{np+s+j}(\tau)\geq x_j$ for $1\le j\le L$. The exterior neighbors satisfy
\[
    u_{np+s}(\tau)\ge0=x_0,
    \qquad
    u_{np+s+L+1}(\tau)
    \ge u_{np+s+L+1-p}(\tau)
    \ge x_{L+1-p}=x_{L+1}.
\]
With $k_{np+s+j_*}=k_{s+j_*}$ and $\omega_{np+s+j_*\pm1/2}=\omega_{s+j_*\pm1/2}$, the definition of \(B_{j,s}^-\) gives
\[
\begin{aligned}
    \dot u_{np+s+j_*}(\tau)
    &=d\Big[k_{s+j_*}\omega_{s+j_*+1/2}\bigl(u_{np+s+j_*+1}(\tau)-x_{j_*}\bigr)+\omega_{s+j_*-1/2}\bigl(u_{np+s+j_*-1}(\tau)-x_{j_*}\bigr)\Big]+g(x_{j_*};a)\\
    &\ge d\Big[k_{s+j_*}\omega_{s+j_*+1/2}\bigl(x_{j_*+1}-x_{j_*}\bigr)+\omega_{s+j_*-1/2}\bigl(x_{j_*-1}-x_{j_*}\bigr)\Big]+g(x_{j_*};a)\\
    &=g(x_{j_*};a)-dB_{j_*,s}^-(\mathbf x;\boldsymbol{\omega},\mathbf k)>0,
\end{aligned}
\]
which is a contradiction. The proof of \ref{lemma:invariant-rectangles-2} is analogous. 
\end{proof}

We are now ready to prove our main result, Theorem~\ref{thm:periodic-pinning}.

\begin{proof}[Proof of Theorem \ref{thm:periodic-pinning}] 
We first prove \ref{thm:periodic-pinning-1}. Fix any $a\in(0,1)$ and a phase $\dInt{s}$. Choose $x\in(a,1)$ and $y\in(0,a)$, and set
\[
    \mathbf x=(x,\ldots,x)\in(a,1)^L,
    \qquad
    \mathbf y=(y,\ldots,y)\in(0,a)^L.
\]
Recalling \eqref{eq:B-} and \eqref{eq:B+}, we obtain
\[
    B_{j,s}^-(\mathbf x;\boldsymbol{\omega},\mathbf k)
    =
    \begin{cases}
        \omega_{s+1/2}\,x, & j=1,\\
        0, & 2\le j\le L,
    \end{cases}
    \qquad
    B_{j,s}^+(\mathbf y;\boldsymbol{\omega},\mathbf k)
    =
    \begin{cases}
        0, & 1\le j<L,\\
        k_{s+L}\omega_{s+L+1/2}\,(y-1), & j=L,
    \end{cases}
\]
where $k_{s+L+p}=k_{s+L}$ and $\omega_{s+L+p+1/2}=\omega_{s+L+1/2}$. Thus \eqref{eq:d-minus-periodic} and \eqref{eq:d-plus-periodic} imply
\[
    d^-(a;\boldsymbol{\omega},\mathbf k)
    \ge \frac{g(x;a)}{\omega_{s+1/2}\,x}>0,
    \qquad
    d^+(a;\boldsymbol{\omega},\mathbf k)
    \ge \frac{g(y;a)}{k_{s+L}\omega_{s+L+1/2}\,(y-1)}>0.
\]
The strict positivity follows from \eqref{Hg} since $g(x;a)>0$, $g(y;a)<0$, and $k_{s+L}\omega_{s+L+1/2}(y-1)<0$. Consequently, both $\mathcal P^-(\boldsymbol{\omega},\mathbf k)$ and $\mathcal P^+(\boldsymbol{\omega},\mathbf k)$ are nonempty.  Since both thresholds are positive, choosing 
\[
0<d<\min\{d^-(a;\boldsymbol\omega,\mathbf k), d^+(a;\boldsymbol\omega,\mathbf k)\}
\]
gives $(a,d)\in\mathcal P^0(\boldsymbol\omega,\mathbf k)$, showing that $\mathcal P^0(\boldsymbol{\omega},\mathbf k)$ is nonempty as well.

Turning to \ref{thm:periodic-pinning-2}, 
consider a traveling wave $(c, \boldsymbol{\Phi})$  
as constructed in
Proposition~\ref{prop:periodic-wave-existence}
and write $\mathbf{u}^{\rm tw}$
for the associated 
solution to \eqref{e:main}
defined via \eqref{eq:wave-ansatz}.
The monotonicity and far-field limits of its profiles imply
\[
    0\le u^{\rm tw}_i(t)\le1,
    \qquad
    u^{\rm tw}_i(t)\le u^{\rm tw}_{i+p}(t),
    \qquad i\in\mathbb Z,\quad t\in\mathbb R.
\] 
Assume that $(a,d)\in\mathcal P^-(\boldsymbol{\omega},\mathbf k)$. By the definition of $\mathcal P^-(\boldsymbol{\omega},\mathbf k)$, there exist a phase $\dInt{s}$ and a vector $\mathbf x=(x_1,\ldots,x_L)$ that satisfy the conditions of Lemma~\ref{lemma:invariant-rectangles}\ref{lemma:invariant-rectangles-1}. Since $u^{\rm tw}_i(0)\to1$ as $i\to+\infty$, there exists $N\in\mathbb Z$ such that

$$ \mathbf{u}^{\rm tw}(0)\in\mathcal S_{np+s}^-(\mathbf x),\qquad n\ge N.$$ By Lemma~\ref{lemma:invariant-rectangles}\ref{lemma:invariant-rectangles-1},

$$\mathbf{u}^{\rm tw}(t)\in\mathcal S_{np+s}^-(\mathbf x),\qquad n\ge N,\quad t\ge0.$$
If $c>0$, we have $i-ct\to -\infty$ when $t\to +\infty$ and hence
\[
    u^{\rm tw}_i(t)\longrightarrow0
\]
at every fixed lattice site $i\in\mathbb Z$. Fixing any $n\ge N$, this contradicts
\[
    u^{\rm tw}_{np+s+1}(t)>x_1>0,
    \qquad t>0.
\]
Therefore $c\le0$. The argument for \ref{thm:periodic-pinning-3} is analogous and \ref{thm:periodic-pinning-4} is an immediate consequence of the fact that $\mathcal P^0(\boldsymbol{\omega}, \mathbf k) =
    \mathcal P^-(\boldsymbol{\omega},\mathbf k)\cap\mathcal P^+(\boldsymbol{\omega},\mathbf k)$.
\end{proof}

%%%%%%%%%%%%%%%%%%%%%%%%%%%%%%%%%%%%%%%%%%%%%%%%%%%%%%%%%%%%%%%%%%%
%%%%%%%%%%%%%%%%%%%%%%%%%%%%%%%%%%%%%%%%%%%%%%%%%%%%%%%%%%%%%%%%%%%
% TAIL CONSTRUCTION
%%%%%%%%%%%%%%%%%%%%%%%%%%%%%%%%%%%%%%%%%%%%%%%%%%%%%%%%%%%%%%%%%%%
%%%%%%%%%%%%%%%%%%%%%%%%%%%%%%%%%%%%%%%%%%%%%%%%%%%%%%%%%%%%%%%%%%%

\section{Tail Construction}
\label{sec:construction}
In this section, we construct the left and right tails featured in the sub-solution~\eqref{eq:mr:def:u:i:s} for the cubic nonlinearity~\eqref{eq:cubic}, working with the periodic diffusion system \eqref{e:main}. We start
in {\S}\ref{subsec:periodic-linear-theory} by discussing the characteristic exponents at the equilibria \(0\) and \(1\).
These determine the exponential convergence rate of the tails
and allow us to extract 
the far-field ordering between stationary and moving tails that we require for the comparison argument.
This allows us to define the two sequences \(\mathcal L_{i;s}\) and \(\mathcal R_{i;s+1}\) in {\S}\ref{subsec:LR}, where we also study their basic properties. We move on to the nonlinear corrections in {\S}\ref{subsec:C_i}. In particular, we determine how the
 positive periodic correction coefficients \(C_i\) in~\eqref{eq:mr:def:u:i:s}
 should be chosen to ensure that the sub-solution
 residual has the correct sign and obtain an \emph{a~priori} condition
 on the free parameter $w$ in~\eqref{eq:mr:def:u:i:s} that restricts the range of the sub-solution to $(0,1)$.

\subsection{Periodic Linear Theory}
\label{subsec:periodic-linear-theory}

In this subsection we recall the parts that we require of the principal eigenvalue theory for periodic lattice equations that was developed in \cite{chen2008traveling}. Writing $\boldsymbol{\phi}=(\phi_i)_{i\in\mathbb Z}=\langle \phi_1, \ldots, \phi_p \rangle $ for a $p$-periodic steady state of \eqref{e:main}, the stationary linearization around $\boldsymbol{\phi}$ is given by
\begin{equation}
\label{eq:linearization}
    0=k_i d\,\omega_{i+1/2}(U_{i+1}-U_i)+d\,\omega_{i-1/2}(U_{i-1}-U_i)+g'(\phi_i;a)U_i,\qquad i\in\mathbb Z.
\end{equation}
Let \(\boldsymbol\psi=(\psi_i)_{i\in\mathbb Z} =\langle\psi_1,\ldots,\psi_p\rangle\) be a positive \(p\)-periodic vector. Searching for exponential tails as customary in periodic media, we consider solutions of the form 
\begin{equation}
\label{eq:U}
    U_i=e^{\lambda i}\psi_i,\qquad\psi_{i+p}=\psi_i
\end{equation}
which leads to the system $\mathcal M(\lambda)[\boldsymbol{\psi}] = 0$, with 
\begin{equation}
\label{eq:periodic-eigenvalue-operator}
[\mathcal M(\lambda)\boldsymbol{\psi}]_i
=k_i d\,\omega_{i+1/2}( e^\lambda\psi_{i+1}-\psi_i)+
d\,\omega_{i-1/2}( e^{-\lambda}\psi_{i-1}-\psi_i)+
g'(\phi_i;a)\psi_i,
\qquad i\in\mathbb Z.
\end{equation}
This is also directly connected to the linearization of the traveling wave MFDE \eqref{e:main-mfde} around $\boldsymbol{\phi}$, which reads
\begin{equation}
\label{eq:linearized-mfde}
\begin{aligned}
    -c v_j'(\xi)&= k_jd\,\omega_{j+1/2}\left[ v_{j+1}(\xi+1)-v_j(\xi)\right]+ d\,\omega_{j-1/2}\left[ v_{j-1}(\xi-1)-v_j(\xi)\right]\\
    &\qquad+ g'(\phi_j;a)v_j(\xi),
    \qquad \dInt{j}.
\end{aligned}
\end{equation}
Indeed, searching for a solution
\begin{equation}
    v_i(\xi) = e^{\lambda \xi} \psi_i, \qquad \psi_{i+p } = \psi_i,
\end{equation}
we directly find
\begin{equation}
    -c \lambda \boldsymbol{\psi}  = \mathcal M (\lambda)[\boldsymbol{\psi}].
\end{equation}

These eigenvalue problems are of Perronâ€“Frobenius type, which aids the analysis considerably. For each $\lambda\in\mathbb R$,  \cite[Theorem 1]{chen2008traveling} ensures the existence of a principal eigenvalue $\mu(\lambda)$ for $\mathcal{M}(\lambda)$ that is accompanied by a strictly positive $p$-periodic eigenfunction $\psi(\lambda)$, i.e. $\mathcal{M}(\lambda) \psi(\lambda) = \mu(\lambda) \psi(\lambda)$.
The admissible exponential rates for positive tails are therefore determined by the characteristic equation 
\begin{equation}
\label{eq:characteristic}
    \mu(\lambda)+c\lambda=0,
\end{equation}
where we take $c=0$ for the stationary 
case \eqref{eq:linearization}.

For the left tails, we take $\boldsymbol{\phi} \equiv 0$ and note $g'(0;a) = -a$. Writing $\mu_L$ and $\boldsymbol{\psi}^L$
for the associated principal eigenvalue and eigenvector, we see that
\begin{equation}
\label{eq:left-linearized}
\begin{aligned}
\mu_L(\lambda)\psi_i^L(\lambda)&=k_i d\,\omega_{i+1/2}\bigl[e^\lambda\psi_{i+1}^L(\lambda)-\psi_i^L(\lambda)\bigr]+d\,\omega_{i-1/2}\bigl[ e^{-\lambda}\psi_{i-1}^L(\lambda) - \psi_i^L(\lambda)\bigr]\\
&\qquad-a\psi_i^L(\lambda),\qquad i\in\mathbb Z.
\end{aligned}
\end{equation}
Taking $\boldsymbol{\phi} \equiv 1$ for the right tails, with $g'(1;a) = a - 1$, we arrive at the corresponding equation
\begin{equation}
\label{eq:right-linearized}
\begin{aligned}
\mu_R(\lambda)\psi_i^R(\lambda)&=k_i d\,\omega_{i+1/2}\bigl[e^\lambda\psi_{i+1}^R(\lambda)-\psi_i^R(\lambda)\bigr]+d\,\omega_{i-1/2}\bigl[e^{-\lambda}\psi_{i-1}^R(\lambda)-\psi_i^R(\lambda)\bigr]\\
&\qquad-(1-a)\psi_i^R(\lambda), \qquad i\in\mathbb Z.
\end{aligned}
\end{equation}
for the associated principal eigenpair $\mu_R$ and $\boldsymbol{\psi}^R$.

\begin{lemma}
\label{lemma:characteristic-roots}
Fix \(a\in(0,1)\) and $d > 0$. For any \(c\in\mathbb R\), there exists a unique positive root \(\lambda_L^c\) that satisfies
\[
    \mu_L(\lambda_L^c)+c\lambda_L^c=0,
    \qquad \lambda_L^c>0,
\]
together with a unique negative root
\(\lambda_R^c\) that satisfies
\[
    \mu_R(\lambda_R^c)+c\lambda_R^c=0,
    \qquad \lambda_R^c<0.
\]
For \(c=0\), we write \(\lambda_L:=\lambda_L^0\) and \(\lambda_R:=\lambda_R^0\).

\end{lemma}

\begin{proof}
The existence and uniqueness of the characteristic roots $\lambda_L^c,~\lambda_R^c$ with the stated signs follow from \cite[Theorem~1 and~\S3.1]{chen2008traveling},
applied to the linearizations at \(0\) and \(1\), respectively.
\end{proof}

Positive and negative roots of~\eqref{eq:characteristic} share uniform bounds.

\begin{lemma}
\label{lemma:exponent-gap}
Under the assumptions of Lemma~\ref{lemma:characteristic-roots}, for every \(c>0\), we have
\[
    0<\lambda_L^c<\lambda_L, \qquad \lambda_R^c<\lambda_R<0.
\]
\end{lemma}

\begin{proof}
Define
\[
    P_L^c(\lambda):=\mu_L(\lambda)+c\lambda, \qquad P_R^c(\lambda):=\mu_R(\lambda)+c\lambda.
\]
Since $P_L^c(0)=-a<0$ and $P_L^c(\lambda_L)=c\lambda_L>0$, continuity and the uniqueness of the positive root from Lemma~\ref{lemma:characteristic-roots} imply $0<\lambda_L^c<\lambda_L$. For the right exponent, we have $P_R^c(\lambda_R)=c\lambda_R<0$. Moreover, \cite[Theorem~1]{chen2008traveling} gives $P_R^c(\lambda)\longrightarrow+\infty$ as $\lambda\to-\infty$. Continuity and the uniqueness of the negative root therefore yield $\lambda_R^c<\lambda_R<0$.
\end{proof}

Using the stationary exponents from Lemma~\ref{lemma:characteristic-roots}, define
\begin{equation}
\label{eq:left-tail}
    U_i^L:=e^{\lambda_L i}\psi_i^L(\lambda_L).
\end{equation}
Since \(\mu_L(\lambda_L)=0\), the sequence $U_i^L$ solves \eqref{eq:left-linearized} with \(\phi_i=0\). Moreover, \(U_i^L>0\) and \(U_i^L\to0\) as \(i\to-\infty\). Similarly, define
\begin{equation}
\label{eq:right-tail}
    U_i^R:=e^{\lambda_R i}\psi_i^R(\lambda_R),
\end{equation}
which solves \eqref{eq:right-linearized} with \(\phi_i=1\) and satisfies \(U_i^R>0\) and \(U_i^R\to0\) as \(i\to+\infty\). The exponents \(\lambda_L^c\) and \(\lambda_R^c\) describe the far-field behavior of actual traveling waves, as stated in the following result.

\begin{proposition}[\cite{chen2008traveling} Theorem 2, Lemma 6]
\label{prop:CGW-tail}
Let the pair $(c, \boldsymbol{\Phi})$ be a solution to~\eqref{e:main-mfde},~\eqref{e:main-mfde-bc}
and write $\mathbf{u}^{\mathrm{tw}}$ for the associated 
traveling wave solution to \eqref{e:main}
defined via \eqref{eq:wave-ansatz}.
Then there exist constants $h_L,h_R>0$ such that
\[
\lim_{i-ct\to-\infty}
\frac{u^{\mathrm{tw}}_i(t)}
{\psi_i^L(\lambda_L^c)e^{\lambda_L^c(i-ct)}}
=h_L,
\]
and
\[
\lim_{i-ct\to+\infty}
\frac{1-u^{\mathrm{tw}}_i(t)}
{\psi_i^R(\lambda_R^c)e^{\lambda_R^c(i-ct)}}
=h_R.
\]
Moreover, if $c>0$, then
\[
\dot u^{\mathrm{tw}}_i(t)<0,
\qquad i\in\mathbb Z,\quad t\in\mathbb R.
\]
\end{proposition}
\begin{proof}
    This conclusion is a direct consequence of \cite[Theorem~2 and Lemma~6]{chen2008traveling}.
\end{proof}

\subsection{Construction of \texorpdfstring{$\mathcal L$}{L} and \texorpdfstring{$\mathcal R$}{R}}
\label{subsec:LR}

Since our main system \eqref{e:main} only connects
the neighboring values \(u_{i-1},u_i,u_{i+1}\), it is more convenient to record how each tail changes from one lattice site to the next. We therefore define the left-tail ratios by
\begin{equation}
\ell_i
:=
\frac{U_{i-1}^L}{U_i^L}
=
e^{-\lambda_L}
\frac{\psi_{i-1}^L(\lambda_L)}{\psi_i^L(\lambda_L)}.
\label{eq:left-ratio}
\end{equation}
The sequence $\boldsymbol{\ell}=(\ell_i)_{i\in\mathbb Z}=\langle \ell_1, \ldots, \ell_p \rangle$ is positive and $p$-periodic. Moreover, substituting
$$
U_{i-1}^L=\ell_i U_i^L,
\qquad
U_{i+1}^L=\frac{U_i^L}{\ell_{i+1}}
$$
into \eqref{eq:linearization} with \(\phi_i=0\) and dividing by \(U_i^L>0\) gives%into \eqref{eq:left-linearized} gives
\begin{equation}
\label{eq:left-ratio-relation}
    d\,\omega_{i-1/2}\ell_i+\frac{k_id\,\omega_{i+1/2}}{\ell_{i+1}} = d\,\omega_{i-1/2} + k_id\,\omega_{i+1/2} +a.
\end{equation}
Analogously, we define the right-tail ratios by
\begin{equation}
r_i:=\frac{U_{i+1}^R}{U_i^R}=e^{\lambda_R}
\frac{\psi_{i+1}^R(\lambda_R)}{\psi_i^R(\lambda_R)}.
\label{eq:right-ratio}
\end{equation}
The sequence $\mathbf{r}=(r_i)_{i\in\mathbb Z}=\langle r_1, \ldots, r_p \rangle$ is positive and $p$-periodic. Substituting
$$
U_{i+1}^R=r_iU_i^R,
\qquad
U_{i-1}^R=\frac{U_i^R}{r_{i-1}}
$$
into \eqref{eq:linearization} with \(\phi_i=1\) and dividing by \(U_i^R>0\) gives %into \eqref{eq:right-linearized} gives
\begin{equation}
\label{eq:right-ratio-relation}
    k_id\,\omega_{i+1/2}r_i+\frac{d\,\omega_{i-1/2}}{r_{i-1}}
    = k_id\,\omega_{i+1/2}+ d\,\omega_{i-1/2}+(1-a).
\end{equation}
Our next result provides useful bounds for ratios $\boldsymbol{\ell}$ and $\mathbf{r}$, both aggregate and pointwise. 

\begin{lemma}
\label{lemma:products}
Let $\boldsymbol{\ell}$ and $\mathbf{r}$ be the positive \(p\)-periodic stationary tail ratios defined by \eqref{eq:left-ratio} and \eqref{eq:right-ratio}, respectively. Then the following statements hold.
\begin{enumerate}[label=(\roman*), ref=(\roman*)]
    \item \label{lemma:products-1}
    The one-period products satisfy
    \[
        0<\Lambda_L:=\prod_{j=0}^{p-1}\ell_j=e^{-p\lambda_L}<1,\qquad0<\Lambda_R:=\prod_{j=0}^{p-1}r_j=e^{p\lambda_R}<1.
    \]

    \item \label{lemma:products-2}
    We have the pointwise bounds
    \[
        0<\ell_i<1,\qquad 0<r_i<1,  \qquad i\in\mathbb Z.
    \]
\end{enumerate}
\end{lemma}

\begin{proof}
We first prove item~\ref{lemma:products-1}. By \eqref{eq:left-ratio}, \eqref{eq:right-ratio}, and the $p$-periodicity of the eigenvectors $\boldsymbol{\psi}^L$ and $\boldsymbol{\psi}^R$, we have
\[
\begin{aligned}
    \prod_{j=0}^{p-1}\ell_j&=e^{-p\lambda_L}\prod_{j=0}^{p-1}\frac{\psi_{j-1}^L}{\psi_j^L}= e^{-p\lambda_L}, \\
    \prod_{j=0}^{p-1}r_j&=e^{p\lambda_R}\prod_{j=0}^{p-1}\frac{\psi_{j+1}^R}{\psi_j^R} =e^{p\lambda_R}.
\end{aligned}
\]
The claimed bounds follow from \(\lambda_L>0\) and \(\lambda_R<0\).

For item~\ref{lemma:products-2}, positivity follows from the definitions. Suppose that \(r_j\ge1\) for some \(j\in\mathbb Z\). Applying \eqref{eq:right-ratio-relation} at \(j+1\) gives
$$
k_{j+1}d\,\omega_{j+3/2}(r_{j+1}-1)=(1-a)+d\,\omega_{j+1/2}
\left(1-\frac1{r_j}\right)>0.
$$
Thus \(r_{j+1}>1\). Applying the same relation successively at \(j+2,\ldots,j+p\) yields
\[
    r_{j+m}>1,\qquad 1\le m\le p.
\]
By periodicity, this implies
\[
    \Lambda_R=\prod_{m=1}^{p}r_{j+m}>1,
\]
contrary to item~\ref{lemma:products-1}. Similarly, if \(\ell_j\ge1\), applying \eqref{eq:left-ratio-relation} at \(j-1\) gives
\[
d\,\omega_{j-3/2}(\ell_{j-1}-1)=a+k_{j-1}d\,\omega_{j-1/2}\left( 1-\frac1{\ell_j}\right)
>0.
\]
Iterating gives \(\ell_{j-m}>1\) for \(1\le m\le p\), contradicting \(\Lambda_L<1\). This proves item~\ref{lemma:products-2}.
\end{proof}

We are now ready to define the exponential tails $\mathcal{L}$ and $\mathcal{R}$ used in the construction of the sub-solution \eqref{eq:mr:def:u:i:s}. Pick a phase $\dInt{s}$. For $i\le s$, we define
\begin{equation}
\label{eq:L}
\mathcal L_{i;s}:=\prod_{j=i+1}^{s}\ell_j,
\end{equation}
while for $i\ge s+1$ we define
\begin{equation}
\label{eq:R}
\mathcal R_{i;s+1}:=\prod_{\nu=s+1}^{i-1}r_{\nu}.
\end{equation}
Setting empty products to one, we observe the normalizations
\begin{equation}
   \mathcal L_{s;s} = 1, \qquad \qquad   \mathcal R_{s+1;s+1}=1.
\end{equation}
Lemma~\ref{lemma:products} shows that \(\mathcal L_{i;s}>0\), \(\mathcal R_{i;s+1}>0\), and
\[
\mathcal L_{i-p;s}=\Lambda_L\mathcal L_{i;s},\qquad
\mathcal R_{i+p;s+1}=\Lambda_R\mathcal R_{i;s+1},
\]
with \(0<\Lambda_L,\Lambda_R<1\). In particular,
\(\mathcal L_{i;s}\to0\) as \(i\to-\infty\) and \(\mathcal R_{i;s+1}\to0\) as \(i\to+\infty\). 
By construction, $(\mathcal L_{i;s})_{i \le s}$ satisfies \eqref{eq:linearization} with \(\boldsymbol{\phi}\equiv0\) for \(i\le s-1\), while \((\mathcal R_{i;s+1})_{i\ge s+1}\) satisfies the same equation with \(\boldsymbol{\phi}\equiv1\) for \(i\ge s+2\).

\subsection{Definition of \texorpdfstring{$C_i$}{Ci}}\label{subsec:C_i}
In this section we construct the two tails of the stationary sub-solution \eqref{eq:mr:def:u:i:s} from the exponential factors \(\mathcal L_{i;s}\) and
\(\mathcal R_{i;s+1}\) defined in \eqref{eq:L}-\eqref{eq:R}. In particular, we construct a suitable periodic sequence $\mathbf C=(C_i)_{i\in\mathbb Z}=\langle C_1,\dots,C_p\rangle$ that allows us to extract the correct sign from the residual for the right tails. The left tails do not need such a correction as the leading order term already provides the correct sign. We recall that
\(F_i\) denotes the stationary residual defined in \eqref{eq:residual} for a time-independent sequence \(\mathbf{u}=(u_i)_{i\in\mathbb Z}\). 

\begin{figure}[t]
    \centering
    \includegraphics[width=1\linewidth]{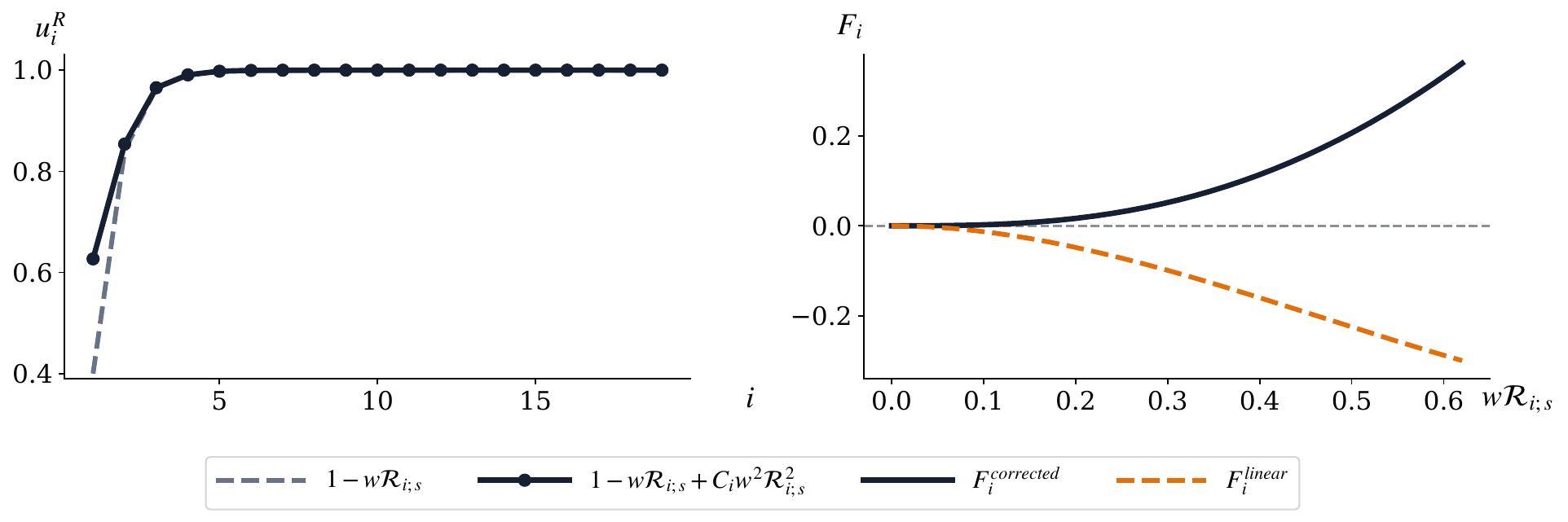}
    \caption{Effect of the quadratic correction on the right tail of profile \eqref{eq:mr:def:u:i:s} and the stationary residual \(F_i\) defined in \eqref{eq:residual}. The parameters are \(\mathbf k=\langle 1,2\rangle\), \(s=0\), \(a=0.6\), and \(d=0.2\), with right tail amplitude \(w=0.6\). The left panel compares the uncorrected right tail \(1-w\mathcal R_{i;s}\) with the quadratically corrected right tail \(1-w\mathcal R_{i;s+1}+C_iw^2(\mathcal R_{i;s+1})^2\). The right panel compares the corresponding right-tail residuals plotted as functions of \(w\mathcal R_{i;s+1}\).}
    \label{fig:correction}
\end{figure}

Let us denote the left and right tails in~\eqref{eq:mr:def:u:i:s} $\mathbf{u}^L$ and $\mathbf{u}^R$ for simplicity.  For the left tail, we introduce an amplitude \(0<v<1\) and set
\begin{equation}
\label{eq:UL}
    u_i^L=v\mathcal L_{i;s}, \qquad i\le s.
\end{equation}
For the right tail, we first write
\begin{equation}
\label{eq:def:c:def:q:i}
    q_i:=w\mathcal R_{i;s+1}, \qquad i\ge s+1,
\end{equation}
using \(w>0\) to control the right-tail amplitude. We then introduce the corrected form
\begin{equation}
\label{eq:UR}
    u_i^R=1-q_i+C_iq_i^2 = 1 - q_i(1 - C_i q_i), \qquad i\ge s+1.
\end{equation}
The coefficients \(C_i\) are chosen so that the second-order terms in the right-tail residual cancel; see Figure~\ref{fig:correction}. More precisely, we recall the right-tail ratios $r$ from \eqref{eq:right-ratio} and note that $r_i = q_{i+1}/q_i$.

We now define \(\mathbf C=(C_i)_{i\in
\mathbb Z}=\langle C_1,\dots,C_p\rangle\) as the \(p\)-periodic solution of
\begin{equation}
\label{eq:C}
k_id\,\omega_{i+1/2}r_i^2C_{i+1}-\bigl(k_id\,\omega_{i+1/2}+d\,\omega_{i-1/2}\bigr)C_i+\frac{d\,\omega_{i-1/2}}{r_{i-1}^2}C_{i-1}
-(1-a)C_i=2-a.
\end{equation}

\begin{proposition}[Existence, uniqueness, positivity, and ordering of \(C_i\)]
\label{prop:C}
The periodic linear system \eqref{eq:C} has a unique \(p\)-periodic solution \(\mathbf C=(C_i)_{i\in\mathbb Z}=\langle C_1,\dots,C_p\rangle\). Moreover, $C_i>0$ and $C_{i+1}r_i<C_i$ for all $i\in\mathbb Z$. 
\end{proposition}
\begin{proof}
Pick two $p$-periodic sequences 
\(\lambda_i>0\) and \(b_i\) and consider the first-order recursion
\begin{equation}
\label{eq:tail:periodic-recursion}
x_i=\lambda_i\,x_{i+1}+b_i,\quad i\in\mathbb Z.
\end{equation}
Then if  \(\Lambda:=\prod_{j=1}^p\lambda_j<1\),
this recursion has a unique
\(p\)-periodic solution
\begin{equation}
\label{eq:tail:periodic-recursion:sol}
x_i=\frac{1}{1-\Lambda}\sum_{j=0}^{p-1}\Bigl(\prod_{l=0}^{j-1}\lambda_{i+l}\Bigr)b_{i+j},
\end{equation}
which is positive when \(b>0\).
This can be readily verified by iterating the recursion $p$ times. The same holds for the backward recursion \(x_{i-1}=\lambda_ix_i+b_{i-1}\).

To apply this, we put \(Z_i:=C_i-r_iC_{i+1}\). Eliminating \(C_{i+1}\) and \(C_{i-1}\) from \eqref{eq:C} with the help of the right-tail relation \eqref{eq:right-ratio-relation} shows that \eqref{eq:C} is equivalent to the pair of first-order recursions
\begin{equation}
\label{eq:tail:C-factorised}
C_i=r_i\,C_{i+1}+Z_i,
\qquad
Z_{i-1}=\kappa_i\,Z_i+\frac{(2-a)\,r_{i-1}^2}{d \omega_{i-1/2}},
\qquad
\kappa_i:=k_i\,\frac{\omega_{i+1/2}}{\omega_{i-1/2}}\,r_{i-1}^2\,r_i .
\end{equation}
Both have positive coefficients and positive forcing, so by \eqref{eq:tail:periodic-recursion:sol} it suffices to check that the two products over a period are less than one. For the first, \(\prod_jr_j=\Lambda_R<1\). For the second, we start by computing \(\prod_j\kappa_j=\bigl(\prod_jk_j\bigr)\Lambda_R^3\). The right-tail relation yields
\[
k_i\frac{\omega_{i+1/2}}{\omega_{i-1/2}}(1-r_i)=\frac{1-r_{i-1}}{r_{i-1}}-\frac{1-a}{d \omega_{i-1/2}}<\frac{1-r_{i-1}}{r_{i-1}}.
\]
Multiplying this over a period and canceling the telescoping factors \((1-r_j)\) we obtain 
\(\prod_jk_j<1/\Lambda_R\) and hence \(\prod_j\kappa_j<\Lambda_R^2<1\). Thus the sequence \(\mathbf{Z} =(Z_i)_{i\in\mathbb Z}\) is uniquely determined and positive, and then so is \(\mathbf{C}\). Finally, \(C_i-r_iC_{i+1}=Z_i>0\), proving the ordering property.
\end{proof}

In order to ensure that $\mathbf{u}^R$ in~\eqref{eq:UR} is strictly positive,
we need to inspect the sign of quadratic
polynomials \(q\mapsto1-q+Cq^2\). In particular, 
writing
\begin{equation}
\label{eq:qstar}
q^\ast(C):=\begin{cases}\dfrac{1-\sqrt{1-4C}}{2C}, & 0<C\le\tfrac14,\\[2mm] 2, & C>\tfrac14,\end{cases}
\end{equation}
we see that $q^*(C) \in (1,2]$ and that
$1 - q + C q^2 > 0$ for $0 \le q < q^\ast(C)$.
Recalling \eqref{eq:def:c:def:q:i}, we now introduce the \emph{positivity threshold}
\begin{equation}
\label{eq:wstar}
w_s^*:=\min_{\dInt{j}}\frac{q^\ast(C_{s+j})}{\mathcal R_{s+j;s+1}} .
\end{equation}

\begin{corollary}
\label{cor:tail:pos:u:r}
We have \(w_s^*>1\), 
together with $u^R_i(w)>0$ for all $i\ge s+1 $
and all \(0<w<w_s^*\).
\end{corollary}

\begin{proof}
Since \(q^\ast>1\) and \(\mathcal R_{s+j;s+1}\le1\), each term in \eqref{eq:wstar} exceeds \(1\), so \(w_s^*>1\). Now let \(0<w<w_*^s\) and \(i\ge s+1\). Write \(i=s+j+np\) with \(\dInt{j}\) and \(n\in\mathbb N_0\); then \(C_i=C_{s+j}\) and \(q_i=w\mathcal R_{i;s+1}=w\Lambda_R^n\mathcal R_{s+j;s+1}\le w\mathcal R_{s+j;s+1}<q^\ast(C_{s+j})=q^\ast(C_i)\), using \(0<\Lambda_R<1\) and the definition of \(w_*^s\). Hence \(u^R_i=1-q_i+C_iq_i^2>0\) by the defining property of \(q^\ast(C)\) in~\eqref{eq:qstar}.
\end{proof}

\begin{remark}
The threshold \(w_*^s\) is sharp whenever the minimum in \eqref{eq:wstar} is attained at an index \(j\) with \(C_{s+j}\le\frac14\): then \(u^R_{s+j}(w_s^*)=0\). If \(C_{s+j}>\frac14\) for the minimizing \(j\), the tail value at that site is positive for all \(w>0\) and the cap \(q^\ast=2\) is merely a convenient finite bound.
\end{remark}

We now turn to the upper bound $\mathbf{u}^R < 1$.
Inspecting \eqref{eq:UR}, it suffices to ensure that $C_i q_i < 1$
for all $i \ge s+1$. The ordering property in Proposition~\ref{prop:C} allows us to restrict attention to the boundary case $i=s+1$.

The upper bound $\mathbf{u}^R<1$ follows immediately.

\begin{corollary}
\label{cor:tail:ci:qi}
Pick \(w>0\). Then for any $j \ge s+1$, we have
\begin{equation}
    C_{j+1} q_{j+1} < C_j q_j.
\end{equation}
In particular,
if \(C_{s+1}q_{s+1}<1\), then \(C_iq_i<1\) and hence $u^R_{i} < 1$ for all \(i\ge s+1\). 
\end{corollary}
\begin{proof}
Applying Proposition \ref{prop:C},
we obtain
\(C_{j+1}q_{j+1}=C_{j+1}r_jq_j<C_jq_j\).
\end{proof}

We now turn our attention to the sub-solution residual of the right tail. It is interesting to note that the product $C_i q_i$ reappears here. Indeed, we can achieve $F_i(\mathbf{u}^R) >0$ by ensuring that $C_i q_i <1$.

\begin{lemma}[Right-tail correction]
\label{lemma:Ci}
For every \(i\ge s+2\), the corrected right tail \eqref{eq:UR} satisfies
\begin{equation}
\label{eq:tail:comp:res:u:r}
    F_i(\mathbf{u}^R)
    =
    q_i^3
    \left[
        (1-C_iq_i)^3
        +(2-a)C_i(2-C_iq_i)
    \right].
\end{equation}
\end{lemma}

\begin{proof}
Pick \(i\geq s+2\). Substituting \eqref{eq:UR} into \(F_i\) and using \eqref{eq:right-ratio-relation}, we obtain
\begin{align*}
    F_i(\mathbf{u}^R)
={}&
k_id\,\omega_{i+1/2}
\Bigg[
    \left(
        1-r_iq_i+C_{i+1}r_i^2q_i^2
    \right)
    -
    \left(
        1-q_i+C_iq_i^2
    \right)
\Bigg]
\\
&+
d\,\omega_{i-1/2}
\Bigg[
    \left(
        1-\frac{q_i}{r_{i-1}}
        +\frac{C_{i-1}}{r_{i-1}^2}q_i^2
    \right)
    -
    \left(
        1-q_i+C_iq_i^2
    \right)
\Bigg]
\\
&+
\left(1-q_i+C_iq_i^2\right)
\left(q_i-C_iq_i^2\right)
\left(1-a-q_i+C_iq_i^2\right).
\end{align*}
Inserting \eqref{eq:C}, we find
\begin{align*}
F_i(\mathbf{u}^R)
={}&
\left[
    1+2(2-a)C_i
\right]q_i^3
-
\left[
    3C_i+(2-a)C_i^2
\right]q_i^4
+3C_i^2q_i^5-C_i^3q_i^6
\\
={}&
q_i^3
\left[
    (1-C_iq_i)^3
    +(2-a)C_i(2-C_iq_i)
\right],
\end{align*}
as claimed.
\end{proof}

Applying Corollary~\ref{cor:tail:ci:qi} once more, 
it hence suffices to require \(C_{s+1}w<1\) to ensure that the sub-solution residual for $\mathbf{u}^R$ has the correct sign. 
This single threshold, combined with \eqref{eq:wstar}, motivates
the definition
\begin{equation}
\label{eq:Ys0:sharp}
W_s^0:=\min\Bigl\{w_s^*,\ \frac1{C_{s+1}}\Bigr\}.
\end{equation}

\begin{corollary}
\label{cor:tail:res:u:r}
If \(0<w<W_s^0\), then
\begin{equation}
\label{eq:q-bounds}
0<C_iq_i<1,\qquad 0<q_i\bigl(1-C_iq_i\bigr)<1,\qquad i\ge s+1 .
\end{equation}
Consequently \(u^R_i\in(0,1)\) for every \(i\ge s+1\) and \(F_i(\mathbf{u}^R)>0\) for every \(i\ge s+2\).
\end{corollary}

\begin{proof}
From \(w<1/C_{s+1}\) and Corollary~\ref{cor:tail:ci:qi}, \(0<C_iq_i<1\) for all \(i\ge s+1\); from \(w<w_s^*\) and Corollary~\ref{cor:tail:pos:u:r}, \(0<q_i(1-C_iq_i)=1-u^R_i<1\) for all \(i\ge s+1\). This yields \eqref{eq:q-bounds}, and the residual sign follows from \eqref{eq:tail:comp:res:u:r}, both factors there being positive when \(C_iq_i<1\).
\end{proof}

We now turn to the left tail $\mathbf{u}^L$ from \eqref{eq:UL}.
A short computation shows that it already satisfies the sub-solution residual inequalities away from the interface
without the need for corrections.

\begin{lemma}
\label{lemma:tail:res:u:l}
Pick $0 < v < 1$. Then
 the left tail \eqref{eq:UL} satisfies $0 < u_i^L \le v < 1$
for every $ i \le s$.
In addition, for every \(i\le s-1\), the residual satisfies
\[
    F_i(\mathbf{u}^L)=(u_i^L)^2(1+a-u_i^L)>0.
\]
\end{lemma}

\begin{proof}
The inclusion $u_i^L \in (0, v]$ follows directly
from Lemma~\ref{lemma:products}\ref{lemma:products-2}.
By \eqref{eq:left-ratio-relation}, we have
\begin{align*}
F_i(\mathbf{u}^L)&=u_i^L\left[k_id\,\omega_{i+1/2}\left(\frac1{\ell_{i+1}}-1\right)+d\,\omega_{i-1/2}(\ell_i-1)\right]+g(u_i^L;a)\\
&=a u_i^L+g(u_i^L;a)\\
&=(u_i^L)^2(1+a-u_i^L)>0.
\end{align*}
which is strictly positive.
\end{proof}

We conclude this section by summarizing the properties of the two tails $\mathbf{u}^L$ and $\mathbf{u}^R$ that will be needed for the construction of the sub-solution $\underline{\mathbf{u}}^{(s)}(v,w)$ from~\eqref{eq:mr:def:u:i:s}.

\begin{corollary}[Properties of the tails]
\label{prop:tail-geometry}
Fix an interface phase \(\dInt{s}\), let \(\mathbf C=\langle C_1,\dots,C_p\rangle\) be the positive \(p\)-periodic solution of \eqref{eq:C}, and let \(\mathbf{u}^L\) and \(\mathbf{u}^R\) be the tails \eqref{eq:UL} and \eqref{eq:UR} with amplitudes
\[
0<v<1,\qquad 0<w<W_s^0 .
\]
Then the following statements hold.
\begin{enumerate}[label=(\roman*), ref=(\roman*)]
\item \label{prop:tail-geometry-1}
We have the bounds \(0<u_i^L<1\) for all \(i\le s\), and \(0<u_i^R<1\) for all \(i\ge s+1\).
\item \label{prop:tail-geometry-2}
We have the limits \(u_i^L\to0\) as \(i\to-\infty\), and \(u_i^R\to1\) as \(i\to+\infty\).
\item \label{prop:tail-geometry-3}
We have the residual inequalities
\[
F_i(\mathbf{u}^L)>0 \hbox{ for } i\le s-1,
\qquad \qquad 
F_i(\mathbf{u}^R)>0 \hbox{ for } i\ge s+2.
\]
\end{enumerate}
\end{corollary}

\begin{proof}
Item~\ref{prop:tail-geometry-1} is Lemma~\ref{lemma:tail:res:u:l} for the left tail and Corollary~\ref{cor:tail:res:u:r} for the right tail. For item~\ref{prop:tail-geometry-2}, the one-period contractions of Lemma~\ref{lemma:products} give \(\mathcal L_{i;s}\to0\) as \(i\to-\infty\) and \(\mathcal R_{i;s+1}\to0\) as \(i\to+\infty\); hence \(u_i^L=v\mathcal L_{i;s}\to0\), and \(q_i=w\mathcal R_{i;s+1}\to0\) gives \(1-u^R_i=q_i(1-C_iq_i)\to0\). Item~\ref{prop:tail-geometry-3} follows from the residual identity of Lemma~\ref{lemma:tail:res:u:l} for \(i\le s-1\) and the identity \eqref{eq:tail:comp:res:u:r} of Lemma~\ref{lemma:Ci} for \(i\ge s+2\). 
\end{proof}

%%%%%%%%%%%%%%%%%%%%%%%%%%%%%%%%%%%%%%%%%%%%%%%%%%%%%%%%%%%%%%%%%%%
%%%%%%%%%%%%%%%%%%%%%%%%%%%%%%%%%%%%%%%%%%%%%%%%%%%%%%%%%%%%%%%%%%%
% EXISTENCE OF STATIONARY SUB-SOLUTIONS
%%%%%%%%%%%%%%%%%%%%%%%%%%%%%%%%%%%%%%%%%%%%%%%%%%%%%%%%%%%%%%%%%%%
%%%%%%%%%%%%%%%%%%%%%%%%%%%%%%%%%%%%%%%%%%%%%%%%%%%%%%%%%%%%%%%%%%%

\section{Existence of stationary sub-solutions}
\label{sec:subsolution}

\begin{figure}[t]
\centering
\includegraphics[width=1\linewidth]{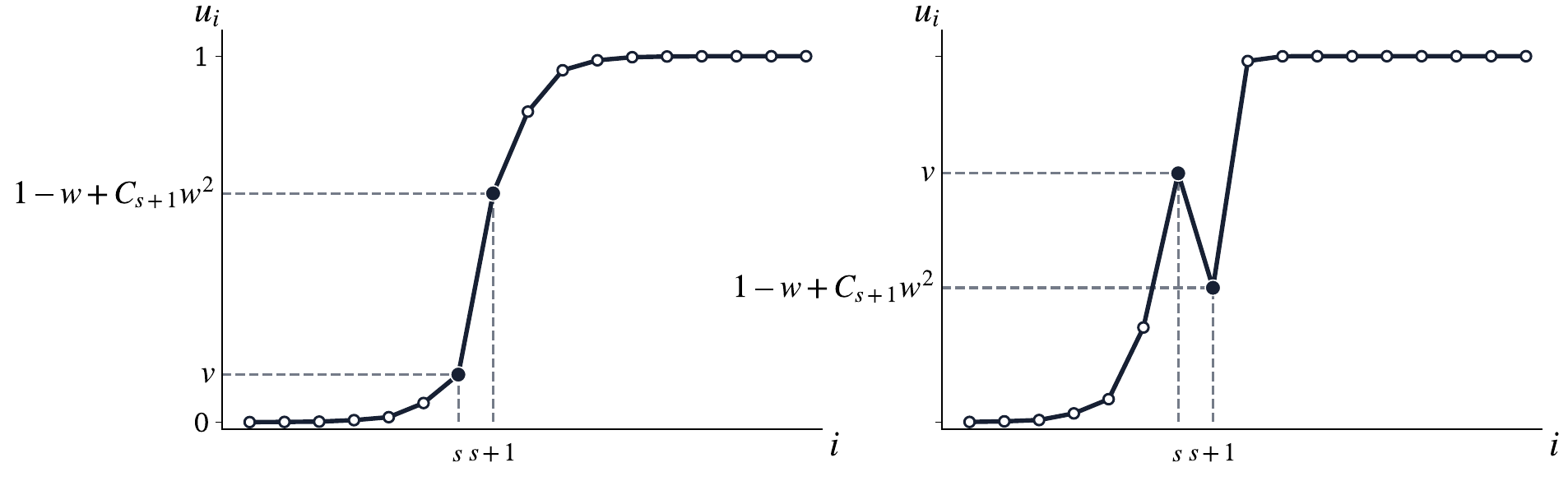}
\caption{Sub-solution profiles $\underline{\mathbf{u}}^{(s)}(v,w)$ from \eqref{eq:mr:def:u:i:s} for \(\mathbf k=\langle 1,2\rangle\) and $\boldsymbol{\omega}=\langle 1 \rangle$. The left panel  uses \(a=0.56\) and \(d=0.25\), with \(v=0.13\) and \(w=0.58\). The right panel uses \(a=0.05\), \(d=0.02\), \(v=0.68\), and \(w=0.65\).}
\label{fig:profile}
\end{figure}

In this section, we  establish the existence of stationary sub-solutions for the periodic diffusion system \eqref{e:main} with the cubic nonlinearity. More precisely, we fix a phase \(\dInt{s}\) and study the sub-solutions $\underline u_i^{(s)}$ introduced in \eqref{eq:mr:def:u:i:s} which combine the left and right tails $\mathbf{u}^L,~\mathbf{u}^R$ constructed in~\S\ref{sec:construction}; see Figure~\ref{fig:profile}.
In view of the preparatory work above, our main task here is to analyze the sub-solution residuals $F_s$ and $F_{s+1}$ at the two interface sites $s$ and $s+1$. During our analysis, we show that the condition $w < 1/C_{s+1}$ is an automatic consequence of the requirement $F_{s+1} \ge 0$.

\paragraph{Residual at $s+1$}

The tail construction incorporates the residual inequalities away from the gluing interfaces $s$ and $s+1$. For the site \(i=s+1\), the residual is

\begin{equation}
\label{eq:F-s+1}
F_{s+1}(v;w)=k_{s+1}d\,\omega_{s+3/2}\bigl(u_{s+2}^R(w)-u_{s+1}^R(w)\bigr)+d\,\omega_{s+1/2}\big(v-u_{s+1}^R(w)\big)+g(u_{s+1}^R(w);a).
\end{equation}
In order to ensure the sign of the residual, we introduce an auxiliary function
\begin{equation}
\label{eq:Vs}
V_s(w):=u_{s+1}^R(w)-\frac{k_{s+1}\omega_{s+3/2}}{\omega_{s+1/2}}\big(u_{s+2}^R(w)-u_{s+1}^R(w)\big)-\frac{g(u_{s+1}^R(w);a)}{d\,\omega_{s+1/2}}.
\end{equation}
%\begin{equation}\label{eq:F-s+1}F_{s+1}(v;w)=d_{s+1}\left(k_{s+1}U_{s,2}(w)-(k_{s+1}+1)U_{s,1}(w)+v\right)+g(U_{s,1}(w);a).\end{equation}Define\begin{equation}\label{eq:Vs}V_s(w):=(k_{s+1}+1)U_{s,1}(w)-k_{s+1}U_{s,2}(w)-\frac{g(U_{s,1}(w);a)}{d_{s+1}}.\end{equation}
%Then $F_{s+1}(v;w)\ge0$ is equivalent to $v\ge V_s(w)$. 
Using \eqref{eq:right-ratio-relation} and \eqref{eq:C} at \(i=s+1\), we can write
\begin{equation}
\label{eq:Vs-polynomial}
\begin{aligned}
V_s(w)
={}&
1-\frac{w}{r_s}+\frac{C_s}{r_s^2}w^2
%\\
%&
-
\frac{w^3}{d\,\omega_{s+1/2}}
\left[
(1-C_{s+1}w)^3
+(2-a)C_{s+1}(2-C_{s+1}w)
\right].
\end{aligned}
\end{equation}
Thus \(V_s\) is a polynomial of degree six in \(w\), with \(V_s(0)=1\) and leading coefficient \(C_{s+1}^3/d\,\omega_{s+1/2}>0\).

\begin{corollary}
\label{cor:sub:res:sp1}
For any pair $(v,w) \in \mathbb R^2$, we have
$F_{s+1}(v,w) \ge 0$ if and only if $v \ge V_s(w)$.
\end{corollary}
\begin{proof}
This follows directly from the computations above,
since $F_{s+1}(v;w) = d\,\omega_{s+1/2} \big(v - V_s(w) \big)$.
\end{proof}

We now show that the restriction $w < 1/C_{s+1}$ encoded in \eqref{eq:Ys0:sharp} can be absorbed by demanding $F_{s+1} \ge 0$.
In particular, the geometric restriction $v \in (0,1)$
is incompatible with this residual requirement whenever $w \ge 1/C_{s+1}$.

\begin{lemma}
\label{lemma:beyond-cap}
Pick any \(w\ge1/C_{s+1}\). Then we have
\begin{equation}
\label{eq:tail:mono:beyond:w:starr}
u^R_{s+2}(w)<u^R_{s+1}(w)
\qquad\text{and hence}\qquad
V_s(w)>u^R_{s+1}(w)\ge1 .
\end{equation}
In particular, no \(v\in(0,1)\) satisfies \(F_{s+1}(v;w)\ge0\).
\end{lemma}

\begin{proof}
Writing \(u^R_{s+1}(w)=1-w+C_{s+1}w^2\) and \(u^R_{s+2}(w)=1-r_{s+1} w+C_{s+2}r_{s+1}^2w^2\), we find
\[
u^R_{s+1}(w)-u^R_{s+2}(w)=w\Bigl[\bigl(C_{s+1}-C_{s+2}r_{s+1}^2\bigr)w-(1-r_{s+1})\Bigr].
\]
Since \(C_{s+2}r_{s+1}<C_{s+1}\) (Proposition~\ref{prop:C}) and \(r_{s+1}<1\), the term in the bracket is increasing in \(w\).
Writing $w_* = 1/C_{s+1}$,
a direct computation yields
\begin{equation}
u^R_{s+1}(w_*)-u^R_{s+2}(w_*)=w_* r_{s+1}\bigl(1-\tfrac{C_{s+2}r_{s+1}}{C_{s+1}}\bigr)>0.    
\end{equation}
Together, this establishes the first inequality in \eqref{eq:tail:mono:beyond:w:starr}. 
For the second, observe that \(u^R_{s+1}(w)\ge1\) yields \(g(u^R_{s+1}(w);a)\le0\), so inspecting \eqref{eq:Vs} leads to
\[
V_s(w)\ge u^R_{s+1}(w)+\frac{k_{s+1}\omega_{s+3/2}}{\omega_{s+1/2}}\bigl(u^R_{s+1}(w)-u^R_{s+2}(w)\bigr)>u^R_{s+1}(w)\ge1.
\]
\end{proof}

\paragraph{Residual at $s$}

It remains to impose the residual inequality at \(i=s\). For the cubic nonlinearity \eqref{eq:cubic}, this residual is

\begin{align}
F_s(v;w)
&:=
F_s\bigl(\underline{\mathbf{u}}^{(s)}(v,w)\bigr)\notag
\\
&=k_sd\,\omega_{s+1/2}\big(u_{s+1}^R(w)-v\big)+d\,\omega_{s-1/2}\big(\ell_sv-v\big)+g(v;a)\label{eq:F-s}
\\
&=
-v^3+(1+a)v^2
-\bigl[a+k_sd\,\omega_{s+1/2}+d\,\omega_{s-1/2}(1-\ell_s)\bigr]v
+d\,\omega_{s+1/2}k_s(1-w+C_{s+1}w^2).\notag
\end{align}

In order to phrase this in a convenient fashion,
we introduce the polynomials
\begin{equation}
    h_s(v):=v^3-(1+a)v^2+\big(a+ d\,\omega_{s+1/2} k_s +d\,\omega_{s-1/2}( 1- \ell_s) \big)v,
    \qquad
    \theta_s(w) = d\,\omega_{s+1/2}k_s(1-w+C_{s+1}w^2),
\end{equation}
which are of degree three and two, respectively.

\begin{corollary}
\label{cor:sub:res:s}
For any pair $(v,w) \in \mathbb R^2$, we have
$F_s(v,w) \ge 0$ if and only if $h_s(v) \le \theta_s(w)$.
\end{corollary}
\begin{proof}
This follows directly from the computations above, 
which show that \(F_s(v;w)=\theta_s(w)-h_s(v)\).

\end{proof}

\paragraph{Admissible set}

\begin{figure}[t]
    \centering
    \includegraphics[width=1\linewidth]{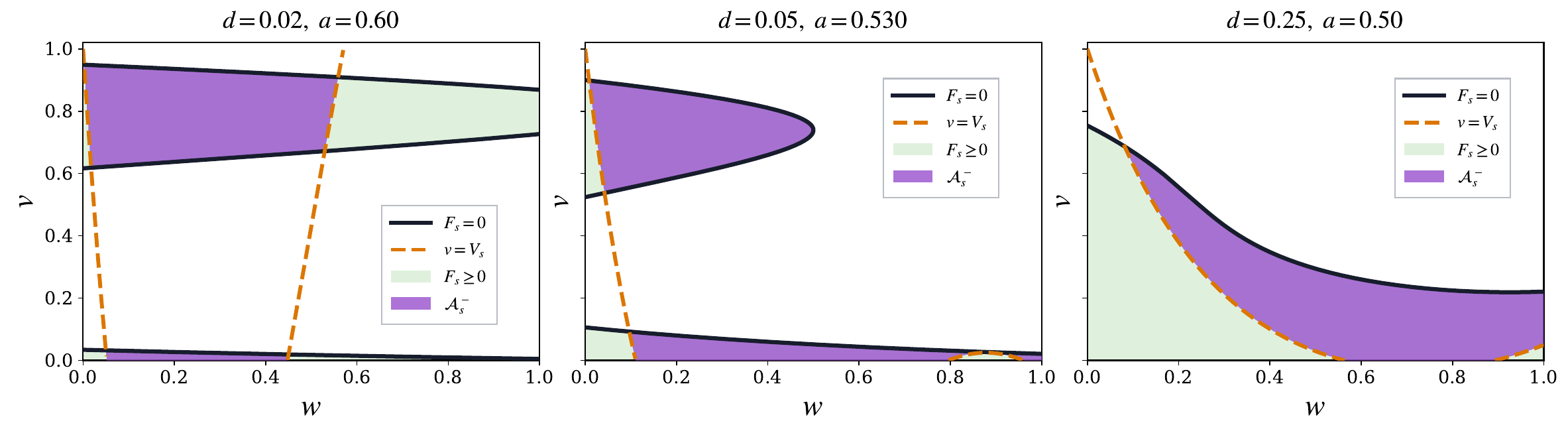}
    \caption{Admissible parameter regions \(\mathcal A_s^-(a,d;\boldsymbol{\omega},\mathbf k)\) defined in \eqref{eq:As-compact}, with \(w\) on the horizontal axis and \(v\) on the vertical axis for \(\mathbf k=\langle1,2\rangle\), \(s=1\), and \(\boldsymbol{\omega} = \langle 1 \rangle\). The black curves are all branches of \(F_s(v;w)=0\), the orange dashed curves are \(v=V_s(w)\). The green regions satisfy \(F_s(v;w)\geq0\), and the purple regions represent \(\mathcal A_s^-\).}
    \label{fig:As}
\end{figure}

\begin{figure}[t]
    \centering
    \includegraphics[width=0.8\linewidth]{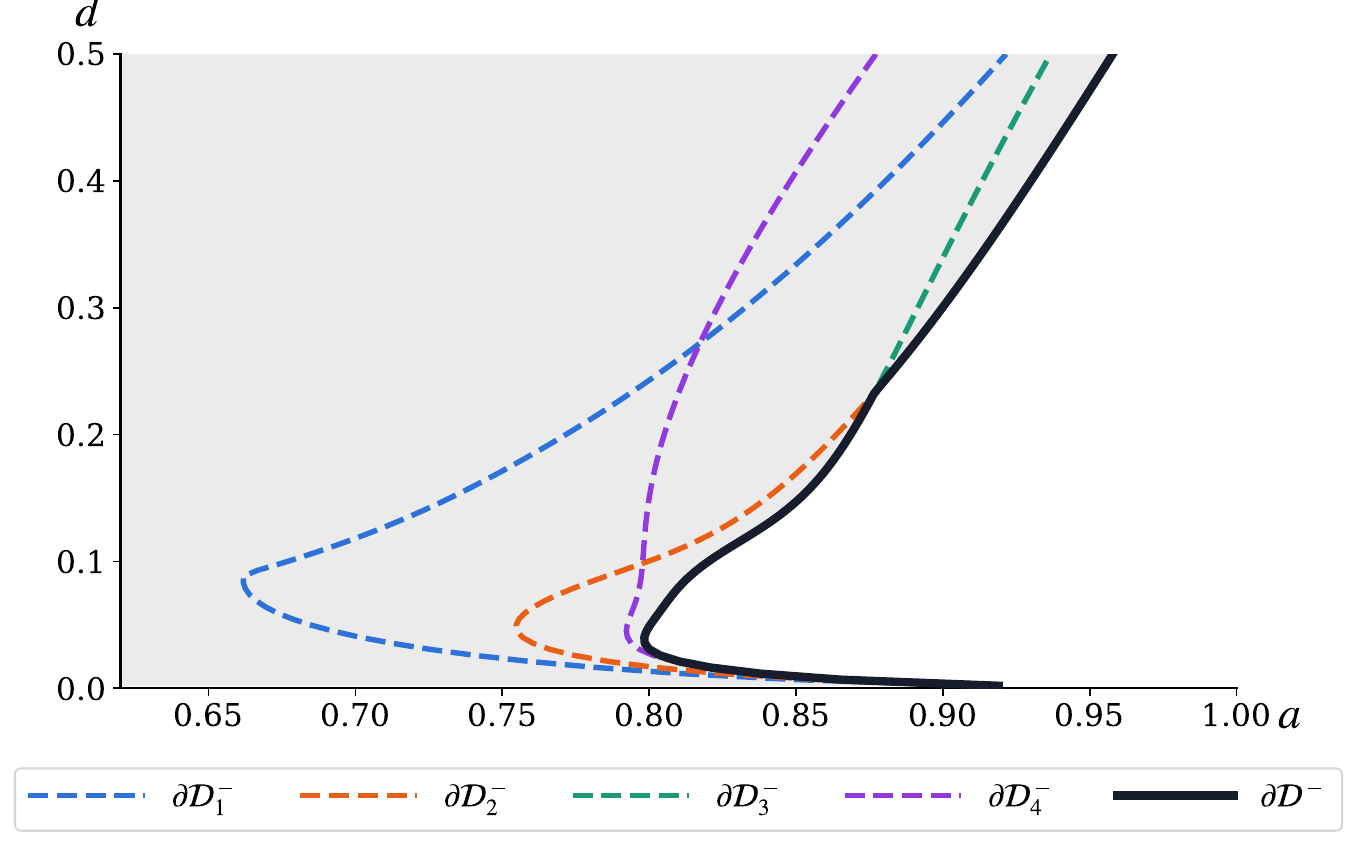}
    \caption{Numerically computed phase-dependent boundaries of the regions \(\mathcal D_s^-(\boldsymbol{\omega},\mathbf k)\) defined in \eqref{eq:Ds-minus}, for \(\mathbf k=\langle 1,2,3,4\rangle\) and $\boldsymbol{\omega} = \langle 1 \rangle$. The dashed curves show the right portions of \(\partial\mathcal D_s^-(\boldsymbol{\omega},\mathbf k)\), \(s\in\{1,2,3,4\}\). At each sampled diffusion value \(d\), the solid black curve selects the largest \(a\) among the four phase-dependent boundaries. The gray shading represents the numerical approximation of the combined region \(\mathcal D^-(\boldsymbol{\omega},\mathbf k)\) defined in \eqref{eq:D-minus} within the displayed window.}
    \label{fig:ad}
\end{figure}

We are now ready to put everything together. 
In view of Lemma \ref{lemma:beyond-cap}
the geometric conditions from Corollary \ref{prop:tail-geometry}
can be reduced to $(v,w) \in (0,1) \times (0, w_s^*)$.
Combining Corollaries~\ref{cor:sub:res:sp1} and \ref{cor:sub:res:s}
leads naturally to the admissible set
\begin{equation}
\label{eq:As-compact}
\mathcal A_s^-(a,d;\boldsymbol{\omega},\mathbf k)
=\Bigl\{(v,w)\in(0,1)\times(0, w_s^*):\ \ h_s(v)\le  \theta_s(w),\ \ V_s(w)\le v\Bigr\}.
\end{equation}
The admissible set is illustrated in Figure~\ref{fig:As}.
The parameter sets \(\mathcal D_s^-(\boldsymbol{\omega}, \mathbf k)\) defined in \eqref{eq:Ds-minus} then select the pairs $(a,d)$ for which $\mathcal A_s^-$ is non-empty. Taking their union over the \(p\) phases gives \(\mathcal D^-(\boldsymbol{\omega},\mathbf k)\) defined in \eqref{eq:D-minus}. Figure \ref{fig:ad} illustrates \eqref{eq:D-minus} schematically.

\begin{proposition}
\label{prop:subsolution}
Let \(a\in(0,1)\) and $d > 0$. Fix an interface phase \(\dInt{s}\). If \(\mathcal A_s^-(a,d;\boldsymbol{\omega},\mathbf k)\neq\varnothing\), then every pair \((v,w)\) in this set defines a global stationary sub-solution \(\underline{\mathbf{u}}^{(s)}(v,w)\) satisfying
\[
0<\underline u_i^{(s)}(v,w)<1, \qquad F_i\bigl(\underline{\mathbf{u}}^{(s)}(v,w)\bigr)\ge0,
\qquad i\in\mathbb Z,
\]
and
\[
\lim_{i\to-\infty}\underline u_i^{(s)}(v,w)=0,
\qquad
\lim_{i\to+\infty}\underline u_i^{(s)}(v,w)=1.
\]
\end{proposition}
\begin{proof}
Choose \((v,w)\in\mathcal A_s^-(a,d;\boldsymbol{\omega},\mathbf k)\). 
Since $V_s(w) \le v < 1$ we have \(0<w<W_s^0\)
by Lemma~\ref{lemma:beyond-cap}.
In particular, Corollary \ref{prop:tail-geometry}
provides the value bounds together with the far-field limits and shows that $F_i(\underline{\mathbf{u}}^{(s)}(v,w) ) \ge 0$ for $i \in \mathbb Z \setminus \{s, s+1\}$. Corollaries~\ref{cor:sub:res:sp1} and \ref{cor:sub:res:s} provide this inequality at the interface sites $i \in \{s, s+1\}$. Thus the profile is a global stationary sub-solution.
\end{proof}

\section{Proof of Theorem~\ref{thm:pinning}}\label{sec:speed-sign}
In this section, we use the properties of sub-solutions from \S\ref{sec:construction},\ref{sec:subsolution} and prove Theorem~\ref{thm:pinning} for the cubic nonlinearity. We first establish non-emptiness of the admissible sets for sufficiently small diffusion.

\begin{lemma}
\label{lemma:nonempty-admissible}
There exists \(d_0>0\) such that for every $0 < d < d_0$ and every phase \(\dInt{s}\), we have
\begin{equation}
\label{eq:general-fixed-a-admissible}
    \left(\frac23,\frac14\right)
    \in\mathcal A_s^-
    \left(\frac12,d;\boldsymbol{\omega},\mathbf k\right).
\end{equation}
Consequently,
\(\mathcal D^-(\boldsymbol{\omega},\mathbf k)\cap\mathcal D^+(\boldsymbol{\omega},\mathbf k) \neq\varnothing\).
\end{lemma}

\begin{proof}
Fix \(a=\frac12\), an interface phase \(\dInt{s}\), and \(v_*=\frac23\), \(w_*=\frac14\). By \eqref{eq:left-ratio-relation}, \eqref{eq:right-ratio-relation}, and Lemma~\ref{lemma:products}\ref{lemma:products-2},
\[
    0<\ell_{i+1}\le2k_i d\,\omega_{i+1/2},
    \qquad
    0<r_{i-1}\le2d\,\omega_{i-1/2}.
\]
In particular, $\ell_i,r_i$ tend uniformly to zero as \(d\to0\). Multiplying \eqref{eq:C} by \(r_{i-1}^2/d\,\omega_{i-1/2}\), and using the periodicity and positivity of \(C_i\) from Proposition~\ref{prop:C} gives \(C_i\to0\) as $d\to 0$. Inspecting \eqref{eq:Ys0:sharp} then gives \(W_s^0>1\) for sufficiently small \(d\).

By \eqref{eq:UR}, \(u_{s+1}^R(w_*)\to3/4\) and \(u_{s+2}^R(w_*)\to1\). Formulas \eqref{eq:F-s} and \eqref{eq:F-s+1} therefore yield
\[
\begin{aligned}
    F_s(v_*;w_*)&\longrightarrow
    g\left(\frac23;\frac12\right)=\frac1{27}>0,
    \\
    F_{s+1}(v_*;w_*)&\longrightarrow
    g\left(\frac34;\frac12\right)=\frac3{64}>0.
\end{aligned}
\]
Hence \(0<w_*<W_s^0\), \(F_s(v_*;w_*)>0\) and \(F_{s+1}(v_*;w_*)>0\) when \(d\) is sufficiently small. Inspecting definition~\eqref{eq:As-compact} now yields
\eqref{eq:general-fixed-a-admissible}. Applying this conclusion first to \((\boldsymbol{\omega},\mathbf k)\), and then to the dual pair \((\widetilde{\boldsymbol{\omega}},\widetilde{\mathbf k})\) from \eqref{eq:mr:refl:dual:transf}, gives
\[
    \left(\frac12,d\right)\in\mathcal D^-(\boldsymbol{\omega},\mathbf k)\cap\mathcal D^+(\boldsymbol{\omega},\mathbf k)
\]
for all sufficiently small \(d>0\).
\end{proof}

We now establish the crucial fact that an admissible construction excludes positive wave speeds. The ordering of the exponents in Lemma \ref{lemma:exponent-gap} is a key ingredient, as it allows us to compare the tails of traveling waves with those of our stationary sub-solution.

\begin{proposition}
\label{prop:admissible-speed}
Fix \(a\in(0,1)\) together with $d > 0$ and consider a traveling wave $(c, \boldsymbol{\Phi})$  
as constructed in Proposition~\ref{prop:periodic-wave-existence}. 

If \(\mathcal A_s^-(a,d;\boldsymbol{\omega},\mathbf k)\neq\varnothing\) for some phase \(\dInt{s}\), then \(c\le0\).
\end{proposition}

\begin{proof}
Write $\mathbf{u}^{\mathrm{tw}}$ for the associated 
traveling wave solution to \eqref{e:main}
defined via \eqref{eq:wave-ansatz}
and
choose \((v,w)\in\mathcal A_s^-(a,d;\boldsymbol{\omega},\mathbf k)\).
Proposition~\ref{prop:subsolution} gives a global stationary sub-solution with values in \((0,1)\) and limits \(0\) and \(1\) at the two ends. Suppose that \(c>0\). The tail construction, Lemma~\ref{lemma:exponent-gap}, and Proposition~\ref{prop:CGW-tail} give
\[
    \frac{\underline u_i^{(s)}}{u^{\rm tw}_i(0)}
    =O\!\left(e^{(\lambda_L-\lambda_L^c)i}\right)
    \longrightarrow0
    \qquad\text{as }i\to-\infty.
\]
On the right,
\[
    1-\underline u_i^{(s)}
    =w\mathcal R_{i;s+1}
      \bigl(1-C_iw\mathcal R_{i;s+1}\bigr)
    \sim w\mathcal R_{i;s+1},
\]
and therefore
\[
    \frac{1-u^{\rm tw}_i(0)}{1-\underline u_i^{(s)}}
    =O\!\left(e^{(\lambda_R^c-\lambda_R)i}\right)
    \longrightarrow0
    \qquad\text{as }i\to+\infty.
\]
Thus there exists \(N\ge1\) such that
\[
    \underline u_i^{(s)}\le u^{\rm tw}_i(0),
    \qquad |i|\ge N.
\]
The comparison argument of Corollary~\ref{cor:speed} now yields a contradiction, showing that indeed \(c\le0\).
\end{proof}

We are now ready to prove our second main result, Theorem~\ref{thm:pinning}.

\begin{proof}[Proof of Theorem~\ref{thm:pinning}]
Lemma~\ref{lemma:nonempty-admissible} proves item~\ref{thm:pinning-1} and the non-emptiness assertion in item~\ref{thm:pinning-3}. 

If \((a,d)\in\mathcal D^-(\boldsymbol{\omega}, \mathbf k)\), then \(\mathcal A_s^-(a,d ;\boldsymbol{\omega},\mathbf k)\neq\varnothing\) for some phase \(\dInt{s}\). Proposition~\ref{prop:admissible-speed} gives \(c(a,d;\boldsymbol{\omega},\mathbf k)\le0\). If \((a,d)\in\mathcal D^+(\boldsymbol{\omega},\mathbf k)\), \eqref{eq:Ds-plus} instead gives
\[
    \mathcal A_s^-
    \bigl(1-a,d; \widetilde{\boldsymbol{\omega}},\widetilde{\mathbf k}\bigr)
    \neq\varnothing
\]
for some phase \(\dInt{s}\). Proposition~\ref{prop:admissible-speed} and Lemma~\ref{lem:mr:reflection} yield \(c(a,d;\boldsymbol{\omega},\mathbf k)\ge0\). This proves item~\ref{thm:pinning-2}, and the pinning assertion in item~\ref{thm:pinning-3} follows from the intersection.

Finally, let \((a,d)\in \overline{\mathcal D^-(\boldsymbol{\omega},\mathbf k)}\cap\mathcal H\). By continuity of the tail construction, compactness of the phase set 
and the bounds $0<v<1$ and $0<w<W_s^0\leq2$, there exists a profile \(\underline{\boldsymbol{u}}\) of \eqref{eq:mr:def:u:i:s} with nonnegative stationary residual at \((a,d)\). Fix \(a'\in(0,a)\). The cubic identity
\[
    g(u;a')-g(u;a)
    =(a-a')u(1-u)>0,
    \qquad 0<u<1,
\]
together with \(g'(0;a')<0,~g'(1;a')<0\), allows us to choose a sufficiently small \(\epsilon\in(0,1)\) such that
\[
    g(u-\epsilon ;a')>g(u;a), \qquad 0\leq u\leq1.
\]
Hence, \(\underline{\boldsymbol{u}}-\langle \epsilon\rangle \) is a stationary sub-solution at \((a',d)\). Applying the comparison principle therefore gives \(c(a',d;\boldsymbol{\omega},\mathbf k)\leq0\). The same argument applies to positive periodic sequence $\boldsymbol{\omega}$. For \((a,d)\in \overline{\mathcal D^+(\boldsymbol{\omega},\mathbf k)}\cap\mathcal H\) and \(a'\in(a,1)\), apply it to the reflected system, noting that \(1-a'<1-a\). Lemma~\ref{lem:mr:reflection} then yields \(c(a',d;\boldsymbol{\omega},\mathbf k)\geq0\), completing the proof of item~\ref{thm:pinning-4}.

\end{proof}

%%%%%%%%%%%%%%%%%%%%%%%%%%%%%%%%%%%%%%%%%%%%%%%%%%%%%%%%%%%%%%%%%%%
%%%%%%%%%%%%%%%%%%%%%%%%%%%%%%%%%%%%%%%%%%%%%%%%%%%%%%%%%%%%%%%%%%%
% SMALL DIFFUSION BOUNDARY EXPANSIONS
%%%%%%%%%%%%%%%%%%%%%%%%%%%%%%%%%%%%%%%%%%%%%%%%%%%%%%%%%%%%%%%%%%%
%%%%%%%%%%%%%%%%%%%%%%%%%%%%%%%%%%%%%%%%%%%%%%%%%%%%%%%%%%%%%%%%%%%

\section{Small diffusion boundary expansions}
\label{sec:estimate}

In this section, in order to gain additional insight into the structure of the sufficient pinning regions $\mathcal{D}^\pm(\boldsymbol{\omega},\mathbf{k})$ from Theorem~\ref{thm:pinning}, we construct explicit curves $(d,a^\pm(d))$ that
(together) describe the edge of these regions in the anti-continuum regime $0 < d \ll 1$. Our approach is motivated by the weak-coupling approximation for propagation failure in the discrete Nagumo equation \cite{erneux1993propagating}, which we adapt to periodic media. 

We first restrict to a corner of the \((a,d)\)-plane near \((a,d)=(1,0)\), where the tail data are under quantitative control and the admissible set splits into at most four cells. Only two of these are critical for the boundary of $\mathcal D^-$, and on them the admissible pairs stay uniformly away from all remaining constraints. These cells can therefore only disappear through a tangency between the curves \(\{F_s=0\}\) and \(\{F_{s+1}=0\}\).

We then solve the resulting tangency system by an expansion in powers of \(\sqrt d\). The key point is that a cell-dependent scaling of the unknowns and of the tail data turns this system into a regular perturbation of an explicit non-degenerate problem. Away from the corner the same tangency conditions can be continued numerically, which captures the boundary of the sufficient pinning region well beyond the reach of the expansions.

\paragraph{Preliminaries}
We consider \(a\) close to \(1\), with the aim of constructing admissible parameter values near the right boundary of \(\mathcal D^-(\boldsymbol{\omega},\mathbf k)\). Fix an interface phase \(\dInt{s}\), and write
\begin{equation}
    \epsilon:=1-a,
    \qquad
    \omega_{\min}:=\min_i\omega_{i+1/2},
    \qquad
    \omega_{\max}:=\max_i\omega_{i+1/2},
    \qquad
    k_{\min}:=\min_i k_i,
    \qquad
    k_{\max}:=\max_i k_i,
\end{equation}
together with
\begin{equation}
\label{eq:M-constant}
    K:=\max_i\bigl(k_i\omega_{i+1/2}+\omega_{i-1/2}\bigr),
    \qquad
    M:=\Bigl(1+\frac{K}{\omega_{\min}}\Bigr)\bigl(1+\sqrt{6K}\bigr)+\tfrac92K .
\end{equation}
For the profile \eqref{eq:mr:def:u:i:s}, the first two right-tail values are
\begin{equation}
\label{eq:asymp:idss:for:usp1}
    u_{s+1}=1-w+C_{s+1}w^2,\qquad u_{s+2}=1-r_{s+1}w+C_{s+2}r_{s+1}^2w^2 .
\end{equation}
Recalling \eqref{eq:residual}, the two interface residuals of the glued profile \(\underline{\mathbf{u}}^{(s)}(v,w)\) take the form
\begin{align}
\label{eq:Fs-edge}
F_s(v;w)&=d\bigl[k_s\omega_{s+1/2}(u_{s+1}-v)+\omega_{s-1/2}(\ell_s-1)v\bigr]+g(v;a)=\theta_s(w)-h_s(v),\\
\label{eq:Fs1-edge}
F_{s+1}(v;w)&=d\bigl[k_{s+1}\omega_{s+3/2}(u_{s+2}-u_{s+1})+\omega_{s+1/2}(v-u_{s+1})\bigr]+g(u_{s+1};a)=d\,\omega_{s+1/2}\bigl(v-V_s(w)\bigr).
\end{align}

Under the natural restriction $u_{s+2} \le 1 $ and $v\le1$,
we see that
\begin{equation}
\label{eq:Fs1-upper}
F_{s+1}(v;w)\ \le\ (1-u_{s+1})\Bigl[d\bigl(k_{s+1}\omega_{s+3/2}+\omega_{s+1/2}\bigr)-u_{s+1}\bigl(a-u_{s+1}\bigr)\Bigr].
\end{equation}

Writing $u_{s+1} = \tau_\pm$ for the two values where
the bracketed term in \eqref{eq:Fs1-upper} vanishes,
we obtain
\begin{equation}
\label{eq:theta-pm}
\tau_\pm:=\tfrac12\bigl(a\pm\sqrt{D_{s+1}}\bigr),
\qquad
D_{s+1}:=a^2-4d\bigl(k_{s+1}\omega_{s+3/2}+\omega_{s+1/2}\bigr).
\end{equation}
In particular, when $D_{s+1} > 0$ we see that
$F_{s+1} < 0$ whenever $u_{s+1} \in (\tau_-, \tau_+)$,
which is inadmissible.
As a final preparation, we provide essential bounds on the tail data \(r_i\), \(C_i\) of the sub-solution profile, which involve the threshold
\begin{equation}
\label{eq:rho0}
\rho_0:=\min\Bigl\{\frac14,\ \frac{\omega_{\min}}{4\max_ik_i\omega_{i+1/2}}\Bigr\}.
\end{equation}

\begin{lemma}
\label{lemma:tail-bounds}
Pick \(\epsilon=1-a\in(0,1)\) and \(d>0\) and write 
\(\rho:=d\,\omega_{\max}/\epsilon\). Then for all 
\(\dInt{i}\)
\begin{equation}
\label{eq:tail-bounds-r}
\frac{d\,\omega_{i+1/2}}{\epsilon+Kd}\ \le\ r_i\ \le\ \frac{d\,\omega_{i+1/2}}{\epsilon+d\,\omega_{i+1/2}}\ \le\ \rho,
\qquad
\ell_{i+1}\le\frac{d\,k_i\omega_{i+1/2}}{a}\le\frac{Kd}{a},
\qquad
C_i\ >\ \frac{(1+\epsilon)\,r_i^2}{d\,\omega_{i+1/2}} .
\end{equation}
Moreover, if \(\rho\le\rho_0\) then we also have
\begin{equation}
\label{eq:tail-bounds-C}
C_i\ \le\ \frac{8(1+\epsilon)\,\rho^2}{3d\,\omega_{\min}}\qquad\text{for all } \dInt{i} .
\end{equation}
\end{lemma}

\begin{proof}
Solving the ratio relations \eqref{eq:right-ratio-relation} and \eqref{eq:left-ratio-relation} for \(r_{i-1}\), resp.\ \(\ell_{i+1}\), gives
\[
r_{i-1}=\frac{d\,\omega_{i-1/2}}{\epsilon+d\,\omega_{i-1/2}+d\,k_i\omega_{i+1/2}(1-r_i)},
\qquad
\ell_{i+1}=\frac{d\,k_i\omega_{i+1/2}}{a+d\,k_i\omega_{i+1/2}+d\,\omega_{i-1/2}(1-\ell_i)} .
\]
Since \(0<r_i<1\), the first denominator is at least \(\epsilon+d\,\omega_{i-1/2}\), which gives the upper bound on \(r_{i-1}\), and at most \(\epsilon+Kd\), which gives the lower bound. Similarly,  the second denominator is at least \(a\) since \(0<\ell_i<1\), which gives the bound on \(\ell_{i+1}\). Writing \eqref{eq:C} with the help of \eqref{eq:right-ratio-relation} as
\begin{equation}
\label{eq:exp:c:rewritten}
d\,k_i\omega_{i+1/2}\,r_i\bigl(r_iC_{i+1}-C_i\bigr)+\frac{d\,\omega_{i-1/2}}{r_{i-1}}\Bigl(\frac{C_{i-1}}{r_{i-1}}-C_i\Bigr)=1+\epsilon
\end{equation}
and using \(r_iC_{i+1}<C_i\) (Proposition~\ref{prop:C}), we obtain the lower bound on \(C_i\). Writing \(\delta_i:=C_i-r_iC_{i+1}>0\), we may also read off the recursion
\[
\delta_{i-1}=\frac{(1+\epsilon)\,r_{i-1}^2}{d\,\omega_{i-1/2}}+\frac{k_i\omega_{i+1/2}}{\omega_{i-1/2}}\,r_ir_{i-1}^2\,\delta_i .
\]
Taking \(r_i\le\rho \le \rho_0 \) turns the recursion into \begin{equation}\delta_{i-1}\le\frac{(1+\epsilon)\rho^2}{d\omega_{\min}}+\lambda\,\delta_i,
\qquad \qquad \hbox{ with }
\lambda:=\frac{\max_ik_i\omega_{i+1/2}}{\omega_{\min}}\rho^3\le\frac{\rho^2}4\le\frac1{64}.
\end{equation}
Taking the maximum over the period, \(\bar\delta:=\max_i\delta_i\) satisfies 
\begin{equation}
\bar\delta\le\frac{(1+\epsilon)\rho^2}{d\omega_{\min}}+\frac1{64}\bar\delta \qquad \hbox{ i.e. }  \qquad \bar\delta\le\frac{64}{63}\frac{(1+\epsilon)\rho^2}{d\omega_{\min}}.
\end{equation}
Finally, 
\(C_i=\delta_i+r_iC_{i+1}\le\bar\delta+\rho\max_jC_j\) gives \(\max_jC_j\le\frac{\bar\delta}{1-\rho}\le\frac43\bar\delta\), which is \eqref{eq:tail-bounds-C}.
\end{proof}

\begin{figure}[t]
\centering
\includegraphics[width=0.75\linewidth]{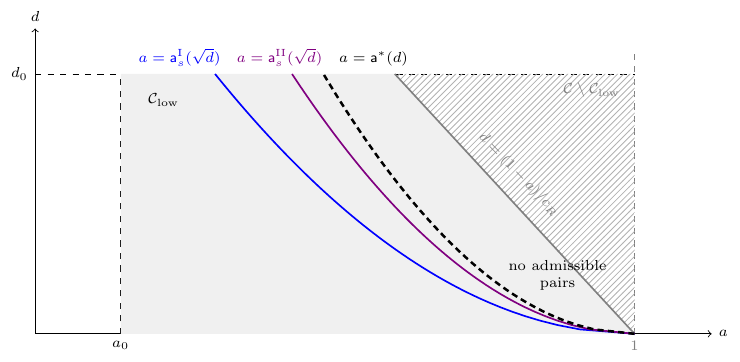}
\caption{Schematic of the corner \(\mathcal C=\{a_0<a<1,\ d<d_0\}\) in the \((a,d)\)-plane, zoomed in near \((a,d)=(1,0)\). In the triangle \(\mathcal C\setminus\mathcal C_{\rm low}\) (hatched) nothing is admissible. The edge curves \(a=\mathsf a_s^{\mathrm I}(\sqrt d)\) and \(a=\mathsf a_s^{\mathrm{II}}(\sqrt d)\) of phase \(s\) describe tangencies of \(\{F_s=0\}\) and \(\{v=V_s(w)\}\). Their upper envelope \(a=\mathsf a^\ast(d)\) (black dashed) is the right boundary of \(\mathcal D^-(\boldsymbol{\omega},\mathbf k)\) in the corner; see Corollary~\ref{cor:envelope}. To its left \(c\le0\) is guaranteed; to its right no admissible pairs exist for any phase in the corner.}
\label{fig:corner}
\end{figure}

\paragraph{The corner.}
With these preparations, we are ready to define the corner area
around $(a,d) = (1, 0)$ where our analysis here is valid.
In particular, we fix  the constants
\begin{equation}
\label{eq:constants}
\begin{aligned}
c_R&:=\max\Bigl\{\frac{\omega_{\max}}{\rho_0},\ 4K,\ 8\omega_{\max},\ \frac{9\,\omega_{\max}^4}{\omega_{\min}^3},\ \frac{86\,\omega_{\max}^2}{\omega_{\min}}\Bigr\},\\[2pt]
d_0&:=\min\Bigl\{\frac{\omega_{\min}}{3c_R^2},\ \frac1{16M^2},\ \frac{k_{\min}\,\omega_{\min}^4}{2^{19}K^2k_{\max}^2\omega_{\max}^3}\Bigr\},\\[2pt]
a_0&:=1-\min\Bigl\{\frac14,\ \frac1{2M}-\sqrt{d_0}\Bigr\},
\end{aligned}
\end{equation}
with \(\rho_0\) from \eqref{eq:rho0} and \(M\) from \eqref{eq:M-constant}. Note that \(a_0\ge\frac34\) and, since \(M\ge\frac92K\), also \(Kd_0<\frac1{16}\) and \(\sqrt{d_0}\le\frac1{4M}\). With these constants we define the corner rectangle $\mathcal{C}$ and a lower region $\mathcal{C}_{\rm low}$ that contains the boundary segments by writing
\begin{equation}
\label{eq:corner}
\mathcal C:=\bigl\{(a,d):\ a_0<a<1,\ 0<d<d_0\bigr\},
\qquad
\mathcal C_{\rm low}:=\bigl\{(a,d)\in\mathcal C:\ 1-a\ge c_R\,d\bigr\},
\end{equation}
see Figure~\ref{fig:corner}. The restriction to this corner will allow us to exploit several crucial inequalities that determine the geometry of the admissible set.

\begin{lemma}
\label{lemma:corner-constants}
For every phase \( 
\dInt{s}
\), the inequalities  
\makeatletter\newcommand{\nameditem}[1]{\item[(#1)]\protected@edef\@currentlabel{(#1)}\label{it:#1}}\makeatother
\begin{enumerate}[leftmargin=*]
\nameditem{S} \(d\,k_s\omega_{s+1/2}<\tfrac4{27}\,a^2\),
\nameditem{M} \(D_{s+1}>0\) and \(\tau_+>\tfrac12\)
\end{enumerate}
hold at every \((a,d)\in\mathcal C\), and the inequalities
\begin{enumerate}[leftmargin=*]
\nameditem{R} \(r_i\le\tfrac14\),\quad \(2\,r_{s+1}C_{s+2}\le C_{s+1}\),\quad \(8\,C_{s+1}(\epsilon+\tau_-)\le1\)
\end{enumerate}
hold at every \((a,d)\in\mathcal C_{\rm low}\).
\end{lemma}

\begin{proof}
\ref{it:S} holds since \(d\,k_s\omega_{s+1/2}\le Kd_0<\frac1{16}<\frac4{27}a_0^2\), and \ref{it:M} holds since \(\tau_-\le2Kd/a\le4Kd_0\le\frac14\) gives \(a>\frac34\ge\frac12+\tau_-\), which is equivalent to \(\tau_+>\frac12\).
Turning to \ref{it:R}, let \((a,d)\in\mathcal C_{\rm low}\), so that \(\epsilon\ge c_Rd\), and put \(\rho:=d\,\omega_{\max}/\epsilon\le\omega_{\max}/c_R\). The definition of \(c_R\) in \eqref{eq:constants} gives \(\rho\le\rho_0\), so that both estimates of Lemma~\ref{lemma:tail-bounds} apply, together with
\begin{equation}
    Kd\le\epsilon/4, \qquad \qquad 
    \epsilon\ge\frac{25\omega_{\max}^4}{3\omega_{\min}^3}d.
\end{equation}
The three inequalities in \ref{it:R} now follow readily:
\begin{itemize}
\item \(r_i\le\rho\le\omega_{\max}/c_R\le\frac14\) by \(c_R\ge4\omega_{\max}\);
\item The lower bounds \eqref{eq:tail-bounds-r} on \(r_{s+1}\) and on \(C_{s+1}\) yield
\[
2\,r_{s+1}C_{s+2}\ \le\ 2\rho\cdot\frac{8(1+\epsilon)\rho^2}{3d\,\omega_{\min}}=\frac{16(1+\epsilon)d^2\omega_{\max}^3}{3\,\epsilon^3\,\omega_{\min}}
\ \le\ \frac{16(1+\epsilon)d\,\omega_{\min}^2}{25\,\epsilon^2\,\omega_{\max}}\ \le\ \frac{(1+\epsilon)r_{s+1}^2}{d\,\omega_{s+3/2}}\ <\ C_{s+1};
\]
\item \(\tau_-\le2Kd/a\le4Kd\le\epsilon\) and \eqref{eq:tail-bounds-C} give \(8C_{s+1}(\epsilon+\tau_-)\le\frac{128(1+\epsilon)\omega_{\max}^2}{3\omega_{\min}}\frac d\epsilon\le\frac{256\omega_{\max}^2}{3\omega_{\min}c_R}\le1\).
\end{itemize}

\end{proof}

\begin{figure}[t]
\centering
\includegraphics[height=0.23\textheight]{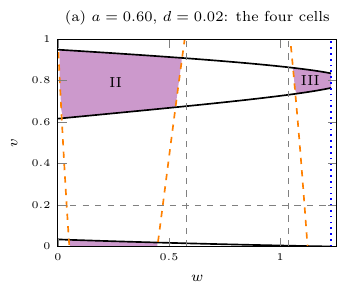}\hfill\includegraphics[height=0.23\textheight]{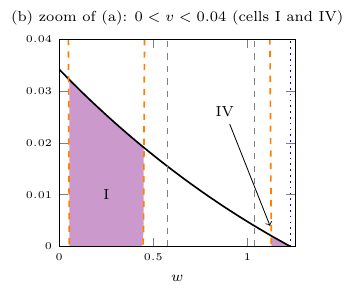}
\caption{The admissible set \(\mathcal A_s^-\) (shaded), computed for \(\mathbf k=\langle1,2\rangle\), \(\omega=\langle 1\rangle\), \(s=1\), \(a=0.60\), \(d=0.02\); (b) zooms into \(0<v<0.04\). Black: the branch \(v_1(w)\) (bottom) and the endpoints \(v_2(w),v_3(w)\) of the upper interval \(\mathcal I_s(w)\) (top) of Lemma~\ref{lemma:zero-set}; orange dashed: \(v=V_s(w)\) (clipped); gray dashed: \(v=a/3\) and the lines \(w=w(\tau_+)\), \(w=w(\tau_-)\); blue dotted: the positivity threshold \(w=w_s^*\). Across the strip \(w(\tau_+)<w<w(\tau_-)\) the lower bound \(V_s\) exceeds \(1\), which is what empties the strip. The cells are marked in the panels: II and III in (a), I and IV in (b).}
\label{fig:cells}
\end{figure}

\paragraph{Splitting the admissible set.}
With the corner in place, we set out to provide a detailed characterization of the admissible set \(\mathcal A_s^-\). The first observation is that two one-line inequalities split \(\mathcal A_s^-\) into four disjoint cells $\mathcal A_s^{\mathrm I}$ through $\mathcal A_s^{\mathrm IV}$ which we refer to as cells I--IV; see Figure~\ref{fig:cells}.

\begin{lemma}
\label{lemma:four-cells}
For \((a,d)\in\mathcal C\) we have the disjoint union
\[
\mathcal A_s^-=\underbrace{\mathcal A_s^{\mathrm I}}_{v<a/3,\ u_{s+1}\ge\tau_+}\ \dot\cup\ \underbrace{\mathcal A_s^{\mathrm{II}}}_{v>a/3,\ u_{s+1}\ge\tau_+}\ \dot\cup\ \underbrace{\mathcal A_s^{\mathrm{III}}}_{v>a/3,\ u_{s+1}\le\tau_-}\ \dot\cup\ \underbrace{\mathcal A_s^{\mathrm{IV}}}_{v<a/3,\ u_{s+1}\le\tau_-},
\]
i.e.\ no pair of \(\mathcal A_s^-\) has \(v=a/3\) or \(\tau_-<u_{s+1}<\tau_+\).
\end{lemma}

\begin{proof}
For \(0<v<a\) we have \(h_s(v)\ge-g(v;a)=v(1-v)(a-v)\), 
so $h_s(a/3) \ge \frac4{27}a^2$. Using \ref{it:S}, we find
\begin{equation}
\label{eq:exp:f:a3:neg}
F_s(a/3;w)\ \le\ d\,k_s\omega_{s+1/2}-\tfrac4{27}a^2\ <\ 0 .
\end{equation}
By \eqref{eq:Fs1-upper}, \(F_{s+1}\ge0\) forces
\[
u_{s+1}\,(a-u_{s+1})\ \le\ d\bigl(k_{s+1}\omega_{s+3/2}+\omega_{s+1/2}\bigr),
\]
which implies\ \(u_{s+1}\notin(\tau_-,\tau_+)\) by \ref{it:M}.
\end{proof}

In cell I the first site close to the upper state is \(s+1\), in cell II it is \(s\); cell III is a profile that dips at \(s+1\), and in cell IV the interface sits at \(s+2\) or further to the right.  We now establish two preliminary results for the two groups of cells, which characterize their corresponding admissible sets and provide \emph{a~prioiri} thresholds that guarantee their emptiness. The results for cells III and IV are easier to establish because there the value $u_{s+1} \le \tau_-$ is small. As a preparation, we 
introduce some notation for the zero set of \(F_s\).

\begin{lemma}
\label{lemma:zero-set}
Let \((a,d)\in\mathcal C\) and \(0<w<W_s^0\). There is a unique \(v_1(w)\in(0,a/3)\) and a closed interval \(\mathcal I_s(w)\subset(a/3,1)\), possibly empty or a single point, such that
\begin{equation}
\label{eq:v123}
\{v\in(0,1):\ F_s(v;w)\ge0\}=(0,v_1(w)]\ \cup\ \mathcal I_s(w).
\end{equation}
When \(\mathcal I_s(w)\neq\varnothing\) we write \(\mathcal I_s(w)=[v_2(w),v_3(w)]\) with \(a/3<v_2(w)\le v_3(w)<1\).
\end{lemma}

\begin{proof}
For fixed \(w\), \(F_s(\cdot\,;w)=\theta_s(w)-h_s(\cdot)\) is a cubic polynomial in \(v\).
Recalling \eqref{eq:exp:f:a3:neg}, we have
\begin{equation}
F_s(0;w)=\theta_s(w)>0, \qquad \qquad  F_s(a/3;w)<0 , \qquad \qquad F_s(1;w)\le-d\,\omega_{s-1/2}(1-\ell_s)<0.
\end{equation}
Its local minimum, if any, lies at \(v_-\ge\frac{(1+a)-\sqrt{1-a+a^2}}3\ge\frac a3\), so \(F_s(\cdot\,;w)\) is decreasing on \((0,a/3]\) and has a unique zero \(v_1(w)\) there. Beyond \(a/3\) it has at most one local maximum, so 
\[ \mathcal{I}_s := \{v \in (a/3,1) :F_s(v;w)\ge0\}
\] 
is either empty, a point, or a closed interval \([v_2(w),v_3(w)]\) with \(F_s<0\) on \((v_1,v_2)\cup(v_3,1]\).
\end{proof}

\begin{figure}[t]
\centering
\includegraphics[width=0.95\linewidth]{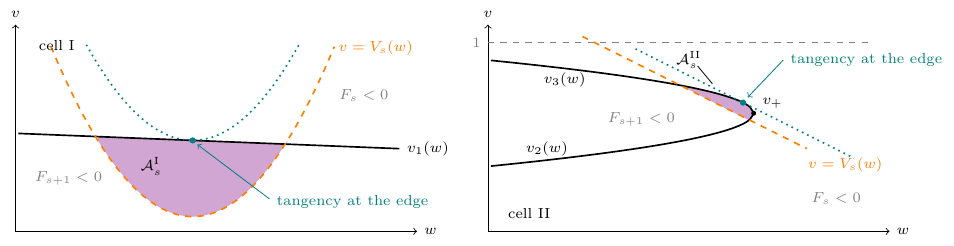}
\caption{Schematic (not to scale) of the bifurcation by which cells I and II empty. 
As \(a\) increases, \(V_s\) moves up and to the right while the band narrows and its tip retracts, until the two curves touch (teal dotted): in cell I at the bottom of the dip of \(V_s\), in cell II at the tip of the band, where \(v_3\) bends down steeply enough to match the slope of \(V_s\).}
\label{fig:bifurcation}
\end{figure}

\begin{lemma}
\label{lemma:cells}
Pick \((a,d)\in\mathcal C\).
\begin{enumerate}[label=(\roman*), ref=(\roman*)]
\item \label{lemma:cells-1} Define the vertical sections by
\[
\mathcal A_s^{\mathrm J}[w]:=\{v\in(0,1):(v,w)\in\mathcal A_s^{\mathrm J}\}.
\]
For \(0<w<W_s^0\) and \(u_{s+1}(w)\ge\tau_+\), they satisfy
\[
\mathcal A_s^{\mathrm I}[w]=(0,v_1(w)]\cap[V_s(w),\infty),
\qquad
\mathcal A_s^{\mathrm{II}}[w]=\mathcal I_s(w)\cap[V_s(w),\infty).
\]
For all other \(w\), these two sections are empty.
\item \label{lemma:cells-2} For \((v,w)\in\mathcal A_s^{\mathrm I}\) we have \(v\le\frac92Kd\). Set
\[
b_{\mathrm I}:=\omega_{s+1/2},\quad b_{\mathrm{II}}:=\omega_{s-1/2},\quad
y_{\mathrm I}:=1-u_{s+1},\quad y_{\mathrm{II}}:=1-v,\quad
E:=M(\epsilon+\sqrt d).
\]
In the respective cell \(\mathcal A_s^{\mathrm J}\), \(\mathrm J\in\{\mathrm I,\mathrm{II}\}\), we have
\[
y_{\mathrm J}(\epsilon-y_{\mathrm J})\ge db_{\mathrm J}(1-E),
\]
and consequently
\[
\frac\epsilon2-\sqrt{\frac{\epsilon^2}{4}-db_{\mathrm J}(1-E)}
\le y_{\mathrm J}\le
\frac\epsilon2+\sqrt{\frac{\epsilon^2}{4}-db_{\mathrm J}(1-E)}.
\]
\item \label{lemma:cells-3}  Recalling the constants \eqref{eq:M-constant}, we have the implications
\begin{align}
\label{eq:thr-I}
\mathcal A_s^{\mathrm I}\neq\varnothing\ &\Rightarrow\ \epsilon^2\ \ge\ 4d\,\omega_{s+1/2}\bigl(1-M(\epsilon+\sqrt d)\bigr)\ \ge\ 2d\,\omega_{s+1/2}\ \ge\ 2d\,\omega_{\min},\\
\label{eq:thr-II}
\mathcal A_s^{\mathrm{II}}\neq\varnothing\ &\Rightarrow\ \epsilon^2\ \ge\ 4d\,\omega_{s-1/2}\bigl(1-M(\epsilon+\sqrt d)\bigr)\ \ge\ 2d\,\omega_{s-1/2}\ \ge\ 2d\,\omega_{\min}.
\end{align}
\end{enumerate}
\end{lemma}

\begin{proof}\emph{\ref{lemma:cells-1}} This is \eqref{eq:v123} combined with 
\eqref{eq:Fs1-edge}. \emph{\ref{lemma:cells-2}, \ref{lemma:cells-3}} Since every coupling term in \(F_i\) is bounded by \(Kd\), the conditions \(F_s\ge0\) and \(F_{s+1}\ge0\) imply \(g(v;a)\ge-Kd\) and \(g(u_{s+1};a)\ge-Kd\). For the cubic \(g\) this leaves only two options for both \(u\in\{v,u_{s+1}\}\). Either \(u<a/3\), which with \(g(u;a)\le-\frac{4a}9u\) yields\begin{equation}\label{eq:u-small}u\le\frac{9Kd}{4a}\le\frac92Kd;\end{equation}or \(u>a/3\). For \(a/3<u\le a\), use \(g(u;a)\le-\frac a3(a-u)^2\), while for \(u>a\) use \(1-u<\epsilon\). Together these give\begin{equation}\label{eq:y-bound}y:=1-u\ \le\ \epsilon+\sqrt{3Kd/a}\ \le\ \epsilon+\sqrt{6Kd}.\end{equation}On cell I we have \(v<a/3\) and \(u_{s+1}\ge\tau_+>a/3\), so \eqref{eq:u-small} applies to \(v\) and \eqref{eq:y-bound} to \(u_{s+1}\); on cell II both \(v\) and \(u_{s+1}\) exceed \(a/3\), so \eqref{eq:y-bound} applies to both. This gives the crude bounds \(v\le\frac92Kd\) in \ref{lemma:cells-2} and \(1-v\le\epsilon+\sqrt{6Kd}\), resp.\ \(1-u_{s+1}\le\epsilon+\sqrt{6Kd}\); it remains to sharpen the latter at the site that carries the reaction, which is \(v\) on cell II and \(u_{s+1}\) on cell I.

On cell II we use \(F_s\ge0\) in the form \eqref{eq:Fs-edge}, bound the reaction by \(g(v;a)\le\epsilon^2/4\) and apply
the bound \(\ell_s\le Kd/a\)  from Lemma~\ref{lemma:tail-bounds}
to find\[\begin{aligned}\frac{\epsilon^2}{4}\ \ge\ g(v;a)\ &\ge\ d\,\omega_{s-1/2}(1-\ell_s)v-d\,k_s\omega_{s+1/2}(1-v)\\&\ge\ d\,\omega_{s-1/2}\Bigl[1-\frac{Kd}{a}-\Bigl(1+\frac{k_s\omega_{s+1/2}}{\omega_{s-1/2}}\Bigr)y\Bigr],\qquad y=1-v .\end{aligned}\]
On cell I we use \(F_{s+1}\ge0\) in the form \eqref{eq:Fs1-edge}, bound the reaction by \(g(u_{s+1};a)\le\epsilon^2/4\) and apply 
\(v\le\frac92Kd\) from \eqref{eq:u-small} to find
\[\begin{aligned}\frac{\epsilon^2}{4}\ \ge\ g(u_{s+1};a)\ &\ge\ d\,\omega_{s+1/2}\bigl(1-y-\tfrac92Kd\bigr)-d\,k_{s+1}\omega_{s+3/2}\,y\\&\ge\ d\,\omega_{s+1/2}\Bigl[1-\tfrac92Kd-\Bigl(1+\frac{k_{s+1}\omega_{s+3/2}}{\omega_{s+1/2}}\Bigr)y\Bigr],\qquad y=1-u_{s+1} .\end{aligned}\]
Inserting \eqref{eq:y-bound} together with \(k_i\omega_{i+1/2}/\omega_{i-1/2}\le K/\omega_{\min}\), \(a\ge\frac12\) and \(d\le\sqrt d\) yields the first inequalities in \eqref{eq:thr-I}--\eqref{eq:thr-II}; the remaining ones are immediate. Keeping the reaction exact in the same estimates gives \(g(1-y_{\mathrm J};a)=(1-y_{\mathrm J})y_{\mathrm J}(\epsilon-y_{\mathrm J})\ge db_{\mathrm J}(1-E)>0\). Since \(0<1-y_{\mathrm J}<1\), this implies the quadratic inequality in \ref{lemma:cells-2}; solving it gives the stated interval. Its radicand is nonnegative by \eqref{eq:thr-I}--\eqref{eq:thr-II}.\end{proof}

\begin{lemma}
\label{lemma:cells34}
Pick \((a,d)\in\mathcal C\).
\begin{enumerate}[label=(\roman*), ref=(\roman*)]
\item \label{lemma:cells34-1} For any \((v,w)\in\mathcal A_s^{\mathrm{III}}\) we have \(v\in(a,1)\). For any \((v,w)\in\mathcal A_s^{\mathrm{IV}}\) we have
\[
C_{s+1}\le\frac{1}{4(1-\tau_-)} .
\]
\item \label{lemma:cells34-2} We have the implications
\begin{align}
\label{eq:thr-III}
\mathcal A_s^{\mathrm{III}}\neq\varnothing\ &\Rightarrow\ \epsilon^2\ \ge\ \tfrac43\,a\,d\,\omega_{s-1/2}\,(1-\ell_s)\ \ge\ \tfrac13\,d\,\omega_{\min},\\
\label{eq:thr-IV}
\mathcal A_s^{\mathrm{IV}}\neq\varnothing\ &\Rightarrow\ (\epsilon+Kd)^2\ \ge\ 2\,d\,\omega_{s+3/2}\ \ge\ 2\,d\,\omega_{\min} .
\end{align}
\end{enumerate}
\end{lemma}

\begin{proof}
On cell III we have \(u_{s+1}\le\tau_-<a/3<v\) since \(\tau_-\le2Kd/a\), so both coupling terms in \eqref{eq:Fs-edge} pull \(F_s\) down and \(F_s\ge0\) forces
\begin{equation}
\label{eq:III-g}
g(v;a)\ \ge\ d\,k_s\omega_{s+1/2}(v-u_{s+1})+d\,\omega_{s-1/2}(1-\ell_s)v\ \ge\ d\,\omega_{s-1/2}(1-\ell_s)\,\frac a3\ >\ 0 .
\end{equation}
Positivity of \(g(v;a)\) means \(v\in(a,1)\), which provides \ref{lemma:cells34-1} for cell III. Combining \eqref{eq:III-g} with the bounds \(g(v;a)\le\epsilon^2/4\), \(a\ge\frac12\) and \(\ell_s\le Kd/a\le\frac12\) gives \eqref{eq:thr-III}.

On cell IV the parabola \(u_{s+1}=1-w(1-C_{s+1}w)\) has minimum \(1-\frac1{4C_{s+1}}\), so \(u_{s+1}\le\tau_-\) is only possible if \(C_{s+1}\le\frac1{4(1-\tau_-)}\), which is \ref{lemma:cells34-1} for cell IV. Combined with the lower bounds on \(C_{s+1}\) and \(r_{s+1}\) in \eqref{eq:tail-bounds-r}, and \((1+\epsilon)(1-\tau_-)\ge\frac12\), this gives
\[
\frac{d^2\,\omega_{s+3/2}^2}{(\epsilon+Kd)^2}\ \le\ r_{s+1}^2\ <\ \frac{d\,\omega_{s+3/2}\,C_{s+1}}{1+\epsilon}\ \le\ \frac{d\,\omega_{s+3/2}}{2},
\]
which is \eqref{eq:thr-IV}.
\end{proof}

Together, the four thresholds above exclude the part of \(\mathcal C\) above the line \(1-a=c_Rd\). In particular, the remainder of our analysis can focus on $\mathcal C_{\rm low}$, where
the additional inequalities \ref{it:R} hold.

\begin{lemma}
\label{lemma:other-triangle}
If \((a,d)\in\mathcal C\setminus\mathcal C_{\rm low}\), then \(\mathcal A_s^-=\varnothing\) for every phase \(\dInt{s}\).
\end{lemma}

\begin{proof}
Let \(\epsilon<c_Rd\). By the thresholds \eqref{eq:thr-I}--\eqref{eq:thr-IV}, non-emptiness of cells I, II and III requires \(\epsilon^2\ge\frac13d\,\omega_{\min}\), and that of cell IV requires \((\epsilon+Kd)^2\ge2d\,\omega_{\min}\). But \(\epsilon^2<c_R^2d^2<\frac13d\,\omega_{\min}\) and \((\epsilon+Kd)^2<(c_R+K)^2d^2\le\frac{25}{16}c_R^2d^2<d\,\omega_{\min}\) by \eqref{eq:constants}. Hence all four cells are empty, and so is \(\mathcal A_s^-\) by Lemma~\ref{lemma:four-cells}.
\end{proof}

\paragraph{Critical cells on $\mathcal C_{\rm low}$}

Focusing on $C_{\rm low}$, we now set out to characterize the boundary of $\mathcal{D}^-$ that is attributable to phase $s$. This will require us to understand the order in which the four components $\mathcal{A}_s^J$ become empty
and the mechanism by which the last component disappears. It will turn out that cells I and II are critical, while cell III can be safely ignored. Cell IV can be absorbed by the analysis at phase $s+1$.

In \(\mathcal C_{\rm low}\) a convenient simplification is available: we can ensure that the quadratic dependence of \(u_{s+1}\) on \(w\) in \eqref{eq:asymp:idss:for:usp1} becomes invertible by removing the increasing branch of the parabola (if present). In particular, we write
\begin{equation}
\label{eq:As-hat}
\widehat{\mathcal A}_s^-:=\mathcal A_s^-\cap\bigl\{C_{s+1}w\le\tfrac12\bigr\},
\qquad
\widehat{\mathcal A}_s^{\mathrm J}:=\mathcal A_s^{\mathrm J}\cap\bigl\{C_{s+1}w\le\tfrac12\bigr\},\quad \mathrm J\in\{\mathrm I,\mathrm{II},\mathrm{III},\mathrm{IV}\},
\end{equation}
and show that nothing is lost.

\begin{lemma}
\label{lemma:restriction}
Pick \((a,d)\in\mathcal C_{\rm low}\), suppose that \((v,w)\in\mathcal A_s^-\) and write $w'=\min\{w,1/C_{s+1}-w\} $.
Then \( (v,w')\in\widehat{\mathcal A}_s^-\), and it lies in the same cell. In particular \(\mathcal A_s^{\mathrm J}\neq\varnothing\Rightarrow\widehat{\mathcal A}_s^{\mathrm J}\neq\varnothing\) for every cell \(\mathrm J\).
\end{lemma}

\begin{proof}
Suppose \(C_{s+1}w>\frac12\) so that \(w':=1/C_{s+1}-w\). This yields \(C_{s+1}w'<\frac12\) and the symmetry of the parabola hence implies \(u_{s+1}(w')=u_{s+1}(w)\). In particular, \(F_s\) is unchanged and the tail stays positive since \(w'<w<w_s^*\). Since \(q_{s+2}=r_{s+1}w\) decreases and \(C_{s+2}q_{s+2}\le\frac12C_{s+1}w<\frac12\) before and after by \ref{it:R}, \(u_{s+2}=1-q_{s+2}(1-C_{s+2}q_{s+2})\) increases, and with it \(F_{s+1}\). The cell is preserved since \(v\) and \(u_{s+1}\) are.
\end{proof}

The key benefit of this restriction is that on
 cells I and II we now have uniform margins away from (almost) all constraints. 
 In particular, the only constraint admissible pairs can approach is \(v=0\) in cell I; see Figure~\ref{fig:cells}.
 This is an essential step for the compactness argument below, allowing us to show that these cells become empty through tangencies rather than boundary effects.

\begin{lemma}
\label{lemma:margins}
Let \((a,d)\in\mathcal C_{\rm low}\). There is \(\gamma=\gamma(d)>0\), independent of \(a\), such that every \((v,w)\in\widehat{\mathcal A}_s^{\mathrm I}\cup\widehat{\mathcal A}_s^{\mathrm{II}}\) satisfies
\begin{equation}
\label{eq:margins}
\begin{aligned}
\gamma\le w\le1-\gamma,&\qquad C_{s+1}w\le\tfrac14, &\qquad v_1(w)\ge \gamma,\\
\mathcal I_s(w)\subset(a/3,1-\gamma],&\qquad |v-a/3|\ge \gamma,&\qquad u_{s+1}-\tau_+\ge \gamma .
\end{aligned}
\end{equation}
Moreover \(w=(1-u_{s+1})(1+O(1-u_{s+1}))\) on \(\widehat{\mathcal A}_s^{\mathrm I}\). 
\end{lemma}

\begin{proof}
We treat the six bounds in \eqref{eq:margins} in turn, and then consider the asymptotics for \(w\).

\emph{Upper bound on \(w\) and on \(C_{s+1}w\).} Combining \(u_{s+1}\ge\tau_+\) with \(C_{s+1}w\le\frac12\) gives
\[
\tfrac12w\ \le\ w(1-C_{s+1}w)\ =\ 1-u_{s+1}\ \le\ 1-\tau_+\ =\ \epsilon+\tau_- .
\]
Hence \(w\le2(\epsilon+\tau_-)\le2(1-a_0)+8Kd_0\), which is at most \(1-\gamma\) by \ref{it:M}, and \(C_{s+1}w\le2C_{s+1}(\epsilon+\tau_-)\le\frac14\) by the third inequality in \ref{it:R}. Since \(w_s^*>1\) by Corollary~\ref{cor:tail:pos:u:r}, the positivity threshold is never approached.

\emph{Bounds on the zero set.} The estimates \(h_s(v)\le(a+dK)v\) and \(u_{s+1}\ge\tau_+\) give
\[
v_1(w)\ \ge\ \frac{\theta_s(w)}{a+dK}\ \ge\ \frac{d\,k_s\omega_{s+1/2}\,\tau_+}{a+dK}\ \ge\ \gamma .
\]
Since \(F_s(1;w)\le-d\,\omega_{s-1/2}(1-\ell_s)<0\) and the slope of \(F_s(\cdot\,;w)\) is bounded uniformly in \(w\), the interval \(\mathcal I_s(w)\) stays a distance \(\gamma\) away from \(v=1\). In addition, the bound 
\eqref{eq:exp:f:a3:neg} keeps both \((0,v_1(w)]\) and \(\mathcal I_s(w)\), hence \(v\), a distance \(\gamma\) away from \(a/3\).

\emph{Lower bound on \(w\).} Since \(V_s(w)\to1\) as \(w\downarrow0\) by \eqref{eq:Vs}, while \(v\le1-\gamma\) in both cells, the constraint \(v\ge V_s(w)\) cannot be met for small \(w\); this gives \(w\ge\gamma\).

\emph{Margin at \(\tau_+\).} 
Inspecting \eqref{eq:Fs1-edge}
and using $u_{s+2} 
\le 1$ and $v 
\le 1 - 
\gamma$, we can sharpen \eqref{eq:Fs1-upper}
to find
\[
F_{s+1}\ \le\ (1-u_{s+1})\bigl[d(k_{s+1}\omega_{s+3/2}+\omega_{s+1/2})-u_{s+1}(a-u_{s+1})\bigr]-\gamma\,d\,\omega_{s+1/2}.
\]
The bracket vanishes at \(u_{s+1}=\tau_+\) and has slope \(\sqrt{D_{s+1}}>0\) there, so \(F_{s+1}\ge0\) forces \(u_{s+1}-\tau_+\ge\gamma'\) for some \(\gamma'>0\) that depends on \(d\) only.

\emph{Asymptotics of \(w\).} On a nonempty cell, \eqref{eq:thr-I}--\eqref{eq:thr-II} give \(\epsilon^2\ge2d\omega_{\min}\). Thus \eqref{eq:tail-bounds-C}, with \(\rho=d\omega_{\max}/\epsilon\), bounds \(C_{s+1}\) uniformly. On cell I, writing \(y=1-u_{s+1}=w(1-C_{s+1}w)\) and using \(C_{s+1}w\le\frac14\) gives \(w\le\frac43y\), so
\[
w=\frac{y}{1-C_{s+1}w}=y\bigl(1+O(y)\bigr).
\]
For fixed \(d>0\), the same threshold inequalities separate the nonempty cells from \(a=1\). All preceding estimates are therefore uniform in \(a\) on the remaining compact range, so a single \(\gamma(d)>0\) can be chosen.
\end{proof}

The following two results quantify that cells I and II are the critical cells for our edge analysis. The geometric intuition can be read off from Figure~\ref{fig:cells}. Indeed, cell III empties as the tip of the band \(\mathcal I_s(w)\), where \(v_2(w)=v_3(w)\), moves leftwards into cell I. On cell IV both $u_s$ and $u_{s+1}$ are small, suggesting that the interface jump towards one can be captured by increasing the phase $s$.

\begin{lemma}
\label{lemma:III}
Let \((a,d)\in\mathcal C_{\rm low}\) and \((v,w)\in\mathcal A_s^{\mathrm{III}}\). Then \((v,w')\in\widehat{\mathcal A}_s^{\mathrm{II}}\), where \(w'\) is the amplitude with \(u_{s+1}(w')=1-\epsilon/2\).
\end{lemma}

\begin{proof}
On cell III, \(F_s\ge0\) and \(u_{s+1}\le\tau_-\le2Kd/a\) give
\[
h_s(v)\ \le\ \theta_s(w)\ =\ d\,k_s\omega_{s+1/2}\,u_{s+1}\ \le\ \frac{2K^2}{a}\,d^2 .
\]
Since \(h_s(v)=-g(v;a)+d\bigl(k_s\omega_{s+1/2}+\omega_{s-1/2}(1-\ell_s)\bigr)v\) and \(v>a/3\), this forces \(g(v;a)>0\), i.e.\ \(v\in(a,1)\), and in particular \(1-v<\epsilon\).

We now replace \(w\) by the amplitude \(w'\) on the decreasing branch with \(u_{s+1}(w')=1-\epsilon/2\); it exists since \(1-\epsilon/2\) lies between \(u_{s+1}(0)=1\) and the minimum \(1-1/(4C_{s+1})\) of the parabola, by \ref{it:R}. Increasing \(u_{s+1}\) increases \(\theta_s\), so
\[
F_s(v;w')\ >\ F_s(v;w)\ \ge\ 0 .
\]
For \(F_{s+1}\) we put \(y:=1-u_{s+1}(w')=\epsilon/2\) and use \eqref{eq:Fs1-edge} together with \(u_{s+2}\ge u_{s+1}\), which holds because \(q_{s+2}=r_{s+1}w'\le\frac14w'\) by \ref{it:R}. This allows us to estimate
\[
F_{s+1}(v;w')\ \ge\ -d\,\omega_{s+1/2}\,(1-v-y)+y(1-y)(\epsilon-y)
\ \ge\ -d\,\omega_{s+1/2}\,\frac\epsilon2+\frac{\epsilon^2}{4}\Bigl(1-\frac\epsilon2\Bigr)\ >\ 0,
\]
where the last inequality uses \(\epsilon\ge c_Rd\ge8\omega_{\max}d\). Finally, \(v>a/3\) and \(u_{s+1}(w')=1-\epsilon/2\ge a\ge\tau_+\), so \((v,w')\in\widehat{\mathcal A}_s^{\mathrm{II}}\).
\end{proof}

\begin{lemma}
\label{lemma:IV}
Let \((a,d)\in\mathcal C_{\rm low}\) and \((v,w)\in\widehat{\mathcal A}_s^{\mathrm{IV}}\). Then there is \(w'=r_{s+1}w\,(1-O(d^3))\) such that \(\bigl(u_{s+1}(w),\,w'\bigr)\in\widehat{\mathcal A}_{s+1}^{\mathrm I}\).
\end{lemma}

\begin{proof}
Write \(\vartheta:=u_{s+1}(w)\le\tau_-\) and consider the phase-\((s+1)\) profile \(\underline{\mathbf{u}}^{(s+1)}(\vartheta,r_{s+1}w)\). It coincides with \(\underline{\mathbf{u}}^{(s)}(v,w)\) at all sites \(i\ge s+1\): at \(i=s+1\) by the choice of \(\vartheta\), and at \(i\ge s+2\) because there both profiles are the right tail with \(q_{s+2}=r_{s+1}w\), i.e.\ \(q_i=w\,\mathcal R_{i;s+1}=r_{s+1}w\,\mathcal R_{i;s+2}\). At the sites \(i\le s\) it equals the left tail \(\vartheta\mathcal L_{i;s+1}\). Consequently its residuals at \(i\le s\) are positive by Lemma~\ref{lemma:tail:res:u:l}, and its residuals at \(i\ge s+2\) are those of the old profile. For the latter, Lemma~\ref{lemma:Ci} and \(C_{s+2}q_{s+2}\le\frac14\) give
\begin{equation}
\label{eq:IV-slack}
F_{s+2}\ \ge\ \tfrac18\,q_{s+2}^3,
\qquad
q_{s+2}\ =\ r_{s+1}w\ \ge\ \tfrac12\,r_{s+1},
\end{equation}
where \(w\ge1-\vartheta\ge\frac12\) was used. Only the residual at \(s+1\) changes, through the value at site \(s\), which is replaced by \(\ell_{s+1}\vartheta\) instead of \(v\le v_1(w)\). Applying \eqref{eq:Fs1-edge}, we find
\begin{equation}
\label{eq:IV-loss}
F_{s+1}^{\rm new}\ \ge\ F_{s+1}^{\rm old}-d\,\omega_{s+1/2}\bigl(v_1(w)-\ell_{s+1}\vartheta\bigr)_+\ \ge\ -\Lambda,
\qquad
\Lambda:=d\,\omega_{s+1/2}\bigl(v_1(w)-\ell_{s+1}\vartheta\bigr)_+ .
\end{equation}

By the left ratio relation \eqref{eq:left-ratio-relation}, \(h_s(v)=\frac{d\,k_s\omega_{s+1/2}}{\ell_{s+1}}v-(1+a)v^2+v^3\), so the defining relation \(h_s(v_1)=\theta_s(w)=d\,k_s\omega_{s+1/2}\vartheta\) of \(v_1\) yields
\[
0\ \le\ v_1-\ell_{s+1}\vartheta\ \le\ \frac{2v_1^2}{a},
\qquad
v_1\ \le\ \frac{2d\,k_s\omega_{s+1/2}\,\vartheta}{a},
\qquad
\vartheta\ \le\ \tau_-\ \le\ \frac{2Kd}{a},
\]
and therefore
\[
\Lambda\ \le\ \frac{32\,K^2k_s^2\,\omega_{\max}^3}{a^5}\,d^5 .
\]

To compensate for this loss, we lower the amplitude to
\[
w':=r_{s+1}w-\delta,\qquad \delta:=\frac{2\Lambda}{d\,k_{s+1}\omega_{s+3/2}}\ \le\ \frac{64\,K^2k_s^2\,\omega_{\max}^3}{a^5\,k_{s+1}\omega_{\min}}\,d^4 .
\]
Put \(q:=r_{s+1}w\). By \eqref{eq:IV-slack} and \eqref{eq:tail-bounds-r}, \(q\ge d\omega_{\min}/4\), whereas \(\delta=O(d^4)\). Hence \(\delta/q=O(d^3)\), so \(w'=q-\delta=q(1-O(d^3))>0\) for small \(d\).
Since \(\partial_{w'}u_{s+2}=-(1-2C_{s+2}w')\le-\frac12\), this raises \(u_{s+2}\) by at least \(\delta/2\), hence \(F_{s+1}\) by at least \(d\,k_{s+1}\omega_{s+3/2}\,\delta/2=\Lambda\), which by \eqref{eq:IV-loss} makes \(F_{s+1}\ge0\). Throughout this change the affected right-tail values have amplitude derivatives of absolute value at most one, since \(C_iq_i\le\frac14\) and the tail ratios are less than one. Since \(|g'(u;a)|\le1\) on \([0,1]\),
\[
|F_{s+2}^{\rm new}-F_{s+2}^{\rm old}|\le(1+2Kd)\delta<2\delta.
\]
By \eqref{eq:IV-slack}, the bounds \(r_{s+1}\ge d\omega_{\min}/(\epsilon+Kd)\ge d\omega_{\min}/2\), and the definition of \(d_0\) in \eqref{eq:constants} (recall \(a\ge a_0\ge\frac34\)),
\[
2\delta\le\frac{d^3\omega_{\min}^3}{2^9}
\le\tfrac18\bigl(\tfrac12r_{s+1}\bigr)^3
\le\tfrac18q^3\le F_{s+2}^{\rm old},
\]
so \(F_{s+2}^{\rm new}>0\). At all sites \(i\ge s+3\), the profile remains a corrected right tail with \(C_iq_i<1\), so its residuals retain the correct sign by Lemma~\ref{lemma:Ci}.

We now claim that the new pair \((\vartheta,w')\) lies in \(\widehat{\mathcal A}_{s+1}^{\mathrm I}\). The three defining conditions are \(\vartheta<a/3\), \(u_{s+2}(w')\ge\tau_+^{(s+1)}\) and \(C_{s+2}w'\le\frac12\). The first holds since \(\vartheta\le\tau_-<a/3\). For the other two we bound \(w'\): since \(w'\le r_{s+1}w\) and \(C_{s+1}w\le\frac12\) on the decreasing branch,
\[
w'\ \le\ r_{s+1}w\ \le\ \frac{r_{s+1}}{2C_{s+1}}\ <\ \frac{d\,\omega_{s+3/2}}{2(1+\epsilon)\,r_{s+1}}\ \le\ \frac{\epsilon+Kd}{2}\ \le\ \epsilon ,
\]
where we used
\eqref{eq:tail-bounds-r}
and \(c_R\ge4K\).
 Hence \[u_{s+2}(w')=1-w'(1-C_{s+2}w')\ge1-w'\ge1-\epsilon=a\ge\tau_+^{(s+1)},\] 
 and \(C_{s+2}w'\le C_{s+2}r_{s+1}w\le\frac12C_{s+1}w\le\frac14\) by \ref{it:R}. Also \(0<w'\le\epsilon<1<w_{s+1}^*\), verifying the shifted positivity threshold.
\end{proof}

We are now ready to characterize the boundaries associated to the two critical cells. 
For \(\mathrm J\in\{\mathrm I,\mathrm{II}\}\) and \(d<d_0\) set
\[
a_s^{\mathrm J}(d):=\sup\bigl\{a:\ (a,d)\in\mathcal C,\ \widehat{\mathcal A}_s^{\mathrm J}(a,d)\neq\varnothing\bigr\}\qquad(\sup\varnothing:=a_0).
\]
Figure~\ref{fig:bifurcation} depicts the two types of bifurcations that can occur at these boundaries. 

\begin{lemma}
\label{lemma:edges}
If \(a^\ast:=a_s^{\mathrm J}(d)>a_0\), then \(\widehat{\mathcal A}_s^{\mathrm J}(a^\ast,d)\neq\varnothing\), and every \((v^\ast,w^\ast)\in\widehat{\mathcal A}_s^{\mathrm J}(a^\ast,d)\) solves the \emph{tangency system}
\begin{equation}
\label{eq:tangency}
F_s=0,\qquad F_{s+1}=0,\qquad 
h_s'(v^\ast)\,V_s'(w^\ast)=\theta_s'(w^\ast)
\qquad\text{at }(a^\ast,v^\ast,w^\ast).
\end{equation}
The contact also satisfies
\begin{equation}
\label{eq:upper-contact}
h_s'(v^\ast)>0,\qquad V_s'(w^\ast)<0.
\end{equation}
Thus the curves \(\{F_s=0\}\) and \(\{v=V_s(w)\}\) touch tangentially at \((v^\ast,w^\ast)\), along the branch \(v_1\) for \(\mathrm J=\mathrm I\) and along the upper endpoint \(v_3\) of \(\mathcal I_s\) for \(\mathrm J=\mathrm{II}\). 
\end{lemma}

\begin{proof}
We first show that the supremum is attained. Pick \((v_n,w_n)\in\widehat{\mathcal A}_s^{\mathrm J}(a_n,d)\) with \(a_n\uparrow a^\ast\). By Lemma~\ref{lemma:margins} these pairs satisfy the uniform bounds
\[
\gamma\le w_n\le1-\gamma,\qquad C_{s+1}w_n\le\tfrac14,\qquad v_n\le1-\gamma,\qquad |v_n-a/3|\ge\gamma,\qquad u_{s+1}(w_n)-\tau_+\ge\gamma,
\]
so after passing to a subsequence we have \((v_n,w_n)\to(v^\ast,w^\ast)\), and the limit satisfies the same bounds. We claim that also \(v^\ast>0\). Indeed, if \(v^\ast=0\), then \(V_s(w^\ast)\le\lim v_n=0\), while \(v_1(w^\ast)\ge\gamma\) by Lemma~\ref{lemma:margins}. Choose a fixed \(0<\widetilde v<\gamma\). Since \(\widetilde v<v_1(w^\ast)\) and \(\widetilde v>V_s(w^\ast)\), both residuals are strictly positive at \((\widetilde v,w^\ast)\). The margin bounds put this point in cell I, and it remains admissible for \(a\) slightly above \(a^\ast\), a contradiction.
The same continuity gives \(F_s\ge0\) and \(F_{s+1}\ge0\) at \((a^\ast,v^\ast,w^\ast)\), so that
\[
(v^\ast,w^\ast)\in\widehat{\mathcal A}_s^{\mathrm J}(a^\ast,d).
\]

We now show that every \((v,w)\in\widehat{\mathcal A}_s^{\mathrm J}(a^\ast,d)\) solves \eqref{eq:tangency}. Suppose not. Then either one of \(F_s\), \(F_{s+1}\) is strictly positive at \((v,w)\), or both vanish with non-parallel gradients. Computing
\[
\nabla_{(v,w)}F_{s}=\bigl(-h_s'(v),\theta_s'(w)\bigr),
\qquad \qquad
\nabla_{(v,w)}F_{s+1}=d\,\omega_{s+1/2}\,\bigl(1,-V_s'(w)\bigr)\ \neq\ 0,
\]
the latter alternative means that the two curves \(\{F_s=0\}\) and \(\{F_{s+1}=0\}\) cross transversally at \((v,w)\). In both alternatives the set \(\{F_s\ge0\}\cap\{F_{s+1}\ge0\}\) contains an open neighborhood of a point arbitrarily close to \((v,w)\); by the margins of Lemma~\ref{lemma:margins} this neighborhood lies in \(\widehat{\mathcal A}_s^{\mathrm J}(a^\ast,d)\), and by continuity of \(F_s,F_{s+1}\) in \(a\) it remains admissible for \(a\) slightly above \(a^\ast\). This contradicts the definition of \(a^\ast\), and hence \eqref{eq:tangency} holds. On the margin set, \(\theta_s'(w)=dk_s\omega_{s+1/2}(-1+2C_{s+1}w)<0\). If \(h_s'(v)\le0\), choose \(\delta w<0\) and a sufficiently large fixed positive ratio \(\delta v/|\delta w|\). The two residual gradients displayed above then have positive scalar product with this direction, so a small perturbation makes both residuals strictly positive. This point persists for \(a>a^\ast\), again a contradiction. Consequently \(h_s'(v^\ast)>0\), and \eqref{eq:tangency} implies \(V_s'(w^\ast)=\theta_s'(w^\ast)/h_s'(v^\ast)<0\). By Lemma~\ref{lemma:cells}\ref{lemma:cells-1}, this selects \(v_1\) in cell I and the upper root \(v_3\) in cell II.
\end{proof}

\paragraph{The tangency system}

We now set out to construct a branch of solutions to each of the tangency systems \eqref{eq:tangency} that bifurcates from the corner point $(a, d) = (1, 0)$. In particular, we use 
\(\sigma:=\sqrt d\) as a parameter and solve the three equations for the three unknowns \((v,w,a)\). This will yield, for each cell, a smooth curve \(a=\mathsf a_s^{\mathrm J}(\sigma)\) through the edge points that can be expanded in powers of \(\sqrt d\).

We start by considering the 
tail data \(\ell_s,r_i,C_i\),
which we scale as
\begin{equation}
\label{eq:exp:tail:data:scaling}
    \ell_i=\sigma^2L_i, \qquad \qquad
    r_i=\sigma R_i,
    \qquad \qquad C_i=C_i.
\end{equation}\
Together with the scaling \(a=1-\sigma A\),  the tail identities \eqref{eq:left-ratio-relation},
\eqref{eq:right-ratio-relation} and \eqref{eq:exp:c:rewritten} 
yield the fixed-point relations
\begin{align}
\label{eq:tail-L}
L_{i+1}&=\frac{k_i\omega_{i+1/2}}{1-\sigma A+\sigma^2\bigl[k_i\omega_{i+1/2}+\omega_{i-1/2}(1-\sigma^2L_i)\bigr]},\\
\label{eq:tail-R}
R_{i-1}&=\frac{\omega_{i-1/2}}{A+\sigma\bigl[\omega_{i-1/2}+k_i\omega_{i+1/2}(1-\sigma R_i)\bigr]},\\
\label{eq:tail-C}
C_{i-1}&=\frac{R_{i-1}^2}{\omega_{i-1/2}}\Bigl[1+\sigma A+\sigma\,\frac{\omega_{i-1/2}}{R_{i-1}}\,C_i+\sigma^3k_i\omega_{i+1/2}R_i\bigl(C_i-\sigma R_iC_{i+1}\bigr)\Bigr].
\end{align}
For $A > 0$ these relations decouple at \(\sigma=0\) 
and the perturbation into $\sigma > 0$ is regular. To quantify this, it is convenient to introduce the lower bound
\begin{equation}
\label{eq:A-range}
A_-:=\tfrac12\sqrt{\omega_{\min}},
\end{equation}
which lies below the leading-order values \(2\sqrt{\omega_{s+1/2}}\), resp.\ \(2\sqrt{\omega_{s-1/2}}\)
that can be read off from item \ref{lemma:cells-3} of Lemma \ref{lemma:cells}.

\begin{lemma}[Scaled tail data]
\label{lemma:scaled-tails}
For every \(\varrho>0\) there is \(\sigma_T(\varrho)>0\) so that for \(A\in[A_-,\varrho]\) and \(0\le\sigma\le\sigma_T(\varrho)\) the relations \eqref{eq:tail-L}--\eqref{eq:tail-C} have a unique $p$-periodic solution \((L_i,R_i,C_i)_{i\in\mathbb Z}\) in a neighborhood of the decoupled solution. It is smooth in \((A,\sigma)\) on \([A_-,\varrho]\times[0,\sigma_T(\varrho)]\), satisfies
\[
L_s=k_{s-1}\omega_{s-1/2}+O(\sigma),\qquad R_i=\frac{\omega_{i+1/2}}{A}+O(\sigma),\qquad C_i=\frac{\omega_{i+1/2}}{A^2}+O(\sigma),
\]
uniformly in \(A\in[A_-,\varrho]\)
and matches
the actual tail data \(\ell_i,r_i,C_i\) at \((a,d)=(1-\sigma A,\sigma^2)\) via \eqref{eq:exp:tail:data:scaling}.
\end{lemma}

\begin{proof}
On \([A_-,\varrho]\) the right-hand sides of \eqref{eq:tail-L}--\eqref{eq:tail-C} are smooth in all their arguments, since the denominators are bounded away from zero for \(\sigma\le\sigma_T(\varrho)\) small. At \(\sigma=0\) the system is triangular and the derivative of the right-hand sides with respect to \((L_i,R_i,C_i)_{i \in \mathbb Z}\) is strictly block-triangular, so that the derivative of the fixed-point map minus the identity is invertible. The implicit function theorem, applied at each \(A\) and with constants uniform on the compact interval \([A_-,\varrho]\), gives the smooth branch and the expansions. To identify this branch with the actual tail data, first use the upper and lower bounds in \eqref{eq:tail-bounds-r} to obtain
\[
R_i=\frac{\omega_{i+1/2}}{A}+O(\sigma)
\]
uniformly on \([A_-,\varrho]\). The estimate \eqref{eq:tail-bounds-C} bounds \(C_i\) uniformly, and \eqref{eq:tail-C} then gives
\[
C_{i-1}=\frac{R_{i-1}^2}{\omega_{i-1/2}}+O(\sigma)
=\frac{\omega_{i-1/2}}{A^2}+O(\sigma).
\]
Likewise, \eqref{eq:tail-L} gives the stated leading value of \(L_i\). Thus the actual tail data lie in the implicit-function neighborhood and satisfy the same relations, so uniqueness identifies them with the smooth branch.
\end{proof}

Returning to \eqref{eq:tangency},
we use \(F_s=\theta_s(w)-h_s(v)\) and \(F_{s+1}=d\,\omega_{s+1/2}(v-V_s(w))\) to write the system in the form \(\mathcal T(v,w,a;\sigma)=0\) with
\begin{equation}
\label{eq:Phi}
\mathcal T:=\begin{pmatrix}\theta_s(w)-h_s(v)\\[1pt] v-V_s(w)\\[1pt] h_s'(v)\,V_s'(w)-\theta_s'(w)\end{pmatrix},
\quad
\partial_{(v,w,a)}\mathcal T=\begin{pmatrix}
-h_s'(v) & \theta_s'(w) & \partial_a\theta_s-\partial_ah_s\\[2pt]
1 & -V_s'(w) & -\partial_aV_s\\[2pt]
h_s''(v)V_s'(w) & h_s'(v)V_s''(w)-\theta_s''(w) & \partial_a\bigl(h_s'V_s'-\theta_s'\bigr)
\end{pmatrix}.
\end{equation}
We emphasize that all quantities also depend on the pair \(a,\sigma\) through the tail data \(\ell_s,r_i,C_i\).

To resolve the upper contacts, we introduce cell-dependent scalings in \(\sigma=\sqrt d\). Lemma~\ref{lemma:apriori} below justifies these scalings for all upper contacts in the corner. We write
\begin{equation}
\label{eq:scaling}
\text{cell I:}\quad v=\sigma^2V,\ \ w=\sigma W; \ 
\qquad\qquad
\text{cell II:}\quad v=1-\sigma Z,\ \ w=\sigma^2W, \ 
\end{equation}
where we recall the common scaling $a = 1 - \sigma A$. 

We now write \(x:=(V,W,A)\) resp.\ \(x:=(Z,W,A)\) for the relevant variables in cells I and II. Dividing out the singular terms of 
\(\mathcal T\) and using the substitutions \eqref{eq:scaling},
we obtain the two nonlinearities
\begin{equation}
\label{eq:Phi-tilde}
\begin{aligned}
\widetilde{\mathcal T}_{\mathrm I}(x;\sigma)&:=
\begin{pmatrix}
\sigma^{-2}[\theta_s(w)-h_s(v)]\\[3pt]
\omega_{s+1/2}[v-V_s(w)]\\[3pt]
-\sigma\omega_{s+1/2}\left[V_s'(w)-\dfrac{\theta_s'(w)}{h_s'(v)}\right]
\end{pmatrix},\\[6pt]
\widetilde{\mathcal T}_{\mathrm{II}}(x;\sigma)&:=
\begin{pmatrix}
\sigma^{-2}[\theta_s(w)-h_s(v)]\\[3pt]
\sigma^{-1}\omega_{s+1/2}[v-V_s(w)]\\[3pt]
\sigma^{-1}\left[h_s'(v)-\dfrac{\theta_s'(w)}{V_s'(w)}\right]
\end{pmatrix}.
\end{aligned}
\end{equation}

which allow us to phrase the tangency system \eqref{eq:tangency} in cell \(\mathrm J\) as
\[
\widetilde{\mathcal T}_{\mathrm J}(x;\sigma)=0.
\]
The maps extend smoothly to \(\sigma=0\) on bounded scaled domains with \(A\ge A_-\), after choosing \(\sigma\) small enough to keep the quotient denominators nonzero; see Lemma~\ref{lemma:R}.
For cell I,
\[
v-V_s(w)=\frac{W(A-W)-\omega_{s+1/2}}{\omega_{s+1/2}}+O(\sigma),
\qquad V_s'(w)=\frac{2W-A}{\sigma\omega_{s+1/2}}+O(1),
\]
whereas for cell II,
\[
v-V_s(w)=\sigma\left[\frac{AW}{\omega_{s+1/2}}-Z\right]+O(\sigma^2).
\]
%(we drop the index \(\mathrm J\) when no confusion can arise). 
Indeed, the leading parts \(\widetilde{\mathcal T}_J^{(0)}:=\widetilde{\mathcal T}_J(\cdot\,;0)\), the leading-order solutions \(x_{0;J}\) and the linearisation \(\mathbf L_J:=D_x\widetilde{\mathcal T}^{(0)}_J(x_{0;J})\) are given by
\begin{equation}
\label{eq:leading-I}
\widetilde{\mathcal T}_I^{(0)}=\begin{pmatrix}k_s\omega_{s+1/2}-V\\ W(A-W)-\omega_{s+1/2}\\ A-2W\end{pmatrix},\qquad
x_{0;I}=\begin{pmatrix}k_s\omega_{s+1/2}\\ \sqrt{\omega_{s+1/2}}\\ 2\sqrt{\omega_{s+1/2}}\end{pmatrix},\qquad
\mathbf L_I=\begin{pmatrix}-1&0&0\\ 0&0&W_0\\ 0&-2&1\end{pmatrix},
\end{equation}
for cell I, with \(\det\mathbf L_I=-2W_0=-2\sqrt{\omega_{s+1/2}}\neq0\); together with
\begin{equation}
\label{eq:leading-II}
\widetilde{\mathcal T}_{II}^{(0)}=\begin{pmatrix}Z(A-Z)-\omega_{s-1/2}\\ WA-\omega_{s+1/2}Z\\ A-2Z\end{pmatrix},\qquad
x_{0;II}=\begin{pmatrix}\sqrt{\omega_{s-1/2}}\\ \tfrac12\omega_{s+1/2}\\ 2\sqrt{\omega_{s-1/2}}\end{pmatrix},\qquad
\mathbf L_{II}=\begin{pmatrix}0&0&Z_0\\ -\omega_{s+1/2}&A_0&W_0\\ -2&0&1\end{pmatrix},
\end{equation}

for cell II, with \(\det\mathbf L_{II} =2A_0Z_0=4\omega_{s-1/2}\neq0\).

We now drop the cell index \(\mathrm J\) whenever the context is clear, and introduce the perturbation variable
\[
y:=x-x_0,\qquad y=(y_V,y_W,y_A)\ \text{ in cell I},\qquad y=(y_Z,y_W,y_A)\ \text{ in cell II}.
\]
Since \(\widetilde{\mathcal T}^{(0)}\) is a quadratic polynomial in \(x\) vanishing at \(x_0\), the scaled system takes the form
\begin{equation}
\label{eq:N}
\widetilde{\mathcal T}(x_0+y;\sigma)=\mathbf L\,y+\mathcal Q(y)+\sigma\,\mathcal R(y;\sigma),
\end{equation}
where \(\mathcal Q\) is the quadratic part of the leading system,
\[
\mathcal Q_{\mathrm I}(y)=\bigl(0,\ y_Wy_A-y_W^2,\ 0\bigr)^{\!\top},
\qquad
\mathcal Q_{\mathrm{II}}(y)=\bigl(y_Zy_A-y_Z^2,\ y_Wy_A,\ 0\bigr)^{\!\top},
\]
and the remainder \(\mathcal R(y;\sigma):=\sigma^{-1}\bigl[\widetilde{\mathcal T}(x_0+y;\sigma)-\widetilde{\mathcal T}^{(0)}(x_0+y)\bigr]\) collects everything that carries a positive power of \(\sigma\), including all tail data.

\begin{lemma}
\label{lemma:R}
Pick any \(\varrho>0\). Choose \(\sigma_\varrho>0\) sufficiently small with \(\sigma_\varrho\le\min\{\sigma_T(|x_0|+\varrho),\ (|x_0|+\varrho)^{-1}\}\).
Then there exists $M_{\varrho} > 0$ so 
that the remainder  \(\mathcal R\) is smooth
on the set
\[
B_\varrho:=\bigl\{(y,\sigma):\ |y|\le\varrho,\ \ A_0+y_A\ge A_-,\ \ 0\le\sigma\le\sigma_\varrho\bigr\}
\]
and admits the bound 
\begin{equation}
\label{eq:R-bound}
|\mathcal R(y;\sigma)|\ \le\ M_\varrho\qquad\text{on }B_\varrho.
\end{equation}
\end{lemma}

\begin{proof}
On \(B_\varrho\), \(A_-\le A\le|x_0|+\varrho\), so Lemma~\ref{lemma:scaled-tails} gives smooth scaled tail data. Uniformly on this bounded set,
\[
h_s'(\sigma^2V)=1+O(\sigma),\qquad
\sigma V_s'(\sigma^2W)=-\frac{A}{\omega_{s+1/2}}+O(\sigma).
\]
Decrease \(\sigma_\varrho\) so that
\[
|h_s'(v)|\ge\tfrac12\quad\text{in cell-I scaling},\qquad
|\sigma V_s'(w)|\ge c_\varrho>0\quad\text{in cell-II scaling}.
\]
Clearing the explicit powers of \(\sigma\) in \eqref{eq:Phi-tilde} using \eqref{eq:tail-L}--\eqref{eq:tail-C} gives smooth extensions of both scaled maps to \(\sigma=0\). Their values there are the leading systems \eqref{eq:leading-I}--\eqref{eq:leading-II}. Consequently the remainder extends smoothly as
\[
\mathcal R(y;\sigma)=\int_0^1\partial_\sigma\widetilde{\mathcal T}(x_0+y;t\sigma)\,dt.
\]
The bound \eqref{eq:R-bound} follows by compactness of \(B_\varrho\).
\end{proof}

The scaled variables are the natural ones for the upper contacts \eqref{eq:upper-contact}: every such solution in a cell is, in these variables, within \(O(\sigma)\) of the leading-order point.

Conversely, small perturbations of \(x_{0;J}\) that solve the scaled system lie in the cell.

\begin{lemma}
\label{lemma:apriori}
There are \(M_1>0\) and \(\sigma_1>0\) such that the following holds for \(0<\sigma<\sigma_1\) and \(\mathrm J\in\{\mathrm I,\mathrm{II}\}\). If \(x=x_{0;J}+y\) solves \(\widetilde{\mathcal T}_{\mathrm J}(x;\sigma)=0\) and the corresponding pair \((v,w)\) from \eqref{eq:scaling} lies in the closure of \(\widehat{\mathcal A}_s^{\mathrm J}(a,\sigma^2)\), \(a=1-\sigma A\), with \((a,\sigma^2)\in\mathcal C\) and \(h_s'(v)>0\), \(V_s'(w)<0\), then
\begin{equation}
\label{eq:apriori}
|y|\ \le\ M_1\,\sigma .
\end{equation}
\end{lemma}

\begin{proof}
Write \(y_R:=1-u_{s+1}\) and \(z:=1-v\). On a nonempty cell, Lemma~\ref{lemma:cells} gives \(d=O(\epsilon^2)\); its preliminary bounds give \(y_R=O(\epsilon)\), and also \(z=O(\epsilon)\) in cell II. Lemma~\ref{lemma:margins} gives \(C_{s+1}w\le\frac14\), while \eqref{eq:tail-bounds-C} bounds \(C_i\) uniformly. In particular, \(w\le\frac43y_R\), and the next right-tail value satisfies \(u_{s+2}\ge u_{s+1}\).

\emph{Cell II.} The sign \(V_s'(w)<0\), differentiated in \eqref{eq:Fs1-edge}, implies \(g'(u_{s+1};a)\le Cd\): here \(\partial_wu_{s+1}<0\) and \(\partial_wu_{s+2}<0\). Expanding the cubic derivative gives
\[
y_R\le\frac\epsilon2+O(\epsilon^2+d)\le\frac{2\epsilon}{3}
\]
after shrinking the corner. Thus \(F_{s+1}=0\) and \(u_{s+2}\ge u_{s+1}\) yield
\[
c_*\epsilon y_R\le g(u_{s+1};a)
\le d\omega_{s+1/2}(1-v)\le Cd\epsilon.
\]
Here \(c_*>0\) and \(C>0\) are independent of \((a,d)\) in the corner. It follows that \(y_R=O(d)\) and \(w=O(d)\). Differentiating \eqref{eq:Fs1-edge} again now gives
\[
V_s'(w)=-\frac{\epsilon}{d\omega_{s+1/2}}\bigl(1+O(\epsilon)\bigr),
\qquad h_s'(v)=\frac{\theta_s'(w)}{V_s'(w)}=O(d^2/\epsilon).
\]
Since \(h_s'(v)=\epsilon-2z+O(z^2+\epsilon z+d)\), we obtain \(z=\epsilon/2+O(\epsilon^2+d)\). Inserting this into \(F_s=0\), using \(\ell_s=O(d)\) and \(y_R=O(d)\), gives
\[
\frac{\epsilon^2}{4}\bigl(1+O(\epsilon)\bigr)
=d\omega_{s-1/2}\bigl(1+O(\epsilon)\bigr).
\]
Hence \(A=\epsilon/\sigma\) is bounded above. Now \(Z=z/\sigma\) and \(W=w/\sigma^2\) are bounded, and \eqref{eq:thr-II} gives \(A\ge\sqrt{2\omega_{s-1/2}}\ge A_-\). Thus \((y,\sigma)\in B_\varrho\) for a fixed \(\varrho\) once \(\sigma\) is small.

These bounds allow us to apply \eqref{eq:N} and \eqref{eq:R-bound}
to conclude \(\mathbf L_{\mathrm{II}}\,y+\mathcal Q_{\mathrm{II}}(y)=O(\sigma)\). Componentwise, this means
\[
Z_0y_A+y_Zy_A-y_Z^2=O(\sigma),\qquad
-\omega_{s+1/2}y_Z+A_0y_W+W_0y_A+y_Wy_A=O(\sigma),\qquad
y_A-2y_Z=O(\sigma).
\]
The third relation gives \(y_A=2y_Z+O(\sigma)\); the first then reads \(y_Z\,(2Z_0+y_Z)=O(\sigma)\). Since \[2Z_0+y_Z=Z_0+Z\ge Z_0>0\] we obtain \(y_Z=O(\sigma)\), hence \(y_A=O(\sigma)\).
The second identity finally gives \(y_W\,(A_0+y_A)=O(\sigma)\). 
We hence conclude \(y_W=O(\sigma)\) and obtain \eqref{eq:apriori} in cell II.

\emph{Cell I.} Here \(v=O(d)\), so \(h_s'(v)=1+O(\epsilon+d)\) is bounded away from zero. The tangency equation gives \(V_s'(w)=\theta_s'(w)/h_s'(v)=O(d)\). Differentiating \eqref{eq:Fs1-edge} therefore gives \(g'(u_{s+1};a)=O(d)\), whence
\[
y_R=\frac\epsilon2+O(\epsilon^2+d).
\]
The equation \(F_{s+1}=0\) now yields
\[
\frac{\epsilon^2}{4}\bigl(1+O(\epsilon)\bigr)
=d\omega_{s+1/2}\bigl(1+O(\epsilon)\bigr).
\]
Thus \(A\) is bounded, \(V=v/\sigma^2\) and \(W=w/\sigma\) are bounded, and \(A\ge A_-\) by \eqref{eq:thr-I}. Applying \eqref{eq:N}--\eqref{eq:R-bound} in this bounded scaled set gives
\[
V=k_s\omega_{s+1/2}+O(\sigma),\qquad
W(A-W)=\omega_{s+1/2}+O(\sigma),\qquad A-2W=O(\sigma).
\]
Since \(W\ge0\), these imply \(W=\sqrt{\omega_{s+1/2}}+O(\sigma)\), \(A=2\sqrt{\omega_{s+1/2}}+O(\sigma)\), and \eqref{eq:apriori} follows in cell I as well.
\end{proof}

\begin{lemma}
\label{lemma:back-in-cell}
There are \(\varrho_0>0\) and \(\sigma_2>0\) such that for \(0<\sigma<\sigma_2\) every solution \(x=x_{0;J}+y\) of \(\widetilde{\mathcal T}_{\mathrm J}(x;\sigma)=0\) with \(|y|\le\varrho_0\) corresponds, through \eqref{eq:scaling}, to a pair \((v,w)\in\widehat{\mathcal A}_s^{\mathrm J}(a,\sigma^2)\)  with \(a=1-\sigma A\).
\end{lemma}

\begin{proof}
The pair satisfies \(F_s=F_{s+1}=0\), and its tail is positive since \(w<1<w_s^*\) (Corollary~\ref{cor:tail:pos:u:r}); hence \((v,w)\in\mathcal A_s^-(a,\sigma^2)\), and it remains to verify the conditions that single out \(\widehat{\mathcal A}_s^{\mathrm J}\), namely
\[
\text{cell I: } v<\tfrac a3;\qquad \text{cell II: } v>\tfrac a3;\qquad\text{both: } u_{s+1}\ge\tau_+,\quad C_{s+1}w\le\tfrac12 .
\]
Choose \(\varrho_0:=\frac14\min\{V_0,W_0,A_0\}\), resp.\ \(\frac14\min\{Z_0,W_0,A_0\}\). Then \(|y|\le\varrho_0\) means that each coordinate of \(x\) lies between \(\frac34\) and \(\frac54\) times the corresponding coordinate of \(x_0\); in particular all coordinates are positive and bounded, and \(A\ge\frac34A_0\ge A_-\). For \(\sigma\) small the conditions now follow from \eqref{eq:scaling}:
\begin{itemize}
\item In cell I, \(v=\sigma^2V<a/3\), and
\[
1-u_{s+1}=\sigma W\,(1-\sigma C_{s+1}W)\ \le\ \sigma W\ \le\ \tfrac54\sigma W_0\ <\ \tfrac74\sigma W_0\ \le\ \sigma A=\epsilon ,
\]
where we used \(A_0=2W_0\); moreover \(C_{s+1}w=\sigma C_{s+1}W\le\frac14\).
\item In cell II, \(v=1-\sigma Z>a/3\), and
\[1-u_{s+1}=\sigma^2W(1-\sigma^2C_{s+1}W)\le\sigma A=\epsilon;\] 
moreover \(C_{s+1}w=\sigma^2C_{s+1}W\le\frac14\).
\end{itemize}
In both cells \(u_{s+1}\ge1-\epsilon=a\ge\tau_+\), so \((v,w)\in\widehat{\mathcal A}_s^{\mathrm J}(a,\sigma^2)\).
\end{proof}

We are now ready to apply the implicit function to solve \eqref{eq:tangency}
and construct the edge curves for cells I and II.
By restricting $d_0$ we can ensure that no other upper contacts satisfying \eqref{eq:upper-contact} exist in the respective cells.

\begin{proposition}
\label{prop:edge-curves}
After shrinking \(d_0\) and increasing \(a_0\) if necessary, for each phase \(s\) and \(\mathrm J\in\{\mathrm I,\mathrm{II}\}\) there are smooth functions
\[
\mathsf a_s^{\mathrm J}:(0,\sqrt{d_0})\to(a_0,1),\qquad v_s^{\mathrm J},\,w_s^{\mathrm J}:(0,\sqrt{d_0})\to(0,1),
\]
such that for every \(\sigma\in(0,\sqrt{d_0})\), \(\bigl(v_s^{\mathrm J}(\sigma),w_s^{\mathrm J}(\sigma),\mathsf a_s^{\mathrm J}(\sigma)\bigr)\) is the \emph{unique upper contact}, that is, the unique solution of \eqref{eq:tangency} satisfying \eqref{eq:upper-contact}, at \(d=\sigma^2\) with \((v,w)\) in the closure of \(\widehat{\mathcal A}_s^{\mathrm J}\); moreover \((v_s^{\mathrm J},w_s^{\mathrm J})(\sigma)\) lies in cell \(\mathrm J\) itself. In the scaled variables \eqref{eq:scaling}, \(x(\sigma)=x_{0;J}+O(\sigma)\) with \(x_{0;J}\) from \eqref{eq:leading-I}, resp.\ \eqref{eq:leading-II}. In particular,
\begin{equation}
\label{eq:edge-curve-expansion}
\mathsf a_s^{\mathrm J}(\sigma)=1-2\sqrt{\omega^{\mathrm J}}\,\sigma+O(\sigma^2),\qquad
1-u^{\mathrm J}_{\rm first}(\sigma)=\sqrt{\omega^{\mathrm J}}\,\sigma+O(\sigma^2),
\end{equation}
where \(\omega^{\mathrm I}:=\omega_{s+1/2}\), \(\omega^{\mathrm{II}}:=\omega_{s-1/2}\) is the coupling of the first high site to its left neighbor, while \(u_{\rm first}=u_{s+1}\) on cell I and \(u_{\rm first}=v\) on cell II.
\end{proposition}

\begin{proof}
By \eqref{eq:N}--\eqref{eq:R-bound} the map \((y,\sigma)\mapsto\mathbf L\,y+\mathcal Q(y)+\sigma\mathcal R(y;\sigma)\) is smooth up to \(\sigma=0\), vanishes at \((0,0)\), and has the invertible \(y\)-derivative \(\mathbf L\) there. By the implicit function theorem there are \(\sigma_0>0\), \(\varrho_1\in(0,\varrho_0]\) and a smooth map \(\sigma\mapsto y(\sigma)\), \(0\le\sigma<\sigma_0\), \(y(0)=0\), such that for \(\sigma<\sigma_0\) the solutions of \(\widetilde{\mathcal T}(x_0+y;\sigma)=0\) with \(|y|\le\varrho_1\) are exactly \(y=y(\sigma)\); moreover \(y(\sigma)=O(\sigma)\).

Decrease \(d_0\) so that \(\sqrt{d_0}\le\min\{\sigma_0,\sigma_1,\sigma_2\}\) and \(M_1\sqrt{d_0}\le\varrho_1\). Then Lemma~\ref{lemma:apriori} places every upper contact satisfying \eqref{eq:tangency} and \eqref{eq:upper-contact} with \((v,w)\) in the closure of \(\widehat{\mathcal A}_s^{\mathrm J}\) in the ball \(|y|\le\varrho_1\), so it equals \(y(\sigma)\): this is the uniqueness. By Lemma~\ref{lemma:back-in-cell}, \(x_0+y(\sigma)\) corresponds to a pair in \(\widehat{\mathcal A}_s^{\mathrm J}\). This contact satisfies \eqref{eq:upper-contact}: in cell I, \(h_s'(v)=1+O(\sigma)>0\), while in cell II, \(\sigma V_s'(w)=-A_0/\omega_{s+1/2}+O(\sigma)<0\); the remaining sign follows from \(h_s'V_s'=\theta_s'<0\).  Undoing the scaling \eqref{eq:scaling} defines \((v_s^{\mathrm J},w_s^{\mathrm J},\mathsf a_s^{\mathrm J})(\sigma)\), and \eqref{eq:edge-curve-expansion} is \(y(\sigma)=O(\sigma)\) read off in the original variables.
\end{proof}

\begin{corollary}
\label{cor:edge-curves}
For \(d<d_0\) and \(\mathrm J\in\{\mathrm I,\mathrm{II}\}\), the edge cell is the single tangency point,
\[
\widehat{\mathcal A}_s^{\mathrm J}\bigl(a_s^{\mathrm J}(d),d\bigr)=\bigl\{(v_s^{\mathrm J},w_s^{\mathrm J})(\sqrt d)\bigr\},\qquad a_s^{\mathrm J}(d)=\mathsf a_s^{\mathrm J}(\sqrt d),
\]
so the graph of \(\mathsf a_s^{\mathrm J}\) consists of the edge points 
%(in particular it lies in \(\mathcal D^-(\boldsymbol{\omega},\mathbf k)\)), 
and
\begin{equation}
\label{eq:curve-in-D}
c(a,d)\le0\qquad\text{for all }(a,d)\in\mathcal C\text{ with }a\le\mathsf a_s^{\mathrm J}(\sqrt d).
\end{equation}
\end{corollary}

\begin{proof}
The tangency point \((v_s^{\mathrm J},w_s^{\mathrm J})(\sigma)\) is admissible at \(a=\mathsf a_s^{\mathrm J}(\sigma)\), so \(a_s^{\mathrm J}(d)\ge\mathsf a_s^{\mathrm J}(\sqrt d)>a_0\). By Lemma~\ref{lemma:edges} the edge cell \(\mathcal A_s^{\mathrm J}(a_s^{\mathrm J}(d),d)\) is non-empty and consists of upper contacts satisfying \eqref{eq:tangency} and \eqref{eq:upper-contact} in the closure of the cell; by Proposition~\ref{prop:edge-curves} the only such upper contact is \((v_s^{\mathrm J},w_s^{\mathrm J},\mathsf a_s^{\mathrm J})(\sqrt d)\), whence the two identities. Finally, for \eqref{eq:curve-in-D} let \(a\le\mathsf a_s^{\mathrm J}(\sqrt d)=a_s^{\mathrm J}(d)\). By definition of the supremum, or by the non-emptiness of the edge cell when \(a=a_s^{\mathrm J}(d)\), there is \(a_1\in[a,a_s^{\mathrm J}(d)]\) with \(\widehat{\mathcal A}_s^{\mathrm J}(a_1,d)\neq\varnothing\), i.e.\ \((a_1,d)\in\mathcal D^-(\boldsymbol{\omega},\mathbf k)\), and Theorem~\ref{thm:pinning}\ref{thm:pinning-4} gives \(c(a,d)\le0\).
\end{proof}

\paragraph{Recursion for the expansion coefficients.}
The functions $(\mathsf a_s^{\mathrm J}, v_s^{\mathrm J}, w_s^{\mathrm J})$ can be  computed to any order by regular perturbation of the tail relations
\eqref{eq:tail-L}--\eqref{eq:tail-C}
coupled to the scaled system \eqref{eq:N}. 
To see this, we first note
that the denominators appearing in the system
\eqref{eq:tail-L}--\eqref{eq:tail-C}
can be expanded geometrically with respect to $\sigma$ and $A$.
In particular, the coefficients \(L_{i,m},R_{i,m},C_{i,m}\) of \(\sigma^m\) can be computed recursively in \(m\), once the
coefficients \(A_0,\dots,A_m\) are known.
On the other hand, we may write \(x=x_0 + y = \sum_{n\ge0}x_n\sigma^n\), with \(x_0\) from \eqref{eq:leading-I}, resp.\ \eqref{eq:leading-II}) and insert this expression into \(\mathbf L\,y+\mathcal Q(y)+\sigma\mathcal R(y;\sigma)=0\). Collecting powers of \(\sigma\) gives
\begin{equation}
\label{eq:recursion}
\mathbf L x_n=-[\sigma^n]\,\widetilde{\mathcal T}\left(x_0+\sum_{1\le m<n}x_m\sigma^m;\sigma\right),\qquad n\ge1,
\end{equation}
where \([\sigma^n]\) extracts the coefficient of \(\sigma^n\). 
At each step, we compute the Taylor coefficients of \eqref{eq:tail-L}--\eqref{eq:tail-C} that are needed in this expression, which never exceed the order \(n\). The recursion is triangular after eliminating this regular tail system: the coefficient multiplying the new unknown \(x_n\) is the fixed invertible matrix \(\mathbf L\).

\paragraph{Explicit coefficients: the general case.}

Carrying the recursion through \(n=2\) for arbitrary periodic \(\mathbf k\) and \(\boldsymbol{\omega}\) is still short, and already exhibits how branching and diffusion compete. Expanding the scaled maps together with the required tail coefficients gives
\begin{align}
\label{eq:aI-general}
\mathsf a_s^{\mathrm I}(\sqrt d)&=1-2\sqrt{\omega_{s+1/2}\,d}+k_{s+1}\omega_{s+3/2}\,d
+\frac{2k_s\omega_{s+1/2}^2+2k_{s+1}\omega_{s+1/2}\omega_{s+3/2}-k_{s+1}\omega_{s+3/2}^2}{2\sqrt{\omega_{s+1/2}}}\,d^{3/2}\notag\\
&\qquad+O(d^2),\\
\label{eq:aII-general}
\mathsf a_s^{\mathrm{II}}(\sqrt d)&=1-2\sqrt{\omega_{s-1/2}\,d}+k_{s}\omega_{s+1/2}\,d
+\frac{2k_{s-1}\omega_{s-1/2}^2+2k_{s}\omega_{s-1/2}\omega_{s+1/2}-k_{s}\omega_{s+1/2}^2}{2\sqrt{\omega_{s-1/2}}}\,d^{3/2}\notag\\
&\qquad+O(d^2),
\end{align}
with tangency points

\begin{align*}
\text{cell I:}\quad& v=d\,k_s\omega_{s+1/2}\Bigl(1+\sqrt{\omega_{s+1/2}\,d}+O(d)\Bigr),
\qquad
w=\sqrt{\omega_{s+1/2}\,d}\,\Bigl(1+\frac{\omega_{s+3/2}}{4\sqrt{\omega_{s+1/2}}}\sqrt d+O(d)\Bigr);\\
\text{cell II:}\quad& v=1-\sqrt{\omega_{s-1/2}\,d}+O(d^{3/2}),
\qquad
w=\tfrac12\omega_{s+1/2}\,d+O(d^{3/2}).
\end{align*} 
At this order the two formulas \eqref{eq:aI-general} 
and \eqref{eq:aII-general} agree up to the phase shift \(s\mapsto s+1\). It is therefore natural to label the configurations by the edge \(e=i+\tfrac12\) joining the last low site \(i\) to the first high site \(i+1\). To this end, we write
\begin{equation}
\label{eq:a-edge}
\begin{aligned}
\mathsf a_{i+1/2}(d):=1&-2\sqrt{\omega_{i+1/2}\,d}+k_{i+1}\omega_{i+3/2}\,d\\
&+\frac{2k_i\omega_{i+1/2}^2+2k_{i+1}\omega_{i+1/2}\omega_{i+3/2}-k_{i+1}\omega_{i+3/2}^2}{2\sqrt{\omega_{i+1/2}}}\,d^{3/2},
\end{aligned}
\end{equation}
so that \(\mathsf a_s^{\mathrm I}=\mathsf a_{s+1/2}+O(d^2)\) and \(\mathsf a_s^{\mathrm{II}}=\mathsf a_{s-1/2}+O(d^2)\). The three coefficients have a transparent meaning: the leading term is set by the coupling \(\omega_{i+1/2}\) across the interface edge, the \(d\)-term by the coupling \(k_{i+1}\omega_{i+3/2}\) of the first high site to its \(k_{i+1}\) right neighbors  and the \(d^{3/2}\)-term brings in the coupling \(k_i\omega_{i+1/2}\) of the last low site to its right neighbors and a correction to the previous effect.

\paragraph{Explicit coefficients for constant diffusion.}
For \(\omega= \langle1\rangle\) the recursion can be carried out in a somewhat manageable fashion  
for a general branching sequence; we record the results for both cells through order \(d^{5/2}\). 
Keeping \(A=\sum_nA_n\sigma^n\)  in unexpanded form,
the required tail data are given by
\begin{align*}
&L_s=k_{s-1}\Bigl[1+A\sigma+\bigl(A^2-k_{s-1}-1\bigr)\sigma^2+O(\sigma^3)\Bigr],\\
&R_{s+1}=\frac1A\Bigl[1-\frac{1+k_{s+2}}{A}\,\sigma+\frac{k_{s+2}^2+3k_{s+2}+1}{A^2}\,\sigma^2+O(\sigma^3)\Bigr],\\
&C_{s+1}=\frac1{A^2}\Bigl[1+\Bigl(A-\frac{1+2k_{s+2}}{A}\Bigr)\sigma+\Bigl(\frac{3k_{s+2}^2+7k_{s+2}+1-2k_{s+3}}{A^2}-(1+2k_{s+2})\Bigr)\sigma^2+O(\sigma^3)\Bigr],\\
&C_{s+2}=\frac1{A^2}\Bigl[1+\Bigl(A-\frac{1+2k_{s+3}}{A}\Bigr)\sigma+O(\sigma^2)\Bigr],
\end{align*}
obtained from \eqref{eq:tail-L}--\eqref{eq:tail-C} with \(L_{s-1},R_{s+3},C_{s+3}\) needed only to leading order.

\emph{Cell I.} With \(x=(V,W,A)\) and \(x_0=(k_s,1,2)\), solving \eqref{eq:recursion} with \eqref{eq:leading-I} gives
\begin{align*}
&V=k_s\Bigl[1+\sigma+(k_s-k_{s+1}+1)\,\sigma^2+O(\sigma^3)\Bigr],\\
&W=1+\tfrac14\sigma-\tfrac14\bigl(2k_s-3k_{s+1}+k_{s+2}-2\bigr)\sigma^2+O(\sigma^3),\\
&A=2-k_{s+1}\sigma-\bigl(k_s+\tfrac12k_{s+1}\bigr)\sigma^2-\tfrac18\bigl(16k_s-2k_{s+1}^2+2k_{s+1}k_{s+2}+5k_{s+1}\bigr)\sigma^3\\
&\qquad-\tfrac1{16}\bigl(20k_s^2-28k_sk_{s+1}+48k_s-2k_{s+1}^3+4k_{s+1}^2k_{s+2}-k_{s+1}^2\\
&\qquad\qquad\quad-2k_{s+1}k_{s+2}^2+k_{s+1}k_{s+2}+6k_{s+1}\bigr)\sigma^4+O(\sigma^5).
\end{align*}
For the edge curve this implies
\begin{equation}
\label{eq:aI-explicit}
\begin{aligned}
\mathsf a_s^{\mathrm{I}}(\sqrt d)=1&-2\sqrt d+k_{s+1}\,d+\bigl(k_{s}+\tfrac12k_{s+1}\bigr)d^{3/2}
+\tfrac18\bigl(16k_{s}+5k_{s+1}+2k_{s+1}k_{s+2}-2k_{s+1}^2\bigr)d^{2}\\
&+\tfrac1{16}\bigl(48k_s+20k_{s}^2-28k_sk_{s+1}+6k_{s+1}-k_{s+1}^2+k_{s+1}k_{s+2}\\
&\qquad\quad+4k_{s+1}^2k_{s+2}-2k_{s+1}k_{s+2}^2-2k_{s+1}^3\bigr)d^{5/2}+O(d^3),
\end{aligned}
\end{equation}
while the tangency point \(v=\sigma^2V\), \(w=\sigma W\) expands as

\begin{align*}
v&=d\,k_s\bigl[1+\sqrt d+(k_s-k_{s+1}+1)\,d+O(d^{3/2})\bigr],\\
w&=\sqrt d\,\bigl[1+\tfrac14\sqrt d-\tfrac14(2k_s-3k_{s+1}+k_{s+2}-2)\,d+O(d^{3/2})\bigr].
\end{align*}

\emph{Cell II.} With \(x=(Z,W,A)\) and \(x_0=(1,\tfrac12,2)\), solving \eqref{eq:recursion} with \eqref{eq:leading-II} gives
\begin{align*}
&Z=1+\tfrac12(k_s-k_{s-1})\,\sigma^2+\tfrac1{16}(11k_s-16k_{s-1})\,\sigma^3+O(\sigma^4),\\
&W=\tfrac12+\tfrac18(2k_s-2k_{s+1}-1)\,\sigma+\tfrac1{16}\bigl(2k_s^2-4k_sk_{s+1}+2k_{s+1}^2+5k_s+3k_{s+1}+5\bigr)\sigma^2+O(\sigma^3),\\
&A=2-k_s\sigma-\bigl(k_{s-1}+\tfrac12k_s\bigr)\sigma^2+\tfrac18\bigl(2k_s^2-2k_sk_{s+1}-5k_s-16k_{s-1}\bigr)\sigma^3\\
&\qquad+\tfrac1{16}\bigl(2k_s^3-4k_s^2k_{s+1}+2k_sk_{s+1}^2+k_s^2-k_sk_{s+1}\\
&\qquad\qquad\quad+28k_sk_{s-1}-6k_s+12k_{s-1}^2-48k_{s-1}\bigr)\sigma^4+O(\sigma^5).
\end{align*}
For the edge curve this implies
\begin{equation}
\label{eq:aII-explicit}
\begin{aligned}
\mathsf a_s^{\mathrm{II}}(\sqrt d)=1&-2\sqrt d+k_s\,d+\bigl(k_{s-1}+\tfrac12k_s\bigr)d^{3/2}
+\tfrac18\bigl(16k_{s-1}+5k_s+2k_sk_{s+1}-2k_s^2\bigr)d^{2}\\
&+\tfrac1{16}\bigl(48k_{s-1}-12k_{s-1}^2-28k_sk_{s-1}+6k_s-k_s^2+k_sk_{s+1}\\
&\qquad\quad+4k_s^2k_{s+1}-2k_sk_{s+1}^2-2k_s^3\bigr)d^{5/2}+O(d^3),
\end{aligned}
\end{equation}
while the tangency point \(v=1-\sqrt d\,Z\), \(w=d\,W\) expands as

\[
v=1-\sqrt d\,\bigl[1+\tfrac12(k_s-k_{s-1})\,d+O(d^{3/2})\bigr],\qquad
w=d\,\bigl[\tfrac12+\tfrac18(2k_s-2k_{s+1}-1)\sqrt d+O(d)\bigr].
\]

It is interesting to note that
\begin{equation}
\label{eq:cellI-wins}
    \mathsf a_s^{\mathrm{I}}(\sqrt d)
    - \mathsf a_{s+1}^{\mathrm{II}}(\sqrt d)
    = 2k_s^2\,d^{5/2}+O(d^3)\ >\ 0 .
\end{equation}

In particular, for equal diffusion  the envelope of all edge curves close to the cornerpoint is always attained by a cell-I curve.
Symbolic manipulation suggests that this difference becomes
$2k_s^2\,\omega_{s+1/2}^{5/2}\,d^{5/2}+O(d^3)$ for the general case,
but we do not pursue this further.

With these expansions in hand we are ready to state the final result of this section, which concerns the boundary $\mathcal{D}^-$ in the corner region. In particular,
the preceding analysis allows us to simply take the right-most edge resulting from all the phases and the two cells to trace the right boundary of \(\mathcal D^-(\boldsymbol{\omega},\mathbf k)\). In addition, our expansions above tell us (in many cases) which configuration attains it.

\begin{corollary}
\label{cor:envelope}
For \(0<d<d_0\) put
\[
\mathsf a^\ast(d):=\max_{s\in[1,p]_{\mathbb N}}\ \max_{\mathrm J\in\{\mathrm I,\mathrm{II}\}}\ \mathsf a_s^{\mathrm J}(\sqrt d).
\]
Then
\begin{enumerate}[label=(\roman*), ref=(\roman*)]
\item \label{cor:envelope-1} \(\mathcal D^-(\boldsymbol{\omega},\mathbf k)\cap\mathcal C\subset\{(a,d):\ a\le\mathsf a^\ast(d)\}\) and \((\mathsf a^\ast(d),d)\in\mathcal D^-(\boldsymbol{\omega},\mathbf k)\); that is,
\[
\sup\{a:\ (a,d)\in\mathcal D^-(\boldsymbol{\omega},\mathbf k)\}=\mathsf a^\ast(d)\qquad\text{for }0<d<d_0 ,
\]
and the supremum is attained.
\item \label{cor:envelope-2} \(c(a,d)\le0\) for all \((a,d)\in\mathcal C\) with \(a\le\mathsf a^\ast(d)\).
\item \label{cor:envelope-3} \(\mathsf a^\ast(d)=\max_{i\in[1,p]_{\mathbb N}}\mathsf a_{i+1/2}(d)+O(d^2)\), with \(\mathsf a_{i+1/2}\) from \eqref{eq:a-edge}. The maximum at this order is decided lexicographically: among all edges, those with the smallest diffusion \(\omega_{i+1/2}=\omega_{\min}\); among these, those with the largest right coupling \(k_{i+1}\omega_{i+3/2}\) of the first high site; among these, those with the largest \(d^{3/2}\)-coefficient in \eqref{eq:a-edge}. In particular
\[
\mathsf a^\ast(d)=1-2\sqrt{\omega_{\min}\,d}+\Bigl(\max_{i:\ \omega_{i+1/2}=\omega_{\min}}k_{i+1}\omega_{i+3/2}\Bigr)\,d+O(d^{3/2}).
\]
\item \label{cor:envelope-4} For constant diffusion \(\omega=\langle1\rangle\),
\[
\mathsf a^\ast(d)=1-2\sqrt d+k_{\max}\,d+\Bigl(\max_{i:\ k_{i+1}=k_{\max}}k_i+\tfrac12k_{\max}\Bigr)\,d^{3/2}+O(d^2).
\]
The boundary is hence set by an interface placed just before a site with the maximal branching number \(k_{\max}\), and, if there are several, before the one whose predecessor branches most. 
\end{enumerate}
\end{corollary}

\begin{proof}
\emph{\ref{cor:envelope-1}} If \((a,d)\in\mathcal D^-(\boldsymbol{\omega},\mathbf k)\cap\mathcal C\), then \(\mathcal A_s^-(a,d)\neq\varnothing\) for some \(s\), and by Lemmas~\ref{lemma:four-cells}, \ref{lemma:restriction}, \ref{lemma:III} and \ref{lemma:IV} one of \(\widehat{\mathcal A}_s^{\mathrm I}\), \(\widehat{\mathcal A}_s^{\mathrm{II}}\), \(\widehat{\mathcal A}_{s+1}^{\mathrm I}\) is non-empty at \((a,d)\); hence \(a\le a_{s'}^{\mathrm J}(d)=\mathsf a_{s'}^{\mathrm J}(\sqrt d)\le\mathsf a^\ast(d)\) for some \((s',\mathrm J)\) by Corollary~\ref{cor:edge-curves}. Conversely, the edge cell of the maximizing pair \((s,\mathrm J)\) is non-empty at \(a=\mathsf a^\ast(d)\), so \((\mathsf a^\ast(d),d)\in\mathcal D^-(\boldsymbol{\omega},\mathbf k)\). \emph{\ref{cor:envelope-2}} is \eqref{eq:curve-in-D} for the maximizing pair. \emph{\ref{cor:envelope-3}} By \eqref{eq:aI-general}--\eqref{eq:a-edge} every \(\mathsf a_s^{\mathrm J}\) equals some \(\mathsf a_{i+1/2}\) up to \(O(d^2)\), and the maximum of finitely many expansions in powers of \(\sqrt d\) is, for small \(d\), attained by the lexicographically largest coefficient sequence; the stated readings of the coefficients follow from \eqref{eq:a-edge}. \emph{\ref{cor:envelope-4}} is the special case \(\omega=\langle 1\rangle \) of \ref{cor:envelope-3}, with the \(d^{3/2}\)-coefficient \(k_i+\tfrac12k_{i+1}\) read off from \eqref{eq:aI-explicit}--\eqref{eq:aII-explicit}.
\end{proof}

\paragraph{Beyond the corner: numerical continuation.}
The edge curves are defined through the tangency system \eqref{eq:tangency}, which is a general smooth nonlinear three-component system. In particular, the solution branches through the points \((v_s^{\mathrm J},w_s^{\mathrm J},\mathsf a_s^{\mathrm J})(\sqrt d)\) can  be followed numerically by standard pseudo-arclength continuation to values of \(d\) that can potentially lie far outside \(\mathcal C\). Outside the corner the structure described by Lemma~\ref{lemma:four-cells} gradually dissolves. For example, when \(F_s(a/3;w)\) turns positive property \ref{it:S} is lost. This ends the four-cell description and tangencies can disappear or emerge. We follow each branch by pseudo-arclength continuation until it either reaches \(d=0.5\) or leaves the admissible region $(v,w) \in (0,1) \times (0, w_s^*)$. Some branches fold back in \(d\) and then descend all the way to the origin \((a,d)=(0,0)\), where \(v,w\to0\). The fold appears as a cusp in the \((a,d)\)-plane, but is smooth in the full \((v,w,a,d)\)-space.

\begin{figure}[t]
\centering
\includegraphics[height=0.3\textheight]{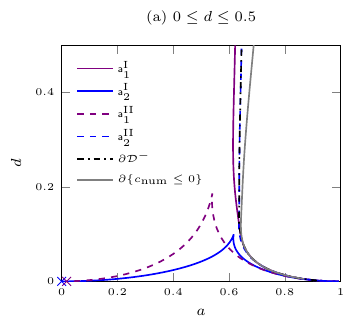}\hfill\includegraphics[height=0.3\textheight]{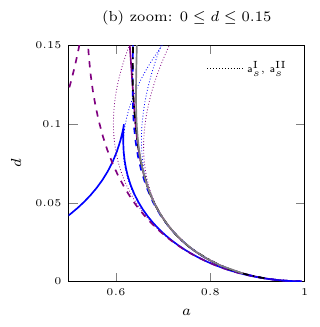}
\caption{Numerical continuation of the tangency curves for \(\mathbf k=\langle1,2\rangle\), \(\omega=\langle 1 \rangle\). (a) The branches of solutions of the tangency system \eqref{eq:tangency} followed in \(\sqrt d\) by pseudo-arclength continuation until $d$ leaves $[0,0.5]$ or the branch leaves the admissible region (crosses). (b) Zoom into \(0\le d\le0.15\), with the truncations of \eqref{eq:aI-explicit}--\eqref{eq:aII-explicit} through \(d^{5/2}\) added (dotted). For small $d$ the right boundary of \(\mathcal D^-(\boldsymbol\omega,\mathbf k)\) is carried by the cell-I branch of phase \(s=1\) and, from \(d\approx0.14\) on, by the cell-II branch of phase \(s=2\). Gray: the boundary \(\partial\{c_{\rm num}\le0\}\) of the actual pinning region, computed by direct simulation. It is tracked closely by \(\partial\mathcal D^-\) up to \(d\approx0.05\) with the relative accuracy $|a_{\rm err}| \le 10^{-4} (1-a)$.}
\label{fig:expansion-check}
\end{figure}

\begin{figure}[t]
\centering
\includegraphics[height=0.3\textheight]{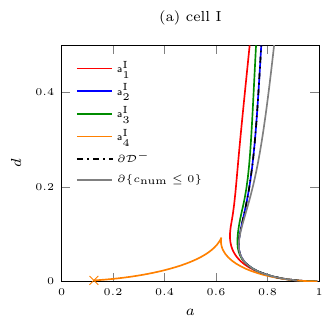}\hfill\includegraphics[height=0.3\textheight]{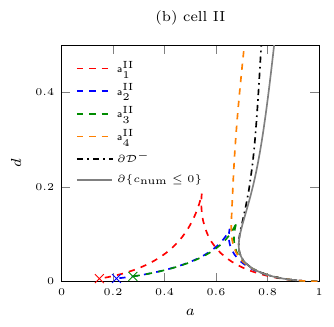}

\caption{The same continuation as in Figure~\ref{fig:expansion-check}(a), for \(\mathbf k=\langle1,2,2,2\rangle\), \(\omega=\langle 1 \rangle\), with cell I in (a) and cell II in (b); the four phases \(s=1,\ldots,4\) are drawn in red, blue, green and orange. Branches are followed to \(d=0.5\) or until they leave the admissible region (crosses). The black dash-dotted curve is \(\partial\mathcal D^-\), which is carried entirely by \(\mathsf a_2^{\rm I}\).
Here
\(\partial\mathcal D^-\) tracks the boundary of the numerical pinning region with the relative accuracy $|a_{\rm err}| \le 3 \cdot 10^{-5} (1-a)$ up to $d = 0.05$.}
\label{fig:continuation-k1222}
\end{figure}

Figures~\ref{fig:expansion-check} and 
\ref{fig:continuation-k1222} explore this
for \(\mathbf k=\langle1,2\rangle\)
and \(\mathbf k=\langle1,2,2,2\rangle\).
In the latter case the bend-back is more pronounced. In both cases
the right boundary of \(\mathcal D^-(\boldsymbol\omega,\mathbf k)\) is fully
described by the envelope of our continuation curves. In particular, the description
of Corollary \ref{cor:envelope} extends
far beyond the corner region.
In both cases the asymptotic expansions \eqref{eq:aI-explicit}--\eqref{eq:aII-explicit} follow the branches closely up to \(d\approx0.05\) and then depart; the continued branches, not the truncated series, are the right objects to compare with the numerically computed pinning boundaries at intermediate values of \(d\).

%%%%%%%%%%%%%%%%%%%%%%%%%%%%%%%%%%%%%%%%%%%%%%%%%%%%%%%%%%%%%%%%%%%
%%%%%%%%%%%%%%%%%%%%%%%%%%%%%%%%%%%%%%%%%%%%%%%%%%%%%%%%%%%%%%%%%%%
% LEFT PROPAGATION FOR LARGE DIFFUSION
%%%%%%%%%%%%%%%%%%%%%%%%%%%%%%%%%%%%%%%%%%%%%%%%%%%%%%%%%%%%%%%%%%%
%%%%%%%%%%%%%%%%%%%%%%%%%%%%%%%%%%%%%%%%%%%%%%%%%%%%%%%%%%%%%%%%%%%

\section{Left-propagation \texorpdfstring{$c<0$}{} for large diffusion}
\label{sec:large-d}

In this section we develop the machinery needed to establish Theorem~\ref{t:left-travel}. The cornerstone of our approach is the construction of a parameter-dependent sub-solution family $\boldsymbol\Psi_{\mathbf{m},A}(\xi;\boldsymbol\omega)$ to~\eqref{e:main-mfde} which satisfies the assumptions of Corollary~\ref{cor:comparison-MFDE}, and thus ensures the leftward propagation $c<0$ of the actual solution to~\eqref{e:main-mfde}. We construct this family in a way that $\mathfrak{D}_{\boldsymbol\omega,\mathbf{k}}[\boldsymbol\Psi_{\mathbf{m},A}(\, \cdot \,;\boldsymbol\omega)](\xi) \ge 0 $ is nonnegative. The non-positivity of the MFDE residual $\boldsymbol{\mathcal{I}}_{a,d,\boldsymbol\omega,\mathbf{k}}[c,\boldsymbol\Psi_{\mathbf{m},A}(\, \cdot \,;\boldsymbol\omega)](\xi)$ can be achieved by choosing $d>0$ large enough while leaving the other parameters fixed. Our construction requires $\widehat k = \prod_{i=1}^p k_i >1$. 

The second part of this section is devoted to the analysis of the quantities $\mu(\boldsymbol\omega,\mathbf{k})$ and $\nu(a)$, given by~\eqref{e:mu}--\eqref{e:nu}, both of which are essential parts of the left-traveling threshold $d^*(a;\boldsymbol\omega,\mathbf{k})$,~\eqref{e:dStar}. The expressions result from maximization and/or minimization of some quantities and finding their value (especially of $\mu(\boldsymbol\omega,\mathbf{k})$,~\eqref{e:mu}) is not a straightforward task for $p > 2$.

\subsection{Sub-solution construction and proof of Theorem~\ref{t:left-travel}}
We now start with the description of the sub-solution form. Each element of the sub-solution consists of a constant part and an exponential part, together with a cubic part that forms a smooth and monotone connection between them. The exponential part is denoted by $\boldsymbol\lambda_{\mathbf{m}}(\xi;\boldsymbol\omega) = \langle \lambda_{j,\mathbf{m}}(\xi;\boldsymbol\omega) \rangle_\dInt{j}$ with 
\begin{equation} \label{e:lambda}
    \lambda_{j,\mathbf{m},A}(\xi;\boldsymbol\omega) := A\left(1- \alpha_{j,\mathbf{m}}(\boldsymbol\omega)  \beta_{j,\mathbf{m}}(\xi) \right) ,
\end{equation}
where $\boldsymbol\alpha_\mathbf{m}(\boldsymbol\omega) = \langle \alpha_{j,\mathbf{m}}(\boldsymbol\omega) \rangle_\dInt{j}$ and $\boldsymbol\beta_\mathbf{m}(\xi) = \langle \beta_{j,\mathbf{m}}(\xi) \rangle_\dInt{j}$ are defined by
\begin{align} 
    \alpha_{j,\mathbf{m}}(\boldsymbol\omega) &:=\frac{1}{p} \left(\sum_{n=1}^{p} \frac{\widehat{\omega}}{\omega_{j+n-1/2}} \prod_{i=1}^{p-n} m_{j-i}  \right), \label{e:alpha} \\
    \beta_{j,\mathbf{m}}(\xi) &:= \left( \prod_{i=1}^{p} \sqrt[p]{m_{j-i}^{i-p}} \right) \sqrt[p]{\widehat{{m}}}^{-\xi}. \label{e:beta}
\end{align}

Note that~\eqref{e:alpha} is precisely~\eqref{e:alpha-intro} from the statement of Theorem~\ref{t:left-travel}. We include the expression here for the convenience of the reader. 

We summarize properties of the functions $\boldsymbol\alpha_{\mathbf{m}}(\boldsymbol\omega), \boldsymbol\beta_{\mathbf{m}}(\xi)$, and $\boldsymbol\lambda_{\mathbf{m},A}(\xi;\boldsymbol\omega)$ in the following technical lemma. Recall the definition \eqref{e:omega}
\[
\Omega_\mathbf{k} := \left\{ \mathbf{m} \in \mathbb{R}^p \, | \, 0<m_j < k_j \text{ for all } \dInt{j}, \; \widehat{m}>1 \right\}. 
\]

\begin{lemma} \label{l:lambda}
The following properties hold:
   \begin{enumerate}[label=(\roman*)]
       \item \label{l:alpha-rotation}
            $ 
            \alpha_{j+1,\mathbf{m}}(\boldsymbol\omega) = m_{j} \alpha_{j,\mathbf{m}}(\boldsymbol\omega) + \dfrac{\widehat{\omega}(1- \widehat{m})}{p \omega_{j+1/2}}.
            $ 
        \item \label{l:beta-rotation} 
            $
                m_j \beta_{j+1,\mathbf{m}}(\xi) = \beta_{j,\mathbf{m}}(\xi)\sqrt[p]{\widehat{m}}.
            $
        \item \label{l:psi-zero}
            Let $A \in (0,1]$, and $\mathbf{m} \in \Omega_\mathbf{k}$, then each $\lambda_{j,\mathbf{m},A}(\xi;\boldsymbol\omega)$ has a unique zero 
            
            \begin{equation}\label{e:xi-j}
                \xi^*_j := \frac{p \log{(\alpha_{j,\mathbf{m}}(\boldsymbol\omega) \beta_{j,\mathbf{m}}(0))}}{\log{(\widehat{m})}}. 
            \end{equation}
        \item \label{l:lambda-zeros-distance}
        Under the assumptions of~\ref{l:psi-zero} we have
        $
             \xi_{j+1}^* - \xi_j^*  < 1.
        $

        \item \label{l:exp-identity}
        \begin{equation} \label{e:exp-identity}
            (\omega_{j+1/2} m_j+\omega_{j-1/2}) \sqrt[p]{\widehat m}^{\xi_j^*} = \omega_{j-1/2}\sqrt[p]{\widehat m}^{\xi_{j-1}^*+1} + \omega_{j+1/2}m_j \sqrt[p]{\widehat m}^{\xi_{j+1}^*-1}. 
        \end{equation}

        \item \label{l:d-lambda}
        
        $
        \mathfrak{D}_{\boldsymbol\omega,\mathbf{k}}[\boldsymbol\lambda_{\mathbf{m},A}(\, \cdot \,;\boldsymbol\omega)]_j(\xi) = \frac{A\widehat{\omega}}{m_j p} (\widehat{{m}}-1)(k_j - m_j) \beta_{j,\mathbf{m}}(\xi). 
        $ 
   \end{enumerate} 
\end{lemma}

\begin{proof}
    \begin{enumerate}[label=(\roman*)]
       \item 
           The proof follows directly from the definition of $\boldsymbol\alpha_{\mathbf{m}}(\boldsymbol\omega)$, \eqref{e:alpha}:
           \begin{align*}
               \frac{p}{\widehat\omega} \alpha_{j+1,\mathbf{m}}(\boldsymbol\omega)& = \sum_{n=1}^{p} \frac{1}{\omega_{j+n+1/2}} \prod_{i=1}^{p-n} m_{j-i+1}   = \frac{1}{\omega_{j+1/2}} + \sum_{n=1}^{p-1} \frac{1}{\omega_{j+n+1/2}} \prod_{i=1}^{p-n} m_{j-i+1},  \\
               \frac{p}{\widehat\omega} m_j \alpha_{j,\mathbf{m}}(\boldsymbol\omega) &= \sum_{n=1}^{p} \frac{m_j}{\omega_{j+n-1/2}} \prod_{i=1}^{p-n} m_{j-i} = \frac{\widehat{m}}{\omega_{j+1/2}} + \sum_{n=2}^{p} \frac{m_j}{\omega_{j+n-1/2}} \prod_{i=1}^{p-n} m_{j-i} \\
               &= \frac{\widehat{m}}{\omega_{j+1/2}} + \sum_{n=2}^{p} \frac{1}{\omega_{j+n-1/2}} \prod_{i=0}^{p-n} m_{j-i} 
               = \frac{\widehat{m}}{\omega_{j+1/2}} + \sum_{n=1}^{p-1} \frac{1}{\omega_{j+n+1/2}} \prod_{i=0}^{p-n-1} m_{j-i} \\
               &= \frac{\widehat{m}}{\omega_{j+1/2}} + \sum_{n=1}^{p-1} \frac{1}{\omega_{j+n+1/2}} \prod_{i=1}^{p-n} m_{j-i+1}.
           \end{align*}
           
        \item
            The proof is carried out using the definition of $\boldsymbol\beta_{\mathbf{m}}(\xi)$, \eqref{e:beta}:
            \begin{align*}
                \frac{\beta_{j+1,\mathbf{m}}(\xi)}{\beta_{j,\mathbf{m}}(\xi)} 
                &= \frac{\prod_{i=1}^{p} \sqrt[p]{m_{j+1-i}^{i-p}}}{\prod_{i=1}^{p} \sqrt[p]{m_{j-i}^{i-p}}}
                =  \frac{\prod_{i=1}^{p} \sqrt[p]{m_{j+1-i}^{i-p}}}{\prod_{i=2}^{p+1} \sqrt[p]{m_{j-i+1}^{i-1-p}}}
                = \frac{\sqrt[p]{m_j^{1-p}}}{\sqrt[p]{1}} \prod_{i=2}^p \sqrt[p]{m_{j-i+1}^{(i-p)-(i-p-1)}}  \\
                &= \frac{\sqrt[p]{m_j}}{m_j} \prod_{i=2}^p \sqrt[p]{m_{j-i+1}} = \frac{\sqrt[p]{\widehat{m}}}{m_j}. 
            \end{align*}

        \item 
            The uniqueness is a consequence of the monotonicity of the exponential function. 
            In order to prove the identity in~\eqref{e:xi-j} the crucial insight  is that 
            \[
                \beta_{j,\mathbf{m}}(\xi)  = \beta_{j,\mathbf{m}}(0) \sqrt[p]{\widehat{m}}^{-\xi}, 
            \]
            cf.~\eqref{e:beta}. The rest follows from direct manipulation 
            \begin{align*}
                A(1-\alpha_{j,\mathbf{m}}(\boldsymbol\omega) \beta_{j,\mathbf{m}}(\xi)) &= 0 \\
                \alpha_{j,\mathbf{m}}(\boldsymbol\omega) \beta_{j,\mathbf{m}}(\xi) &= 1 \\
                \log(\alpha_{j,\mathbf{m}}(\boldsymbol\omega) \beta_{j,\mathbf{m}}(0)) - \frac{\xi}{p} \log(\widehat{m}) &= 0,
            \end{align*}
            and~\eqref{e:xi-j} follows. 

        \item
            Using the definition of $\xi_j^*$ and logarithmic identities, we get
            \begin{align*}
                \xi_{j+1}^* - \xi_j^* &= \frac{p \log(\alpha_{j+1,\mathbf{m}}(\boldsymbol\omega) \beta_{j+1,\mathbf{m}}(0))}{\log(\widehat{m})} - \frac{p \log(\alpha_{j,\mathbf{m}}(\boldsymbol\omega) \beta_{j,\mathbf{m}}(0))}{\log(\widehat{m})}= \frac{\log \left( \frac{\alpha_{j+1,\mathbf{m}}(\boldsymbol\omega) \beta_{j+1,\mathbf{m}}(0)}{\alpha_{j,\mathbf{m}}(\boldsymbol\omega) \beta_{j,\mathbf{m}}(0)} \right)}{\log(\sqrt[p]{\widehat{m}})}. 
            \end{align*}
            It is now sufficient to show that 
            \[
            \frac{\alpha_{j+1,\mathbf{m}}(\boldsymbol\omega) \beta_{j+1,\mathbf{m}}(0)}{\alpha_{j,\mathbf{m}}(\boldsymbol\omega) \beta_{j,\mathbf{m}}(0)} < \sqrt[p]{\widehat{m}}. 
            \]
            One has

            \[
                \begin{aligned}
                & \frac{\alpha_{j+1,\mathbf{m}}(\boldsymbol\omega) \beta_{j+1,\mathbf{m}}(0)}{\alpha_{j,\mathbf{m}}(\boldsymbol\omega)\beta_{j,\mathbf{m}}(0)} = \frac{\alpha_{j+1,\mathbf{m}}(\boldsymbol\omega)\sqrt[p]{\widehat{m}}}{\alpha_{j,\mathbf{m}}(\boldsymbol\omega) m_j} =\frac{\left(\alpha_{j,\mathbf{m}}(\boldsymbol\omega) m_j + \frac{\widehat{\omega
                }(1-\widehat{m})}{p \omega_{j+1/2}}\right)\sqrt[p]{\widehat{m}}}{\alpha_{j,\mathbf{m}}(\boldsymbol\omega) m_j} \\
                &\qquad = \left(1 + \frac{\widehat{\omega}(1-\widehat{m})}{p \omega_{j+1/2} \alpha_{j,\mathbf{m}}(\boldsymbol\omega) m_j} \right)\sqrt[p]{\widehat{m}} < \sqrt[p]{\widehat{m}},
                \end{aligned}
            \]
            in which we used Lemmas~\ref{l:lambda}\ref{l:alpha-rotation},~\ref{l:lambda}\ref{l:beta-rotation}. The statement follows from the assumption $\widehat{m}>1$. 
    
            \item 
            We can alternatively express $\xi_j^*$ as 
            \[
                \xi_j^* = \log_{\sqrt[p]{\widehat m}}(\alpha_{j,\mathbf{m}}(\boldsymbol\omega) \beta_{j,\mathbf{m}}(0)). 
            \]
            From the proof of Lemma~\ref{l:lambda}\ref{l:lambda-zeros-distance} and using the  combination of Lemmas~\ref{l:lambda}\ref{l:alpha-rotation},~\ref{l:lambda}\ref{l:beta-rotation} analogously we have 
            \begin{align*}
                \frac{\alpha_{j+1,\mathbf{m}}(\boldsymbol\omega) \beta_{j+1,\mathbf{m}}(0)}{\alpha_{j,\mathbf{m}}(\boldsymbol\omega) \beta_{j,\mathbf{m}}(0)} &= \frac{(\alpha_{j,\mathbf{m}}(\boldsymbol\omega) m_j + \frac{\widehat{\omega}(1-\widehat{m})}{p \omega_{j+1/2}})\sqrt[p]{\widehat{m}}}{\alpha_{j,\mathbf{m}}(\boldsymbol\omega) m_j} , \\
                \frac{\alpha_{j-1,\mathbf{m}}(\boldsymbol\omega) \beta_{j-1,\mathbf{m}}(0)}{\alpha_{j,\mathbf{m}}(\boldsymbol\omega) \beta_{j,\mathbf{m}}(0)} &= 
                \frac{\alpha_{j,\mathbf{m}}(\boldsymbol\omega)- \frac{\widehat{\omega}(1-\widehat m)}{p\omega_{j-1/2}}}{\alpha_{j,\mathbf{m}}(\boldsymbol\omega) \sqrt[p]{\widehat m}}.
            \end{align*}
            The right-hand side of~\eqref{e:exp-identity} is then 
           \begin{align*}
                &\omega_{j-1/2}\sqrt[p]{\widehat m}^{\xi_{j-1}^*+1} + \omega_{j+1/2}m_j \sqrt[p]{\widehat m}^{\xi_{j+1}^*-1}  \\ 
                &\qquad =  \omega_{j-1/2}\sqrt[p]{\widehat m} \, \alpha_{j-1,\mathbf{m}}(\boldsymbol\omega) \beta_{j-1,\mathbf{m}}(0) + \omega_{j+1/2}\frac{m_j}{\sqrt[p]{\widehat m}}\alpha_{j+1,\mathbf{m}}(\boldsymbol\omega)( \beta_{j+1,\mathbf{m}}(0) \\
                &\qquad= \left( \frac{\omega_{j-1/2}\left(\alpha_{j,\mathbf{m}}(\boldsymbol\omega)- \frac{\widehat{\omega}(1-\widehat m)}{p\omega_{j-1/2}} \right) + \omega_{j+1/2} \left( \alpha_{j,\mathbf{m}}(\boldsymbol\omega) m_j + \frac{\widehat{\omega}(1-\widehat{m})}{p\omega_{j+1/2}} \right)}{\alpha_{j,\mathbf{m}}(\boldsymbol\omega)}\right)\alpha_{j,\mathbf{m}}(\boldsymbol\omega) \beta_{j,\mathbf{m}}(0) \\
                &\qquad = (\omega_{j+1/2} m_j +\omega_{j-1/2} ) \alpha_{j,\mathbf{m}}(\boldsymbol\omega) \beta_{j,\mathbf{m}}(0) = (\omega_{j+1/2} m_j +\omega_{j-1/2} ) \sqrt[p]{\widehat m}^{\xi_j^*}. 
            \end{align*}  

        \item 
            The proof of this claim is strictly computational. Start with  
           \begin{equation} \label{e:d-cal-lem-1} 
                \begin{aligned}
            &\mathfrak{D}_{\boldsymbol\omega,\mathbf{k}}[\boldsymbol\lambda_{\mathbf{m},A}(\, \cdot \,;\boldsymbol\omega)]_j(\xi) =  -A \big( 
                    k_j \omega_{j+1/2} \alpha_{j+1,\mathbf{m}}(\boldsymbol\omega)\beta_{j+1,\mathbf{m}}(\xi+1)  
                    \\ &\qquad - (k_j \omega_{j+1/2} + \omega_{j-1/2} ) \alpha_{j,\mathbf{m}}(\boldsymbol\omega) \beta_{j,\mathbf{m}}(\xi)
                     + \omega_{j-1/2}\alpha_{j-1,\mathbf{m}}(\boldsymbol\omega) \beta_{j-1,\mathbf{m}}(\xi-1)\big),
                \end{aligned}
            \end{equation} 
            which directly follows from the definition of the operator $\mathfrak{D}_{\boldsymbol\omega,\mathbf{k}}$,~\eqref{e:lapl-cont}, and the definition of $\boldsymbol\lambda_{\mathbf{m},A}(\xi;\boldsymbol\omega)$,~\eqref{e:lambda}. From Lemma~\ref{l:lambda}\ref{l:beta-rotation} and the definition of $\boldsymbol\beta_{\mathbf{m}}(\xi)$ we have 
            \begin{align*}
                \beta_{j-1,\mathbf{m}}(\xi-1) &= \beta_{j-1,\mathbf{m}}(\xi) \sqrt[p]{\widehat{m}} =  m_{j-1} \beta_{j,\mathbf{m}}(\xi), \\
                \beta_{j+1,\mathbf{m}}(\xi+1) &= \beta_{j+1,\mathbf{m}}(\xi) \frac{1}{\sqrt[p]{\widehat{m}}} =  \frac{1}{m_{j}} \beta_{j,\mathbf{m}}(\xi).
            \end{align*}
            Plugging these two expressions into~\eqref{e:d-cal-lem-1} yields 
            \begin{equation}\label{e:d-cal-lem-2}
                \begin{aligned}
                & \mathfrak{D}_{\boldsymbol\omega,\mathbf{k}}[\boldsymbol\lambda_{\mathbf{m},A}(\, \cdot \,;\boldsymbol\omega)](\xi) 
                    = - \frac{1}{m_j} A \beta_{j,\mathbf{m}}(\xi) \big(
                        \omega_{j+1/2} k_j(\alpha_{j+1,\mathbf{m}}(\boldsymbol\omega) 
                        -m_j \alpha_{j,\mathbf{m}}(\boldsymbol\omega) )
                         \\
                         &\qquad + \omega_{j-1/2} m_j (m_{j-1} \alpha_{j-1,\mathbf{m}}(\boldsymbol\omega) 
                        - \alpha_{j,\mathbf{m}}(\boldsymbol\omega)
                        ) 
                    \big).
                \end{aligned}
            \end{equation}
            From Lemma~\ref{l:lambda}\ref{l:alpha-rotation} we have 
            
             \[
                 \alpha_{j+1,\mathbf{m}}(\boldsymbol\omega) - m_{j} \alpha_{j,\mathbf{m}}(\boldsymbol\omega) = \frac{\widehat{\omega} (1 - \widehat{m})}{p \omega_{j+1/2}}, \qquad \qquad
                m_{j-1} \alpha_{j-1,\mathbf{m}}(\boldsymbol\omega) - \alpha_{j,\mathbf{m}}(\boldsymbol\omega) = \frac{\widehat{\omega}(\widehat{m} -1)}{p \omega_{j-1/2}}.
            \]
            Plugging both into~\eqref{e:d-cal-lem-2} concludes the proof.
   \end{enumerate} 
\end{proof}

\begin{figure}[t]
    \centering
    \includegraphics[width=0.8\linewidth]{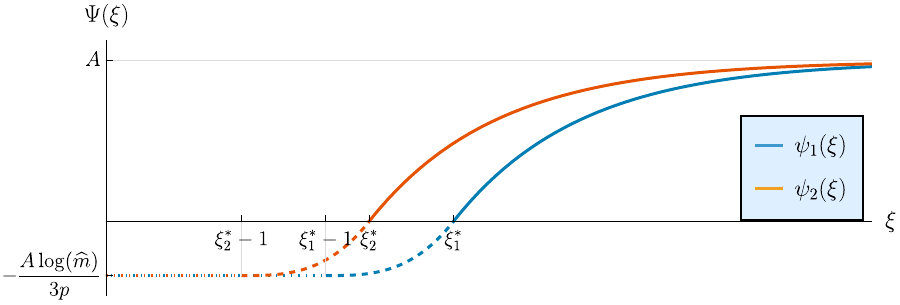}
    \caption{Schematic representation of profiles $\Psi_{j,\mathbf{m},A}$ defined by Eq.~\eqref{e:psi}. Each of the profiles consists of a constant part (dotted), a cubic transition part $\pi_i$, \eqref{e:pi} (dashed), and an exponential part $\lambda_i$, \eqref{e:lambda} (full). The values $\xi^*_i$ are zeros of the exponential functions defined in Lemma~\ref{l:lambda}\ref{l:psi-zero}. They are not further than one unit apart, see Lemma~\ref{l:lambda}\ref{l:lambda-zeros-distance}. Finally, the monotonicity of the profiles, the common base of the exponentials $\lambda_i$, and Lemma~\ref{l:psi}\ref{l:psi-at-cub} ensure that the profiles are ordered.}    
    \label{f:left-travel-psi}
\end{figure}

With these properties of the exponentials $\boldsymbol\lambda_{\mathbf{m},A}(\xi;\boldsymbol\omega)$ 
in hand,
we can finally construct the bounded sub-solution profiles $\boldsymbol\Psi_{\mathbf{m},A}(\xi;\boldsymbol\omega)$ by
\begin{align} \label{e:psi}
\Psi_{j,\mathbf{m},A}(\xi;\boldsymbol\omega) = \begin{cases}
        -A   \dfrac{ \log(\widehat{m})}{3 p} , & \xi \le \xi^*_j - 1, \\
        A \pi_{j,\mathbf{m}}(\xi), &\xi^*_j-1 < \xi \le \xi_j^*, \\
        \lambda_{j,\mathbf{m},A}(\xi;\boldsymbol\omega), &  \xi_j^*< \xi,
    \end{cases}
\end{align}
where $\pi_{j,\mathbf{m}}(\xi)$ is a third-order polynomial
\begin{equation} \label{e:pi}
    \pi_{j,\mathbf{m}}(\xi) :=  \dfrac{\log{\widehat{m}}}{3p} 
    \left[
        \left( 
            \xi - (\xi_j^*-1)
        \right)^3
        -1
    \right ], 
\end{equation}
see Figure~\ref{f:left-travel-psi}. Note that the polynomial $\pi_j$ is implicitly dependent on $\boldsymbol\omega$ through $\xi_j^*$. Since this dependence does not influence further calculations and considerations, we omit it for the sake of clarity. 

Based on the value of the variable $\xi$ with respect to the zero $\xi_j^*$ we speak about $\Psi_{j,\mathbf{m},A}(\xi;\boldsymbol\omega)$ being in the constant, cubic, or exponential regime. Also note that the use of a cubic as a transition between the constant and exponential regime originates purely from the fact that it is the simplest polynomial to be able to connect them smoothly. 

We formally prove two properties of the profiles $\Psi_{j,\mathbf{m},A}(\xi;\boldsymbol\omega)$, namely their smoothness and their ordering in the cubic regime. 

\begin{lemma} \label{l:psi}
    Let $\mathbf{m} \in \Omega_\mathbf{k}$, then
    \begin{enumerate}[label=({\roman*})]
        \item \label{l:psi-smooth}
            The profiles $\Psi_{j,\mathbf{m},A}(\xi;\boldsymbol\omega)$,~\eqref{e:psi}, are well defined and $\Psi_{j,\mathbf{m},A}(\xi;\boldsymbol\omega) \in C^1(\mathbb{R})$. 
        \item \label{l:psi-at-cub}
            Let $\xi_j^*-1 < \xi < \xi_j^*$, then 
            \[
                \Psi_{j+1,\mathbf{m},A}(\xi+1;\boldsymbol\omega) - \Psi_{j,\mathbf{m},A}(\xi;\boldsymbol\omega) > 0. 
            \]
    \end{enumerate}
\end{lemma}
\begin{proof}
    \begin{enumerate}[label=({\roman*})]
    \item
        Since any $\Psi_{j,\mathbf{m},A}(\xi;\boldsymbol\omega)$ is defined piecewise from analytic functions, we focus only on the behavior at the points $\xi_j^*-1$ and $\xi_j^*$. The continuity of the profile follows from the identities
        \[
            - A\pi_{j,\mathbf{m}}(\xi_j^*-1) = \dfrac{-A \log{\widehat{m}}}{3p}, \qquad A\pi_{j,\mathbf{m}}(\xi_j^*) = 0 = \lambda_{j,\mathbf{m},A}(\xi_j^*;\boldsymbol\omega). 
        \]
        
        To show smoothness it is sufficient to note that the derivative of $\pi_j$ vanishes at $\xi_j^*-1$ and that the derivatives of $A \pi_{j,\mathbf{m}}(\xi)$ and $\lambda_{j,\mathbf{m},A}(\xi;\boldsymbol\omega)$ coincide at $\xi_j^*$. 
        Compute
        \[ 
        \pi_{j,\mathbf{m}}'(\xi) = \dfrac{\log{\widehat{m}}}{p}(\xi-(\xi_j^*-1))^2,
        \]
        thus 
        \begin{equation}\label{e:l:psi:pider}
        \pi_{j,\mathbf{m}}'(\xi_j^*-1) =0, \qquad \pi_{j,\mathbf{m}}'(\xi_j^*) = \dfrac{\log{\widehat{m}}}{p}.
        \end{equation}
        Further,
        \[
        \lambda'_{j,\mathbf{m},A}(\xi;\boldsymbol\omega) = A \dfrac{\log{\widehat m}}{p} \alpha_{j,\mathbf{m}}(\boldsymbol\omega) \beta_{j,\mathbf{m}}(\xi)
        \]
        and 
        \begin{equation}\label{e:l:psi:lder}
        \lambda'_{j,\mathbf{m},A}(\xi_j^*;\boldsymbol\omega) = A \dfrac{\log{\widehat m}}{p}
        \end{equation}
        since $\alpha_{j,\mathbf{m}}(\boldsymbol\omega) \beta_{j,\mathbf{m}}(\xi_j^*) = 1$ from Lemma~\ref{l:lambda}\ref{l:psi-zero}. Comparing \eqref{e:l:psi:pider} and \eqref{e:l:psi:lder} concludes the proof. 
    \item 
        The assumption on $\xi$ ensure that the profile $\Psi_{j,\mathbf{m},A}(\xi;\boldsymbol\omega)$ is in the cubic regime. On the other hand, Lemma~\ref{l:lambda}\ref{l:lambda-zeros-distance} yields that $\Psi_{j+1,\mathbf{m},A}(\xi+1;\boldsymbol\omega)$ can be either in the cubic or exponential regime. In particular, the assumption $\xi_{j+1}^* - \xi_j^* < 1$ ensures that
        \begin{equation} \label{e:psi-at-cub-tmp-1}
            \xi+1 \ge \xi_j^* > \xi_{j+1}^*-1. 
        \end{equation}
        Therefore, we consider two cases:
    
        \noindent \textbf{Case $\xi > \xi_{j+1}^* -1$:} Since $\Psi_{j,\mathbf{m},A}(\xi;\boldsymbol\omega)$ is in the cubic regime and $\Psi_{j+1,\mathbf{m},A}(\xi+1;\boldsymbol\omega)$ is in the exponential regime, we have $\Psi_{j,\mathbf{m},A}(\xi;\boldsymbol\omega) <0$ and $\Psi_{j+1,\mathbf{m},A}(\xi+1;\boldsymbol\omega)>0$, respectively. This directly leads to the conclusion. 
    
        \noindent \textbf{Case $\xi \le \xi_{j+1}^*-1$:} Now, we have 
        \begin{align*}
            \Psi_{j+1,\mathbf{m},A}(\xi+1;\boldsymbol\omega) &= \frac{ A\log(\widehat{m})}{3p}\left( (\xi+1-(\xi_{j+1}^*-1))^3-1 \right), \\
            \Psi_{j,\mathbf{m},A}(\xi;\boldsymbol\omega) &= \frac{A \log(\widehat{m})}{3p}\left( (\xi-(\xi_j^*-1))^3-1 \right).
        \end{align*}
        Expressing the difference we can compute
        \begin{align*}
            &\frac{3p}{A \log(\widehat m)}\left( \Psi_{j+1,\mathbf{m},A}(\xi+1;\boldsymbol\omega) - \Psi_{j,\mathbf{m},A}(\xi;\boldsymbol\omega) \right)  \\
            &=\left( (\xi+1-(\xi_{j+1}^*-1))^3-1 \right) - \left( (\xi-(\xi_j^*-1))^3-1 \right) 
            = (\xi+1-(\xi_{j+1}^*-1))^3 - (\xi-(\xi_{j}^*-1))^3  \\
            &= (\xi_j^* - \xi_{j+1}^* +1)\left((\xi+1-(\xi_{j+1}^*-1))^2 + (\xi-(\xi_{j}^*-1))^2 +(\xi+1-(\xi_{j+1}^*-1))(\xi-(\xi_{j}^*-1)) \right)\\ 
            &>0. 
        \end{align*}
        Indeed, as a direct consequence of Lemma~\ref{l:lambda}\ref{l:lambda-zeros-distance} we have $\xi_j^* - \xi_{j+1}^* +1>0$. The inequality $\xi+1-(\xi_{j+1}^*-1)>0$ follows from~\eqref{e:psi-at-cub-tmp-1} and $\xi-(\xi_{j}^*-1)>0$ is implied by the assumption of this lemma.
    \end{enumerate}
\end{proof}

The next task is to show the non-negativity  $\mathfrak{D}_{\boldsymbol\omega,\mathbf{k}}[\boldsymbol\Psi_{\mathbf{m},A}(\, \cdot \,;\boldsymbol\omega)](\xi) \ge 0$, with strict positivity on $\xi \in [\xi_j^*-1, \xi_j^*]$ for $j \in [1,p]_\mathbb{N}$. In order to be able to handle the convex cubic transition part we apply Jensen's inequality. 

\begin{lemma}[Jensen's inequality]
Let $n \in \mathbb{N}, n \ge 2$, $\varphi: [a,b] \to \mathbb{R}$ be a convex function where $-\infty \le a < b \le \infty$. Furthermore let $w_1, \ldots, w_n > 0$ be positive weights and $x_1, \ldots, x_n \in [a,b]$. Then %the following holds
we have 
\[
    \varphi \left( \frac{ \sum_{i=1}^n w_i x_i }{\sum_{i=1}^n w_i}\right) \le \frac{\sum_{i=1}^n w_i \varphi (x_i)}{\sum_{i=1}^n w_i}.
\]
\end{lemma}

\begin{figure}[t]
    \centering
    \includegraphics[width=.8\textwidth]{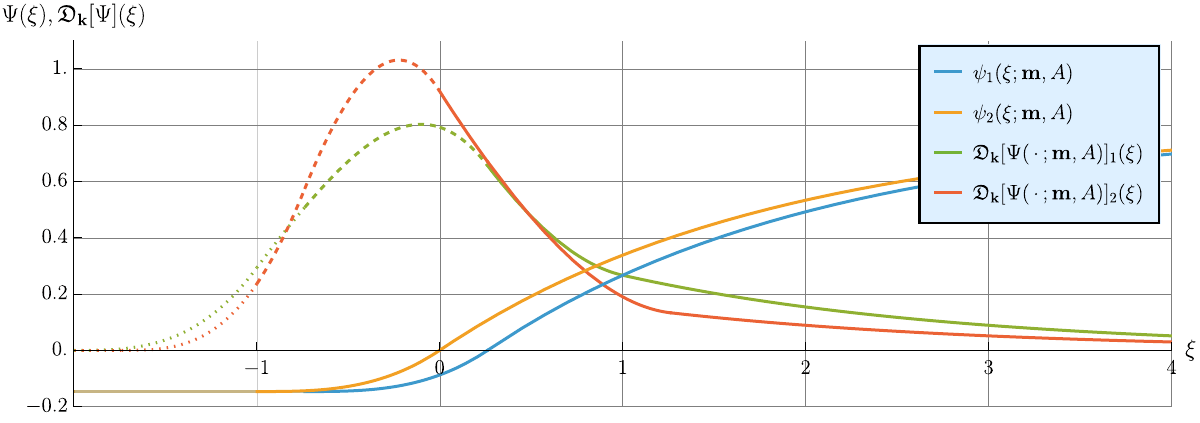}
    \caption{Illustration of Lemma~\ref{l:Dcal}. The green and red curves depict $\mathfrak{D}_{\boldsymbol\omega,\mathbf{k}}[\boldsymbol\Psi_{\mathbf{m},A}(\, \cdot \,;\boldsymbol\omega)]_j(\xi)$ for $j = 1,2$. The used parameters are $\mathbf{k}=\langle 2,4 \rangle$, $\mathbf{m} = \langle 1,3 \rangle$, $A = 0.8$. The diffusion terms are split into parts corresponding to the constant regime (dotted), cubic regime (dashed), and exponential regime (full) of the respective underlying profiles,~\eqref{e:psi}. Note that this division is still crude since the regime of the diffusion terms $\mathfrak{D}_{\boldsymbol\omega,\mathbf{k}}[\boldsymbol\Psi_{\mathbf{m},A}(\xi;\boldsymbol\omega)]_j$ depends also on the regime change of the neighboring profiles $\Psi_{j,\mathbf{m},A}(\xi;\boldsymbol\omega)$. For details see proofs of separate claims of Lemma~\ref{l:Dcal}. }
    \label{f:left-travel-diff}
\end{figure}

In the following lemma, the claims could be written in a more concise form. The current structure however mirrors various types of techniques used to prove them. We encourage the reader to consult Figure~\ref{f:left-travel-diff} for illustration.

\begin{lemma}\label{l:Dcal}
Let $\widehat{k} >1$, $\mathbf{m} \in \Omega_\mathbf{k}$ and $A \in (0,1)$, then 
\begin{enumerate}[label=({\roman*})]
    \item $\mathfrak{D}_{\boldsymbol\omega,\mathbf{k}}[\boldsymbol\Psi_{\mathbf{m},A}(\, \cdot \,;\boldsymbol\omega)](\xi) \ge 0$ for $\xi \le \xi_j^*-1$,\label{l:Dcal-1}
    \item $\mathfrak{D}_{\boldsymbol\omega,\mathbf{k}}[\boldsymbol\Psi_{\mathbf{m},A}(\, \cdot \,;\boldsymbol\omega)](\xi) > 0$ for $\xi = \xi_j^*-1$,
    \item $\mathfrak{D}_{\boldsymbol\omega,\mathbf{k}}[\boldsymbol\Psi_{\mathbf{m},A}(\, \cdot \,;\boldsymbol\omega)](\xi) > 0$ for $\xi \in (\xi_j^*-1, \xi_j^*)$,
    \item $\mathfrak{D}_{\boldsymbol\omega,\mathbf{k}}[\boldsymbol\Psi_{\mathbf{m},A}(\, \cdot \,;\boldsymbol\omega)](\xi) > 0$ for $\xi = \xi_j^*$ \label{l:Dcal-4},
    \item $\mathfrak{D}_{\boldsymbol\omega,\mathbf{k}}[\boldsymbol\Psi_{\mathbf{m},A}(\, \cdot \,;\boldsymbol\omega)](\xi) \ge \mathfrak{D}_{\boldsymbol\omega,\mathbf{k}}[\boldsymbol\lambda_{\mathbf{m},A}(\, \cdot \,;\boldsymbol\omega)](\xi) > 0$ for $\xi > \xi_j^*$,
    \item $\max\limits_{\xi \in \mathbb{R}} \Psi'_{j,\mathbf{m},A}(\xi; \mathbf{\boldsymbol{\omega}}) = \Psi'_{j,\mathbf{m},A}(\xi_j^*;\mathbf{\boldsymbol{\omega}}) = A \frac{\log \widehat m}{p}.$
\end{enumerate}
\end{lemma}
\begin{proof}
In this proof we repeatedly use Lemma~\ref{l:lambda}\ref{l:lambda-zeros-distance}. In particular, we have 
\begin{align}
    \xi_{j+1}^* &< \xi_j^* + 1, \label{e:xi_jp1}\\
    \xi_{j-1}^* &> \xi_j^* - 1. \label{e:xi_jm1}
\end{align}
For each of the items we first determine which regime  the profiles $\Psi_{j,\mathbf{m},A}(\xi;\boldsymbol\omega)$ belong to and use \eqref{e:lapl-cont}, i.e.,
\begin{align*}
    &\mathfrak{D}_{\boldsymbol\omega,\mathbf{k}} [\boldsymbol\Psi_{\mathbf{m},A}(\, \cdot \,;\boldsymbol\omega)]_j(\xi) \\ 
    &=\omega_{j+1/2} k_j \Psi_{j+1,\mathbf{m},A}(\xi+1;\boldsymbol\omega) - (\omega_{j+1/2} k_j+ \omega_{j-1/2}) \Psi_{j,\mathbf{m},A}(\xi;\boldsymbol\omega) +\omega_{j-1/2 }\Psi_{j-1,\mathbf{m},A}(\xi-1;\boldsymbol\omega) \\
    &= \omega_{j+1/2} k_j (\Psi_{j+1,\mathbf{m},A}(\xi+1;\boldsymbol\omega) - \Psi_{j,\mathbf{m},A}(\xi;\boldsymbol\omega) ) + \omega_{j-1/2 }(\Psi_{j-1,\mathbf{m},A}(\xi-1;\boldsymbol\omega)- \Psi_{j,\mathbf{m},A}(\xi;\boldsymbol\omega)).
\end{align*}
\begin{enumerate}[label=({\roman*})]
    \item Since $\xi \le \xi_j^* -1$, $\Psi_{j,\mathbf{m},A}(\xi;\boldsymbol\omega)= -A \frac{\log(\widehat{m})}{3p}$ is in the constant regime, and we have 
    \[
    \xi - 1 \le \xi_j^* - 2 < \xi_{j-1}^* -1
    \]
    by~\eqref{e:xi_jm1}, i.e., $\Psi_{j-1,\mathbf{m},A}(\xi-1;\boldsymbol\omega) =  -A \frac{\log(\widehat{m})}{3p}$ is also in the constant regime. Regardless of the position of the zero $\xi_{j+1}^*$ the estimate from below $\Psi_{j+1,\mathbf{m},A}(\xi+1;\boldsymbol\omega) \ge -A \frac{\log(\widehat{m})}{3p}$ holds. Putting all together we obtain
   \begin{align*}
  &\mathfrak{D}_{\boldsymbol\omega,\mathbf{k}} [\boldsymbol\Psi_{\mathbf{m},A}(\, \cdot \,;\boldsymbol\omega)]_j(\xi)  \\&=\omega_{j+1/2} k_j (\Psi_{j+1,\mathbf{m},A}(\xi+1;\boldsymbol\omega) - \Psi_{j,\mathbf{m},A}(\xi;\boldsymbol\omega) ) + \omega_{j-1/2 }(\Psi_{j-1,\mathbf{m},A}(\xi-1;\boldsymbol\omega)- \Psi_{j,\mathbf{m},A}(\xi;\boldsymbol\omega))  \\ 
        &= \omega_{j+1/2}k_j(\Psi_{j+1,\mathbf{m},A}(\xi+1;\boldsymbol\omega) - \Psi_{j,\mathbf{m},A}(\xi;\boldsymbol\omega)) \ge 0. 
    \end{align*} 
    \item The argumentation is very similar to the previous point. We can however 
    use the sharper bound 
    \[
        \xi + 1 = \xi_j^* > \xi_{j+1}^* - 1, 
    \]
    i.e., $\Psi_{j+1,\mathbf{m},A}(\xi+1;\boldsymbol\omega)$ is firmly outside the constant regime and thus $\Psi_{j+1,\mathbf{m},A}(\xi_j^*+1;\boldsymbol\omega) > -A \frac{\log(\widehat{m})}{3p} = \Psi_{j,\mathbf{m},A}(\xi_j^*;\boldsymbol\omega)$. Subsequently 
   \begin{align*}
  &\mathfrak{D}_{\boldsymbol\omega,\mathbf{k}} [\boldsymbol\Psi_{\mathbf{m},A}(\, \cdot \,;\boldsymbol\omega)]_j(\xi)  \\
  &=\omega_{j+1/2} k_j (\Psi_{j+1,\mathbf{m},A}(\xi+1;\boldsymbol\omega) - \Psi_{j,\mathbf{m},A}(\xi;\boldsymbol\omega) ) + \omega_{j-1/2 }(\Psi_{j-1,\mathbf{m},A}(\xi-1;\boldsymbol\omega)- \Psi_{j,\mathbf{m},A}(\xi;\boldsymbol\omega)) \\ 
        &= \omega_{j+1/2}k_j(\Psi_{j+1,\mathbf{m},A}(\xi+1;\boldsymbol\omega)) - \Psi_{j,\mathbf{m},A}(\xi;\boldsymbol\omega)) \ge 0. 
    \end{align*} 
    \item
    Assume $\xi \in ( \xi_j^*-1, \xi_j^*) $, then by Lemma~\ref{l:lambda}\ref{l:lambda-zeros-distance} we have 
    \begin{align*}
        \xi - 1 &< \xi_j^* -1 < \xi_{j-1}^*, \\
        \xi + 1 &> \xi_j^* > \xi_{j+1}^* - 1,
    \end{align*}
    i.e., if the middle profile $\Psi_{j,\mathbf{m},A}(\xi;\boldsymbol\omega)$ in $\mathfrak{D}_{\boldsymbol\omega,\mathbf{k}}$ is in the cubic regime, the previous one $\Psi_{j-1,\mathbf{m},A}(\xi-1;\boldsymbol\omega)$ cannot be in the exponential regime. Analogously, the next one, $\Psi_{j+1,\mathbf{m},A}(\xi+1;\boldsymbol\omega)$, cannot be in the constant regime. We separate the proof in several cases, each of them given by a strict inequality. The switch points are taken care of by the fact that  $\mathfrak{D}_{\boldsymbol\omega,\mathbf{k}}$ is continuously differentiable. 

    In all four cases we use the observation that
    \[
    \begin{aligned}
        &\mathfrak{D}_{\boldsymbol\omega,\mathbf{k}}[\boldsymbol\Psi_{\mathbf{m},A}(\, \cdot \,;\boldsymbol\omega)]_j(\xi) - \mathfrak{D}_{\boldsymbol\omega,\mathbf{m}}[\boldsymbol\Psi_{\mathbf{m},A}(\, \cdot \,;\boldsymbol\omega)]_j(\xi) = \omega_{j+1/2}(k_j - m_j)(\Psi_{j+1,\mathbf{m},A}(\xi+1;\boldsymbol\omega) \\ 
        &\qquad - \Psi_{j,\mathbf{m},A}(\xi;\boldsymbol\omega)) > 0, 
    \end{aligned}
    \]
    since $m_j< k_j$ and the ordering result in Lemma~\ref{l:psi}\ref{l:psi-at-cub}. This enables us to drop the parameter $k_j$ in the computation.

    \textbf{Case $\xi_{j-1}^*-1 \le \xi \le \xi_{j+1}^*$:}
    We now consider the case in which all three profiles are in the cubic regime as a base case. 

    We use Jensen's inequality twice. First, consider a convex function $\tau_1(s):=\sqrt[p]{\widehat m}^s$. Applying the identity~\eqref{e:exp-identity} from Lemma~\ref{l:lambda}\ref{l:exp-identity} we have 
    \begin{align*}
        \sqrt[p]{\widehat m} ^{\xi_j^*} =  \frac{ \omega_{j-1/2}\sqrt[p]{\widehat m}^{\xi_{j-1}^*+1}+ \omega_{j+1/2}  m_j \sqrt[p]{\widehat m}^{\xi_{j+1}^*-1}}{\omega_{j+1/2} m_j + \omega_{j-1/2}} \ge \sqrt[p]{\widehat m}^{\frac{\omega_{j-1/2}(\xi_{j-1}^*+1)+\omega_{j+1/2}m_j(\xi_{j+1}^*-1)}{\omega_{j+1/2} m_j + \omega_{j-1/2}}}
    \end{align*}
    which by monotonicity of the exponential function yields
    \[
       (\omega_{j+1/2} m_j + \omega_{j-1/2})\xi_j^* \ge \omega_{j-1/2}(\xi_{j-1}^*+1)+\omega_{j+1/2} m_j(\xi_{j+1}^*-1).
    \]
    After switching the sign and adding $(\omega_{j+1/2} m_j + \omega_{j-1/2})(\xi+1)$ to both sides of the inequality and rearranging,  we get
    \begin{align} \label{e:jensen-1}
        (\xi-(\xi_j^*-1)) \le \frac{\omega_{j-1/2}((\xi -1)-(\xi_{j-1}^*-1)) + \omega_{j+1/2} m_j ((\xi +1)-(\xi_{j+1}^*-1))}{\omega_{j+1/2} m_j + \omega_{j-1/2}}. 
    \end{align} 
    Now, the assumptions of the current case ensure that both sides of the inequality are positive. Thus, we can cube both sides and use Jensen's inequality again on the right-hand side to find
    \begin{align*}
       (\xi-(\xi_j^*-1))^3 &\le \left( \frac{\omega_{j-1/2}((\xi -1)-(\xi_{j-1}^*-1)) + \omega_{j+1/2} m_j ((\xi +1)-(\xi_{j+1}^*-1))}{\omega_{j+1/2} m_j + \omega_{j-1/2}} \right)^3 \\
        & \le \frac{\omega_{j-1/2}((\xi -1)-(\xi_{j-1}^*-1))^3 + \omega_{j+1/2} m_j ((\xi +1)-(\xi_{j+1}^*-1))^3}{\omega_{j+1/2} m_j + \omega_{j-1/2}}.
    \end{align*}
    Comparing with $\mathfrak{D}_{\boldsymbol\omega,\mathbf{m}}[\boldsymbol\Psi_{\mathbf{m},A}(\, \cdot \,;\boldsymbol\omega)]_j(\xi)$ and using the estimate 
    \[
    \mathfrak{D}_{\boldsymbol\omega,\mathbf{k}}[\boldsymbol\Psi_{\mathbf{m},A}(\, \cdot \,;\boldsymbol\omega)]_j(\xi) > \mathfrak{D}_{\boldsymbol\omega,\mathbf{m}}[\boldsymbol\Psi_{\mathbf{m},A}(\, \cdot \,;\boldsymbol\omega)]_j(\xi)
    \]
    concludes this part of the proof.

    \textbf{Case $(\xi-1 \le \xi_{j-1}^*-1) \wedge (\xi +1 \le \xi_{j+1}^*) $: }Now assume that $\Psi_{j,\mathbf{m},A}(\xi-1;\boldsymbol\omega) = -\frac{A \log(\widehat m)}{3 p}$ while $\Psi_{j+1,\mathbf{m},A}(\xi+1;\boldsymbol\omega)$ is in the cubic regime. First note that we derived~\eqref{e:jensen-1} regardless of the regime assumption. In particular $(\xi -1)-(\xi_{j-1}^*-1)<0$, thus we deduce
    \begin{align*}
        (\xi-(\xi_j^*-1)) \le \frac{\omega_{j+1/2} m_j ((\xi +1)-(\xi_{j+1}^*-1))}{\omega_{j+1/2} m_j + \omega_{j-1/2}} = \frac{1 \cdot 0 + \omega_{j+1/2} m_j ((\xi +1)-(\xi_{j+1}^*-1))}{\omega_{j+1/2} m_j + \omega_{j-1/2}}. 
    \end{align*}
    Next we apply Jensen's inequality with one argument being $0$, yielding 
    \[
       (\omega_{j+1/2} m_j + \omega_{j-1/2})(\xi-(\xi_j^*-1))^3 \le \omega_{j+1/2} m_j ((\xi +1)-(\xi_{j+1}^*-1))^3.
    \]
    Adding $\omega_{j+1/2} m_j + \omega_{j-1/2}$ to both sides of the inequality directly implies $\mathfrak{D}_{\boldsymbol\omega,\mathbf{m}}[\boldsymbol\Psi_{\mathbf{m},A}(\, \cdot \,;\boldsymbol\omega)]_j(\xi) \allowbreak \ge 0$.

    \textbf{Case $(\xi_{j-1}^* < \xi-1) \wedge (\xi_{j+1}^* < \xi + 1 )$:} Now the profile $\Psi_{j+1,\mathbf{m},A}(\xi+1;\boldsymbol\omega)$ is in the exponential regime and $\Psi_{j-1,\mathbf{m},A}(\xi-1;\boldsymbol\omega)$ is in the cubic regime. The idea of the proof differs from the previous ones. In particular, we use the fact that $\mathfrak{D}_{\boldsymbol\omega,\mathbf{m}}[\boldsymbol\Psi_{\mathbf{m},A}(\, \cdot \,;\boldsymbol\omega)]_j(\xi)$ is strictly positive in the endpoints of the interval $(\max\{ \xi_{j-1}^*+1, \xi_{j+1}^*\}, \xi_{j})$. The positivity at the right endpoint follows from Lemma~\ref{l:Dcal}\ref{l:Dcal-4}. The remaining ingredient is to prove concavity of $\mathfrak{D}_{\boldsymbol\omega,\mathbf{m}}[\boldsymbol\Psi_{\mathbf{m},A}(\, \cdot \,;\boldsymbol\omega)]_j(\xi)$ over the sought interval.
    Compute 
    \begin{align*}
       &\mathfrak{D}''_{\boldsymbol\omega,\mathbf{m}}[\boldsymbol\Psi_{\mathbf{m},A}(\, \cdot \,;\boldsymbol\omega)]_j(\xi)  \\
       &\qquad = \omega_{j+1/2} m_j \Psi''_{j+1}(\xi+1)  - (\omega_{j+1/2} m_j+\omega_{j-1/2}) \Psi''_{j}(\xi) + \omega_{j-1/2} \Psi''_{j-1}(\xi-1) \\ 
       &\qquad \le  \omega_{j-1/2} \Psi''_{j-1}(\xi-1) - (\omega_{j+1/2} m_j+\omega_{j-1/2}) \Psi''_{j}(\xi) \\ 
       &\qquad = \frac{2 \log(\widehat m)}{p} \left[ \omega_{j-1/2} (\xi -1 - (\xi_{j-1}^* - 1)) - (\omega_{j+1/2} m_j + \omega_{j-1/2}) (\xi - (\xi_j^*-1)) \right] \\
       &\qquad \le \frac{2 \log(\widehat m)}{p} \omega_{j-1/2} (\xi_j^* - \xi_{j-1}^* -1) < 0.
    \end{align*}
    In the second step, we used the strict concavity of $\Psi_{j+1,\mathbf{m},A}(\xi+1;\boldsymbol\omega)$, resp. $\lambda_{j+1,\mathbf{m},A}(\xi+1;\boldsymbol\omega)$. 
    
    \textbf{Case $(\xi-1 < \xi_{j-1}^*) \wedge (\xi_{j+1}^* < \xi + 1 )$:} This situation is similar to the previous one, but simpler. In particular, the only non-strictly-concave part in the sum $\mathfrak{D}_{\boldsymbol\omega,\mathbf{m}}[\boldsymbol\Psi_{\mathbf{m},A}(\, \cdot \,;\boldsymbol\omega)]_j(\xi)$ is the profile $\Psi_{j-1,\mathbf{m},A}(\xi-1;\boldsymbol\omega)$. However, it is constant and we can thus directly deduce  $\mathfrak{D}''_{\boldsymbol\omega,\mathbf{m}}[\boldsymbol\Psi_{\mathbf{m},A}(\, \cdot \,;\boldsymbol\omega)]_j(\xi) \allowbreak <0$. 
    
    \item 
    The cornerstone observation is that  $\Psi_{j,\mathbf{m},A}(\xi;\boldsymbol\omega) \ge \lambda_{j,\mathbf{m},A}(\xi;\boldsymbol\omega)$, $\xi \in \mathbb{R}$, and furthermore \allowbreak $\Psi_{j,\mathbf{m},A}(\xi;\boldsymbol\omega) \allowbreak > \lambda_{j,\mathbf{m},A}(\xi;\boldsymbol\omega)$, $ \xi < \xi_j^*$. 
    \begin{align*}
        &\mathfrak{D}_{\boldsymbol\omega,\mathbf{k}}[\boldsymbol\Psi_{\mathbf{m},A}(\, \cdot \,;\boldsymbol\omega)]_j(\xi)  \\
        &= \omega_{j+1/2} k_j \Psi_{j+1,\mathbf{m},A}(\xi+1;\boldsymbol\omega)  - (\omega_{j+1/2} k_j+\omega_{j-1/2}) \Psi_{j,\mathbf{m},A}(\xi;\boldsymbol\omega) + \omega_{j-1/2} \Psi_{j-1,\mathbf{m},A}(\xi-1;\boldsymbol\omega)\\
        &= \omega_{j+1/2} k_j \Psi_{j+1,\mathbf{m},A}(\xi+1;\boldsymbol\omega)  - (\omega_{j+1/2} k_j+\omega_{j-1/2}) \lambda_{j,\mathbf{m},A}(\xi;\boldsymbol\omega) + \omega_{j-1/2} \Psi_{j-1,\mathbf{m},A}(\xi-1;\boldsymbol\omega)\\
        &> \omega_{j+1/2} k_j \lambda_{j+1,\mathbf{m},A}(\xi+1;\boldsymbol\omega)  - (\omega_{j+1/2} k_j+\omega_{j-1/2}) \lambda_{j,\mathbf{m},A}(\xi;\boldsymbol\omega) + \omega_{j-1/2} \lambda_{j-1,\mathbf{m},A}(\xi-1;\boldsymbol\omega)\\
        &= \mathfrak{D}_{\boldsymbol\omega,\mathbf{k}}[\boldsymbol\lambda_{\mathbf{m},A}(\, \cdot \,;\boldsymbol\omega)]_j>0.
    \end{align*} 
    The strict inequality follows from the fact that $\xi - 1 = \xi_j^*-1 < \xi_{j-1}^*$, i.e., $\Psi_{j-1,\mathbf{m},A}(\xi-1;\boldsymbol\omega)$ cannot be in the exponential regime.
    \item 
    As in the previous point, we can estimate the $\boldsymbol\Psi_{\mathbf{m},A}(\xi;\boldsymbol\omega)$ profiles from below by the $\boldsymbol\lambda_{\mathbf{m},A}(\xi;\boldsymbol\omega)$ exponential profiles. This leads to the identical chain of inequalities
    \begin{align*}
        & \mathfrak{D}_{\boldsymbol\omega,\mathbf{k}}[\boldsymbol\Psi_{\mathbf{m},A}(\, \cdot \,;\boldsymbol\omega)]_j(\xi)  \\ 
        &= \omega_{j+1/2} k_j \Psi_{j+1,\mathbf{m},A}(\xi+1;\boldsymbol\omega)  - (\omega_{j+1/2} k_j+\omega_{j-1/2}) \Psi_{j,\mathbf{m},A}(\xi;\boldsymbol\omega) + \omega_{j-1/2} \Psi_{j-1,\mathbf{m},A}(\xi-1;\boldsymbol\omega)\\
        &= \omega_{j+1/2} k_j \Psi_{j+1,\mathbf{m},A}(\xi+1;\boldsymbol\omega)  - (\omega_{j+1/2} k_j+\omega_{j-1/2}) \lambda_{j,\mathbf{m},A}(\xi;\boldsymbol\omega) + \omega_{j-1/2} \Psi_{j-1,\mathbf{m},A}(\xi-1;\boldsymbol\omega)\\
        &> \omega_{j+1/2} k_j \lambda_{j+1,\mathbf{m},A}(\xi+1;\boldsymbol\omega)  - (\omega_{j+1/2} k_j+\omega_{j-1/2}) \lambda_{j,\mathbf{m},A}(\xi;\boldsymbol\omega) + \omega_{j-1/2} \lambda_{j-1,\mathbf{m},A}(\xi-1;\boldsymbol\omega)\\
        &= \mathfrak{D}_{\boldsymbol\omega,\mathbf{k}}[\boldsymbol\lambda_{\mathbf{m},A}(\, \cdot \,;\boldsymbol\omega)]_j>0.
    \end{align*} 
    for $\xi > \xi_j^*$. The last inequality directly uses Lemma~\ref{l:lambda}\ref{l:d-lambda} and the strict inequality $m_j < k_j$ that follows from $\mathbf{m} \in \Omega_\mathbf{k}$. 
    \item 
    This result follows from the smoothness result $\Psi_{j,\mathbf{m},A}(\xi;\boldsymbol\omega) \in C^1(\mathbb{R})$ (Lemma~\ref{l:psi}\ref{l:psi-smooth}), direct computation, and properties of elementary functions.
\end{enumerate}
\end{proof}

Note that Lemma~\ref{l:Dcal} implies $\mathfrak{D}_{\boldsymbol\omega,\mathbf{k}}[\boldsymbol\Psi_{\mathbf{m},A}(\, \cdot \,;\boldsymbol\omega)]_j(\xi) \ge 0$, with the inequality being strict for $\xi \ge \xi_j^*-1$.
This positions us towards using the comparison principle for MFDEs, Corollary~\ref{cor:comparison-MFDE}, for 
sufficiently large diffusion strengths $d \gg 1$. Indeed, the influence of the linear part then outweighs the nonlinear part. 

We now generalize the definition of $\nu$ from \eqref{e:nu} by writing
\begin{align} \label{e:nu-tilde}
        \widetilde\nu_A(a) &:= \max\limits_{s \in [1-\frac{a}{A},1]} -\frac{g(A(1-s);a)}{As}, \qquad a < A \le 1.
\end{align}
We also re-introduce the expressions $\mu_{j,\mathbf{m}}(\boldsymbol\omega,\mathbf{k})$ originally defined in~\eqref{e:mu_j-intro}
\begin{align} \label{e:mu_j}
   \mu_{j,\mathbf{m}}(\boldsymbol\omega,\mathbf{k}) = \frac{m_j \alpha_{j,\mathbf{m}}(\boldsymbol\omega) }{(\widehat{m} -1)(k_j-m_j)}, 
\end{align}
and define
\begin{align} \label{e:d-tilde}
   \widetilde{d}_{j,\mathbf{m},A}(a;\boldsymbol\omega,\mathbf{k}) := \frac{p}{\widehat{\omega}} \mu_{j,\mathbf{m}}(\boldsymbol\omega,\mathbf{k}) \widetilde{\nu}_A(a).
\end{align}
In contrast to the formulas in Theorem~\ref{t:left-travel},  \eqref{e:nu-tilde}--\eqref{e:d-tilde} depend on the sub-solution parameters $\mathbf{m}$ and $A$. We now show that the lower threshold $\widetilde{d}_{j,\mathbf{m},A}(a;\boldsymbol\omega,\mathbf{k})$ is sufficient to be able to apply Corollary~\ref{cor:comparison-MFDE} for any $\mathbf{m} \in \Omega_\mathbf{k}$, and $A \in [a, 1)$. In particular, we need to incorporate the strictly negative sub-solution speed $\overline{c}<0$. 

\begin{lemma}\label{l:negative-speed}
Let $\widehat{k}>1$, $\mathbf{m} \in \Omega_\mathbf{k}$, and pick $A \in (a,1)$. Then for each $d>0$ such that $d>\widetilde{d}_{j,\mathbf{m},A}(a;\boldsymbol\omega,\mathbf{k})$ for all $j \in [1,p]_\mathbb{N}$, there exists $\overline{c}<0$ such that the inequality $\boldsymbol{\mathcal{I}}_{a,d,\boldsymbol\omega,\mathbf{k}}[\overline{c},\boldsymbol\Psi_{\mathbf{m},A}(\, \cdot \,;\boldsymbol\omega)]_j(\xi) \allowbreak< 0$ holds for all $\xi \in \mathbb{R}$. 
\end{lemma}
\begin{proof}
    The strict inequality in $d > \widetilde{d}_{j,\mathbf{m},A}(a;\boldsymbol\omega,\mathbf{k})$ ensures the existence of $\varepsilon_j>0$ such that 
    \begin{align}
        d := \widetilde{d}_{j,\mathbf{m},A}(a;\boldsymbol\omega,\mathbf{k}) +\mu_{j,\mathbf{m}}(\boldsymbol\omega,\mathbf{k}) \varepsilon_j. 
    \end{align}
    We further define several auxiliary variables
    \begin{align}
        C_j &:= \inf\limits_{\xi \in [\xi_j^* -1, \xi_j^*]} \mathfrak{D}_{\boldsymbol\omega,\mathbf{k}}[\boldsymbol\Psi_{\mathbf{m},A}(\, \cdot \,;\boldsymbol\omega)]_j(\xi) >0
        \label{e:Cj}, \\
        D_j &:= \frac{A \log(\widehat{m})}{p}\alpha_{j,\mathbf{m}}(\boldsymbol\omega) \beta_{j,\mathbf{m}}(\xi_j^*) = \max\limits_{\xi \in \mathbb{R}} \Psi'_{j,\mathbf{m}}(\xi;\boldsymbol\omega) > 0 \label{e:Dj},
        \end{align}
        and set
        \begin{align}
            0 > \overline{c} := \max\limits_{j \in [1,p]_\mathbb{N}} \left\{ -\frac{\varepsilon_j \widehat{\omega}}{\log(\widehat{m})} ,
            -\frac{C_j}{D_j} \mu_{j,\mathbf{m}}(\boldsymbol\omega,\mathbf{k}) \varepsilon_j 
            \right\}. 
        \end{align}

        We separate the sought inequality into three disjoint cases. 

        \noindent\textbf{Case $ \xi < \xi_j^*-1$:} The profile $\Psi_{j,\mathbf{m},A}(\xi;\boldsymbol\omega)$ is in the constant regime. One can deduce
        \begin{align*}
            \boldsymbol{\mathcal{I}}_{a,d,\boldsymbol\omega,\mathbf{k}}[\overline{c},\boldsymbol\Psi_{\mathbf{m},A}(\, \cdot \,;\boldsymbol\omega)]_j(\xi) &= -\overline{c} \Psi'_{j,\mathbf{m},A}(\xi;\boldsymbol\omega) -  d \mathfrak{D}_{\boldsymbol\omega,\mathbf{k}}[\boldsymbol\Psi_{\mathbf{m},A}(\, \cdot \,;\boldsymbol\omega)]_j(\xi) - g(\Psi_{j,\mathbf{m},A}(\xi;\boldsymbol\omega);a) \\
            &=- d \mathfrak{D}_{\boldsymbol\omega,\mathbf{k}}[\boldsymbol\Psi_{\mathbf{m},A}(\, \cdot \,;\boldsymbol\omega)]_j(\xi) - g(\Psi_{j,\mathbf{m},A}(\xi;\boldsymbol\omega);a)  \\
            &< - d \mathfrak{D}_{\boldsymbol\omega,\mathbf{k}}[\boldsymbol\Psi_{\mathbf{m},A}(\, \cdot \,;\boldsymbol\omega)]_j(\xi) < 0,
        \end{align*}
        where we used the fact that $\Psi_{j,\mathbf{m},A}(\xi;\boldsymbol\omega)<0$ is constant and negative for $\xi <  \xi_j^*-1$, that $g(s;a)>0$ for $s<0$, and that $d>0$, together with Lemma~\ref{l:Dcal}\ref{l:Dcal-1}. 

        \noindent\textbf{Case $\xi \in [\xi_j^*-1, \xi_j]$:} The profile $\Psi_{j,\mathbf{m},A}(\xi;\boldsymbol\omega)$ is in the transition cubic regime.
        \begin{align*}
            \boldsymbol{\mathcal{I}}_{a,d,\boldsymbol\omega,\mathbf{k}}[\overline{c},\boldsymbol\Psi_{\mathbf{m},A}(\, \cdot \,;\boldsymbol\omega)]_j(\xi) &= -\overline{c} \Psi'_{j,\mathbf{m},A}(\xi;\boldsymbol\omega) - d \mathfrak{D}_{\boldsymbol\omega,\mathbf{k}}[\boldsymbol\Psi_{\mathbf{m},A}(\, \cdot \,;\boldsymbol\omega)]_j(\xi) - g(\Psi_{j,\mathbf{m},A}(\xi;\boldsymbol\omega);a) \\
            & \le -\overline{c} \Psi'_{j,\mathbf{m},A}(\xi;\boldsymbol\omega) - d \mathfrak{D}_{\boldsymbol\omega,\mathbf{k}}[\boldsymbol\Psi_{\mathbf{m},A}(\, \cdot \,;\boldsymbol\omega)]_j(\xi) \\
            &\le  \frac{C_j}{D_j}\mu_{j,\mathbf{m}}(\boldsymbol\omega,\mathbf{k}) \varepsilon_j D_j - \left( \frac{p}{\widehat{\omega}} \mu_{j,\mathbf{m}}(\boldsymbol\omega,\mathbf{k}) \widetilde\nu_A(a)+\mu_{j,\mathbf{m}}(\boldsymbol\omega,\mathbf{k}) \varepsilon_j\right) C_j \\
            &= -\frac{p}{\widehat{\omega}}  \mu_{j,\mathbf{m}}(\boldsymbol\omega,\mathbf{k}) \widetilde\nu_A(a) C_j < 0.
            \end{align*}

        \noindent\textbf{Case $\xi > \xi_j^*$:} The profile $\Psi_{j,\mathbf{m},A}(\xi;\boldsymbol\omega)$ is in the exponential regime.   
        \begin{align*}
            &\boldsymbol{\mathcal{I}}_{a,d,\mathbf{k}}[\overline{c},\boldsymbol\Psi_{\mathbf{m},A}(\, \cdot \,;\boldsymbol\omega)]_j(\xi)  \\
           &= -\overline{c} \Psi'_{j,\mathbf{m},A}(\xi;\boldsymbol\omega) - d \mathfrak{D}_{\boldsymbol\omega,\mathbf{k}}[\boldsymbol\Psi_{\mathbf{m},A}(\, \cdot \,;\boldsymbol\omega)]_j(\xi) - g(\Psi_{j,\mathbf{m},A}(\xi;\boldsymbol\omega);a) \\
            & \le -\overline{c} \lambda'_{j,\mathbf{m},A}(\xi;\boldsymbol\omega) - d \mathfrak{D}_{\boldsymbol\omega,\mathbf{k}}[\boldsymbol\lambda_{\mathbf{m},A}(\, \cdot \,;\boldsymbol\omega)]_j(\xi) - g(\lambda_{j,\mathbf{m},A}(\xi;\boldsymbol\omega);a) \\
            & \le \frac{\varepsilon_j \widehat{\omega}}{\log(\widehat{m})} \frac{A \log(\widehat{m})}{p} \alpha_{j,\mathbf{m}}(\boldsymbol\omega) \beta_{j,\mathbf{m}}(\xi) - \left(\frac{p}{\widehat{\omega}} \mu_{j,\mathbf{m}}(\boldsymbol\omega,\mathbf{k}) \widetilde{\nu}_A(a) + \varepsilon_j \mu_{j,\mathbf{m}}(\boldsymbol\omega,\mathbf{k})\right)\frac{A \widehat{\omega}\alpha_{j,\mathbf{m}}(\boldsymbol\omega)\beta_{j,\mathbf{m}}(\xi)}{p \mu_{j,\mathbf{m}}(\boldsymbol\omega,\mathbf{k})}   \\
            & \qquad \qquad -g(A(1-\alpha_{j,\mathbf{m}}(\boldsymbol\omega)\beta_{j,\mathbf{m}}(\xi));a) \\
            & = \widetilde{\nu}_A(a) A \alpha_{j,\mathbf{m}}(\boldsymbol\omega)\beta_{j,\mathbf{m}}(\xi)- g(A(1-\alpha_{j,\mathbf{m}}(\boldsymbol\omega)\beta_{j,\mathbf{m}}(\xi));a) \\
            & = A \alpha_{j,\mathbf{m}}(\boldsymbol\omega)\beta_{j,\mathbf{m}}(\xi) \left(-\widetilde{\nu}_A(a) - \frac{g(A(1-\alpha_{j,\mathbf{m}}(\boldsymbol\omega)\beta_{j,\mathbf{m}}(\xi));a)}{A \alpha_{j,\mathbf{m}}(\boldsymbol\omega)\beta_{j,\mathbf{m}}(\xi)} \right) \le 0.
        \end{align*}
        In the considerations above, we used the formula from Lemma~\ref{l:Dcal} and the definition of $\widetilde\mu_j(\mathbf{m})$. In particular, we have 
        \[
            \mathfrak{D}_{\boldsymbol\omega,\mathbf{k}}[\boldsymbol\lambda_{\mathbf{m},A}(\, \cdot \,;\boldsymbol\omega)]_j(\xi) = \frac{A \widehat{\omega} \alpha_{j,\mathbf{m}}(\boldsymbol\omega)\beta_{j,\mathbf{m}}(\xi)}{p \mu_{j,\mathbf{m}}(\boldsymbol\omega,\mathbf{k})}. 
        \]
        The final inequality follows from the properties of the nonlinearity $g$ and the definition of $\widetilde{\nu}_A(a)$. 
\end{proof}

\begin{proof}[Proof of Theorem~\ref{t:left-travel}]
    The existence of a solution to~\eqref{e:main-mfde} satisfying~\eqref{e:main-mfde-bc} follows from Proposition~\ref{prop:periodic-wave-existence}. The uniqueness is a consequence of $c<0$ which is the focus of the rest of the proof. 

    We start with the observation that both $\mu$, $\nu$, Eq.~\eqref{e:mu},~\eqref{e:nu}, resp., are well-defined. The function 
    \[
        \max_{\dInt{j}} \frac{m_j \alpha_{j,\mathbf{m}}(\boldsymbol\omega)}{(\widehat{{m}}-1)(k_j-m_j)}, 
    \]
    is continuous and positive in the open domain $\Omega_\mathbf{k}$,\eqref{e:omega}.
    \[
        \Omega_\mathbf{k} = \left\{ \mathbf{m} \in \mathbb{R}^p \, | \, 0<m_j < k_j, \; \widehat{m}>1 \right\}.
    \]
    It is however straightforward to see that its value diverges to infinity for $\mathbf{m}$ approaching the boundary $\partial \Omega_\mathbf{k}$, and thus the minimum in $\mu$ is well-defined. The correctness of $\nu$ and its positivity follows from the condition \eqref{Hg}.

    Now, let us pick any $d := d^*(a;\boldsymbol\omega, \mathbf{k}) + \varepsilon > d^*(a;\boldsymbol\omega,\mathbf{k})$. Choose 
    \[
        \widetilde{\mathbf{m}} \in \mathrm{argmin}_{m \in \Omega_\mathbf{k}} \max_{\dInt{j}} \frac{m_j \alpha_{j,\mathbf{m}}(\boldsymbol\omega)}{(\widehat{{m}}-1)(k_j-m_j)}, 
    \]
    and by definition we have 
    \[
        \mu(\boldsymbol\omega,\mathbf{k}) = \min_{\mathbf{m} \in \Omega_\mathbf{k}} \max_{j \in [1,p]_\mathbb{N}} \mu_{j,\mathbf{m}}(\boldsymbol\omega,\mathbf{k}) = \max_{j \in [1,p]_\mathbb{N}} \mu_{j,\widetilde{\mathbf{m}}}(\boldsymbol\omega,\mathbf{k}).
    \]
    At the same time, by continuity of $\widetilde{\nu}_A(a)$ in the parameter $A$, we can find $a < \widetilde{A}<1$ possibly close to $1$ such that 
    \[
        d = \frac{p}{\widehat\omega} \mu(\boldsymbol\omega,\mathbf{k}) \nu(a) + \varepsilon \ge \frac{p}{\widehat\omega} \mu_{j,\widetilde{\mathbf{m}}}(\boldsymbol\omega,\mathbf{k}) \widetilde\nu_1(a) + \varepsilon \ge \frac{p}{\widehat\omega} \mu_{j,\widetilde{\mathbf{m}}}(\boldsymbol\omega, \mathbf{k}) \widetilde{\nu}_{\widetilde{A}}(a)+ \frac{\varepsilon}{2} > \widetilde{d}_{j,\widetilde{\mathbf{m}},\widetilde{A}}(a;\boldsymbol\omega,\mathbf{k}). 
    \]
    Finally, we can apply Lemma~\ref{l:negative-speed} and Corollary~\ref{cor:comparison-MFDE} to conclude that $c<0$. 

    Subsequently, we exploit the uniqueness of the profile from Proposition~\ref{prop:periodic-wave-existence}.
\end{proof}

\subsection{Evaluating \texorpdfstring{$\mu(\boldsymbol\omega,\mathbf{k})$, $\nu(a)$}{mu,nu}}
\label{sec:large-d-analysis}
It has been shown in the previous section in Proof of Theorem~\ref{t:left-travel} that both $\mu(\boldsymbol\omega,\mathbf{k})$ and $\nu(a)$ are well defined. Here we elaborate on how to obtain their value. We want to remind the reader that the quantities are crucial for evaluation of the diffusion threshold $d^*(a;\boldsymbol\omega,\mathbf{k})$,~\eqref{e:dStar}, above which the solutions to~\eqref{e:main-mfde},~\eqref{e:main-mfde-bc} have negative speed $c<0$ for all $a \in (0,1)$. 

\subsubsection{Explicit form and algebraic system for minimizers of \texorpdfstring{$\mu$}{mu}}
Here we show a system derived from Karush-Kuhn-Tucker (KKT) conditions which is satisfied by minimizers to $\mu$. For the sake of clarity, we recall the definition~\eqref{e:mu_j}, resp.~\eqref{e:mu_j-intro},
    \[
        \mu(\boldsymbol\omega,\textbf{k}) =  \min_{\mathbf{m} \in \Omega_\mathbf{k}} \max_{\dInt{j}} \mu_{j,\mathbf{m}}(\boldsymbol\omega,\mathbf{k})
    \]
    where 
    \[
    \mu_{j,\mathbf{m}}(\boldsymbol\omega,\mathbf{k})= \frac{m_j \alpha_{j,\mathbf{m}}(\boldsymbol\omega)}{(\widehat{{m}}-1)(k_j-m_j)},  
    \]

The important observation here is that the optimal $\mathbf{m}$ lies inside the region $\Omega_{\mathbf{k}}$, since at least one of $\mu_{j,\mathbf{m}}(\boldsymbol\omega,\mathbf{k})$ diverges to infinity when $\mathbf{m}$ approaches $\partial \Omega_{\mathbf{k}}$. This leads us to the KKT conditions for a minmax problem in an open domain 
\[
\begin{gathered}
\sum_{j=1}^p\lambda_j\nabla_{\mathbf m}\mu_{j,\mathbf m^*}=0,
\qquad \lambda_j\ge0,\qquad \sum_{j=1}^p\lambda_j=1,\\
\lambda_j(\mu_{j,\mathbf m^*}-\mu^*(\boldsymbol\omega,\mathbf{k}))=0,\qquad \dInt{j}.
\end{gathered}
\]
where $\mu^*(\boldsymbol\omega,\mathbf{k}) := \max\limits_{\dInt{j}} \mu_{j,\mathbf{m}^*}$. 
Our aim now is to show, that at optimum $\mathbf{m}^*$ all $\lambda_j>0$ are positive and thus all constraints are active. Subsequently, all $\mu_{j,\mathbf{m}^*}(\boldsymbol\omega,\mathbf{k}) = \mu^*(\boldsymbol\omega,\mathbf{k})$ are equal. 

The KKT condition can be equivalently rewritten in the matrix form
\begin{equation}
(\nabla_{\!\!\mathbf{m}\,} \boldsymbol\mu_{\mathbf{m}}(\boldsymbol\omega,\mathbf{k}))^\top \boldsymbol{\lambda} = 0, \label{e:mu-kkt}
\end{equation}
where $\nabla_{\!\!\mathbf{m}\,} \boldsymbol\mu_{\mathbf{m}}(\boldsymbol\omega,\mathbf{k})$ is a standard Jacobian matrix and $\boldsymbol{\lambda}$ is a column vector of Lagrange multipliers.

Next, Lemma~\ref{l:lambda}\ref{l:alpha-rotation} grants the shift identity for the entries of $\boldsymbol{\mu}_\mathbf{m}(\boldsymbol\omega,\mathbf{k})$, namely
\begin{align} 
(k_{j-1} - m_{j-1}) \mu_{j-1,\mathbf{m}}(\boldsymbol\omega,\mathbf{k}) = \frac{m_{j-1} \alpha_{j-1,\mathbf{m}}(\boldsymbol\omega)}{\widehat{m} -1} = \left( \frac{k_j}{m_j}-1 \right) \mu_{j,\mathbf{m}}(\boldsymbol\omega,\mathbf{k}) + \frac{\widehat\omega}{p \omega_{j-1/2}}.
\label{e:mu-shift}
\end{align}
Applying the gradient operator in the variable $\mathbf{m}$ on both sides and regrouping yields the identity
\[
 (k_{j-1}-m_{j-1})\nabla_{\!\!\mathbf{m}\,} \mu_{j-1,\mathbf{m}}(\boldsymbol\omega,\mathbf{k})  - \left( \frac{k_j}{m_j}-1 \right)\nabla_{\!\!\mathbf{m}\,} \mu_{j,\mathbf{m}}(\boldsymbol\omega,\mathbf{k}) = \mu_{j-1,\mathbf{m}}(\boldsymbol\omega,\mathbf{k}) \mathbf{e}_{j-1}- \mu_{j,\mathbf{m}}(\boldsymbol\omega,\mathbf{k}) \frac{k_j}{m_j^2} \mathbf{e}_{j}. 
\]
where $\mathbf{e}_j$ is the $j$-th Euclidean basis vector. Again, the equality can be rewritten in matrix form as 
\begin{equation}
A_\mathbf{m}(\mathbf{k}) \nabla_{\!\!\mathbf{m}\,}\boldsymbol\mu_{\mathbf{m}}(\boldsymbol\omega,\mathbf{k}) = B_\mathbf{m}(\boldsymbol\mu_{\mathbf{m}}(\boldsymbol\omega,\mathbf{k}),\mathbf{k}), \label{e:mu-shift-grad}
\end{equation}
where
\[
A_\mathbf{m}(\mathbf{k}) = 
\begin{bmatrix}
   1-\frac{k_1}{m_1} & 0 & 0 & \ldots & k_p - m_p \\ 
   k_1 - m_1 & 1-\frac{k_2}{m_2} & 0 & \ldots & 0 \\ 
   0 & k_2 - m_2 &  1-\frac{k_3}{m_3} & \vdots & 0 \\
   \vdots & \vdots & \ddots & \ddots & \vdots  \\
   0 & 0 & \ldots & k_{p-1} - m_{p-1} & 1-\frac{k_p}{m_p}  
\end{bmatrix},
\]
and 
\[
\setlength{\arraycolsep}{0pt}
 B_\mathbf{m}(\boldsymbol\mu_{\mathbf{m}}(\boldsymbol\omega,\mathbf{k}),\mathbf{k}) = 
\begin{bmatrix}
   -\mu_{1,\mathbf{m}}(\boldsymbol\omega,\mathbf{k})\frac{k_1}{m_1^2} & 0 & 0 & \ldots & \mu_{p,\mathbf{m}}(\boldsymbol\omega,\mathbf{k}) \\ 
   \mu_{1,\mathbf{m}}(\boldsymbol\omega,\mathbf{k}) & -\mu_{2,\mathbf{m}}(\boldsymbol\omega,\mathbf{k})\frac{k_2}{m_2^2} & 0 & \ldots & 0 \\ 
   0 & \mu_{2,\mathbf{m}}(\boldsymbol\omega,\mathbf{k}) & - \mu_{3,\mathbf{m}}(\boldsymbol\omega,\mathbf{k})\frac{k_3}{m_3^2} & \vdots & 0 \\
   \vdots & \vdots & \ddots & \ddots & \vdots  \\
   0 & 0 & \ldots & \mu_{p-1,\mathbf{m}}(\boldsymbol\omega,\mathbf{k}) & -\mu_{p,\mathbf{m}}(\boldsymbol\omega,\mathbf{k})\frac{k_p}{m_p^2}  
\end{bmatrix}.
\setlength{\arraycolsep}{1pt}
\]

Transposing~\eqref{e:mu-shift-grad} and expressing the gradient $\nabla_{\!\!\mathbf{m}\,}\boldsymbol\mu_{\mathbf{m}}(\boldsymbol\omega,\mathbf{k})$, we get 
\begin{equation} \label{e:mu-shift-grad-transposed}
(\nabla_{\!\!\mathbf{m}\,}\boldsymbol\mu_{\mathbf{m}}(\boldsymbol\omega,\mathbf{k}))^\top = B^\top_\mathbf{m}(\boldsymbol\mu_{\mathbf{m}}(\boldsymbol\omega,\mathbf{k}),\mathbf{k}) (A^\top_\mathbf{m}(\mathbf{k}))^{-1}.
\end{equation}
It is important to notice that $A_\mathbf{m}(\mathbf{k})$ can be inverted since $\det A_\mathbf{m}(\mathbf{k})$ can be computed by expanding along the first row; the minors are then computed from the triangular matrices
\begin{align*}
    \det A_\mathbf{m}(\mathbf{k}) & = (-1)^2 \left(1-\frac{k_1}{m_1}\right) \prod_{j=2}^p \left(1- \frac{k_j}{m_j} \right) + (-1)^{p+1} (k_p-m_p) \prod_{j=1}^{p-1} (k_j-m_j) \\
    &= \prod_{j=1}^p (k_j-m_j) \left[- \frac{1}{\widehat{m}} + (-1)^{p+1}  \right]\neq 0, 
\end{align*}
noting that $m_j < k_j$ and $\widehat{m}>1$ for $\mathbf{m} \in \Omega_\mathbf{k}$. 

Plugging~\eqref{e:mu-shift-grad-transposed} into~\eqref{e:mu-kkt} and introducing transformed Lagrange multipliers 
$\boldsymbol{\lambda}'$ we obtain two simultaneous matrix equations
\begin{align}
    B^\top_\mathbf{m}(\boldsymbol\mu_{\mathbf{m}}(\boldsymbol\omega,\mathbf{k}),\mathbf{k}) \boldsymbol{\lambda}' &= 0, \label{e:mu-sys-1} \\
    (A^\top_\mathbf{m}(\mathbf{k}))^{-1} \boldsymbol{\lambda} &= \boldsymbol{\lambda}'. \label{e:mu-sys-2} 
\end{align}
Using the fact that $\boldsymbol\mu_{\mathbf{m}}(\boldsymbol\omega,\mathbf{k})>0$,~\eqref{e:mu-sys-1}, leads to a circular system of equalities
\begin{equation}
\frac{k_{j-1}}{m_{j-1}^2} \lambda_{j-1}' = \lambda_j'. \label{e:mu-lambda-prime}
\end{equation}
The multipliers $\boldsymbol{\lambda}$ are nonzero since $\sum_{j=1}^p \lambda_j = 1$, subsequently $\boldsymbol{\lambda}'$ cannot be identically zero since $(A^\top_\mathbf{m}(\mathbf{k}))^{-1}$ is regular. This enforces 
\begin{equation}
\prod_{j=1}^p \frac{k_j}{m_j^2} = 1,\label{eq:mu-cond1}
\end{equation}
and all $\lambda_j'>0$ to be positive. 

Finally, the equation~\eqref{e:mu-sys-2} yields
\begin{align*}
    \lambda_j = - \left( \frac{k_j}{m_j}-1 \right) \lambda_{j}' + (k_j-m_j) \lambda_{j+1}' = \frac{(k_j-m_j)^2}{m_j^2} \lambda'_j >0, 
\end{align*}
where we used~\eqref{e:mu-lambda-prime}. Since all multipliers $\lambda_j>0$ are positive we know that at the optimum $\mathbf{m}^*$, all $\mu_{j,\mathbf{m}}(\boldsymbol\omega,\mathbf{k}) =  \mu^*(\boldsymbol{\omega})$ have equal value. Using this fact, one can simplify~\eqref{e:mu-shift} to 
\[
(k_{j-1} - m_{j-1}) - \left( \frac{k_j}{m_j}-1 \right)  = \frac{\widehat\omega}{p \omega_{j-1/2}\mu^*(\boldsymbol\omega,\mathbf{k})}
\]

Finally, in order to find candidate solutions to the optimization problem included in $\boldsymbol\mu_{\mathbf{m}}(\boldsymbol\omega,\mathbf{k})$ it suffices to solve the system of $p$ equations
\begin{align}
    %\begin{split}
       \omega_{j-1/2}\left( (k_{j-1}-m_{j-1}) - \left( \frac{k_j}{m_j} -1 \right) \right) &= \omega_{j+1/2}\left( (k_{j}-m_{j}) - \left( \frac{k_{j+1}}{m_{j+1}} -1 \right) \right), \qquad j \in [1, p-1]_\mathbb{N},\notag \\
        \prod_{j=1}^p k_j &= \prod_{j=1}^p m_j^2.\label{e:mu-sys-final}
    %\end{split}
\end{align}

Unfortunately, we were not able to obtain an analytical solution to the system~\eqref{e:mu-sys-final} with the exception of the cases $p=1$ and $p=2$ with $\boldsymbol\omega = \langle 1 \rangle$. In particular, the system for $p=1$ reduces to the second equation only yielding $m = \sqrt{k}$. This was treated in~\cite{hupkes2023propagation}. 

If $p=2$, i.e., $\mathbf{k} = \langle k_1, k_2 \rangle$, then the only feasible solution is
\begin{align*}
    m_1 &= \frac{1}{2} \sqrt{k_1}\left( \sqrt{k_1} - \sqrt{k_2} +\sqrt{4+\left(\sqrt{k_1}-\sqrt{k_2}\right)^2} \right) , \\
    m_2 &= \frac{1}{2} \sqrt{k_2}\left( \sqrt{k_2} - \sqrt{k_1} +\sqrt{4+\left(\sqrt{k_2}-\sqrt{k_1}\right)^2} \right), 
\end{align*}
leading to
\[
\mu^*(\langle 1 \rangle,\mathbf{k}) = \left(2 + k_1 + k_2 - 
   \sqrt{4 + \left(\sqrt{k_1} - \sqrt{k_2}\right)^2}\cdot\left(\sqrt{k_1} + \sqrt{k_2}\right) \right)^{-1}. 
\]

It is however important to note that~\eqref{e:mu-sys-final} is a finite system of equations which can be treated via numerical means to obtain an approximate value of $\mu(\boldsymbol\omega,\mathbf{k})$. 

\subsubsection{Form of \texorpdfstring{$\nu$}{nu}}
Recalling~\eqref{e:nu}, we see that
\[
        \nu(a) = \max\limits_{s \in [1-a,1]} -\frac{g(1-s;a)}{s} 
\]
depends only on the form of $g$ and the parameter $a \in (0,1)$. The function in the maximum operator is positive inside the interval $s \in (1-a,1)$ and zero at the endpoints. We thus have a (possibly nonunique) maximum inside the interval, and due to the condition \eqref{Hg} differential calculus can be used. In particular, we are interested in roots of 
\[
g'(1-s;a) s + g(1-s;a) = 0, \qquad s \in (1-a;1).
\]
The situation is much simpler for $g_\mathrm{cub}$ since 
\[
-\frac{g_{\rm cub}(1-s;a)}{s} = -(1-s)(1-s-a),
\]
which is maximized at the midpoint between the roots, i.e., $s^* = 1- a/2$. One then arrives at 
\[
\nu_\mathrm{cub}(a) = \frac{a^2}{4}. 
\]
Similar computations can be performed for the nonlinearity $g_\mathrm{exp}(a)$ to get
\[
\nu_\mathrm{exp}(a) = \frac{1}{2} \left(\sqrt{a^2+1}-1\right) e^{\sqrt{a^2+1}+a-1}.
\]
%%%%%%%%%%%%%%%%%%%%%%%%%%%%%%%%%%%%%%%%%%%%%%%%%%%%%%%%%%%%%%%%%%%
%%%%%%%%%%%%%%%%%%%%%%%%%%%%%%%%%%%%%%%%%%%%%%%%%%%%%%%%%%%%%%%%%%%
% FINITE CORE OPTIMIZATON
%%%%%%%%%%%%%%%%%%%%%%%%%%%%%%%%%%%%%%%%%%%%%%%%%%%%%%%%%%%%%%%%%%%
%%%%%%%%%%%%%%%%%%%%%%%%%%%%%%%%%%%%%%%%%%%%%%%%%%%%%%%%%%%%%%%%%%%

\section{Finite core optimization}
\label{sec:numerics-sharp}

In this section we describe how the regions \(\mathcal D^\pm(\boldsymbol{\omega},\mathbf k)\) are computed. In addition, we show how a mild generalization of the construction \eqref{eq:mr:def:u:i:s} that allows the interface to spread over several sites  narrows the gap between these regions and the observed pinning region. 

\begin{figure}[t]
    \centering
    \includegraphics[width=0.95\linewidth]{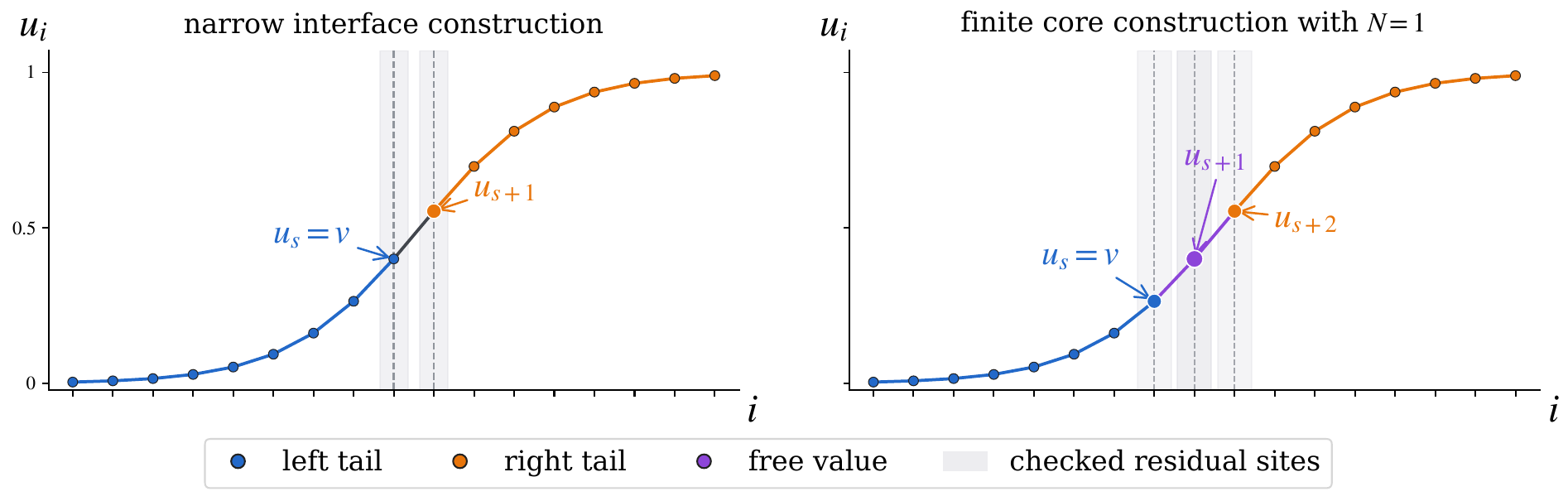}
    \caption{Schematic comparison between the narrow-interface construction \eqref{eq:mr:def:u:i:s} and the finite-core construction \eqref{eq:fc-profile}. In the left panel, the two tails meet at the adjacent values \(u_s=v\) and \(u_{s+1}\). The right panel shows \(N=1\), with one free core value \(h_1\) at \(s+1\) and the corrected right tail starting at \(s+2\).}
    \label{fig:method_explanation}
\end{figure}

\paragraph{Finite-core profiles.}
The profile \eqref{eq:mr:def:u:i:s} glues the left tail at site \(s\) to the corrected right tail at site \(s+1\), so the transition is compressed into two adjacent sites. We relax this by inserting \(N\in\mathbb N_0\) free values between the tails. 
In particular, we fix a phase \(s\dInt{}\), put \(m:=s+N\), and recall that \(\mathcal R_{i;m+1}:=\prod_{\nu=m+1}^{i-1}r_\nu\) is the right-tail product normalized at \(m+1\). For \(v>0\), \(\mathbf h=(h_1,\ldots,h_N)\), and \(w>0\) define $\underline{\mathbf{u}}^{(s, N)}$ (and generalize $\underline{\mathbf{u}}^{(s)}$ from~\eqref{eq:mr:def:u:i:s})
\begin{equation}
\label{eq:fc-profile}
\underline{u}_i^{(s,N)}(v,\mathbf h,w):=
\begin{cases}
v\,\mathcal L_{i;s}, & i\le s,\\
h_j, & i=s+j,\ 1\le j\le N,\\
1-w\,\mathcal R_{i;m+1}+C_iw^2\mathcal R_{i;m+1}^2, & i\ge m+1,
\end{cases}
\end{equation}
where \(\mathbf C\) is the positive periodic solution of \eqref{eq:C}. For \(N=0\) this is \eqref{eq:mr:def:u:i:s}; see Figure~\ref{fig:method_explanation} for a sketch. 

Exactly as in \S\ref{sec:subsolution}, the tails take values in \((0,1)\) and have positive residuals at the sites \(i\le s-1\) and \(i\ge m+2\) as soon as \(0<v<1\) and \(0<w<w_m^*\), where \(w_m^*\) is the positivity threshold \eqref{eq:wstar}. The only residuals that are not automatically non-negative are those at the \(N+2\) interface sites \(s,\ldots,m+1\). In analogy with \eqref{eq:As-compact},
we therefore define the generalized admissible set
\begin{equation}
\label{eq:fc-admissible}
\mathcal A_{s,N}^-(a,d;\boldsymbol{\omega},\mathbf k):=\Bigl\{(v,\mathbf h,w)\in(0,1)\times(0,1)^N\times(0,w_m^*):\ F_i\bigl(\underline{\mathbf{u}}^{(s,N)}(v,\mathbf h,w)\bigr)\ge0\ \text{ for } s\le i\le m+1\Bigr\},
\end{equation}
so that \(\mathcal A_{s,0}^-=\mathcal A_s^-\). Every element of \(\mathcal A_{s,N}^-\) yields a stationary sub-solution, and this forces \(c\le0\) by Corollary~\ref{cor:speed}: its tail-ordering hypothesis is verified exactly as in the proof of Proposition~\ref{prop:admissible-speed}, since \(\underline{\mathbf{u}}^{(s,N)}\) has the same tails as \eqref{eq:mr:def:u:i:s}.

We now pick a subset $\mathcal{N} \subset \mathbb N_0$ with $0 \in \mathcal{N}$ and introduce the finite-core region
\begin{equation}
\label{eq:fc-region}
\mathcal D_{\rm fc}^-(\boldsymbol{\omega},\mathbf k; \mathcal N):=\bigl\{(a,d)\in\mathcal H:\ \mathcal A_{s,N}^-(a,d;\boldsymbol{\omega},\mathbf k)\neq\varnothing\ \text{ for some } s\dInt{},\ N\in \mathcal N\bigr\}\ \supset\ \mathcal D^-(\boldsymbol{\omega},\mathbf k).
\end{equation}
By construction, we have \(c\le0\) on \(\mathcal D_{\rm fc}^-(\boldsymbol{\omega},\mathbf k;\mathcal N)\). The companion region is defined through the reflection principle of Lemma~\ref{lem:mr:reflection}, with the transformed sequences \eqref{eq:mr:refl:dual:transf}:
\begin{equation}
\label{eq:fc-region-plus}
\mathcal D_{\rm fc}^+(\boldsymbol{\omega},\mathbf k;\mathcal N):=\bigl\{(a,d)\in\mathcal H:\ (1-a,d)\in\mathcal D_{\rm fc}^-(\widetilde{\boldsymbol{\omega}},\widetilde{\mathbf k};\mathcal N)\bigr\}\ \supset\ \mathcal D^+(\boldsymbol{\omega},\mathbf k).
\end{equation}
By \eqref{eq:mr:refl:c} we have \(c\ge0\) on \(\mathcal D_{\rm fc}^+(\boldsymbol{\omega},\mathbf k;\mathcal N)\), and hence \(c=0\) on the intersection
\begin{equation}
\label{eq:fc:def:fc:0}
\mathcal D_{\rm fc}^0(\boldsymbol{\omega},\mathbf k;\mathcal N)
= 
\mathcal D_{\rm fc}^-(\boldsymbol{\omega},\mathbf k;\mathcal N)\cap\mathcal D_{\rm fc}^+(\boldsymbol{\omega},\mathbf k;\mathcal N).
\end{equation}

For each $d>0$, we define the rightmost threshold of
$\mathcal D_{\rm fc}^-$ and the leftmost threshold of
$\mathcal D_{\rm fc}^+$ by
\begin{align}
a_{\rm fc}^-(d;\boldsymbol{\omega},\mathbf k;\mathcal N)
&:=\sup\bigl\{a\in(0,1):
(a,d)\in\mathcal D_{\rm fc}^-
(\boldsymbol{\omega},\mathbf k;\mathcal N)\bigr\},\\
a_{\rm fc}^+(d;\boldsymbol{\omega},\mathbf k;\mathcal N)
&:=\inf\bigl\{a\in(0,1):
(a,d)\in\mathcal D_{\rm fc}^+
(\boldsymbol{\omega},\mathbf k;\mathcal N)\bigr\},
\end{align}
with the conventions $\sup\varnothing:=0$ and
$\inf\varnothing:=1$.
The reflection identity \eqref{eq:fc-region-plus} gives
\begin{equation}
\label{eq:num:bnd:a:pm:fc}
a_{\rm fc}^+(d;\boldsymbol{\omega},\mathbf k;\mathcal N)
=
1-a_{\rm fc}^-(d;\widetilde{\boldsymbol{\omega}},
\widetilde{\mathbf k};\mathcal N).
\end{equation}
For $\mathcal N=\{0\}$, these definitions recover the
corresponding thresholds for $\mathcal D^\pm$.

\paragraph{Numerical procedure.}
Non-emptiness of \(\mathcal A_{s,N}^-\) is decided through the sign of the scalar quantity
\begin{equation}
\label{eq:fc-tau}
\tau_{s,N}(a,d;\boldsymbol{\omega},\mathbf k):=\sup_{(v,\mathbf h,w)\in(0,1)\times(0,1)^N\times(0,w_m^*)}\ \min_{s\le i\le m+1}F_i\bigl(\underline{\mathbf{u}}^{(s,N)}(v,\mathbf h,w)\bigr).
\end{equation}
Indeed, \(\tau_{s,N}>0\) implies \(\mathcal A_{s,N}^-\neq\varnothing\), \(\tau_{s,N}<0\) implies \(\mathcal A_{s,N}^-=\varnothing\), and \(\tau_{s,N}=0\) is the boundary case, which is never encountered numerically. 
For $N=0$, we have the equivalent sign criterion
\begin{equation}
\label{eq:fc-tau0}
\tau_{s,0}(a,d;\boldsymbol\omega,\mathbf k)>0
\quad\Longleftrightarrow\quad
\sup_{0<w<w_s^*}
\sup_{\substack{0<v<1\\v\ge V_s(w)}}F_s(v;w)>0.
\end{equation}
This is a one-dimensional optimization in \(w\), with the inner supremum found by comparing the admissible endpoint and the real critical points of the cubic. 

For $N\ge1$, we numerically maximize \eqref{eq:fc-tau}
over the $N+2$ variables using a standard constrained optimizer.
For each phase $s$ and core length $N$, we introduce the numerical proxy
\begin{equation}
\label{eq:fc-threshold}
\widehat a_{s,N}^-(d;\boldsymbol{\omega},\mathbf k)
:=
\sup\bigl\{a\in(0,1):
\tau_{s,N}(a,d;\boldsymbol{\omega},\mathbf k)>0\bigr\},
\end{equation}
which we trace by continuation in $d$.
At each step, the right-most sign change 
of $\tau_{s,N}$
is bracketed and refined by bisection,
with all tail data recomputed at the current $a$ from
\eqref{eq:left-ratio-relation}, \eqref{eq:right-ratio-relation},
and \eqref{eq:C}.

The associated approximations for the boundaries
\eqref{eq:num:bnd:a:pm:fc} are then given by
\[
\widehat a_{\rm fc}^-
(d;\boldsymbol{\omega},\mathbf k;\mathcal N)
:=
\max_{s\in[1,p]_{\mathbb N}}\max_{N\in\mathcal N}
\widehat a_{s,N}^-(d;\boldsymbol{\omega},\mathbf k),
\qquad
\widehat a_{\rm fc}^+
(d;\boldsymbol{\omega},\mathbf k;\mathcal N)
:=
1-\widehat a_{\rm fc}^-
(d;\widetilde{\boldsymbol{\omega}},
\widetilde{\mathbf k};\mathcal N).
\]
For \(\omega=\langle 1 \rangle\) and $\mathcal N = \{ 0 \}$, these produce the regions \(\mathcal D^\pm(\boldsymbol{\omega},\mathbf k)\) in Figures~\ref{fig:pin-region-k12} and \ref{fig:ad}; both are compared with the numerically observed pinning boundaries in Figure~\ref{fig:comparison_region}. Allowing the interface to spread over a few sites visibly reduces the conservativeness of the two-tail construction, in particular for larger \(d\), while in the corner \((a,d)\to(1,0)\) the boundaries merge.

%%%%%%%%%%%%%%%%%%%%%%%%%%%%%%%%%%%%%%%%%%%%%%%%%%%%%%%%%%%%%%%%%%%
%%%%%%%%%%%%%%%%%%%%%%%%%%%%%%%%%%%%%%%%%%%%%%%%%%%%%%%%%%%%%%%%%%%
% BIBLIOGRAPHY
%%%%%%%%%%%%%%%%%%%%%%%%%%%%%%%%%%%%%%%%%%%%%%%%%%%%%%%%%%%%%%%%%%%
%%%%%%%%%%%%%%%%%%%%%%%%%%%%%%%%%%%%%%%%%%%%%%%%%%%%%%%%%%%%%%%%%%%

\bibliographystyle{klunumHJ}
\bibliography{ref}

\end{document}